\documentclass[11pt]{amsart}
\usepackage[utf8]{inputenc}
\usepackage[T1]{fontenc}
\usepackage{amsmath,amssymb,amsthm,mathtools}
\usepackage{enumitem}
\usepackage{array}
\usepackage[margin=2.6cm]{geometry}
\usepackage{xcolor}
\usepackage[colorlinks=true,linkcolor=blue!50!black,citecolor=blue!50!black,urlcolor=blue!50!black]{hyperref}

\theoremstyle{plain}
\newtheorem{theorem}{Theorem}[section]
\newtheorem{proposition}[theorem]{Proposition}
\newtheorem{lemma}[theorem]{Lemma}
\newtheorem{corollary}[theorem]{Corollary}

\theoremstyle{definition}

\newtheorem{example}[theorem]{Example}
\theoremstyle{remark}
\newtheorem{remark}[theorem]{Remark}
\newtheorem*{remark*}{Remark}
\newcommand{\kk}{K}
\newcommand{\C}{\mathbb{C}}
\newcommand{\Q}{\mathbb{Q}}
\newcommand{\R}{\mathbb{R}}

\newcommand{\PP}{\mathbb{P}}
\newcommand{\F}{\mathbb{F}}
\newcommand{\A}{\mathbb{A}}
\newcommand{\Hess}{\operatorname{Hess}}
\newcommand{\BH}{\operatorname{BH}}
\newcommand{\adj}{\operatorname{adj}}
\newcommand{\rk}{\operatorname{rank}}
\newcommand{\grad}{\nabla}
\newcommand{\vv}[1]{\langle #1\rangle}
\newcommand{\JC}{\mathrm{JC}}
\newcommand{\HC}{\mathrm{HC}}
\newcommand{\Rc}{\mathcal R}
\newcommand{\T}{^{\mathsf T}}
\DeclareMathOperator{\supp}{supp}
\DeclareMathOperator{\spn}{span}
\DeclareMathOperator{\tr}{tr}

\begin{document}

\title[Monge--Amp\`ere equation in four variables]{Polynomial solutions of the Monge--Amp\`ere equation in four variables: affine directions and the Legendre transform}

\author{El\.zbieta Adamus}
\address{Faculty of Applied Mathematics, AGH University of Krakow, al.~Mickiewicza 30, 30-059 Krak\'ow, Poland}
\email{esowa@agh.edu.pl}
\urladdr{https://orcid.org/0000-0003-3533-6279}

\author{Zbigniew Hajto}
\address{Faculty of Mathematics and Computer Science, Jagiellonian University, ul.~\L{}ojasiewicza 6, 30-348 Krak\'ow, Poland}
\email{zbigniew.hajto@uj.edu.pl}
\urladdr{https://orcid.org/0000-0002-1103-4404}


\begin{abstract}
Let $\phi$ be a polynomial in four variables over an algebraically closed field of characteristic zero whose Hessian determinant is a nonzero constant, that is, a polynomial solution of the Monge--Amp\`ere equation. Recently Q.~Tang stated in a Hexagon preprint that the gradient of $\phi$ is a polynomial automorphism unless $\phi$ is, in suitable linear coordinates, the doubling $x_3C(x_1,x_2)+x_4A(x_1,x_2)+E(x_1,x_2)$ of a non-injective Keller map $(C,A)$ of the plane, and proposed a strategy of proof. Following this strategy, we give a complete proof of the stronger statement that $\phi$ has an affine direction, a vector $v\neq0$ with $D_v^2\phi=0$, so that a single vector is isotropic for all the quadratic forms $\Hess\phi(x)$. A partial Legendre transform then gives the dichotomy. If the gradient of $\phi$ is a polynomial automorphism, the Legendre transform of $\phi$ is again a polynomial. The dichotomy shows that the Hessian conjecture in four variables is equivalent to the planar Jacobian conjecture. Every real solution with negative Hessian determinant has a polynomial Legendre transform. For solutions with a direction of degree at most two, in particular in degree at most six, we give a second proof by homogeneous layers. We also show that a polynomial in three variables with vanishing bordered Hessian depends on two linear forms, and that the five-variable potentials of Meng and Yang have neither an affine direction nor a constant pivot.
\end{abstract}

\maketitle

\section{Introduction}\label{sec:intro}

\subsection{Polynomial solutions and the Legendre transform}
Let $\kk$ be a field of characteristic zero. We study polynomials $\phi\in\kk[x_1,\dots,x_n]$ with
\begin{equation}\label{eq:MA}
\det\Hess\phi=c,\qquad c\in\kk^\times ,
\end{equation}
where $\Hess\phi=(\partial_i\partial_j\phi)$. Over $\R$ this is the Monge--Amp\`ere equation of an improper affine hypersphere. By the theorem of J\"orgens, Calabi and Pogorelov, a convex solution defined on all of $\R^n$ is a quadratic polynomial \cite{Jor54,Cal58,Pog72}.

Without convexity the equation has many polynomial solutions, and the questions about them are algebraic. The classical symmetry of the Monge--Amp\`ere equation is the Legendre transform. If $\grad\phi$ is a diffeomorphism, then $\phi^*(y)=\vv{y,x}-\phi(x)$ at $y=\grad\phi(x)$ satisfies $\grad\phi^*=(\grad\phi)^{-1}$ and $\det\Hess\phi^*(\grad\phi(x))\cdot\det\Hess\phi(x)=1$. For a polynomial solution of \eqref{eq:MA} the gradient $\grad\phi\colon\kk^n\to\kk^n$ is a Keller map, that is, a polynomial map with constant nonzero Jacobian determinant, and $\phi^*$ is a polynomial if and only if $\grad\phi$ is a polynomial automorphism (Proposition~\ref{prop:leg-global}). We call a polynomial solution of \eqref{eq:MA} a \emph{Keller potential}. The Hessian conjecture $\HC_n$, formulated by Meng \cite{Men06} through the formal Legendre transform and independently, without the Legendre transform, by de~Bondt and van~den~Essen \cite{dBvdE05}, see \cite[Remark~4.5]{AHL}, asserts that the gradient of every Keller potential in $n$ variables is a polynomial automorphism. Equivalently, the polynomial solutions of \eqref{eq:MA} in $n$ variables form a class which is stable under the Legendre transform \cite{AHL}. Both the theorem of J\"orgens, Calabi and Pogorelov and the Hessian conjecture concern entire solutions. The former classifies the entire convex solutions, whereas the latter asks whether the polynomial solutions are preserved by the Legendre duality. Geometrically, the graph $\Lambda_\phi=\{(x,\grad\phi(x))\}$ of a Keller potential is an exact Lagrangian subvariety of $T^*\A^n$ whose projection to the momentum space is \'etale, and $\phi^*$ is a polynomial if and only if $\Lambda_\phi$ is also a graph over the momentum space (Remark~\ref{rem:lagrangian}).

In two variables $x_1x_2+g(x_1)$ is a solution for every $g\in\kk[x_1]$. If $F=(F_1,F_2)$ is a Keller map of the plane, then
\begin{equation}\label{eq:doubling}
\phi_F=x_3F_1(x_1,x_2)+x_4F_2(x_1,x_2)
\end{equation}
is a solution in four variables, with $\det\Hess\phi_F=(\det JF)^2$ (Lemma~\ref{lem:doubling-det}). Following \cite{MY26} we call $\phi_F$ the \emph{doubling} of $F$. It is a generating function of the cotangent lift of $F$ \cite{Arn89} (Remark~\ref{rem:lagrangian}). More generally, a polynomial of the form $E(x_1,x_2)+x_3C(x_1,x_2)+x_4A(x_1,x_2)$ in some linear coordinates will be called a doubling of $(C,A)$. The gradient of $\phi_F$ is injective if and only if $F$ is. Thus $\HC_n$ is the symmetric case of the Jacobian conjecture $\JC_n$ of Keller \cite{Kel39}, one has $\JC_n\Rightarrow\HC_n$, and the doubling gives $\HC_{2n}\Rightarrow\JC_n$ \cite{Men06}, \cite{dBvdE05}. Dillen proved $\HC_2$ \cite{Dil91} and de Bondt proved $\HC_3$ \cite{dB15}. Starting from a non-injective Keller map of $\kk^3$ announced by Alp\"oge,\footnote{The announcement was made in July 2026, see \cite[Section~1]{MY26}. Hence $\JC_n$ fails for $n\ge3$, and the plane is the only open case of the Jacobian conjecture.} Meng and Yang obtained a Keller potential in five variables with non-injective gradient, hence $\HC_n$ fails for $n\ge5$ \cite{MY26}. This leaves $\JC_2$ and $\HC_4$, with $\HC_4\Rightarrow\JC_2$. In \cite[Remark~4.5]{AHL} we observed that in a fixed dimension the argument of Meng gives $\HC_{2n}\Rightarrow\JC_n\Rightarrow\HC_n$, and in particular does not give the converse implication $\JC_2\Rightarrow\HC_4$. This converse is a consequence of the main result of the present paper (Corollary~\ref{cor:equiv}).

\subsection{Affine directions}
For $v\in\kk^n\setminus\{0\}$ let $\deg_v\phi$ be the degree of $t\mapsto\phi(x+tv)$ over $\kk[x]$, so that $\deg_v\phi\le m$ if and only if $D_v^{m+1}\phi=0$. We call $v$ an \emph{affine direction} of $\phi$ if $D_v^2\phi=0$, that is, if $\phi$ is affine on every line of direction $v$, and a \emph{constant pivot} if $D_v^2\phi\in\kk^\times$. Equivalently, $v\T\Hess\phi(x)\,v=0$ for every $x$, so $v$ is isotropic for all the quadratic forms $\Hess\phi(x)$ at once. Over an algebraically closed field every quadratic form in four variables has nonzero isotropic vectors, and an affine direction is a single vector which is isotropic for all of them. The doubling \eqref{eq:doubling} has the affine directions $e_3$ and $e_4$. Our main result says that every Keller potential in four variables shares this property of the doubling.

\begin{theorem}[affine directions]\label{thm:main-aff}
Let $\kk$ be algebraically closed. Every Keller potential $\phi\in\kk[x_1,\dots,x_4]$ has an affine direction.
\end{theorem}

The forms $q_x(v)=v\T\Hess\phi(x)\,v$, $x\in\kk^4$, span a linear system of quadrics in $\PP^3$, generated by the coefficients $q_\alpha$ of $v\T\Hess\phi(x)\,v=\sum_\alpha q_\alpha(v)\,x^\alpha$. Theorem~\ref{thm:main-aff} says that its base locus, the zero set of the ideal $(q_\alpha)_\alpha$ of $\kk[v_1,\dots,v_4]$, is nonempty. For a doubling the base locus contains a line, the projectivized plane of affine directions (Remark~\ref{rem:commonplane}). For the potential of Example~\ref{ex:case-ii} it is a single point, so in general no more than a point can be expected. If $\grad\phi$ is not a polynomial automorphism, then $\phi$ is a doubling over $\overline\kk$ (Corollary~\ref{cor:equiv}), and the base locus contains a line. Conversely, the base locus contains a line if and only if $\phi$ is a doubling, since for a plane $U$ of affine directions polarization gives $D_uD_{u'}\phi=0$ for all $u,u'\in U$.

An affine direction reduces the problem to the plane. If $e_4$ is an affine direction, then $\phi=x_4A+B$ with $A,B\in\kk[x_1,x_2,x_3]$. The coefficient of $x_4^2$ in $\det\Hess\phi$ forces $A$ to depend on two linear forms only, and then either $\phi$ is a doubling, or a partial Legendre transform reduces the injectivity of $\grad\phi$ to $\HC_2$ (Theorem~\ref{thm:affine}). We say that a Keller potential $\phi$ in four variables satisfies the \emph{dichotomy}~\textup{(D)} if $\grad\phi$ is a polynomial automorphism or $\phi$ is a doubling of a planar Keller map. By Theorem~\ref{thm:main-aff} every Keller potential in four variables satisfies \textup{(D)} over an algebraic closure of $\kk$. Consequently $\HC_4$ is equivalent to $\JC_2$ over every algebraically closed field of characteristic zero, and a counterexample to $\HC_4$ is the doubling of a counterexample to $\JC_2$ (Corollary~\ref{cor:equiv}). By the theorem of Moh that $\JC_2$ holds in degree at most $100$ \cite{Moh83}, a Keller potential in four variables whose gradient is not an automorphism has degree at least $102$. The bound of Guccione, Guccione, Horruitiner and Valqui on counterexamples to $\JC_2$, from a preprint \cite[Theorem~2.1]{GGHV22}, raises this to $109$, and below degree $126$ the planar Keller map of the doubling then has component degrees $72$ and $108$.

Over $\R$, complex conjugation of the plane of affine directions of a doubling shows that a real Keller potential whose gradient is not a polynomial automorphism is a doubling in real linear coordinates (Corollary~\ref{cor:realdoubling}). Its Hessian then has signature $(2,2)$, and every real Keller potential whose Hessian has signature $(3,1)$ or $(1,3)$, equivalently negative Hessian determinant, has a polynomial Legendre transform (Corollaries~\ref{cor:lorentz} and~\ref{cor:sign}). A real Keller potential in four variables has a real affine direction unless it is quadratic with definite Hessian (Corollary~\ref{cor:realaff}). Together with the theorem of J\"orgens, Calabi and Pogorelov this leaves the signature $(2,2)$ of the doubling as the only possible signature of a real counterexample to $\HC_4$ (Corollary~\ref{cor:sign}). In contrast with the definite case, in signature $(2,2)$ neither positivity of the Hessian determinant nor analyticity of an entire solution ensures injectivity of the gradient (Remark~\ref{rem:22}). The real case is treated in Section~\ref{ssec:real}.

Recently Q.~Tang stated in a Hexagon preprint \cite{Tan26}, with reference to \cite[Theorem~4.17 and Remark~4.5]{AHL}, that every Keller potential in four variables satisfies \textup{(D)}, which gives the converse implication $\JC_2\Rightarrow\HC_4$. The proof of Theorem~\ref{thm:main-aff} in Sections~\ref{sec:translation} to~\ref{sec:tree} follows his strategy. The Newton polytope of $\phi$ is explored along paths of weights, and at the walls of these paths the leading weighted forms are described by translation directions of their gradient images (Section~\ref{sec:translation}). Every branch of the weighted analysis can be brought to an end in an affine direction or in a contradiction, without inverting the gradient, and this yields the stronger statement of Theorem~\ref{thm:main-aff}. We give complete proofs of the statements taken from that preprint, and the correspondence, together with the differences, is listed in Remark~\ref{rem:attribution}. The key intermediate statement is Theorem~\ref{thm:cong}, due to Tang \cite[Theorem~2.1]{Tan26}. If a rational map of $\kk^3$ has nilpotent Jacobian matrix $N$ and $HN$ is symmetric for the Hessian matrix $H$ of a polynomial with $\det H\neq0$, then the components of the map satisfy a nontrivial affine linear relation. Tang's proof rests on two classifications, that of congruences of lines with one-dimensional focal locus by De~Poi \cite{DP04} and that of polynomial Hessians of rank two and three by de~Bondt \cite{dB16}. Our account uses neither of them. They are replaced by the construction of Gordan and Noether, by quasi-translations, and by a study of the field of invariants of the kernel field of a Jacobian matrix (Lemmas~\ref{lem:cong-dir}, \ref{lem:E6} and~\ref{lem:E7}, Appendix~\ref{app:dB62}). Some paths reach a weight vector at which the leading form has constant Hessian determinant. They end by a direct argument (Theorem~\ref{thm:touch}). Invertibility enters only through Theorem~\ref{thm:affine}.

\subsection{A proof by homogeneous layers}
For potentials with a direction of degree at most two we give a second proof, which uses no paths of weights and stays close to the classical treatment of the Monge--Amp\`ere equation. It has three tools. The first is the partial Legendre transform (Lemma~\ref{lem:legendre}). It turns a potential which is quadratic in some of the variables into a family of potentials in the remaining variables, to which $\HC_2$ and $\HC_3$ apply. The second consists of identities between the homogeneous layers of $\det\Hess\phi$, and between its coefficients along a direction, which are governed by the bordered Hessian
\[
\BH(A)=\grad A^{\mathsf T}\adj(\Hess A)\grad A ,
\]
where $\adj$ denotes the adjugate matrix. If $\phi=x_4A(x_1,x_2,x_3)+B(x_1,x_2,x_3)$, the coefficient of $x_4^2$ in $\det\Hess\phi$ is $-\BH(A)$. The third is the construction of Gordan and Noether \cite{GN1876}, which we apply to rational and algebraic functions (Lemma~\ref{lem:graph}). The first two tools exploit the block structure of $\Hess\phi$. An invertible block in the variables in which $\phi$ is quadratic is removed by a Schur complement, which is the partial Legendre transform, while the coefficient of the highest power of a variable in $\det\Hess\phi$ is governed by the bordered Hessian of the leading coefficient. The first result concerns potentials with a direction of degree at most two.

\begin{theorem}[directions of degree at most two]\label{thm:main-dir}
Let $\phi\in\kk[x_1,\dots,x_4]$ be a Keller potential with $\deg_v\phi\le2$ for some $v\neq0$. Then $\grad\phi$ is a polynomial automorphism, or over $\overline\kk$ the potential $\phi$ is a doubling of a non-injective planar Keller map of degree at most $\deg\phi-1$.
\end{theorem}

Besides the three tools the proof uses the theorems $\HC_2$ and $\HC_3$ and the classification of polynomials in three variables with vanishing Hessian determinant \cite{dBvdE04}, \cite{dB18}, which we prove again in Appendix~\ref{app:dB62}. In low degree a direction of degree at most two always exists.

\begin{theorem}[low degree]\label{thm:main-low}
Let $\kk$ be algebraically closed and let $\phi\in\kk[x_1,\dots,x_4]$ be a Keller potential.
\begin{enumerate}[label=\textup{(\roman*)}]
\item If $\deg\phi\le4$, then $\phi$ has an affine direction or a constant pivot.
\item If $\deg\phi\le5$, then $\deg_v\phi\le2$ for some $v\neq0$.
\end{enumerate}
\end{theorem}

Part (i) has the same content as the structure theorem of Ni for quartic potentials \cite{Ni26}, and is proved here by a weight identity for the Hessian determinant. Part (ii) is proved by the homogeneous layers of $\det\Hess\phi$. With the theorem of Moh that $\JC_2$ holds in degree at most $100$ \cite{Moh83}, Theorems~\ref{thm:main-dir} and~\ref{thm:main-low} give $\HC_4$ in degree at most five. In this degree $\HC_4$ was obtained by Niedbala Giraudin with certificates of computer algebra in an unrefereed preprint \cite{G2}, and for real potentials of Lorentzian signature by Liu \cite{Liu26}. Our proof does not use these computations. In degree six the homogeneous layers, together with the homogeneous parts for the single weight vector $(1,1,1,2)$ which the layers determine, give a direction of degree at most two (Theorems~\ref{thm:deg6}, \ref{thm:deg6r2}, \ref{thm:deg6r1} and Corollary~\ref{cor:deg6}). So the dichotomy \textup{(D)} for potentials of degree at most six is proved in Sections~\ref{sec:prelim} to~\ref{sec:low}, independently of Sections~\ref{sec:translation} to~\ref{sec:tree}. The two approaches are complementary. The weighted proof applies in every degree, while the homogeneous layers give a direct treatment of all potentials with a direction of degree at most two and of all potentials of degree at most six, without paths of weights.

For a form $A$ of degree $e\ge2$ the Euler identity gives $\BH(A)=\frac e{e-1}A\det\Hess A$, and the theorem of Gordan and Noether describes the forms in at most four variables with vanishing Hessian determinant. Part of this extends to arbitrary polynomials.

\begin{theorem}[vanishing bordered Hessian]\label{thm:main-BH}
Let $\kk$ be a field of characteristic zero.
\begin{enumerate}[label=\textup{(\roman*)}]
\item If $A\in\kk[x_1,\dots,x_n]$ and $\BH(A)=0$, then $\det\Hess A=0$.
\item If $A\in\kk[x_1,x_2,x_3]$ and $\BH(A)=0$, then $A\in\kk[\ell_1,\ell_2]$ for two linear forms $\ell_1,\ell_2$.
\item If $\phi=\sum_{j=0}^mt^jc_j(x_1,x_2,x_3)$ is a Keller potential with $m\ge1$ and $c_m\neq0$, then $c_m$ depends on two linear forms only.
\end{enumerate}
\end{theorem}

Part (i) is proved by a formal Legendre transform. Parts (ii) and (iii) use the construction of Gordan and Noether, applied to rational and algebraic functions. Part (ii) fails in four variables, and the polynomials in four variables with vanishing bordered Hessian are described in Proposition~\ref{prop:BH4}. Besides the polynomials in three linear forms, they are, after a linear change of coordinates, the polynomials $g(x_1,L)$ with $L=\ell_2(x_1)x_2+\ell_3(x_1)x_3+\ell_4(x_1)x_4$, a linear form in $x_2,x_3,x_4$ whose coefficients are polynomials in $x_1$.

\begin{remark*}
Over $\R$ the quotient $\BH(A)/|\grad A|^4$ is the Gaussian curvature of the level surfaces of a function $A$ of three variables \cite{Gol05}, \cite[Proposition~7.2.3]{Rus16}. Part (ii) thus says that a polynomial in three variables whose level surfaces are all flat is a function of two linear forms. Its level surfaces are then cylinders over plane curves. For a single irreducible projective surface $V(f)\subset\PP^3$, vanishing curvature means that $f$ divides its Hessian determinant, and the surface is a cone or developable \cite[Theorem~7.2.5]{Rus16}.
\end{remark*}

\subsection{Pascal finite maps}
For a polynomial map $F$ let $\sigma_F(g)=g\circ F$ and $\Delta_F=\sigma_F-\mathrm{Id}$. The map $F$ is Pascal finite if $\Delta_F^m$ annihilates the coordinate functions for some $m$ \cite{ABCH19}. In characteristic zero the Pascal finite automorphisms are the exponentials of locally nilpotent derivations, and their inverses are computed by a finite algorithm of substitutions \cite{ABCH19}. The construction of Gordan and Noether attaches to a polynomial with singular Hessian a quasi-translation, that is, a polynomial automorphism $X+H$ with inverse $X-H$ \cite{dB18}. The proofs in Section~\ref{sec:translation}, of Theorem~\ref{thm:main-BH} and in Appendix~\ref{app:dB62} use only the identity $H(X+sH)=H$ of quasi-translations. Behind this identity is a single mechanism. For a polynomial $h$, a polynomial relation $R(\grad h)=0$ gives the quasi-translation $X+H$ with $H=(\grad R)(\grad h)$, a Pascal finite map of depth at most two, defined by the condition $\Delta_{X+H}^2(X)=0$ over every field. In characteristic zero $X+H$ is the exponential of the locally nilpotent derivation $\sum_iH_i\partial_i$, whose ring of constants contains the components of $\grad h$ (Section~\ref{ssec:qt}). The quasi-translations are exactly the Pascal finite maps of Pascal depth at most two (Proposition~\ref{prop:qt-pf}). For such a map $F=X+H$ the inverse $X-H$ is obtained by a single substitution. So the maps produced by the construction of Gordan and Noether form the first level of the class of Pascal finite maps. By Theorem~\ref{thm:apex}, due to Tang \cite[Theorem~2.3]{Tan26}, the closure of the gradient image of a polynomial in four variables with Hessian of rank three has a translation direction, although the kernel of the Hessian need not contain a nonzero constant vector (Example~\ref{ex:apex}). The quasi-translation of Proposition~\ref{prop:qt-hess} links the varying kernel direction in the source with the constant translation direction in the target. As observed by Tang \cite[Corollary~2.4]{Tan26}, this gives an affirmative answer to Problem~1 of de~Bondt \cite[Section~6]{dB18} on polynomials in four variables with vanishing Hessian determinant (Corollary~\ref{cor:dBproblem}).

\subsection{Five variables}
In five variables affine directions need not exist. Let $\Rc_n$ denote the set of Keller potentials in $n$ variables with neither an affine direction nor a constant pivot. By Theorem~\ref{thm:main-aff}, $\Rc_4=\varnothing$.

\begin{theorem}[Meng--Yang potentials]\label{thm:main-MY}
The Keller potentials $\psi_{\lambda,\mu}$, $\lambda\neq0$, of \cite{MY26}, among them the counterexample to $\HC_5$, lie in $\Rc_5$.
\end{theorem}

The potentials $\psi_{\lambda,\mu}$ have degree at most two exactly along the directions of a three-dimensional space (Proposition~\ref{prop:psi-L0}). Hence in five variables a direction of degree at most two implies neither that the gradient is an automorphism nor that an affine direction exists, and Theorems~\ref{thm:main-aff} and~\ref{thm:main-dir} have no analog there. The proof of Theorem~\ref{thm:main-MY} reads off five coefficients of $D_\xi^2\psi_{\lambda,\mu}$, regarded as a polynomial in the coordinates. No computer calculation is used in the proofs of this paper.

The results of the paper are ordered by the directional degree. A vector $v\neq0$ is an affine direction if $\deg_v\phi\le1$, and a constant pivot if $\deg_v\phi=2$ and $D_v^2\phi\in\kk^\times$. Theorem~\ref{thm:main-aff} gives a direction of degree at most one for every Keller potential in four variables. Theorem~\ref{thm:main-dir} shows that a direction of degree at most two already gives the dichotomy~\textup{(D)}, and Theorem~\ref{thm:main-low} and Corollary~\ref{cor:deg6} give such a direction in degree at most six. In five variables the potentials $\psi_{\lambda,\mu}$ have a three-dimensional space of directions of degree exactly two, and none of them is an affine direction or a constant pivot.

\subsection{Organization}
The paper has three parts. Sections~\ref{sec:prelim} to~\ref{sec:low} prove Theorems~\ref{thm:main-dir}, \ref{thm:main-low} and~\ref{thm:main-BH} by partial Legendre transforms and homogeneous layers. Section~\ref{sec:prelim} collects the bordered Hessian, the Legendre transform and its partial version, the construction of Gordan and Noether for rational and algebraic functions, and its relation to Pascal finite maps. Section~\ref{sec:BH} proves Theorem~\ref{thm:main-BH}. Section~\ref{sec:affine} treats affine directions and constant pivots, and Sections~\ref{sec:block} and~\ref{sec:quadratic} treat potentials which are quadratic in two variables, or in one variable with nonconstant coefficient. Section~\ref{sec:thmA} proves Theorem~\ref{thm:main-dir}, and Section~\ref{sec:low} proves Theorem~\ref{thm:main-low} and the results in degree six. These results also follow from Theorem~\ref{thm:main-aff} and Corollary~\ref{cor:equiv}. Sections~\ref{sec:block} to~\ref{sec:low} prove them without Sections~\ref{sec:translation} to~\ref{sec:tree}, by the homogeneous layers of $\det\Hess\phi$. Sections~\ref{sec:translation} to~\ref{sec:tree} prove Theorem~\ref{thm:main-aff}. Section~\ref{sec:translation} attaches quasi-translations to Hessians of rank three and proves that their components are linearly dependent, Section~\ref{sec:walls} develops weighted leading forms and the exit lemmas, and Section~\ref{sec:tree} treats the highest homogeneous part. These three sections use only Sections~\ref{sec:prelim} to~\ref{sec:affine} of the first part. Section~\ref{sec:tang} draws the consequences of Theorem~\ref{thm:main-aff}, Section~\ref{sec:MY} concerns five variables, and Section~\ref{sec:open} contains concluding remarks. Appendix~\ref{app:dB62} gives a proof of the classification of polynomials in three variables with vanishing Hessian determinant by the methods of the paper.

From Section~\ref{sec:prelim} to Section~\ref{sec:MY} the field $\kk$ is algebraically closed of characteristic zero. Statements on injectivity descend to an arbitrary field of characteristic zero by Lemma~\ref{lem:descent}.

\section{Preliminaries}\label{sec:prelim}

\subsection{The bordered Hessian}
For $A\in\kk[x_1,\dots,x_n]$ put
\[
\BH(A):=\grad A^{\mathsf T}\adj(\Hess A)\grad A=-\det\begin{pmatrix}\Hess A&\grad A\\ \grad A^{\mathsf T}&0\end{pmatrix}.
\]
For $n=2$ one has $\BH(A)=A_1^2A_{22}-2A_1A_2A_{12}+A_2^2A_{11}$, where $A_i=\partial_iA$ and $A_{ij}=\partial_i\partial_jA$.

\begin{proposition}\label{prop:BHdet}%
Let $A\in\kk[x_1,\dots,x_n]$ with $\BH(A)=0$. Then $\det\Hess A=0$.
\end{proposition}
The idea is that the formal Legendre transform $u$ of $A$ at a point satisfies $Y\T\Hess u\,Y=\BH(A)/\det\Hess A$. So $\BH(A)=0$ makes $u$ radially affine, that is $u(\lambda Y)=a(Y)+\lambda b(Y)$ with $a$ homogeneous of degree $0$ and $b$ homogeneous of degree $1$, and this is incompatible with $A$ being a polynomial.
\begin{proof}
Suppose $\Delta=\det\Hess A\neq0$. Let $d=\deg A$, choose $x_0\in\kk^n$ with $\Delta(x_0)\neq0$, put $y_0=\grad A(x_0)$, and let $\eta=(\eta_1,\dots,\eta_n)$ be indeterminates. By the formal inverse function theorem there is $G\in\kk[[\eta]]^n$ with $G(0)=x_0$ and $\grad A(G)=Y$, where $Y=y_0+\eta$. Differentiation with respect to $\eta$ gives $\Hess A(G)\cdot J_\eta G=I$, where $J_\eta G=(\partial G_i/\partial\eta_k)$ is the Jacobian matrix of $G$, and $\Delta(G)$ is a unit of $\kk[[\eta]]$. Put $u=\vv{G,Y}-A(G)\in\kk[[\eta]]$. Then
\[
\frac{\partial u}{\partial\eta_k}=G_k+\sum_i\frac{\partial G_i}{\partial\eta_k}\bigl(Y_i-(\partial_iA)(G)\bigr)=G_k,
\]
hence $\Hess_\eta u=J_\eta G=\Hess A(G)^{-1}$ and
\[
Y^{\mathsf T}\Hess_\eta u\;Y=\grad A(G)^{\mathsf T}\Hess A(G)^{-1}\grad A(G)=\frac{\BH(A)(G)}{\Delta(G)}=0 .
\]
Let $E=\sum_jY_j\,\partial/\partial\eta_j$, a derivation of $\kk[[\eta]]$ with $EY_k=Y_k$ and $\partial_{\eta_k}E-E\partial_{\eta_k}=\partial_{\eta_k}$. Then $E^2u=Eu+Y^{\mathsf T}\Hess_\eta u\,Y=Eu$. Hence $f=Eu-u$ satisfies $Ef=0$, and
\[
EG_k=E\,\partial_{\eta_k}u=\partial_{\eta_k}f,\qquad E^2G_k=\partial_{\eta_k}Ef-\partial_{\eta_k}f=-EG_k .
\]
For $k\ge0$ let $T_k\in\kk[[\eta]]^n$ have the components $\sum_{j_1,\dots,j_k}(\partial_{j_1}\cdots\partial_{j_k}\partial_iA)(G)\,EG_{j_1}\cdots EG_{j_k}$. Then $T_0=\grad A(G)=Y$, and the chain rule together with $E^2G=-EG$ gives $ET_k=T_{k+1}-kT_k$. Since $EY=Y$, induction on $k$ gives $T_k=k!\,Y$ for every $k$. On the other hand $T_d=0$, since all partial derivatives of order $d+1$ of $A$ vanish. As $Y\neq0$ and $\operatorname{char}\kk=0$, this is a contradiction.
\end{proof}

\begin{remark*}
For a form $A$ of degree $e\ge2$ the proposition follows directly from the identity $\BH(A)=\frac e{e-1}A\det\Hess A$. The proof above is needed because $A$ is not assumed to be homogeneous.
\end{remark*}

\subsection{Results from the literature}
We use the following results.
\begin{enumerate}[label=\textup{(E\arabic*)}]
\item\label{E:dB} $\HC_n$ holds for $n\le3$, by \cite{Dil91} for $n=2$ and by \cite{dB15} for $n=3$.
\item\label{E:dB62} Let $F$ be a field of characteristic zero and $A\in F[x_1,x_2,x_3]$ with $\det\Hess A=0$. Then $\partial_cA\in F$ for some $c\in F^3\setminus\{0\}$, or there is $T\in\mathrm{GL}_3(F)$ such that $A(Tx)=a_1(x_1)+a_2(x_1)x_2+a_3(x_1)x_3$ with $a_1,a_2,a_3\in F[x_1]$ \cite[Theorem~3.3]{dBvdE04}, \cite[Theorem~6.2]{dB18}. A proof by the methods of this paper is given in Appendix~\ref{app:dB62}.%
\item\label{E:GN} A form in at most four variables with vanishing Hessian determinant has a nonzero vertex direction $\kappa$, that is $\partial_\kappa f=0$ \cite{GN1876}, \cite{Los04}, \cite[Theorem~3.1]{dB18}, over an arbitrary field of characteristic zero \cite[Theorem~A]{BFP23}, see also \cite[Chapter~7]{Rus16}. The statement fails in five variables \cite{Per1900}, and the irreducible forms in five variables with vanishing Hessian determinant which are not cones are classified in \cite[Theorem~7.4.14]{Rus16}.
\item\label{E:Zaks} If $h_1,\dots,h_r\in\kk[u_1,\dots,u_n]$ generate a $\kk$-algebra of transcendence degree $1$, then $h_1,\dots,h_r\in\kk[g]$ for some $g\in\kk[u]$ \cite{Zak71}, \cite[Chapter~1]{Now94}.
\item\label{E:dBvdE} Let $F$ be a field of characteristic zero and $f\in F[u_1,u_2]$ with $\det\Hess f=0$. Then there are $\lambda\in F^2\setminus\{0\}$ and $\gamma\in F$ with $\lambda_1\partial_1f+\lambda_2\partial_2f=\gamma$ \cite[Theorem~3.1]{dBvdE04}, \cite[Theorem~3.1]{dB18}.%
\item\label{E:Dil} Let $f\in\C[u_1,u_2]$ with $\det\Hess f\in\C^\times$ and $\deg f\ge3$. Then after a linear change of coordinates $\partial_2^2f=0$ \cite{Dil91}, \cite[Section~1]{dB15}.%
\end{enumerate}
Results \ref{E:dB62} and \ref{E:dBvdE} are proved in \cite{dBvdE04} over an arbitrary field of characteristic zero. In \ref{E:dB62} the condition $\partial_cA\in F$ cannot be replaced by $\partial_cA=0$, since $A=x_1^2x_2^2+x_3$ has $\det\Hess A=0$ and $\partial_3A=1$, but $\partial_cA\neq0$ for every $c\neq0$, and $\Hess A$ vanishes on no two-dimensional subspace. In \ref{E:dBvdE} the constant $\gamma$ cannot be omitted, since $f=(u_1+u_2)^3+u_1$ satisfies $(\partial_1-\partial_2)f=1$, and the statement fails in three variables, since $x_1^2x_2+x_1^3x_3$ has vanishing Hessian determinant and no such $\lambda$.

\begin{lemma}\label{lem:descent}%
Let $\kk_0$ be a field of characteristic zero and let $F\colon\kk_0^n\to\kk_0^n$ be a Keller map whose extension to $\overline{\kk_0}^{\,n}$ is injective. Then $F$ is a polynomial automorphism of $\kk_0^n$.
\end{lemma}
\begin{proof}
By \cite{BBR62} and \cite[Theorem~4.2.1]{vdE00} the map $F$ has an inverse $G\in\overline{\kk_0}[x]^n$. For $\sigma\in\operatorname{Gal}(\overline{\kk_0}/\kk_0)$ the map $\sigma(G)$ is also an inverse of $F$, so $\sigma(G)=G$ and $G\in\kk_0[x]^n$.%
\end{proof}

\begin{lemma}\label{lem:lindescent}%
Let $\kk_0\subseteq L$ be fields and let $V\subseteq L^n$ be the space of solutions of a system of homogeneous linear equations with coefficients in $\kk_0$. Then $V$ is spanned over $L$ by $V\cap\kk_0^n$. In particular, for $A\in\kk_0[x_1,\dots,x_n]$ this applies to each of the spaces of vectors $c\in L^n$ with $D_cA=0$, with $D_cA\in L$, and with $\Hess A\cdot c=0$.
\end{lemma}
\begin{proof}
Gaussian elimination over $\kk_0$ gives a basis of the solutions in $\kk_0^n$, and the rank of a matrix with entries in $\kk_0$ does not change under the extension to $L$. Writing $D_cA=\sum_ic_i\partial_iA$ and comparing the coefficients of all monomials, respectively of the nonconstant monomials, gives linear equations in $c$ with coefficients in $\kk_0$, and the same holds for the entries of $\Hess A\cdot c$.
\end{proof}

\begin{lemma}[linear components of fibers]\label{lem:linfib}%
Let $\kk$ be algebraically closed, let $f\colon\kk^n\to\kk^m$ be a polynomial map whose general fibers have dimension $k$, and let $D\subseteq\kk^n$ be dense. For $x\in D$ let $U(x)\subseteq\kk^n$ be a linear subspace of dimension $k$ with $x+U(x)\subseteq f^{-1}(f(x))$. Then there is a dense open subset $\Omega\subseteq\kk^n$ such that for $x,x'\in D\cap\Omega$ the affine subspaces $x+U(x)$ and $x'+U(x')$ are irreducible components of fibers of $f$, and they are equal or disjoint. If they are disjoint, then $U(x)+U(x')\neq\kk^n$. In particular, for $k=n-1$ one has $U(x)=U(x')$.
\end{lemma}
\begin{proof}
In characteristic zero there is a dense open subset of the closure of the image of $f$ over which the fibers are smooth of pure dimension $k$. Let $\Omega$ be its inverse image. For $x\in D\cap\Omega$ the closed irreducible set $x+U(x)$ of dimension $k$ is an irreducible component of the fiber through $x$. Distinct irreducible components of a smooth variety are disjoint, and components of distinct fibers are disjoint. If $U(x)+U(x')=\kk^n$, then $x'-x=u-u'$ with $u\in U(x)$, $u'\in U(x')$, and $x+u=x'+u'$ is a common point. For $k=n-1$ the condition $U(x)+U(x')\neq\kk^n$ means $U(x)=U(x')$.
\end{proof}

\begin{lemma}\label{lem:luroth}%
Let $\kk$ be algebraically closed and let $L$ be a field with $\kk\subseteq L\subseteq\kk(x_1,\dots,x_n)$ and $\operatorname{trdeg}_\kk L=1$. Then $L=\kk(s)$ for some $s\in L$.
\end{lemma}
\begin{proof}
The field $L$ is finitely generated over $\kk$, as a subfield of a finitely generated extension, say $L=\kk(g_1,\dots,g_r)$ with $g_1\notin\kk$. On a general line $\tau\mapsto a+\tau b$ of $\kk^n$ all $g_i$ are defined and $g_1$ is nonconstant. Restriction to the line gives a $\kk$-algebra homomorphism $\kk[g_1,\dots,g_r]\to\kk(\tau)$. Its kernel is a prime ideal of a domain of dimension one, and it is not maximal, since the image is not contained in $\kk$. Hence the kernel is zero, and restriction extends to an embedding of $L$ into $\kk(\tau)$. By the theorem of L\"uroth its image is $\kk(s_0)$ for some $s_0$, and $L=\kk(s)$ for the preimage $s$ of $s_0$.
\end{proof}

\subsection{The Legendre transform}
\begin{proposition}\label{prop:leg-global}%
Let $\phi\in\kk[x_1,\dots,x_n]$ be a Keller potential with $\det\Hess\phi=c$. Then $\grad\phi$ is a polynomial automorphism if and only if there is $\psi\in\kk[y_1,\dots,y_n]$ with $\grad\psi(\grad\phi(x))=x$. In this case $G=(\grad\phi)^{-1}$ and
\[
\phi^*=\vv{y,G(y)}-\phi(G(y))\in\kk[y]
\]
satisfy $\grad\phi^*=G$, $\det\Hess\phi^*=c^{-1}$ and $\phi^{**}=\phi$, and every such $\psi$ differs from $\phi^*$ by a constant.
\end{proposition}
\begin{proof}
Let $G$ be the inverse of $\grad\phi$. Then
\[
\frac{\partial\phi^*}{\partial y_k}=G_k+\sum_i\frac{\partial G_i}{\partial y_k}\bigl(y_i-(\partial_i\phi)(G)\bigr)=G_k ,
\]
so $\Hess\phi^*=JG=\Hess\phi(G)^{-1}$ and $\det\Hess\phi^*=c^{-1}$. Since $\grad\phi^*=G$ has the inverse $\grad\phi$, one gets $\phi^{**}=\vv{x,\grad\phi}-\phi^*(\grad\phi)=\phi$. Conversely, if $\grad\psi\circ\grad\phi=x$, then $\grad\phi$ is injective, hence an automorphism by Lemma~\ref{lem:descent}, and $\grad\psi=G=\grad\phi^*$.
\end{proof}

\begin{remark*}
Thus $\HC_n$ says that the polynomial solutions of \eqref{eq:MA} in $n$ variables form a set stable under the Legendre transform. In \cite{AHL} the transform $\phi^*$ is computed from the Legendre potential $\vv{x,\grad\phi}-\phi$. Over $\R$, for a convex solution on $\R^n$, the theorem of J\"orgens, Calabi and Pogorelov says that $\phi$ and $\phi^*$ are quadratic.
\end{remark*}

\subsection{The partial Legendre transform}
\begin{lemma}\label{lem:legendre}
Let $x=(x_1,\dots,x_n)$, $y=(y_1,\dots,y_m)$ and $\phi=B(x)+\vv{y,\Lambda(x)}+\frac12y^{\mathsf T}S(x)y$ with $B\in\kk[x]$, $\Lambda\in\kk[x]^m$, $S=S^{\mathsf T}\in\kk[x]^{m\times m}$ and $\delta=\det S\neq0$. Put $W=S^{-1}$ and, for indeterminates $p=(p_1,\dots,p_m)$,
\[
\psi_p=-B+\tfrac12(p-\Lambda)^{\mathsf T}W(p-\Lambda)\in\kk(x)[p].
\]
\begin{enumerate}[label=\textup{(\alph*)}]
\item $\grad_x\psi_p(x)=-\grad_x\phi(x,y)$ at $y=W(x)(p-\Lambda(x))$.
\item $\det\Hess\phi(x,y)=(-1)^n\delta(x)\det\Hess_x\psi_p(x)$ at $p=\Lambda(x)+S(x)y$.
\item $\det\Hess\phi=c\in\kk$ if and only if $\det\Hess_x\psi_p=(-1)^nc/\delta$ in $\kk(x)[p]$.
\end{enumerate}
\end{lemma}
\begin{proof}
Write $\Phi=\Hess_x\phi$ and $N=\partial_x(\grad_y\phi)$, an $m\times n$ matrix, so that $\Hess\phi=\begin{psmallmatrix}\Phi&N^{\mathsf T}\\ N&S\end{psmallmatrix}$ and $\det\Hess\phi=\delta\det(\Phi-N^{\mathsf T}WN)$. The vector $y(x,p)=W(p-\Lambda)$ satisfies $\grad_y\phi(x,y(x,p))=p$. By the chain rule $\grad_x\psi_p=(\partial_xy)^{\mathsf T}p-\grad_x\phi-(\partial_xy)^{\mathsf T}\grad_y\phi=-\grad_x\phi$, which is (a). Differentiating $\grad_y\phi(x,y(x,p))=p$ gives $N+S\,\partial_xy=0$, so $\partial_xy=-WN$. Differentiating (a) gives $\Hess_x\psi_p=-\Phi-N^{\mathsf T}\partial_xy=-(\Phi-N^{\mathsf T}WN)$, which is (b). The substitution $p=\Lambda(x)+S(x)y$ is an automorphism of $\kk(x)[y]$ over $\kk(x)$, so (c) follows from (b).
\end{proof}

\begin{lemma}\label{lem:doubling-det}%
Let $M$ be a symmetric $m\times m$ matrix and $N$ an $m\times m$ matrix over a commutative ring. Then $\det\begin{psmallmatrix}M&N\T\\ N&0\end{psmallmatrix}=(-1)^m(\det N)^2$. In particular, for $x=(x_1,\dots,x_m)$, $y=(y_1,\dots,y_m)$ and $\phi=E(x)+\vv{y,F(x)}$ with $E\in\kk[x]$, $F\in\kk[x]^m$, one has $\det\Hess\phi=(-1)^m(\det JF)^2$.
\end{lemma}
\begin{proof}
$\begin{psmallmatrix}M&N\T\\ N&0\end{psmallmatrix}=\begin{psmallmatrix}N\T&M\\ 0&N\end{psmallmatrix}\begin{psmallmatrix}0&I\\ I&0\end{psmallmatrix}$, and the last factor has determinant $(-1)^m$. For $\phi$ one has $\partial_{y_i}\partial_{x_j}\phi=\partial_jF_i$ and $\partial_{y_i}\partial_{y_k}\phi=0$.
\end{proof}

\begin{remark}
For $m=1$ and $S$ constant, Lemma~\ref{lem:legendre} is the mechanism of Theorem~\ref{thm:pivot} and of the Schur descent of \cite{MY26}. For $n=m=2$ and $\det S\in\kk^\times$ it reduces the injectivity of $\grad\phi$ to $\HC_2$ (Theorem~\ref{thm:block}).
\end{remark}

\subsection{The construction of Gordan and Noether for rational functions}

\begin{lemma}\label{lem:graph}%
\textup{(a)} Let $h\in\kk(x_1,\dots,x_n)$ with $\det\Hess h=0$. Then there is $H\in\kk(x)^n\setminus\{0\}$ with $\Hess h\cdot H=0$ such that the derivation $D=\sum_iH_i\partial_i$ of $\kk(x)$ satisfies $D(H_i)=0$ for every $i$.

\textup{(b)} Let $a\in\kk[x_1,\dots,x_n]$ be nonconstant with $\det\Hess(1/a)=0$, and let $H$ be as in \textup{(a)} for $h=1/a$. Then $\Hess a\cdot H=0$ and $\vv{\grad a,H}=0$. In particular $\det\Hess a=0$ and $\BH(a)=0$.
\end{lemma}
\begin{proof}
(a) Put $f=\grad h$. The components of $f$ are algebraically dependent over $\kk$. Otherwise $\kk(x)$ is algebraic over $\kk(f)$, and each $x_j$ is a root of a polynomial $P_j(f,T)$, $P_j\in\kk[\Phi_1,\dots,\Phi_n,T]$, with $(\partial_TP_j)(f,x_j)\neq0$, since $\operatorname{char}\kk=0$. Applying $\partial_k$ to $P_j(f,x_j)=0$ gives $Q\cdot\Hess h=-\mathrm{diag}\bigl((\partial_TP_j)(f,x_j)\bigr)$ with $Q_{ji}=(\partial_{\Phi_i}P_j)(f,x_j)$, so $\det\Hess h\neq0$, a contradiction. Choose $R\in\kk[\Phi]\setminus\{0\}$ of minimal degree with $R(f)=0$ and put $H=(\grad_\Phi R)(f)$. Since $R$ is nonconstant and $\operatorname{char}\kk=0$, some $\partial_{\Phi_i}R$ is nonzero, of smaller degree, so $H_i\neq0$. Applying $\partial_j$ to $R(f)=0$ and using the symmetry of $\Hess h$ gives $\Hess h\cdot H=0$. Finally $D(H_k)=\sum_l(\partial_{\Phi_l}\partial_{\Phi_k}R)(f)\,(\Hess h\cdot H)_l=0$.

(b) By (a), $D(\partial_kh)=(\Hess h\cdot H)_k=0$, so $c=Dh=\vv{\grad h,H}$ satisfies $Dc=0$. Since $D(H_i)=0$, induction gives $D^ma=\sum_{i_1,\dots,i_m}H_{i_1}\cdots H_{i_m}\,\partial_{i_1}\cdots\partial_{i_m}a$, so $D^ma=0$ for $m>\deg a$. On the other hand $a=h^{-1}$, $Dh=c$ and $Dc=0$ give $D^ma=(-1)^mm!\,c^m\,h^{-1-m}$. Hence $c=0$. Then $Da=-a^2c=0$, that is $\vv{\grad a,H}=0$, and $\grad a=-a^2\grad h$ gives $D(\partial_ka)=-2a\,(Da)\,\partial_kh-a^2D(\partial_kh)=0$, that is $\Hess a\cdot H=0$. Since $H\neq0$, $\det\Hess a=0$. If $\rk\Hess a\le n-2$, then $\adj\Hess a=0$. If $\rk\Hess a=n-1$, the columns of the symmetric matrix $\adj\Hess a$ lie in the kernel $\kk(x)H$, so $\adj\Hess a=\kappa HH^{\mathsf T}$ with $\kappa\in\kk(x)$. In both cases $\BH(a)=\kappa\vv{\grad a,H}^2=0$.
\end{proof}

\begin{remark*}
For a polynomial $h$, part (a) is the construction of Gordan and Noether, and $x+H$ is a quasi-translation \cite[Proposition~2.1]{dB18}, that is, by Proposition~\ref{prop:qt-pf}, a Pascal finite map of Pascal depth at most two. With $H$ polynomial the same construction appears in Lemma~\ref{lem:slice}. Since $D(H)=0$, one has $p(x+sH(x))=\sum_ms^mD^mp/m!$ for every polynomial $p$, so in (b) the polynomials $a$ and $\partial_ka$ are constant on the lines $s\mapsto x+sH(x)$.
\end{remark*}

\subsection{Quasi-translations and Pascal finite maps}\label{ssec:qt}
For a polynomial map $F\in\kk[X]^n$, $X=(X_1,\dots,X_n)$, let $\sigma_F$ be the $\kk$-algebra endomorphism $g\mapsto g\circ F$ of $\kk[X]$ and $\Delta_F=\sigma_F-\mathrm{Id}$. Following \cite{ABCH19}, $F$ is called \emph{Pascal finite} if there is $m\ge1$ with $\Delta_F^m(X_i)=0$ for $i=1,\dots,n$, and the least such $m$ is the \emph{Pascal depth} of $F$. A polynomial map $X+H$ is a \emph{quasi-translation} if it is a polynomial automorphism with inverse $X-H$ \cite[Section~1]{dB18}.

\begin{proposition}\label{prop:qt-pf}%
Let $H\in\kk[X]^n$. The following are equivalent.
\begin{enumerate}[label=\textup{(\alph*)}]
\item $X+H$ is Pascal finite of Pascal depth at most two.
\item $H(X+H)=H$.
\item $JH\cdot H=0$.
\item $X+H$ is a quasi-translation.
\end{enumerate}
If they hold, then $D=\sum_iH_i\partial_i$ is a locally nilpotent derivation of $\kk[X]$ with $\sigma_{X+H}=\exp(D)$, and $H(X+sH)=H$ for every $s\in\kk$.
\end{proposition}
The equivalence of \textup{(b)}, \textup{(c)} and \textup{(d)} is \cite[Proposition~1.3]{dB18}, and we include the argument.
\begin{proof}
Put $F=X+H$ and $\Delta=\Delta_F$. Since $\Delta_F(X_i)=H_i$ and $\Delta_F(H_i)=H_i(X+H)-H_i$, (a) and (b) are equivalent. If $X-H$ is the inverse of $X+H$, then $X=(X-H)\circ(X+H)=X+H-H(X+H)$, so (d) implies (b).

(c)$\Rightarrow$(d). The derivation $D$ satisfies $D(X_i)=H_i$ and $D(H_i)=(JH\cdot H)_i=0$, so $D^2(X_i)=0$ and $D$ is locally nilpotent. Hence $\exp(sD)$ is a $\kk$-algebra homomorphism $\kk[X]\to\kk[X][s]$ with $\exp(sD)(X_i)=X_i+sH_i$, that is $\exp(sD)(g)=g(X+sH)$. From $D(H_i)=0$ one gets $H(X+sH)=H$. For $s=1$ and $s=-1$ this gives $\sigma_{X+H}=\exp(D)$ and $\sigma_{X-H}=\exp(-D)$, and these are mutually inverse.

(a)$\Rightarrow$(c). From $\sigma_F(gh)=\sigma_F(g)\sigma_F(h)$ one gets $\Delta(gh)=\Delta(g)h+g\Delta(h)+\Delta(g)\Delta(h)$, and induction gives $\Delta^{a+b-1}(gh)=0$ whenever $\Delta^ag=\Delta^bh=0$. Hence $\Delta$ is locally nilpotent on $\kk[X]$. For $t\in\kk$ put $\sigma^t=\sum_{k\ge0}\binom tk\Delta^k$. For fixed $g,h$ the identity $\sigma^t(gh)=\sigma^t(g)\sigma^t(h)$ is polynomial in $t$ and holds for $t\in\{0,1,2,\dots\}$, hence for all $t$. Its derivative at $t=0$ shows that $D'=\sum_{k\ge1}(-1)^{k+1}\Delta^k/k$ is a derivation of $\kk[X]$. Since $\Delta^2(X_i)=0$, one has $D'(X_i)=H_i$, so $D'=D$. As $D$ commutes with $\Delta$, $D(H_i)=D(\Delta X_i)=\Delta(D X_i)=\Delta^2(X_i)=0$, which is (c).
\end{proof}

\begin{remark*}
The implication (a)$\Rightarrow$(c) is the case of depth two of the theorem of \cite{ABCH19} that the Pascal finite automorphisms in characteristic zero are the exponentials of locally nilpotent derivations.
\end{remark*}

\emph{Quasi-translations and relations among the components of the gradient.} The constructions of Gordan and Noether in this paper rest on one algebraic mechanism. Let $h\in\kk[X]$ be a polynomial with singular Hessian, and let $R\neq0$ be a polynomial of least degree with $R(\grad h)=0$. Then $H=(\grad R)(\grad h)$ is a nonzero polynomial vector, and $X+H$ is a quasi-translation which leaves $\grad h$ invariant, that is, a Pascal finite map of Pascal depth at most two (Proposition~\ref{prop:qt-pf} and the remark above). This notion is defined by the finite-difference condition $\Delta_{X+H}^2(X)=0$, equivalently $H(X+H)=H$, which is meaningful over every field and involves no exponentials. In characteristic zero it has the following interpretation. The derivation $D=\sum_iH_i\partial_i$ is a locally nilpotent derivation of $\kk[X]$ with $\sigma_{X+H}=\exp(D)$, and its ring of constants $\kk[X]^D$ contains the components of $H$ and of $\grad h$, since $D(\partial_kh)=(\Hess h\cdot H)_k=0$. If $\rk\Hess h=n-1$, the ring of constants of $D$ has transcendence degree $n-1$ over $\kk$ \cite{Now94}, and so does $\kk(\grad h)$. Hence the field of fractions of $\kk[X]^D$ is a finite algebraic extension of $\kk(\grad h)$. For a rational function $h$, Lemma~\ref{lem:graph}(a) gives only a derivation of the field $\kk(x)$ with these properties, which in general does not come from a locally nilpotent derivation of $\kk[X]$. For a Hessian of rank $n-2$ two independent relations give two such derivations $D_{G_1},D_{G_2}$ with $D_{G_1}(G_2)=D_{G_2}(G_1)=0$, which therefore commute (proof of Lemma~\ref{lem:E7}), and for a Hessian of rank at most one there are $n-1$ of them (Lemma~\ref{lem:app-low}). For a Hessian of rank $n-1$ the relation $R$ generates the ideal of relations, and a constant vector $d$ is a translation direction of the closure of the gradient image if and only if $D_dR=0$, that is, if and only if $R$ does not involve one coordinate after a linear change of the target coordinates (Proposition~\ref{prop:qt-hess}). Since $d\T H=(D_dR)(\grad h)$, this means that the linear form $\vv{d,X}$ is invariant under the quasi-translation $X+H$, because $\vv{d,X+H}=\vv{d,X}+d\T H$. In characteristic zero it lies in the ring of constants of $D$. In Example~\ref{ex:apex} this linear form is $Z_4$. In four variables such a vector exists by Theorem~\ref{thm:apex}, due to Tang \cite[Theorem~2.3]{Tan26}.

\begin{proposition}\label{prop:gamma-alpha}%
Let $B\in\kk[x_1,\dots,x_n]$ be nonconstant, and let $\alpha\in\Q\setminus\{0\}$ with $1/\alpha\notin\{1,2,3,\dots\}$. If
\[
B\det\Hess B+(\alpha-1)\BH(B)=0,
\]
then $\BH(B)=0$ and $\det\Hess B=0$.
\end{proposition}
\begin{proof}
Write $\alpha=p/q$ with $q\ge1$ and $\gcd(p,q)=1$. Let $\mathbb L=\kk(x)[\vartheta]$, where $\vartheta$ is a root of an irreducible factor of $T^q-B$, put $h=\vartheta^p\in\mathbb L$, and extend $\partial_1,\dots,\partial_n$ to $\mathbb L$, which is possible in a unique way since $\mathbb L/\kk(x)$ is separable. Then $\partial_ih=\alpha h\,\partial_iB/B$, and
\[
\det\Hess h=\alpha^nh^nB^{-n-1}\bigl(B\det\Hess B+(\alpha-1)\BH(B)\bigr)=0 .
\]
The proof of Lemma~\ref{lem:graph}(a) applies verbatim in $\mathbb L$, since $\mathbb L$ is algebraic over $\kk(x)$. In particular the components of $\grad h$ are algebraically dependent over $\kk$. Otherwise $\mathbb L$ would be algebraic over the subfield they generate, a derivation of $\mathbb L$ vanishing on them would be zero, and $\det\Hess h\ne0$. It gives $H\in\mathbb L^n\setminus\{0\}$ with $\Hess h\cdot H=0$ and $D(H_i)=0$ for $D=\sum_iH_i\partial_i$. Hence $D(\partial_kh)=0$, and $c=Dh$ satisfies $Dc=0$. Put $\psi=c/(\alpha h)$. Then $DB=\psi B$ and $D\psi=-\alpha\psi^2$, so induction gives
\[
D^mB=\prod_{i=0}^{m-1}(1-i\alpha)\,\psi^mB .
\]
The product does not vanish, since $1/\alpha\notin\{1,2,\dots\}$. As in the proof of Lemma~\ref{lem:graph}(b), $D^mB=0$ for $m>\deg B$, because $D(H_i)=0$. Therefore $\psi=0$, $c=0$, $\vv{\grad B,H}=DB=0$, and $D(\partial_kB)=D\bigl(B\,\partial_kh/(\alpha h)\bigr)=0$, that is $\Hess B\cdot H=0$. The conclusion follows as at the end of the proof of Lemma~\ref{lem:graph}(b).
\end{proof}

\begin{remark}
The hypothesis on $\alpha$ is necessary. For $\alpha=\frac12$ and $B=x_1^2+x_2^2$ one has $\det\Hess B=4$, $\BH(B)=8B$ and $B\det\Hess B-\frac12\BH(B)=0$, while $\BH(B)\neq0$. More generally, for $m\ge1$, $\alpha=\frac1m$ and $B=(x_1^2+x_2)^m$ one has $\det\Hess(B^\alpha)=0$, hence $B\det\Hess B+(\alpha-1)\BH(B)=0$ by the identity in the proof, while $\BH(B)=2m^3(x_1^2+x_2)^{3m-3}\neq0$. So every excluded value of $\alpha$ is necessary. The case $\alpha=-1$ is Lemma~\ref{lem:graph}(b), since $a^{2n+1}\det\Hess(1/a)=(-1)^n\bigl(a\det\Hess a-2\BH(a)\bigr)$.
\end{remark}

\subsection{Further tools}
\begin{lemma}\label{lem:det2}
Let $R$ be a commutative ring and $J=\begin{psmallmatrix}0&-1\\1&0\end{psmallmatrix}$. For $X,Y\in R^{2\times2}$ put $\{X,Y\}=\operatorname{tr}(\adj(X)Y)$, a symmetric bilinear form. Then $\det(X+Y)=\det X+\{X,Y\}+\det Y$. For $X$ symmetric and $u\in R^2$ one has $u^{\mathsf T}\adj(X)u=(Ju)^{\mathsf T}X(Ju)$ and $\adj(uu^{\mathsf T})=(Ju)(Ju)^{\mathsf T}$.
\end{lemma}
\begin{proof}
Direct computation, $\adj$ being linear on $2\times2$ matrices.
\end{proof}

\begin{lemma}\label{lem:hess2F}
Let $F$ be a field of characteristic zero and $f\in F[u_1,u_2]$ with $\det\Hess f=0$. Then, after possibly exchanging $u_1$ and $u_2$, $f=Q(u_1+\alpha u_2)+m_1u_1+m_2u_2+m_0$ with $\alpha,m_i\in F$ and $Q\in F[T]$, and $\Hess f=Q''(\ell)\grad\ell\grad\ell^{\mathsf T}$ for $\ell=u_1+\alpha u_2$.
\end{lemma}
\begin{proof}
Take $\lambda,\gamma$ from \ref{E:dBvdE}, put $\ell=\lambda_2u_1-\lambda_1u_2$ and choose a linear form $r$ with $\vv{\lambda,\grad r}=1$. Then $\sum_i\lambda_i\partial_i(f-\gamma r)=0$, so $f-\gamma r\in F[\ell]$. It remains to normalize $\ell$.
\end{proof}

\begin{lemma}\label{lem:devel}%
Let $f\in\kk(x_1,x_2)$ with $\det\Hess f=0$ and $\Hess f\neq0$.
\begin{enumerate}[label=\textup{(\alph*)}]
\item There are $\rho\in\kk(x)^\times$ and $n\in\kk(x)^2$ with $\Hess f=\rho\,nn^{\mathsf T}$. Put $u=Jn$. The integral curves of $u$, called the \emph{generators} of $f$, are straight lines, $\grad f$ is constant along each generator, and the restriction of $f$ to a generic generator is an affine function of the parameter of the line. For every $g$ one has $\{\Hess f,\Hess g\}=\rho\,u^{\mathsf T}\Hess g\,u$.
\item If $f$ has a pole along an irreducible curve $L$, then either all generators are parallel and $L$ is a line parallel to them, or all generators pass through one point $x^*\in L$, that is $u\in\kk(x)\,(x-x^*)$.
\item If all generators pass through $x^*$, then $f=c+F(x-x^*)$ with $c\in\kk$ and $F$ rational and homogeneous of degree one.
\end{enumerate}
\end{lemma}
\begin{proof}
(a) A symmetric $2\times2$ matrix of rank one is of the form $\rho nn^{\mathsf T}$. Lemma~\ref{lem:graph}(a) gives $H\in\kk(x)^2\setminus\{0\}$ with $\Hess f\cdot H=0$ and $D(H)=0$ for $D=H_1\partial_1+H_2\partial_2$, so $H\in\kk(x)u$. Then $D(\partial_kf)=0$, hence $D^2f=0$, and $\exp(sD)$, a homomorphism $\kk(x)\to\kk(x)[[s]]$ agreeing with the substitution $x\mapsto x+sH$ on polynomials, gives $f(x+sH)=f+s\,Df$ and $\grad f(x+sH)=\grad f$. The last formula follows from Lemma~\ref{lem:det2}, since $\adj(\rho nn^{\mathsf T})=\rho uu^{\mathsf T}$.

(b) Let $W\in\kk[x]$ be an irreducible equation of $L$, and let $H$ and $D$ be as in the proof of (a). Write $H=G/a$ with $a\in\kk[x]\setminus\{0\}$ and $G=(G_1,G_2)\in\kk[x]^2$, $\gcd(G_1,G_2)=1$, and put $E=aD=G_1\partial_1+G_2\partial_2$. From $D(H_i)=0$ one gets $aE(G_i)=G_iE(a)$, so $E(G_i)=cG_i$ with $c=E(a)/a$. If $c=c_1/c_2$ with $\gcd(c_1,c_2)=1$, then $c_2E(G_i)=c_1G_i$, so $c_2$ divides $G_1$ and $G_2$ and is constant. Hence $c\in\kk[x]$. By the proof of (a), $D^2f=0$, and $D^2f=a^{-2}(E^2f-cEf)$ gives $E^2f=cEf$.

Let $\nu$ be the valuation of $\kk(x)$ defined by $W$. Suppose $W\nmid E(W)$. If $h=W^mh_0$ with $m\neq0$ and $\nu(h_0)=0$, then $E(h)=W^{m-1}\bigl(mh_0E(W)+WE(h_0)\bigr)$ with $\nu(E(h_0))\ge0$. Since $m\neq0$ and $W$ divides neither $h_0$ nor $E(W)$, the first summand has valuation $0$, so $E(h)\neq0$ and $\nu(E(h))=m-1$. Since $\nu(f)<0$, this gives $\nu(E^2f)=\nu(f)-2<\nu(f)-1\le\nu(cEf)$, contrary to $E^2f=cEf$. Hence $E(W)=rW$ with $r\in\kk[x]$.

We show that $L$ is a line. Induction using $E(x_i)=G_i$ and $E(gG_i)=(E(g)+cg)G_i$ gives $E^m(x_i)=c_mG_i$ with $c_m\in\kk[x]$ and $c_1=1$. The homomorphism $\exp(sE)\colon\kk[x]\to\kk[x][[s]]$ is therefore the substitution $x\mapsto x+\sigma G$ with $\sigma=\sum_{m\ge1}c_ms^m/m!$. As $W$ divides every $E^m(W)$, one has $W(x+\sigma G)\in W\,\kk[x][[s]]$. Since $G$ is primitive, $W$ does not divide both $G_1$ and $G_2$, and there is $q\in L$ with $G(q)\neq0$. The polynomial $g(\tau)=W(q+\tau G(q))\in\kk[\tau]$ satisfies $g(\sigma(q))=0$ in $\kk[[s]]$, where $\sigma(q)$ has order one in $s$. Hence $g=0$, the equation of the line through $q$ with direction $G(q)$ divides $W$, and $L$ is this line.

After a linear change of coordinates, $L=\{x_1=0\}$ and $f=P/(x_1^eQ)$ with $e\ge1$, $P,Q\in\kk[x]$ and $x_1\nmid PQ$. By the proof of (a), $f(x+sH)=f+s\,Df$ in $\kk(x)[[s]]$, and the substitution $s\mapsto as$ gives $f(x+sG)=f+s\,Ef$. Hence
\[
P(x+sG)=(x_1+sG_1)^e\,Q(x+sG)\,(f+s\,Ef)\qquad\text{in }\kk(x)[s].
\]
If $G_1=0$, then $G_2\in\kk^\times$, and the generators are the lines parallel to $L$. Let $G_1\neq0$. Evaluation at $s=-x_1/G_1$ gives $P(0,t)=0$ with $t=x_2-x_1G_2/G_1\in\kk(x)$. Since $x_1\nmid P$, the polynomial $P(0,T)\in\kk[T]$ is nonzero, so $t$ is algebraic over $\kk$, and $t=t_0\in\kk$ because $\kk$ is algebraically closed in $\kk(x)$. Then $(x_2-t_0)G_1=x_1G_2$, so $G_1=\rho x_1$ and $G_2=\rho(x_2-t_0)$ with $\rho\in\kk[x]$, and $\rho\in\kk^\times$ by primitivity of $G$. Thus $u\in\kk(x)\,(x-x^*)$ with $x^*=(0,t_0)\in L$.

(c) Put $w=x-x^*$ and $F=\vv{w,\grad f}$. Since $w$ is tangent to the generator through $x$, $\Hess f\cdot w=0$, so $\grad F=\grad f$ and $f-F$ is constant. By the proof of (a), $D(\partial_kf)=0$ for every $k$, and $H$ is a nonzero multiple of $w$. Hence $E(\partial_kf)=0$ for the Euler operator $E=\sum_iw_i\partial_i$. Write $\partial_kf=P/R$ with coprime polynomials $P,R$ in $w$. Then $R\,E(P)=P\,E(R)$, so $R$ divides $E(R)$. As $\deg E(R)\le\deg R$, one has $E(R)=cR$ with $c\in\kk$, and then $E(P)=cP$. So $P$ and $R$ are homogeneous of the same degree, $\partial_kf$ is homogeneous of degree zero in $w$, and $F=\vv{w,\grad f}$ is homogeneous of degree one in $w$.
\end{proof}

\begin{remark*}
Geometrically, the first part of the proof of (b) shows that $L$ is invariant under the vector field $G$, hence $L$ is a generator. The second part shows that a generator meeting $L$ at a point where the numerator $P$ of $f$ does not vanish would carry a pole of $f$, which is excluded since $f$ is affine along the generators.
\end{remark*}

\section{Polynomials with vanishing bordered Hessian}\label{sec:BH}

In this section we prove Theorem~\ref{thm:main-BH}. Up to sign, the bordered Hessian of $A$ is the coefficient of $x_4^2$ in the Hessian determinant of $x_4A+B$, so its vanishing is the first condition imposed by an affine direction. By Proposition~\ref{prop:BHdet} it forces $\det\Hess A=0$, and in three variables the classification \ref{E:dB62} then leaves only the polynomials in two linear forms. The leading coefficient $c_m$ of a potential in a variable $t$ satisfies an equation of the same type, which Proposition~\ref{prop:gamma-alpha} reduces to $\BH(c_m)=0$.

\begin{lemma}\label{lem:2d}%
If $A\in\kk[x_1,x_2]$ and $\BH(A)=0$, then $A=f(\ell)$ with $\ell$ a linear form and $f\in\kk[T]$.
\end{lemma}
\begin{proof}
By Proposition~\ref{prop:BHdet}, $\det\Hess A=0$, and \ref{E:dBvdE} gives $\lambda\in\kk^2\setminus\{0\}$ and $\gamma\in\kk$ with $\lambda_1\partial_1A+\lambda_2\partial_2A=\gamma$. After a linear change of coordinates $\lambda=e_2$, so $A=g(x_1)+\gamma x_2$ with $g\in\kk[x_1]$. Then $\BH(A)=\gamma^2g''$. If $\gamma=0$, then $A=g(x_1)$. Otherwise $g''=0$, and $A$ is an affine function of the linear form $\ell=g'(0)x_1+\gamma x_2$.
\end{proof}

\begin{theorem}[reduction to two variables]\label{thm:cylinder}%
If $A\in\kk[x_1,x_2,x_3]$ and $\BH(A)=0$, then $A\in\kk[x_1,x_2]$ after a linear change of coordinates.
\end{theorem}
\begin{proof}
By Proposition~\ref{prop:BHdet}, $\det\Hess A=0$. Suppose first that $\partial_cA=\gamma\in\kk$ for some $c\neq0$. After a linear change with $c=e_3$ one has $A=g(x_1,x_2)+\gamma x_3$. If $\gamma=0$, this is the claim. If $\gamma\neq0$, then $\adj\Hess A=\det\Hess_{(x_1,x_2)}g\cdot e_3e_3^{\mathsf T}$ and $\BH(A)=\gamma^2\det\Hess_{(x_1,x_2)}g$, so $\det\Hess g=0$. By Lemma~\ref{lem:hess2F}, $g=Q(\ell)+m$ with $\ell$ a linear form and $m$ an affine function of $(x_1,x_2)$, and $A=Q(\ell)+m+\gamma x_3$ is a function of the two linear forms $\ell$ and $m+\gamma x_3$, up to constants. Otherwise \ref{E:dB62} gives, after a linear change, $A=a_1(x_1)+a_2(x_1)x_2+a_3(x_1)x_3$ with $a_i\in\kk[x_1]$. Then $\grad A=(a_1'+a_2'x_2+a_3'x_3,\,a_2,\,a_3)$, the second and third rows of $\Hess A$ are $(a_2',0,0)$ and $(a_3',0,0)$, and
\[
\adj\Hess A=\begin{pmatrix}0&0&0\\0&-a_3'^2&a_2'a_3'\\0&a_2'a_3'&-a_2'^2\end{pmatrix},\qquad \BH(A)=-\bigl(a_2a_3'-a_3a_2'\bigr)^2 .
\]
Hence $a_2a_3'=a_3a_2'$. If $a_3=0$, then $A\in\kk[x_1,x_2]$. Otherwise $(a_2/a_3)'=0$, so $a_2=\kappa a_3$ with $\kappa\in\kk$, and $A=a_1(x_1)+a_3(x_1)(\kappa x_2+x_3)$ lies in $\kk[x_1,\kappa x_2+x_3]$.
\end{proof}

\begin{remark}\label{rem:cyl}
Theorem~\ref{thm:cylinder} does not hold for $\det\Hess A=0$ alone, since $A=x_1x_2+x_3$ has $\det\Hess A=0$ and linearly independent partial derivatives, while $\BH(A)=-1$. It does not hold in four variables either, since $A=\sum_{i=1}^3x_{i+1}g_i(x_1)$ has $\BH(A)=0$ and, for $g_1,g_2,g_3$ linearly independent over $\kk$, depends on all four variables. Its Hessian has rank two. Hessian rank three occurs as well. Let $p(t)=(1,t,t^2)$, $s=\langle p(x_4),(x_1,x_2,x_3)\rangle=x_1+x_4x_2+x_4^2x_3$ and $A=g(x_4,s)$ with $g\in\kk[t,s]$. The vector field $u=(t^2,-2t,1,0)$, $t=x_4$, satisfies $D_ux_4=D_us=D_u(x_2+2x_4x_3)=0$. The entries of $\grad A$ are polynomials in $x_4$, $s$ and $x_2+2x_4x_3$, hence $\Hess A\cdot u=0$ and $\langle\grad A,u\rangle=0$. If $\Hess A$ has rank three, then $\adj\Hess A$ is a multiple of $uu^{\mathsf T}$ and $\BH(A)=0$. For $A=s^2$ the rank is three, and there is no $c\neq0$ with $\partial_cA=0$. All polynomials in four variables with $\BH(A)=0$ are described in Proposition~\ref{prop:BH4}.%
Geometrically, for $\kk=\C$ and a general value $c$, the condition $\BH(A)=0$ says that the Gauss map of the surface $\{A=c\}$ has rank at most one, and Theorem~\ref{thm:cylinder} says that this surface is a cylinder. A proof along these lines can be based on the classification of surfaces with degenerate Gauss map \cite[\S2]{GH79}, \cite[Chapter~3]{FP01}.
\end{remark}

\begin{theorem}\label{thm:gamma}%
If $a\in\kk[x_1,x_2,x_3]$ and $\det\Hess(1/a)=0$, then $a\in\kk[x_1,x_2]$ after a linear change of coordinates. If moreover $a\in\kk[x_1,x_2]$ and $a\det\Hess a-2\BH(a)=0$, then $a=f(\ell)$ for a linear form $\ell$.
\end{theorem}
\begin{proof}
By Lemma~\ref{lem:graph}(b), $\BH(a)=0$, and Theorem~\ref{thm:cylinder} applies. In two variables one has $a^5\det\Hess(1/a)=a\det\Hess a-2\BH(a)$, so Lemma~\ref{lem:graph}(b) gives $\BH(a)=0$, and Lemma~\ref{lem:2d} applies.
\end{proof}

\begin{theorem}\label{thm:lead}%
Let $\phi=\sum_{j=0}^mt^jc_j(x)$ with $x=(x_1,x_2,x_3)$, $m\ge2$ and $c=c_m\neq0$, and assume $\det\Hess\phi\in\kk^\times$. Then $\det\Hess_{(x,t)}(t^mc)=0$, and after a linear change of $x$ one has $c\in\kk[x_1,x_2]$.
\end{theorem}
\begin{proof}
In $\Hess\phi$ the entries indexed by $(x,x)$, $(x,t)$ and $(t,t)$ have degree in $t$ at most $m$, $m-1$ and $m-2$, with leading parts $t^m\Hess c$, $mt^{m-1}\grad c$ and $m(m-1)t^{m-2}c$. Every term of the Leibniz expansion has degree at most $4m-2$ in $t$, and the coefficient of $t^{4m-2}$ is
\[
N(c)=\det\begin{pmatrix}\Hess c&m\grad c\\ m\grad c^{\mathsf T}&m(m-1)c\end{pmatrix}=m(m-1)\,c\det\Hess c-m^2\BH(c).
\]
The same matrix, up to powers of $t$, is $\Hess(t^mc)$. Since $4m-2>0$, $N(c)=0$. If $c$ is constant there is nothing to prove. Otherwise $N(c)=m(m-1)\bigl(c\det\Hess c+(\alpha-1)\BH(c)\bigr)$ with $\alpha=-1/(m-1)$, and $1/\alpha=1-m$ is not a positive integer. Proposition~\ref{prop:gamma-alpha} gives $\BH(c)=0$, and Theorem~\ref{thm:cylinder} applies.
\end{proof}

\begin{remark}
For $m=2$ the theorem also follows from Theorem~\ref{thm:gamma}, since then $N(c)=2(c\det\Hess c-2\BH(c))$. For $m=1$ the top coefficient is $N(c)=-\BH(c)$, and Theorem~\ref{thm:cylinder} is needed in full.
\end{remark}

\begin{proof}[Proof of Theorem~\ref{thm:main-BH}]
Over an algebraically closed field, part (i) is Proposition~\ref{prop:BHdet}, part (ii) is Theorem~\ref{thm:cylinder}, and part (iii) is Theorem~\ref{thm:lead} for $m\ge2$ and, for $m=1$, follows from part (ii), since the coefficient of $t^2$ in $\det\Hess\phi$ is $-\BH(c_1)$ as in the proof of Theorem~\ref{thm:affine}. Over an arbitrary field $\kk$ of characteristic zero, part (i) is an identity between polynomials, and in parts (ii) and (iii) the space of vectors $c$ with $D_cA=0$, respectively $D_cc_m=0$, is nonzero over $\overline\kk$ and hence contains a nonzero vector of $\kk^3$ by Lemma~\ref{lem:lindescent}.
\end{proof}
\section{Affine directions and constant pivots}\label{sec:affine}

If $v$ is an affine direction, then $\grad\phi$ is affine along every line of direction $v$ (Remark~\ref{rem:lines}). With $v=e_4$ the potential is $x_4A+B$, and the coefficient of $x_4^2$ in $\det\Hess\phi$, which is $-\BH(A)$, forces $A$ to depend on two variables. The coefficient of $x_4$ then gives either $B_{33}=0$, which is the doubling, or an affine $A$, and in the second case the injectivity of $\grad\phi$ is reduced to $\HC_2$. A constant pivot is treated by a partial Legendre transform in one variable and $\HC_3$.

\begin{theorem}\label{thm:affine}%
Let $\phi=x_4A(x_1,x_2,x_3)+B(x_1,x_2,x_3)$ be a Keller potential. After a linear change of $(x_1,x_2,x_3)$ one has $A=A(x_1,x_2)$ and exactly one of the following holds.
\begin{enumerate}[label=\textup{(\roman*)}]
\item $B_{33}=0$. Then $B=x_3C(x_1,x_2)+E(x_1,x_2)$, $(C,A)$ is a Keller map of the plane, and $\grad\phi$ is injective if and only if $(C,A)$ is injective.
\item $B_{33}\neq0$. Then after a further linear change $A=\alpha x_1+\beta$ with $\alpha\neq0$, $\det\Hess_{(x_2,x_3)}B\in\kk^\times$, and $\grad\phi$ is injective.
\end{enumerate}
\end{theorem}
\begin{proof}
Since $\Hess\phi=\begin{psmallmatrix}x_4\Hess A+\Hess B&\grad A\\ \grad A^{\mathsf T}&0\end{psmallmatrix}$, one has $\det\Hess\phi=-\grad A^{\mathsf T}\adj(x_4\Hess A+\Hess B)\grad A$, and its coefficient of $x_4^2$ is $-\BH(A)$. So $\BH(A)=0$ and Theorem~\ref{thm:cylinder} gives $A=A(x_1,x_2)$. Then
\[
\det\Hess\phi=-B_{33}x_4\BH(A)-B_{33}\bigl(B_{11}A_2^2-2B_{12}A_1A_2+B_{22}A_1^2\bigr)+(A_1B_{23}-A_2B_{13})^2 .
\]
The coefficient of $x_4$ vanishes, so $B_{33}=0$ or $\BH(A)=0$ in two variables.
(i) $B_{33}=0$ gives $B=x_3C+E$ and $\det\Hess\phi=\det J(C,A)^2$ by Lemma~\ref{lem:doubling-det}. Since $\grad\phi=\bigl(J(C,A)^{\mathsf T}(x_3,x_4)^{\mathsf T}+\grad E,\,C,\,A\bigr)$, a collision of $(C,A)$ at $a\neq b$ lifts to a collision of $\grad\phi$ at $(a,0,0)$ and $(b,x_3',x_4')$ with $(x_3',x_4')^{\mathsf T}=J(C,A)(b)^{-\mathsf T}(\grad E(a)-\grad E(b))$. Conversely injectivity of $(C,A)$ forces equality of the first two coordinates and then of the last two.
(ii) Lemma~\ref{lem:2d} gives $A=f(x_1)$ after a linear change, and $\det\Hess\phi=-f'(x_1)^2\det\Hess_{(x_2,x_3)}B$. Both factors are constants, so $f$ is affine. Equality of gradients at two points gives $x_1=x_1'$ from the last coordinate, then $(x_2,x_3)=(x_2',x_3')$ by \ref{E:dB} for $n=2$ applied to $B(x_1,\cdot,\cdot)$, and finally $x_4=x_4'$.
\end{proof}

\begin{remark}\label{rem:lines}
An affine direction is isotropic for every $\Hess\phi(x)$, but it is not a kernel vector, since $\det\Hess\phi=c\neq0$ gives $\Hess\phi(x)\,v\neq0$. Let $v$ be an affine direction of a Keller potential $\phi$ in $n$ variables. Since $\partial_t^2\grad\phi(x+tv)=\grad(D_v^2\phi)(x+tv)=0$,
\[
\grad\phi(x+tv)=\grad\phi(x)+t\,\Hess\phi(x)\,v ,
\]
and $\vv{\grad\phi(x+tv),v}=D_v\phi(x)$ does not depend on $t$. Thus $\grad\phi$ maps every line $x+\kk v$ by an affine isomorphism onto a line in the hyperplane $\{p:\vv{p,v}=D_v\phi(x)\}$, and the direction $w=\Hess\phi(x)\,v$ of this line is isotropic for $\Hess\phi(x)^{-1}$, since $w\T\Hess\phi(x)^{-1}w=v\T\Hess\phi(x)\,v=0$. It need not be isotropic for $\Hess\phi(x)$. For $\phi=x_1x_4+\frac12(x_1^2+x_2^2+x_3^2)$ and $v=e_4$ one has $\Hess\phi\,v=e_1$ and $e_1\T\Hess\phi\,e_1=1$. For $v=e_4$ and $\phi=x_4A+B$ as in Theorem~\ref{thm:affine} the hyperplane is $\{p_4=A(x)\}$. A collision of $\grad\phi$ can therefore occur only between distinct lines of the family.
\end{remark}

\begin{remark}[the gradient of a doubling]\label{rem:fibers}
Let $\phi=E(x)+\vv{y,F(x)}$ with $x=(x_1,x_2)$, $y=(x_3,x_4)$ and a Keller map $F$ of the plane. Put $\tau(x,y)=\bigl(x,\grad E(x)+JF(x)\T y\bigr)$ and $\pi(u,v)=(v,u)$. Since $JF(x)$ is invertible over $\kk[x]$, the map $\tau$ is a polynomial automorphism of $\kk^4$, and
\[
\grad\phi=\pi\circ(F\times\mathrm{Id}_{\kk^2})\circ\tau .
\]
Thus, up to a polynomial automorphism of the source and a permutation of the target coordinates, $\grad\phi$ is the product of $F$ with the identity of $\kk^2$. Consequently the fiber of $\grad\phi$ over $(p,q)$ is isomorphic to $F^{-1}(q)$, the image of $\grad\phi$ is $\kk^2\times F(\kk^2)$, and $\grad\phi$ is proper if and only if $F$ is. The inclusion $\kk[\grad\phi]\subseteq\kk[x,y]$ arises from $\kk[F]\subseteq\kk[x]$ by adjoining two independent variables, so the extension $\kk(x,y)/\kk(\grad\phi)$ has the degree of $\kk(x)/\kk(F)$ and a normal closure with the same Galois group. This is more than the equivalence of injectivity in Theorem~\ref{thm:affine}\textup{(i)}. In particular a counterexample to $\HC_4$ has the same fibers as the non-injective planar Keller map of which it is a doubling. For the Nambu tower of \cite[Theorem~4.7]{AHL} the doubling gives no condition beyond those for $F$ \cite[Corollary~5.5]{AHL}.
\end{remark}

\begin{example}\label{ex:case-ii}
The potential
\[
\phi=x_1x_4+(x_2+x_1x_3)\bigl(x_1x_2+(1+x_1^2)\,x_3\bigr)
\]
is of type (ii) of Theorem~\ref{thm:affine}, and it is not a doubling. Here $A=x_1$ and $B=x_1x_2^2+(1+2x_1^2)x_2x_3+x_1(1+x_1^2)x_3^2$, so $B_{33}=2x_1(1+x_1^2)\neq0$ and
\[
\det\Hess_{(x_2,x_3)}B=4x_1^2(1+x_1^2)-(1+2x_1^2)^2=-1 .
\]
Hence $\det\Hess\phi=-A_1^2\det\Hess_{(x_2,x_3)}B=1$, and $\grad\phi$ is a polynomial automorphism by Theorem~\ref{thm:affine}(ii). A doubling has a plane of affine directions, while the affine directions of $\phi$ are the nonzero multiples of $e_4$. Indeed, let $D_v^2\phi=0$. The coefficients of $x_1^3$ and of $x_1$ in $D_v^2\phi$ are $2v_3^2$ and $2v_2^2+2v_3^2$, so $v_2=v_3=0$. Then $D_v^2\phi=v_1^2\phi_{11}+2v_1v_4$ with $\phi_{11}=4x_2x_3+6x_1x_3^2$, so $v_1=0$.
\end{example}

\begin{corollary}\label{cor:ceiling}%
The gradient of every Keller potential in four variables with an affine direction is a polynomial automorphism if and only if $\JC_2$ holds.
\end{corollary}
\begin{proof}
An affine direction becomes $e_4$ after a linear change, and then $\phi$ is affine in $x_4$. By Theorem~\ref{thm:affine} and Lemma~\ref{lem:descent}, $\JC_2$ implies that $\grad\phi$ is an automorphism. Conversely, for a Keller map $F$ of the plane the doubling $\phi_F$ of \eqref{eq:doubling} has an affine direction, and $\grad\phi_F$ is injective if and only if $F$ is, by Theorem~\ref{thm:affine}(i) \cite[Proposition~2.3]{MY26}.
\end{proof}

\begin{proposition}\label{prop:block}%
Let $n,m\ge1$, $x=(x_1,\dots,x_n)$, $y=(y_1,\dots,y_m)$, and let $\phi=B(x)+\vv{y,\Lambda(x)}+\frac12y\T Sy$ be a Keller potential with $B\in\kk[x]$, $\Lambda\in\kk[x]^m$ and $S=S\T\in\mathrm{GL}_m(\kk)$ constant. If $\HC_n$ holds, then $\grad\phi$ is a polynomial automorphism. Conversely, if the gradients of all Keller potentials of this form are polynomial automorphisms, then $\HC_n$ holds.
\end{proposition}
\begin{proof}
Put $W=S^{-1}$ and, for $p\in\kk^m$, $\psi_p=-B+\frac12(p-\Lambda)\T W(p-\Lambda)\in\kk[x]$. By Lemma~\ref{lem:legendre}(c), $\det\Hess_x\psi_p=(-1)^n\det\Hess\phi/\det S\in\kk^\times$, so every $\psi_p$ is a Keller potential in $n$ variables, and under $\HC_n$ its gradient is injective. If $\grad\phi(x,y)=\grad\phi(x',y')$, the last $m$ coordinates give $p=\Lambda(x)+Sy=\Lambda(x')+Sy'$, and the first $n$ coordinates give $\grad\psi_p(x)=\grad\psi_p(x')$ by Lemma~\ref{lem:legendre}(a). Hence $x=x'$, then $y=y'$, and $\grad\phi$ is a polynomial automorphism by Lemma~\ref{lem:descent}. Conversely, for a Keller potential $\psi$ in $n$ variables the potential $\psi(x)+\frac12\vv{y,y}$ has the above form, and its gradient $(\grad\psi,y)$ is a polynomial automorphism if and only if $\grad\psi$ is.
\end{proof}

\begin{remark}
The proof is a partial Legendre transform in the variables $y$. For $m=1$ it is the Schur descent of \cite{MY26} read in the opposite direction.
\end{remark}

A vector $v\neq0$ is a constant pivot of $\phi$ if and only if, after a linear change of coordinates in which $v$ is the coordinate vector of a variable $u$, one has $\phi=\frac\kappa2u^2+u\,a(z)+b(z)$ with $\kappa\in\kk^\times$. This is the form of Proposition~\ref{prop:block} with $m=1$ and $S=(\kappa)$.

\begin{theorem}\label{thm:pivot}%
Let $\phi$ be a Keller potential in four variables with a constant pivot. Then $\grad\phi$ is an automorphism.
\end{theorem}
\begin{proof}
Proposition~\ref{prop:block} with $n=3$ and $m=1$, since $\HC_3$ holds by \ref{E:dB}.
\end{proof}

\begin{corollary}\label{cor:hc5}%
$\HC_4$ fails if and only if there is a Keller potential in five variables with a constant pivot and non-injective gradient.
\end{corollary}
\begin{proof}
Proposition~\ref{prop:block} with $n=4$ and $m=1$, and Lemma~\ref{lem:descent}.
\end{proof}

\begin{remark}\label{rem:block}
By $\HC_2$ and $\HC_3$, the gradient of a Keller potential in five variables which has the form of Proposition~\ref{prop:block} with $m\ge2$ in some linear coordinates is a polynomial automorphism. A constant invertible block contains a constant pivot, since $u\T Su\neq0$ for some $u\in\kk^m$. So the potentials $\psi_{\lambda,\mu}$ of Section~\ref{sec:MY} have no such block for any $m\ge1$, by Theorem~\ref{thm:MYcore}.
\end{remark}
\section{Solutions quadratic in two variables}\label{sec:block}

Here $\phi$ is quadratic in two of the variables. The partial Legendre transform of Lemma~\ref{lem:legendre} in these variables turns $\phi$ into a family $\psi_p$ of potentials in the two remaining variables with $\det\Hess\psi_p=\kappa/\det S$. If $\det S$ is constant, $\HC_2$ applies to every member of the family. Otherwise the top layer of $\det\Hess\phi$ in $y$ and the five layers of $\det\Hess\psi_p$ in $p$ restrict $S$ and the other coefficients so strongly that $\phi$ has an affine direction or is quadratic in three variables, where Lemma~\ref{lem:quadblock} applies.
In this section $x=(x_1,x_2)$, $y=(y_1,y_2)$ and
\begin{equation}\label{eq:blockZ}
\phi=B(x)+\vv{y,\Lambda(x)}+\tfrac12y^{\mathsf T}S(x)y,\qquad B\in\kk[x],\ \Lambda\in\kk[x]^2,\ S=S^{\mathsf T}\in\kk[x]^{2\times2}.
\end{equation}
We call such polynomials the \emph{block class}.

\begin{theorem}\label{thm:block}
If $\phi$ is as in \eqref{eq:blockZ} and $\det\Hess\phi=\kappa\in\kk^\times$, then $\grad\phi$ is a polynomial automorphism or $D_u^2\phi=0$ for some $u\in\kk^4\setminus\{0\}$.
\end{theorem}

In the second case Theorem~\ref{thm:affine} applies. As a polynomial in $y$, $\det\Hess\phi$ has degree at most four, and its part of degree four is $\det\Hess(\frac12y^{\mathsf T}S(x)y)$. Hence
\begin{equation}\label{eq:toplayer}
\det\Hess_{(x,y)}\bigl(\tfrac12y^{\mathsf T}S(x)y\bigr)=0 .
\end{equation}

\begin{lemma}\label{lem:quadblock}%
Let $z=(z_1,\dots,z_n)$ and $\phi=\frac12z^{\mathsf T}S(x_1)z+\vv{m(x_1),z}+B_0(x_1)$ with $S=S^{\mathsf T}\in\kk[x_1]^{n\times n}$, $m\in\kk[x_1]^n$, $B_0\in\kk[x_1]$ and $\det\Hess\phi=\kappa\in\kk^\times$. Then either there is $n_0\in\kk^n\setminus\{0\}$ with $S(x_1)n_0=0$, so that $D^2_{(0,n_0)}\phi=0$, or $\det S\in\kk^\times$, $S^{-1}$ is affine in $x_1$, and $\grad\phi$ is a polynomial automorphism.
\end{lemma}
\begin{proof}
Here $\Hess\phi=\begin{psmallmatrix}\phi_{11}&w^{\mathsf T}\\ w&S\end{psmallmatrix}$ with $\phi_{11}=\frac12z^{\mathsf T}S''z+\vv{m'',z}+B_0''$ and $w=S'z+m'$, and $\kappa=\det S\cdot\phi_{11}-w^{\mathsf T}\adj(S)w$. Primes denote derivatives in $x_1$. Suppose $\det S\neq0$. The part of degree two in $z$ gives $\frac12\det S\cdot S''=S'\adj(S)S'$, that is $S''=2S'S^{-1}S'$, since both sides are symmetric matrices. Differentiating $SS^{-1}=I$ twice gives $(S^{-1})''=S^{-1}\bigl(2S'S^{-1}S'-S''\bigr)S^{-1}$, hence $(S^{-1})''=0$. Hence $S^{-1}=T_0+x_1T_1$ with constant $T_i$, and $\det S=1/\det(T_0+x_1T_1)$ is a polynomial only if it is constant. So $\det S\in\kk^\times$ and $S^{-1}$ is polynomial. Let $\psi_p(x_1)$ be the transform of Lemma~\ref{lem:legendre} in $z$. Then $\psi_p''=-\kappa/\det S\in\kk^\times$, so $\psi_p$ is quadratic in $x_1$ and $\psi_p'$ is injective. If $\grad\phi(x_1,z)=\grad\phi(x_1',z')$, then $p=S(x_1)z+m(x_1)=S(x_1')z'+m(x_1')$, and Lemma~\ref{lem:legendre}(a) gives $\psi_p'(x_1)=\psi_p'(x_1')$. Therefore $x_1=x_1'$ and $z=S^{-1}(p-m)=z'$, and Lemma~\ref{lem:descent} applies. Suppose $\det S=0$. If $\rk S\le n-2$, then $\adj S=0$ and $\kappa=0$. So $\adj S=\sigma\nu\nu^{\mathsf T}$ with $\sigma\neq0$ and $\nu\in\kk[x_1]^n$ a primitive generator of $\ker S$, and $\kappa=-\sigma\vv{S'z+m',\nu}^2$. Constancy in $z$ gives $S'\nu=0$. Differentiating $S\nu=0$ gives $S\nu'=0$, so $\nu'=r\nu$ with $r\in\kk(x_1)$. Since the entries of $\nu$ are coprime, $r$ is a polynomial, and comparing degrees in $\nu'=r\nu$ gives $r=0$. Therefore $\nu=n_0$ is constant.
\end{proof}

\subsection{Rank at most one}
If $S=0$, every $u=(0,w)$ works. If $\rk S=1$, then $S_{11}S_{22}=S_{12}^2$, and unique factorization in $\kk[x]$ gives
\begin{equation}\label{eq:rk1}
S=d\,ff^{\mathsf T},\qquad d\in\kk[x]\setminus\{0\},\quad f=(f_1,f_2)\in\kk[x]^2,\quad\gcd(f_1,f_2)=1.
\end{equation}
Indeed, with $e=\gcd(S_{11},S_{22})$, $S_{11}=e\mu$, $S_{22}=e\nu$, one has $e\mid S_{12}$ and $\mu\nu=(S_{12}/e)^2$ with $\gcd(\mu,\nu)=1$, so $\mu=\varepsilon f_1^2$, $\nu=\varepsilon^{-1}f_2^2$ for a unit $\varepsilon$. If $f_1/f_2\in\kk\cup\{\infty\}$, then $f$ is constant and $u=(0,Jf)$ satisfies $D_u^2\phi=0$. It remains to treat $f_1/f_2\notin\kk$.

\begin{proposition}\label{prop:rk1top}%
Let $d\in\kk[x]\setminus\{0\}$ and $f\in\kk[x]^2$ with $f_1/f_2\notin\kk$, and put $\Theta=\frac12d\vv{y,f}^2$. Then $\det\Hess_{(x,y)}\Theta=0$ if and only if there is a linear form $\ell$ in $x$ with $d\,ff^{\mathsf T}\in\kk[\ell]^{2\times2}$.
\end{proposition}
\begin{proof}
If $d\,ff^{\mathsf T}\in\kk[\ell]^{2\times2}$, then $\frac12d\vv{y,f}^2$ is a function of the three linear forms $\ell,y_1,y_2$ and its Hessian is singular.

Conversely, let $\mathbb L$ be the field obtained from $\kk(x)$ by adjoining a square root of $d$, put $F=\sqrt d\,f\in\mathbb L^2$ and $A=\vv{y,F}$, so that $\frac12d\vv{y,f}^2=\frac12A^2$. Let $JF$ be the Jacobian matrix, $Y=\sum_ky_k\Hess F_k$ and $v=\adj(JF)F$. For arbitrary $F$
\begin{equation}\label{eq:rk1det}
\det\Hess\bigl(\tfrac12A^2\bigr)=A^3\bigl(3A(\det JF)^2-v^{\mathsf T}Yv\bigr).
\end{equation}
\emph{Case $\det JF\neq0$.} Since $A\neq0$ and both sides of $3A(\det JF)^2=v^{\mathsf T}Yv$ are linear in $y$, the vanishing of \eqref{eq:rk1det} means $w^{\mathsf T}\Hess F_k\,w=3F_k$ for $k=1,2$, with $w=JF^{-1}F=v/\det JF$. Differentiating $JF\,w=F$ along $w$ gives $w^{\mathsf T}\Hess F_k\,w+\vv{\grad F_k,(Jw)w}=F_k$, so the condition is $JF\,(Jw)w+2w=0$, that is $(Jw)w=-2w$. Choose $x_0\in\kk^2$ with $\det JF(x_0)\neq0$ and $d(x_0)f_1(x_0)\neq0$, and fix a square root of $d(x_0)$, so that $F$ and $w$ are given by power series at $x_0$. Let $x(\tau)\in\kk[[\tau]]^2$ be the formal solution of $\dot x=w(x)$, $x(0)=x_0$. Then $\ddot x=(Jw)\dot x=-2\dot x$, so $x(\tau)=x_0+\sigma c$ with $c=\frac12w(x_0)$ and $\sigma=1-e^{-2\tau}$. Moreover $\frac{d}{d\tau}F(x(\tau))=JF\,w=F$, so $F(x_0+\sigma c)=(1-\sigma)^{-1/2}F(x_0)$, and squaring the first component,
\[
d(x_0+\sigma c)\,f_1(x_0+\sigma c)^2=(1-\sigma)^{-1}d(x_0)f_1(x_0)^2\qquad\text{in }\kk[[\sigma]].
\]
Since $\sigma\in\tau\kk[[\tau]]$ has order one, this is an identity in $\kk[[\sigma]]$. Its left side is a polynomial in $\sigma$ and its right side is not, a contradiction.%

\emph{Case $\det JF=0$.} Put $P_1=df_1^2$, $P_2=df_2^2$, $P_{12}=df_1f_2\in\kk[x]$. Then $\det J(P_1,P_2)=4F_1F_2\det JF=0$, so $P_1,P_2$ are algebraically dependent and $P_{12}^2=P_1P_2$. If all three are constant, $d\,ff^{\mathsf T}$ is constant. Otherwise the subalgebra $\kk[P_1,P_2,P_{12}]$ has transcendence degree one, since $P_1,P_2$ are algebraically dependent and $P_{12}$ is algebraic over $\kk(P_1,P_2)$, and \ref{E:Zaks} gives $P_1=\gamma_1(g)$, $P_2=\gamma_2(g)$, $P_{12}=\gamma_{12}(g)$ for some $g\in\kk[x]$, so $d\,ff^{\mathsf T}=\Gamma\Gamma^{\mathsf T}(g)$ with $\Gamma=(\sqrt{\gamma_1},\gamma_{12}/\sqrt{\gamma_1})$, a vector of algebraic functions of one variable, that is, of elements of an algebraic closure of $\kk(s)$, evaluated at $s=g$. For arbitrary $\Gamma(s)$ and $g(x)$
\begin{equation}\label{eq:rk1Gamma}
\det\Hess\tfrac12\vv{y,\Gamma(g)}^2=-A^3\vv{y,\Gamma'(g)}\,\mathcal W(\Gamma)(g)^2\,\BH(g),\qquad A=\vv{y,\Gamma(g)},\quad\mathcal W(\Gamma)=\Gamma_1\Gamma_2'-\Gamma_1'\Gamma_2 .
\end{equation}
If $\mathcal W(\Gamma)=0$, the direction of $\Gamma$, hence of $f$, is constant, which is excluded. If $\Gamma'=0$, then $d\,ff^{\mathsf T}$ is constant. Otherwise $\BH(g)=0$, so $g=h(\ell)$ by Lemma~\ref{lem:2d}, and $d\,ff^{\mathsf T}\in\kk[\ell]^{2\times2}$.
\end{proof}

Thus the vanishing of the Hessian determinant forces a family of quadratic forms of rank one on the fiber to depend on a single linear coordinate of the base.

\begin{proposition}\label{prop:rk1}%
Let $\phi$ be as in \eqref{eq:blockZ} with $S=d\,ff^{\mathsf T}\in\kk[x_1]^{2\times2}$ of rank one, $f_1/f_2\notin\kk$, and $\det\Hess\phi=\kappa\in\kk^\times$. Then $\grad\phi$ is a polynomial automorphism.
\end{proposition}
\begin{proof}
Factoring in $\kk[x_1]$ we may assume $d,f_1,f_2\in\kk[x_1]$, and $\omega=f_2f_1'-f_1f_2'\neq0$. Put $A=\vv{y,f}$ and let $[y^j]$ denote the part of degree $j$ in $y$. For arbitrary $B,\Lambda$
\begin{equation}\label{eq:rk1layers}
[y^3]\det\Hess\phi=-d^3\omega^2A^2\vv{y,\partial_2^2\Lambda},
\end{equation}
and, when $\partial_2^2\Lambda=0$, $[y^2]\det\Hess\phi=-d^3\omega^2A^2\,\partial_2^2B+R(x_1,y)$ with $R$ independent of $x_2$ and of $B$. Since $\det\Hess\phi$ is constant, the first identity gives $\Lambda=\Lambda_0(x_1)+x_2\Lambda_1(x_1)$, and the second gives $\partial_2^2B=\beta(x_1)$. Hence $\phi$ has degree at most two in $z=(x_2,y_1,y_2)$ with coefficients in $\kk[x_1]$, and its quadratic part is $\frac12z^{\mathsf T}S_zz$ with
\[
S_z=\begin{pmatrix}\beta&\Lambda_1^{\mathsf T}\\ \Lambda_1&d\,ff^{\mathsf T}\end{pmatrix},\qquad\det S_z=-d\,\varepsilon^2,\quad\varepsilon=\Lambda_1^{\mathsf T}Jf .
\]
Apply Lemma~\ref{lem:quadblock}. If $\det S_z=0$, then $\varepsilon=0$ and $(0,Jf)$ spans $\ker S_z$, the rank being two since otherwise $\kappa=0$. A constant kernel vector would then be parallel to $(0,Jf(x_1))$, contrary to $f_1/f_2\notin\kk$. Therefore the second alternative of Lemma~\ref{lem:quadblock} holds.
\end{proof}

\begin{proof}[Proof of Theorem~\ref{thm:block} when $\rk S\le1$]
We may assume \eqref{eq:rk1} with $f_1/f_2\notin\kk$. By \eqref{eq:toplayer} and Proposition~\ref{prop:rk1top}, after a linear change of $x$ one has $d\,ff^{\mathsf T}\in\kk[x_1]^{2\times2}$, and Proposition~\ref{prop:rk1} applies.
\end{proof}

\subsection{Rank two. The Legendre system}
Now $\delta=\det S\neq0$. Let $W=S^{-1}$ and, for $p\in\kk^2$, let $\psi_p=-B+\frac12(p-\Lambda)^{\mathsf T}W(p-\Lambda)$ be the transform of Lemma~\ref{lem:legendre}. Write
\[
\psi_p=W_p-p^{\mathsf T}W\Lambda+u_0,\qquad W_p=\tfrac12p^{\mathsf T}Wp,\quad u_0=\tfrac12\Lambda^{\mathsf T}W\Lambda-B .
\]
By Lemma~\ref{lem:legendre}(c), $\det\Hess_x\psi_p=\kappa/\delta$ for all $p$. Lemma~\ref{lem:det2} and sorting by degree in $p$ give the five layers
\begin{equation}\label{eq:layers}
\begin{aligned}
&(p^4)\ \det\Hess W_p=0,\qquad (p^3)\ \{\Hess W_p,\Hess(p^{\mathsf T}W\Lambda)\}=0,\\
&(p^2)\ \det\Hess(p^{\mathsf T}W\Lambda)+\{\Hess W_p,\Hess u_0\}=0,\qquad (p^1)\ \{\Hess(p^{\mathsf T}W\Lambda),\Hess u_0\}=0,\\
&(p^0)\ \det\Hess u_0=\kappa/\delta .
\end{aligned}
\end{equation}
Thus the partial Legendre transform produces a family of rational functions $\psi_p$ of two variables, indexed by the parameter $p$, whose Hessian determinants all equal $\kappa/\delta$, a function which does not depend on $p$, and the five equations \eqref{eq:layers} are the conditions for this family to be compatible with one determinant equation.

\begin{proof}[Proof of Theorem~\ref{thm:block} when $\delta\in\kk^\times$]
Then $W$ is polynomial and $\psi_p\in\kk[x]$ has $\det\Hess\psi_p=\kappa/\delta\in\kk^\times$ for every $p$, so $\grad\psi_p$ is injective by \ref{E:dB}. If $\grad\phi(x,y)=\grad\phi(x',y')$, then $p=\Lambda(x)+S(x)y=\Lambda(x')+S(x')y'$ and, by Lemma~\ref{lem:legendre}(a), $\grad\psi_p(x)=\grad\psi_p(x')$. Hence $x=x'$ and $y=y'=W(x)(p-\Lambda(x))$, and Lemma~\ref{lem:descent} applies.
\end{proof}

From now on $\delta\notin\kk$.

\begin{lemma}\label{lem:Wpole}
$W$ has a pole along every irreducible component $L$ of $\{\delta=0\}$, and so does $W_p$ for generic $p$.
\end{lemma}
\begin{proof}
$\det W=1/\delta$ has a pole along $L$, so some entry of $W$ does. Write the Laurent expansion $W=\sum_{j\ge-e}W_j\ell_L^j$ along $L$ with $e\ge1$ and $W_{-e}\not\equiv0$ on $L$. Then $p^{\mathsf T}W_{-e}p\not\equiv0$ on $L$ for generic $p$.
\end{proof}

\begin{proposition}\label{prop:types}%
Assume $\delta\notin\kk$ and the layer $(p^4)$. After a linear change of $x$, a linear change of $y$ and a rescaling, exactly one of the following holds.
\begin{enumerate}[label=\textup{(T\arabic*)}]
\item $S=S(x_1)$.
\item $S=\begin{psmallmatrix}a(x_1)-x_2b(x_1)^2&b(x_1)\\ b(x_1)&0\end{psmallmatrix}$ with $b\notin\kk$, so that $\delta=-b^2$.
\item $W=W_{\mathrm{aff}}(x)-\frac k{x_1}vv^{\mathsf T}$ with $v=(1,-x_2)^{\mathsf T}$, $k\in\kk^\times$ and $W_{\mathrm{aff}}$ a matrix of affine functions. Moreover $S=c_0^{-1}\adj(x_1W)$ and $\delta=x_1/c_0$ with $c_0\in\kk^\times$.
\end{enumerate}
\end{proposition}
\begin{proof}
For $\xi\in\operatorname{Sym}_2(\kk)$ put $\vv{W,\xi}=\xi_{11}w_{11}+2\xi_{12}w_{12}+\xi_{22}w_{22}$ and $Q(\xi)=\det\Hess_x\vv{W,\xi}$, a quadratic form on the three-dimensional space $\operatorname{Sym}_2(\kk)$ with coefficients in $\kk(x)$. The layer $(p^4)$ says $Q(pp^{\mathsf T})=0$ for all $p$, so $Q$ vanishes on the irreducible conic $\{\det\xi=0\}$ and $Q=\lambda(x)\det\xi$.

\emph{Step 1. $\lambda=0$.} Then every element of $V=\operatorname{span}_\kk(w_{11},w_{12},w_{22})$ has zero Hessian determinant. Not all of them are affine, since otherwise $W$ would be polynomial and $\det W\cdot\delta=1$ would force $\delta\in\kk$. If $f,g\in V$ have $\Hess f=\rho_1n_1n_1^{\mathsf T}\neq0$ and $\Hess g=\rho_2n_2n_2^{\mathsf T}\neq0$, then $0=\det(\Hess f+\Hess g)=\rho_1\rho_2\det(n_1,n_2)^2$. So all non-affine elements of $V$ have the same generators (Lemma~\ref{lem:devel}(a)). By Lemmas~\ref{lem:devel}(b) and~\ref{lem:Wpole} applied to $W_p\in V$ for generic $p$ there are two cases.

(1a) The generators are parallel. Then every component of $\{\delta=0\}$ is a line parallel to them. Choose coordinates in which the generators are the lines $x_1=\mathrm{const}$. Then $\delta\in\kk[x_1]$, every entry of $W$ is affine along the generators, $W=W_0(x_1)+x_2W_1(x_1)$, and every $\Hess w_{ij}$ is a multiple of $e_1e_1^{\mathsf T}$, so $W_1=C$ is constant. Since $\det W=1/\delta$ does not depend on $x_2$, $\det C=0$ and $\{W_0,C\}=0$. If $C=0$, this is (T1). Otherwise $C=\sigma kk^{\mathsf T}$ with $k^{\mathsf T}\adj(W_0)k=0$. After a linear change of $y$, $k=e_2$, so $(W_0)_{11}=0$. Then $S=\delta\adj(W)=\delta\adj W_0+x_2\sigma\delta e_1e_1^{\mathsf T}$ and $S_0=\delta\adj W_0=\begin{psmallmatrix}a&b\\ b&0\end{psmallmatrix}$ with $a,b\in\kk[x_1]$. Hence $\delta=\det S=-b^2$, $S=S_0-x_2\sigma b^2e_1e_1^{\mathsf T}$, and a rescaling of $x_2$ gives (T2), with $b\notin\kk$ because $\delta\notin\kk$.

(1b) The generators pass through one point $x^*$. By Lemma~\ref{lem:devel}(c) applied to every element of $V$, $W=C+H(x-x^*)$ with $C$ constant and $H$ rational and homogeneous of degree one. Expand $S=\sum_{j\ge0}S_j$ into homogeneous components in $w=x-x^*$. Then $S(C+H)=I$ gives $S_0C=I$ and $S_jC+S_{j-1}H=0$, so $S_j=(-1)^jC^{-1}(HC^{-1})^j$. As $S_j=0$ for large $j$, $HC^{-1}$ is nilpotent, and $\delta=\det C^{-1}\det(I+HC^{-1})^{-1}=\det C^{-1}\in\kk$, a contradiction.

\emph{Step 2. $\lambda\neq0$.} Write $H_{ij}=\Hess w_{ij}$. After a generic linear change of $y$, $H_{11}\neq0\neq H_{22}$, since otherwise $\Hess(p^{\mathsf T}Wp)=0$ for all $p$ and $W$ would be affine. Comparing coefficients in $Q(\xi)=\lambda\det\xi$ gives
\[
\det H_{11}=\det H_{22}=0,\qquad\{H_{11},H_{12}\}=\{H_{12},H_{22}\}=0,\qquad\lambda=\{H_{11},H_{22}\}=-4\det H_{12}.
\]
Write $H_{ii}=\rho_ia_ia_i^{\mathsf T}$ with $a_i\in\kk(x)^2$, $\rho_i\in\kk(x)$. Then $\lambda=\rho_1\rho_2\det(a_1,a_2)^2\neq0$, so $a_1,a_2$ are independent and $\det H_{12}\neq0$. By Lemma~\ref{lem:det2} the conditions $\{H_{ii},H_{12}\}=0$ say that $Ja_1$ and $Ja_2$ are isotropic for the nondegenerate form $H_{12}$. Hence $H_{12}=\sigma(a_1a_2^{\mathsf T}+a_2a_1^{\mathsf T})$ with $\sigma\in\kk(x)$ and $4\sigma^2=\rho_1\rho_2$. It follows that
\[
\Hess_x(p^{\mathsf T}Wp)=\sum_{i,j}p_ip_jH_{ij}=\rho_1^{-1}(Ap)(Ap)^{\mathsf T},\qquad A=(\rho_1a_1\ \ 2\sigma a_2)\in\mathrm{GL}_2(\kk(x)).
\]
So the generators of $W_p$ have direction $JA(x)p$. Apply Lemma~\ref{lem:devel}(b) to $W_p$ for generic $p$.

(2a) For generic $p$ the generators of $W_p$ are parallel. Then $A(x)p=s(x,p)v(p)$. For generic $p,p'$ the vectors $v(p),v(p')$ are independent, since $\rk A=2$, and $A(p+p')=sv(p)+s'v(p')$ also has constant direction. So $s/s'$ does not depend on $x$, and $A(x)=\rho(x)A_0$. After the linear change $p\mapsto A_0^{-1}p$ of $y$ one gets $\Hess(p^{\mathsf T}Wp)=\nu\,pp^{\mathsf T}$, that is $\Hess w_{11}=\nu e_1e_1^{\mathsf T}$, $\Hess w_{22}=\nu e_2e_2^{\mathsf T}$ and $2\Hess w_{12}=\nu(e_1e_2^{\mathsf T}+e_2e_1^{\mathsf T})$. The first equation gives $\partial_2\partial_jw_{11}=0$ and $\nu=\partial_1^2w_{11}$, hence $\partial_2\nu=0$, and the second gives $\partial_1\nu=0$ in the same way. So $\nu$ is constant, and the three equations make $w_{11},w_{22},w_{12}$ polynomials of degree at most two. So $W$ is polynomial and $\delta\in\kk$, a contradiction.

(2b) For generic $p$ the generators of $W_p$ pass through a point $x^*(p)$. By Lemma~\ref{lem:devel}(b) every polar curve of $W_p$ contains $x^*(p)$, and by Lemma~\ref{lem:devel}(c) it is a line through $x^*(p)$. By Lemma~\ref{lem:Wpole} the polar curves of $W_p$ include every component of $\{\delta=0\}$. If $\{\delta=0\}$ contains two distinct lines, or if $x^*$ does not depend on $p$, then $x^*$ is fixed and the argument of (1b) gives a contradiction. Hence $\{\delta=0\}$ is a single line, say $\{x_1=0\}$, and $x^*(p)=(0,t(p))$ with $t$ nonconstant, rational and homogeneous of degree zero. The generator of $W_p$ through $x$ is the line through $x$ and $(0,t(p))$, so $A(x)p=s(x,p)(t(p)-x_2,x_1)$. Linearity in $p$ forces $s=\ell(x)^{\mathsf T}p/x_1$ and $(\ell(x)^{\mathsf T}p)\,t(p)=m(x)^{\mathsf T}p$. Since $t(p)=m(x)^{\mathsf T}p/\ell(x)^{\mathsf T}p$ is nonconstant and independent of $x$, the pair $(m(x),\ell(x))$ is proportional to a fixed pair of independent vectors. After a linear change of $y$ one gets $t=p_1/p_2$ and
\[
\Hess_x(p^{\mathsf T}Wp)=\nu(x)\,\eta\eta^{\mathsf T},\qquad\eta=(p_1-x_2p_2,\ x_1p_2)^{\mathsf T}.
\]
Comparing the coefficients of $p_1^2$, $p_1p_2$, $p_2^2$ gives $\Hess w_{11}=\nu e_1e_1^{\mathsf T}$, $2\Hess w_{12}=\nu\bigl(e_1(-x_2,x_1)+(-x_2,x_1)^{\mathsf T}e_1^{\mathsf T}\bigr)$ and $\Hess w_{22}=\nu(-x_2,x_1)^{\mathsf T}(-x_2,x_1)$. From the first, $w_{11}=f(x_1)+c_1x_2$ with $f''=\nu$. From the second, $w_{12}=g(x_1)+x_2g_1(x_1)$ with $g_1'=\frac12x_1\nu$, $g''=0$ and $g_1''=-\nu$. Therefore $x_1\nu'=-3\nu$ and $\nu=-2k/x_1^3$ with $k\in\kk^\times$. Integrating, $f=-k/x_1+$ affine and $g_1=k/x_1+$ const. The third equation gives $w_{22}=-kx_2^2/x_1+$ affine. Altogether $W=W_{\mathrm{aff}}-\frac k{x_1}vv^{\mathsf T}$ with $v=(1,-x_2)^{\mathsf T}$. Put $Q_1=x_1W=x_1W_{\mathrm{aff}}-kvv^{\mathsf T}\in\kk[x]^{2\times2}$. Then $\det Q_1|_{x_1=0}=\det(-kvv^{\mathsf T})=0$, so $\det Q_1=x_1D$ with $D\in\kk[x]$, and $\delta=x_1^2/\det Q_1=x_1/D$. Since $\delta$ is a nonconstant polynomial, $D=c_0\in\kk^\times$. Then $S=x_1Q_1^{-1}=c_0^{-1}\adj Q_1$, which is (T3).
\end{proof}

\begin{proposition}\label{prop:T1}%
In type \textup{(T1)}, $\grad\phi$ is a polynomial automorphism or $D_u^2\phi=0$ for some $u\neq0$.
\end{proposition}
\begin{proof}
Here $W=W(x_1)$, so $\Hess W_p=\operatorname{diag}(\frac12p^{\mathsf T}W''p,0)$, and $W''\neq0$, since otherwise $W$ would be affine, hence polynomial, and $\delta\in\kk$. For every $F$, $\{\Hess W_p,\Hess F\}=\frac12(p^{\mathsf T}W''p)\,\partial_2^2F$. The layer $(p^3)$ gives $\partial_2^2(W\Lambda)=0$, hence $\partial_2^2\Lambda=0$ and $\Lambda=\Lambda_0(x_1)+x_2\Lambda_1(x_1)$ with $\Lambda_i\in\kk[x_1]^2$. Then $p^{\mathsf T}W\Lambda=f_0(x_1)+x_2f_1(x_1)$ with $f_i=p^{\mathsf T}W\Lambda_i$, and $\det\Hess(p^{\mathsf T}W\Lambda)=-(f_1')^2$. The layer $(p^2)$ reads $\frac12(p^{\mathsf T}W''p)\,\partial_2^2u_0=\bigl(p^{\mathsf T}(W\Lambda_1)'\bigr)^2$. The right side and $W''$ depend only on $x_1$, so $\partial_2^2u_0=h(x_1)$. Thus $u_0$, and with it $B=\frac12\Lambda^{\mathsf T}W\Lambda-u_0$, is quadratic in $x_2$. Consequently $\phi$ has degree at most two in $z=(x_2,y_1,y_2)$ with coefficients in $\kk[x_1]$, and Lemma~\ref{lem:quadblock} applies.
\end{proof}

\begin{remark}
Type \textup{(T1)} occurs. For $\phi=-\frac\kappa2x_1^2-x_2y_1+\frac12x_1y_1^2+\frac12y_2^2$ one has $\det\Hess\phi=\kappa$ and $S=\operatorname{diag}(x_1,1)$, $\delta=x_1$. Lemma~\ref{lem:quadblock} applies with $z=(x_2,y_1,y_2)$, $S_z=\begin{psmallmatrix}0&-1&0\\-1&x_1&0\\0&0&1\end{psmallmatrix}$, $\det S_z=-1$, and $\grad\phi$ is an automorphism.
\end{remark}

\begin{proposition}\label{prop:T2}%
There is no $\phi$ of type \textup{(T2)} with $\det\Hess\phi\in\kk^\times$.
\end{proposition}
\begin{proof}
Here $W=\begin{psmallmatrix}0&1/b\\ 1/b&x_2-a/b^2\end{psmallmatrix}$ and $W_p=p_1p_2/b+\frac12p_2^2(x_2-a/b^2)$, so $\Hess W_p=\omega e_1e_1^{\mathsf T}$ with $\omega=p_1p_2(1/b)''-\frac12p_2^2(a/b^2)''$. Since $b\notin\kk$, $(1/b)''\neq0$. For every $F$, $\{\Hess W_p,\Hess F\}=\omega\,\partial_2^2F$. The layer $(p^3)$ gives $W\Lambda=m_0(x_1)+x_2m_1(x_1)$, and the layer $(p^2)$ reads $\omega\,\partial_2^2u_0=(p^{\mathsf T}m_1')^2$. Comparing the coefficients of $p_1^2$, $p_1p_2$ and $p_2^2$, and using that $\omega$ has no term in $p_1^2$, gives $m_1'=0$ and $\partial_2^2u_0=0$. So $m_1=(d_1,d_2)$ is constant and $u_0=u_{00}(x_1)+x_2u_{01}(x_1)$. The layer $(p^1)$ then holds, and the layer $(p^0)$ reads $-(u_{01}')^2=\det\Hess u_0=\kappa/\delta=-\kappa/b^2$.

We show $u_{01}\in\kk[x_1]$. Write $S=S_0-x_2b^2e_1e_1^{\mathsf T}$ with $S_0=\begin{psmallmatrix}a&b\\ b&0\end{psmallmatrix}$. Since $\Lambda=S(m_0+x_2m_1)$ is polynomial, $S_0m_0=(am_{0,1}+bm_{0,2},\,bm_{0,1})$ is polynomial, in particular $bm_{0,1}\in\kk[x_1]$. The coefficient of $x_2$ in $u_0=\frac12(m_0+x_2m_1)^{\mathsf T}S(m_0+x_2m_1)-B$ is $u_{01}=m_1^{\mathsf T}S_0m_0-\frac12(bm_{0,1})^2-B_1$, where $B_1\in\kk[x_1]$ is the coefficient of $x_2$ in $B$. So $u_{01}$ is a polynomial, and $\kappa/b^2=(u_{01}')^2\in\kk[x_1]$ forces $b\in\kk^\times$, a contradiction.
\end{proof}

\begin{proposition}\label{prop:T3}%
There is no $\phi$ of type \textup{(T3)} with $\det\Hess\phi\in\kk^\times$.
\end{proposition}
\begin{proof}
Put $\tau=v^{\mathsf T}p=p_1-x_2p_2$. The affine part $W_{\mathrm{aff}}$ does not contribute to Hessians, so $W_p=\frac12p^{\mathsf T}W_{\mathrm{aff}}p-k\tau^2/(2x_1)$ and
\[
\Hess W_p=-\frac k{x_1^3}\eta\eta^{\mathsf T},\qquad\eta=(\tau,\ x_1p_2)^{\mathsf T},
\]
with generators through $(0,t)$, $t=p_1/p_2$, of direction $J\eta\parallel(x_1,x_2-t)$.

\emph{The layer $(p^3)$.} Write $p^{\mathsf T}W\Lambda=p_1g_1+p_2g_2=p_2(tg_1+g_2)$ with $g=W\Lambda$, and let $D_t=x_1^2\partial_1^2+2x_1(x_2-t)\partial_1\partial_2+(x_2-t)^2\partial_2^2$. The layer is $D_t(tg_1+g_2)=0$. Writing $D_t=D_0-2t(x_1\partial_1\partial_2+x_2\partial_2^2)+t^2\partial_2^2$, this is an identity of degree three in $t$ whose coefficients give
\[
\partial_2^2g_1=0,\quad\partial_2^2g_2=2(x_1\partial_1\partial_2+x_2\partial_2^2)g_1,\quad D_0g_1=2(x_1\partial_1\partial_2+x_2\partial_2^2)g_2,\quad D_0g_2=0 .
\]
The first equation gives $g_1=g_{10}(x_1)+x_2g_{11}(x_1)$, and the second $g_2=x_1g_{11}'x_2^2+x_2g_{21}(x_1)+g_{20}(x_1)$. Separating powers of $x_2$, the third equation gives $x_1g_{11}''+2g_{11}'=0$ and $x_1g_{10}''=2g_{21}'$, and the fourth gives $g_{20}''=0$, $x_1g_{21}''+2g_{21}'=0$, and an equation for $g_{11}$ which follows from the third. Hence $g_{11}=c_1-c/x_1$, $g_{21}=c_2-c'/x_1$, $g_{10}=c'/x_1+a_0+a_1x_1$, and $g_{20}$ is affine, with constants $c,c',c_i,a_i$. This gives
\[
p^{\mathsf T}W\Lambda=p^{\mathsf T}m_{\mathrm{aff}}(x)+\frac{\sigma\tau}{x_1},\qquad\sigma=c'-cx_2,
\]
with $m_{\mathrm{aff}}$ affine.

\emph{The layer $(p^2)$.} For $N=N(x_2)$ one has $\det\Hess(N/x_1)=(2NN''-N'^2)/x_1^4$, and with $N=\sigma\tau$ this gives $\det\Hess(p^{\mathsf T}W\Lambda)=-(cp_1-c'p_2)^2/x_1^4$. By Lemma~\ref{lem:devel}(a), $\{\Hess W_p,\Hess u_0\}=-\frac k{x_1^3}(J\eta)^{\mathsf T}\Hess u_0(J\eta)$. The coefficients of $p_1^2$, $p_1p_2$ and $p_2^2$ in the layer $(p^2)$ form a triangular linear system for the entries $u_{0,22},u_{0,12},u_{0,11}$ of $\Hess u_0$ with nonzero diagonal $(-kx_1^{-3},2kx_1^{-2},-kx_1^{-1})$. Its unique solution is
\[
\Hess u_0=-\frac1{kx_1^3}\begin{pmatrix}\sigma^2&c\sigma x_1\\ c\sigma x_1&c^2x_1^2\end{pmatrix},
\]
a matrix of rank at most one. The layer $(p^0)$ gives $0=\det\Hess u_0=\kappa/\delta=\kappa c_0/x_1$, a contradiction.
\end{proof}

\begin{remark}
Type \textup{(T3)} satisfies the layer $(p^4)$. For $W_{\mathrm{aff}}=\operatorname{diag}(0,1)$ and $k=1$ one gets $S=\begin{psmallmatrix}x_2^2-x_1&x_2\\ x_2&1\end{psmallmatrix}$, with $\det S=-x_1$ and $\det\Hess(\frac12y^{\mathsf T}Sy)=0$. This $S$ neither depends on a single linear form nor is affine in any direction. It is excluded only by the lower layers.
\end{remark}

\begin{proof}[Proof of Theorem~\ref{thm:block}]
Combine the case $\rk S\le1$, the case $\delta\in\kk^\times$, Proposition~\ref{prop:types} and Propositions~\ref{prop:T1}, \ref{prop:T2} and~\ref{prop:T3}.
\end{proof}

\section{One quadratic variable with nonconstant coefficient}\label{sec:quadratic}

Here $\phi$ is quadratic in one variable $t$, and the coefficient $c$ of $\frac12t^2$ is not constant. By Theorem~\ref{thm:lead}, $c$ depends on two variables. The partial Legendre transform in $t$ gives a family $\psi_s$ of potentials in three variables, and the vanishing of the coefficients of $\det\Hess\psi_s$ in $s$ shows that, in suitable linear coordinates, $\phi$ is quadratic in a second variable as well. Then Section~\ref{sec:block} or Lemma~\ref{lem:quadblock} applies.

\begin{theorem}\label{thm:quad}
Let $x=(x_1,x_2,x_3)$ and $\phi=a(x)+b(x)t+\frac12c(x)t^2$ with $a,b,c\in\kk[x]$, $c\notin\kk$, and $\det\Hess\phi=\kappa\in\kk^\times$. Then $\grad\phi$ is a polynomial automorphism or $D_u^2\phi=0$ for some $u\in\kk^4\setminus\{0\}$.
\end{theorem}

By Theorem~\ref{thm:lead} with $m=2$, after a linear change of $x$ one has $c=c(x_1,x_2)$. Apply Lemma~\ref{lem:legendre} in the variable $t$ ($n=3$, $m=1$),
\begin{equation}\label{eq:psis}
\psi_s(x)=\frac{(s-b)^2}{2c}-a=\tfrac12s^2W-s\beta+u_0,\qquad W=\frac1c,\quad\beta=\frac bc,\quad u_0=\frac{b^2}{2c}-a,
\end{equation}
so that $\det\Hess_x\psi_s=-\kappa/c$ for all $s$. Let $L_k$ be the coefficient of $s^k$ in $\det\Hess_x\psi_s$. Then $L_k=0$ for $1\le k\le6$ and $L_0=\det\Hess u_0=-\kappa/c$.

\begin{lemma}\label{lem:blockreach}
Suppose that, after a linear change of $x$ into coordinates $(x',y)$ with $y$ one or two of the new coordinates, $c$ depends only on $x'$, $\beta$ is affine in $y$ and $u_0$ is quadratic in $y$, with coefficients rational in $x'$. Then $\phi$ has degree at most two in $(y,t)$ with coefficients in $\kk[x']$.
\end{lemma}
\begin{proof}
$b=c\beta$ and $a=\frac12c\beta^2-u_0$, and the coefficients of the expansion of the polynomial $\phi$ in $(y,t)$ are polynomials in $x'$.
\end{proof}

If $y$ consists of one coordinate, $\phi$ has the form \eqref{eq:blockZ} in the variables $(x';y,t)$ and Theorem~\ref{thm:block} applies. If $y$ consists of two coordinates, $x'$ is a single variable and Lemma~\ref{lem:quadblock} applies. In both cases the conclusion of Theorem~\ref{thm:quad} follows, and it suffices to reach one of the two situations.

\subsection{\texorpdfstring{Case A, $\det\Hess_{(x_1,x_2)}W\neq0$}{Case A}}
\begin{proposition}\label{prop:caseA}
If $\det\Hess_{(x_1,x_2)}W\neq0$, then $\beta$ is affine and $u_0$ is quadratic in $x_3$.
\end{proposition}
\begin{proof}
Since $W=W(x_1,x_2)$, the matrix $\Hess_xW$ has zero third row and column. A direct expansion gives $L_6=0$ and $L_5=-\frac14\det\Hess_{(x_1,x_2)}W\cdot\partial_3^2\beta$. Hence $\partial_3^2b=c\,\partial_3^2\beta=0$, that is $\beta=\beta_0+x_3\beta_1$ with $\beta_i\in\kk(x_1,x_2)$. Then
\[
L_4=\tfrac14\det\Hess_{(x_1,x_2)}W\cdot\partial_3^2u_0-\tfrac12\grad\beta_1^{\mathsf T}\adj\bigl(\Hess_{(x_1,x_2)}W\bigr)\grad\beta_1,
\]
so $\partial_3^2u_0$ depends only on $(x_1,x_2)$.
\end{proof}

\subsection{\texorpdfstring{Case B, $\det\Hess_{(x_1,x_2)}W=0$}{Case B}}
\begin{lemma}\label{lem:caseB}
If $\det\Hess_{(x_1,x_2)}W=0$, then after a linear change of $(x_1,x_2)$ one has $c=c(x_1)$.
\end{lemma}
\begin{proof}
In two variables $c^5\det\Hess(1/c)=c\det\Hess c-2\BH(c)$, so the hypothesis says $c\det\Hess c-2\BH(c)=0$, and Theorem~\ref{thm:gamma} gives $c=f(\ell)$ for a linear form $\ell$.
\end{proof}

From now on $c=c(x_1)$, $W=W(x_1)$ and $W''\neq0$, primes denoting $\partial_{x_1}$. Write $y=(x_2,x_3)$, $\beta_1=\partial_1\beta$, $u_{0,1}=\partial_1u_0$, $\beta_{11}=\partial_1^2\beta$, $u_{0,11}=\partial_1^2u_0$, and let $\{\cdot,\cdot\}$ be as in Lemma~\ref{lem:det2} for $2\times2$ Hessians in $y$. The Schur complement of the entry $(x_1,x_1)$ of $\Hess_x\psi_s$ gives, for arbitrary $\beta,u_0$,
\begin{equation}\label{eq:Llayers}
\begin{aligned}
L_6&=L_5=0,\qquad L_4=\tfrac12W''\det\Hess_y\beta,\\
L_3&=-\tfrac12W''\{\Hess_y\beta,\Hess_yu_0\}-\beta_{11}\det\Hess_y\beta+\grad_y\beta_1^{\mathsf T}\adj(\Hess_y\beta)\grad_y\beta_1,\\
L_2&=\tfrac12W''\det\Hess_yu_0+\beta_{11}\{\Hess_y\beta,\Hess_yu_0\}+u_{0,11}\det\Hess_y\beta\\
&\quad-\grad_y\beta_1^{\mathsf T}\adj(\Hess_yu_0)\grad_y\beta_1-2\grad_yu_{0,1}^{\mathsf T}\adj(\Hess_y\beta)\grad_y\beta_1,\\
L_1&=-\beta_{11}\det\Hess_yu_0-u_{0,11}\{\Hess_y\beta,\Hess_yu_0\}+2\grad_yu_{0,1}^{\mathsf T}\adj(\Hess_yu_0)\grad_y\beta_1\\&\quad+\grad_yu_{0,1}^{\mathsf T}\adj(\Hess_y\beta)\grad_yu_{0,1},\\
L_0&=u_{0,11}\det\Hess_yu_0-\grad_yu_{0,1}^{\mathsf T}\adj(\Hess_yu_0)\grad_yu_{0,1}=\det\Hess_xu_0 .
\end{aligned}
\end{equation}
Since $L_4=0$, $\beta\in\kk(x_1)[x_2,x_3]$ has zero Hessian determinant in $y$. By Lemma~\ref{lem:hess2F} over $\kk(x_1)$, after possibly exchanging $x_2$ and $x_3$,
\begin{equation}\label{eq:betaform}
\beta=Q(x_1,\ell)+m_1+m_2x_2+m_3x_3,\qquad\ell=x_2+\alpha(x_1)x_3,\quad Q\in\kk(x_1)[\ell],\ \alpha,m_i\in\kk(x_1),
\end{equation}
and $\Hess_y\beta=Q_{\ell\ell}\grad_y\ell\grad_y\ell^{\mathsf T}$.

\begin{proposition}\label{prop:Qll}
If $Q_{\ell\ell}\neq0$, then $\alpha\in\kk$, and in the coordinates $(x_1,\ell,x_3)$ the function $\beta$ is affine and $u_0$ is quadratic in $x_3$.
\end{proposition}
\begin{proof}
Put $v=(\alpha,-1)^{\mathsf T}$, orthogonal to $\grad_y\ell=(1,\alpha)^{\mathsf T}$, so that $\adj(\Hess_y\beta)=Q_{\ell\ell}vv^{\mathsf T}$, and put $G=v^{\mathsf T}\grad_y\beta_1=\alpha m_2'-m_3'-\alpha'Q_\ell$. Then $\{\Hess_y\beta,\Hess_yu_0\}=Q_{\ell\ell}\,\partial_v^2u_0$, and $L_3=0$ reads $Q_{\ell\ell}\bigl(W''\partial_v^2u_0-2G^2\bigr)=0$. Hence $\partial_v^2u_0=2G^2/W''$. In the coordinates $(x_1,\ell,x_3)$ one has $\partial_v=-\partial_{x_3}$ at fixed $(x_1,\ell)$, and integrating twice gives
\[
u_0=x_3^2\frac{G^2}{W''}+x_3U_1(x_1,\ell)+U_0(x_1,\ell),\qquad U_i\in\kk(x_1)[\ell].
\]
Substituting into $L_2$ and passing to the coordinates $(x_1,\ell,x_3)$ gives a polynomial of degree two in $x_3$ with leading coefficient $-12G^2Q_{\ell\ell}^2\alpha'^2/W''$. Since $L_2=0$, $Q_{\ell\ell}\neq0$ and $W''\neq0$, this gives $G\alpha'=0$ in $\kk(x_1)[\ell]$. Suppose $\alpha'\neq0$. Then $G=0$, that is $\alpha'Q_\ell=\alpha m_2'-m_3'$. Here $\alpha,m_2,m_3\in\kk(x_1)$, so $Q_\ell\in\kk(x_1)$ does not depend on $\ell$, and $Q_{\ell\ell}=0$, a contradiction. Hence $\alpha\in\kk$, and $(x_1,\ell,x_3)$ are linear coordinates. In these coordinates $\beta=Q(x_1,\ell)+m_1+m_2\ell+(m_3-\alpha m_2)x_3$ is affine in $x_3$. Since $G\in\kk(x_1)[\ell]$, the formula above shows that $u_0$ is quadratic in $x_3$.%
\end{proof}

\begin{proposition}\label{prop:Qll0}
If $Q_{\ell\ell}=0$, that is $\beta=m_1+m_2x_2+m_3x_3$, then $\beta$ is affine in $y$, and either $u_0$ is quadratic in $y=(x_2,x_3)$, or there are independent linear forms $\tilde\ell,\eta$ in $y$ with constant coefficients such that in the coordinates $(x_1,\tilde\ell,\eta)$ the function $\beta$ is affine and $u_0$ is quadratic in $\eta$.
\end{proposition}
\begin{proof}
Now $\Hess_y\beta=0$ and $L_4=L_3=0$. A generic linear change of $(x_2,x_3)$ with constant coefficients preserves the assumptions, so we may assume $m_2'\neq0$ whenever $(m_2',m_3')\neq0$, and that no exchange of $x_2,x_3$ is needed in Lemma~\ref{lem:hess2F}. Put $g=\grad_y\beta_1=(m_2',m_3')$ and $g^\perp=Jg$. By Lemma~\ref{lem:det2}, \eqref{eq:Llayers} reduces to
\[
L_2=\tfrac12W''\det\Hess_yu_0-\partial_{g^\perp}^2u_0,\qquad L_0=\det\Hess_xu_0,\]\[L_1=-\beta_{11}\det\Hess_yu_0+2g^{\mathsf T}\adj(\Hess_yu_0)\grad_yu_{0,1}.
\]
\emph{Sub-case $g=0$}, that is $m_2,m_3\in\kk$. Then $L_2=0$ gives $\det\Hess_yu_0=0$, and Lemma~\ref{lem:hess2F} over $\kk(x_1)$ gives $u_0=R(x_1,\tilde\ell)+p_1+p_2x_2+p_3x_3$ with $\tilde\ell=x_2+\tilde\alpha(x_1)x_3$, $R\in\kk(x_1)[\tilde\ell]$, $\tilde\alpha,p_i\in\kk(x_1)$. For such $u_0$
\[
\det\Hess_xu_0=-R_{\tilde\ell\tilde\ell}\bigl(D-\tilde\alpha'R_{\tilde\ell}\bigr)^2,\qquad D=\tilde\alpha p_2'-p_3' .
\]
So $L_0=-\kappa/c$ reads $R_{\tilde\ell\tilde\ell}(D-\tilde\alpha'R_{\tilde\ell})^2=\kappa/c$. If $\tilde\alpha'\neq0$, then $\partial_{\tilde\ell}(D-\tilde\alpha'R_{\tilde\ell})^3=-3\tilde\alpha'\kappa/c$, so the cube of a polynomial in $\tilde\ell$ is affine in $\tilde\ell$, the polynomial is independent of $\tilde\ell$, $R_{\tilde\ell\tilde\ell}=0$, and the left side vanishes, a contradiction. Hence $\tilde\alpha\in\kk$ and $R_{\tilde\ell\tilde\ell}=\kappa/(cD^2)$ depends only on $x_1$, so $u_0$ is quadratic in $y$.

\emph{Sub-case $g\neq0$.} Put $\tilde u=u_0-(g^{\mathsf T}y)^2/W''$. Then $\det\Hess_y\tilde u=\det\Hess_yu_0-\frac2{W''}\partial_{g^\perp}^2u_0=\frac2{W''}L_2=0$, so
\begin{equation}\label{eq:u0g}
u_0=R(x_1,\tilde\ell)+\frac{(m_2'x_2+m_3'x_3)^2}{W''}+p_1+p_2x_2+p_3x_3,\qquad\tilde\ell=x_2+\tilde\alpha(x_1)x_3 .
\end{equation}
If $R_{\tilde\ell\tilde\ell}$ depends only on $x_1$, $u_0$ is quadratic in $y$. Assume $R_{\tilde\ell\tilde\ell\tilde\ell}\neq0$ and put $\Delta=\tilde\alpha m_2'-m_3'$, $\tilde\Delta=\tilde\alpha m_2''-m_3''$ and $\mathcal W=m_2'm_3''-m_2''m_3'$. For $u_0$ as in \eqref{eq:u0g}, in the coordinates $(x_1,\tilde\ell,x_3)$, $L_2=0$ identically, $L_1$ has degree at most one in $x_3$ and $L_0$ degree at most two, with
\[
[x_3^1]L_1=-\frac{2\Delta^2}{W''^2}R_{\tilde\ell\tilde\ell}\bigl(3W''\tilde\Delta-2W'''\Delta\bigr),\qquad[x_3^2]L_0=-\frac{2\Delta^2}{W''^3}\bigl(\mathcal A\,R_{\tilde\ell\tilde\ell}+4\mathcal W^2\bigr),
\]
where, with $\Delta_3=\tilde\alpha m_2'''-m_3'''$,
\[
\mathcal A=W''''\Delta^2-4W'''\Delta\tilde\Delta+W''\bigl(6\tilde\Delta^2-2\Delta\Delta_3-4\tilde\alpha'\mathcal W\bigr)\in\kk(x_1).
\]%
These formulas, and those below, follow by expanding the Schur complement \eqref{eq:Llayers} for the given $\beta$ and $u_0$.
Since $L_1=0$ and $L_0=-\kappa/c$ does not depend on $x_3$, both coefficients vanish.

(i) $\Delta=0$. Then $L_1=0$ automatically, and
\[
L_0=-\frac{\bigl(W''R_{\tilde\ell\tilde\ell}+2m_2'^2\bigr)T^2}{W''^3},\qquad T=W''\tilde\alpha'R_{\tilde\ell}+2\tilde\alpha'm_2'^2\tilde\ell+W''(p_3'-\tilde\alpha p_2').
\]
Since $\partial_{\tilde\ell}T=\tilde\alpha'(W''R_{\tilde\ell\tilde\ell}+2m_2'^2)$, if $\tilde\alpha'\neq0$ then $\partial_{\tilde\ell}(T^3)=3\tilde\alpha'\kappa W''^3/c$. As before $T$ is then independent of $\tilde\ell$, so $W''R_{\tilde\ell\tilde\ell}+2m_2'^2=0$ and $L_0=0$, a contradiction. Hence $\tilde\alpha\in\kk$, $T\in\kk(x_1)$, and $L_0=-\kappa/c$ shows that $R_{\tilde\ell\tilde\ell}$ depends only on $x_1$, contrary to the assumption.

(ii) $\Delta\neq0$, $\mathcal A\neq0$. Then $[x_3^2]L_0=0$ gives $R_{\tilde\ell\tilde\ell}=-4\mathcal W^2/\mathcal A$, which depends only on $x_1$, contrary to the assumption.

(iii) $\Delta\neq0$, $\mathcal A=0$. Then $[x_3^2]L_0=0$ gives $\mathcal W=0$. As $m_2'\neq0$, $\mathcal W=m_2'^2(m_3'/m_2')'=0$ gives $m_3=\lambda_0m_2+c_3$ with $\lambda_0,c_3\in\kk$. Then $\Delta=(\tilde\alpha-\lambda_0)m_2'\neq0$ and $\tilde\Delta=(\tilde\alpha-\lambda_0)m_2''$, so $[x_3^1]L_1=0$ gives $3W''m_2''=2W'''m_2'$. Under these relations
\[
[x_3^0]L_1=\frac{2\Delta R_{\tilde\ell\tilde\ell}}{W''^2}\bigl(F(x_1)-W''^2\tilde\alpha'R_{\tilde\ell}\bigr)
\]
with $F=W''\bigl(W''(\tilde\alpha p_2'-p_3')+(\lambda_0-\tilde\alpha)m_1''m_2'\bigr)$.
Since $R_{\tilde\ell\tilde\ell}\neq0$, $\tilde\alpha'\neq0$ would make $R_{\tilde\ell}$ independent of $\tilde\ell$, a contradiction. So $\tilde\alpha\in\kk$ and $\tilde\alpha\neq\lambda_0$. The linear forms $\tilde\ell=x_2+\tilde\alpha x_3$ and $\eta=x_2+\lambda_0x_3$ are independent with constant coefficients, $g^{\mathsf T}y=m_2'\eta$, so by \eqref{eq:u0g} $u_0=R(x_1,\tilde\ell)+\frac{m_2'^2}{W''}\eta^2+(\text{affine in }\tilde\ell,\eta)$ is quadratic in $\eta$, and $\beta=m_1+m_2\eta+c_3x_3$ is affine in $(\tilde\ell,\eta)$.
\end{proof}

\begin{proof}[Proof of Theorem~\ref{thm:quad}]
In Case~A, Proposition~\ref{prop:caseA} and Lemma~\ref{lem:blockreach} give the block $(x_3,t)$. In Case~B, Lemma~\ref{lem:caseB} gives $c=c(x_1)$. Then Proposition~\ref{prop:Qll} gives the block $(x_3,t)$ in the coordinates $(x_1,\ell,x_3)$, and Proposition~\ref{prop:Qll0} gives the block $(\eta,t)$ in the coordinates $(x_1,\tilde\ell,\eta)$ or the block $(x_2,x_3,t)$. In each case Theorem~\ref{thm:block} or Lemma~\ref{lem:quadblock} applies.
\end{proof}

\begin{corollary}\label{cor:IIb}
Let $\phi=\frac12a(x_1)t^2+b\,t+c(x_1,x_2,x_3)$ with $a$ nonconstant and $b=B(x_1,\ell)+m_2x_2+m_3x_3+\delta$, where $\ell=p(x_1)x_2+q(x_1)x_3$, $p,q$ coprime, $[p:q]$ nonconstant, $B_{\ell\ell}\neq0$, and $m_2,m_3,\delta\in\kk[x_1]$. Then $\det\Hess\phi\notin\kk^\times$.
\end{corollary}
\begin{proof}
Suppose $\det\Hess\phi\in\kk^\times$. We are in Case~B with $c$ replaced by $a=a(x_1)$. After exchanging $x_2,x_3$ if necessary, $p\neq0$, and $\ell'=x_2+(q/p)x_3$ satisfies $\ell=p\ell'$. Then $\beta=b/a$ is of the form \eqref{eq:betaform} with $Q(x_1,\ell')=B(x_1,p\ell')/a$, so $Q_{\ell'\ell'}=p^2B_{\ell\ell}/a\neq0$, and $\alpha=q/p$ is nonconstant. This contradicts Proposition~\ref{prop:Qll}.
\end{proof}

\begin{remark}
For the potentials of Corollary~\ref{cor:IIb} the kernel of $\Hess_{(x_2,x_3)}b$ is spanned by $(-q,p)$, whose direction depends on $x_1$. Corollary~\ref{cor:IIb} says that this twist is incompatible with a constant Hessian determinant.%
\end{remark}
\section{Proof of the theorem on directions of degree at most two}\label{sec:thmA}

Let $\phi$ have a direction $v$ of degree at most two. After a linear change $v=e_4$, and with $t=x_4$ and $x=(x_1,x_2,x_3)$ one has $\phi=a(x)+b(x)t+\frac12c(x)t^2$. The proof distinguishes three cases. For $c=0$ it uses Theorem~\ref{thm:affine}, which rests on Theorem~\ref{thm:cylinder}. For $c\in\kk^\times$ it uses Theorem~\ref{thm:pivot}, which rests on $\HC_3$. For $c\notin\kk$ it uses Theorem~\ref{thm:quad}, which rests on Theorems~\ref{thm:lead} and~\ref{thm:block}, the partial Legendre transform and $\HC_2$. In the second alternative of Theorems~\ref{thm:block} and~\ref{thm:quad} the potential has an affine direction, and the first case applies.

\begin{proof}[Proof of Theorem~\ref{thm:main-dir}]%
By Lemma~\ref{lem:descent} it suffices to work over $\overline\kk$. If $c=0$, Theorem~\ref{thm:affine} shows that $\grad\phi$ is an automorphism, or that $\phi=x_3C+x_4A+E$ with $C,A,E\in\kk[x_1,x_2]$ and $(C,A)$ a Keller map which is not injective unless $\grad\phi$ is. If $c\in\kk^\times$, Theorem~\ref{thm:pivot} shows that $\grad\phi$ is an automorphism. If $c\notin\kk$, Theorem~\ref{thm:quad} gives an automorphism or an affine direction, and the case $c=0$ applies to the affine direction. All changes of coordinates in Theorems~\ref{thm:affine} and~\ref{thm:quad} are linear, and the vector $u$ of Theorem~\ref{thm:quad} is constant. Hence in the form $x_3C+x_4A+E$ one has $\deg C,\deg A\le\deg\phi-1$.
\end{proof}

\begin{corollary}\label{cor:thmA}
\textup{(a)} The gradient of every Keller potential in four variables with a direction of degree at most two is a polynomial automorphism if and only if $\JC_2$ holds. \textup{(b)} Let $\phi$ be a Keller potential in four variables with a direction of degree at most two. If $\deg\phi\le101$, then $\grad\phi$ is a polynomial automorphism, and the same holds if $\deg\phi\le108$. If $\deg\phi\le125$ and $\grad\phi$ is not a polynomial automorphism, then over $\overline\kk$ the potential $\phi$ is a doubling of a planar Keller map whose components have the degrees $72$ and $108$.
\end{corollary}
\begin{proof}
(a) follows from Theorem~\ref{thm:main-dir} and Corollary~\ref{cor:ceiling}. (b) In the second alternative of Theorem~\ref{thm:main-dir} the map $(C,A)$ is not injective and $\deg C,\deg A\le\deg\phi-1$. A non-injective Keller map $(P,Q)$ of the plane satisfies $\max(\deg P,\deg Q)\ge101$ by \cite{Moh83}, which excludes $\deg\phi\le101$. By \cite[Theorem~2.1]{GGHV22}, it satisfies $\max(\deg P,\deg Q)\ge125$ or $\{\deg P,\deg Q\}=\{72,108\}$. For $\deg\phi\le108$ both alternatives are excluded, and for $\deg\phi\le125$ only the second remains.
\end{proof}

\begin{remark}\label{rem:liu}
In the second alternative of Theorem~\ref{thm:main-dir} the Hessian has the block form $\begin{psmallmatrix}*&N^{\mathsf T}\\ N&0\end{psmallmatrix}$ with $N=J(C,A)$. For a real potential which is a doubling in real coordinates the Hessian therefore has signature $(2,2)$. Liu \cite{Liu26} studies real Keller potentials in four variables whose Hessian has signature $(3,1)$. In that signature he shows that the gradient is an automorphism whenever $D_\xi^2\phi$ is constant for some $\xi\neq0$ \cite[Proposition~2.5]{Liu26}, and that $\HC_4$ holds for potentials of degree at most five \cite[Corollary~7.9]{Liu26}. By Corollary~\ref{cor:lorentz}, $\HC_4$ holds in signature $(3,1)$ without restriction on the degree.
\end{remark}
\section{Low degree}\label{sec:low}

In low degree we find a direction of degree at most two directly. Since $\det\Hess\phi$ is constant, its homogeneous layers of positive degree vanish, and the highest of them involve only the top forms of $\phi$. Together with the theorem of Gordan and Noether this successively removes variables from the top forms. In degree at most five the homogeneous layers suffice. In degree six we use in addition the homogeneous parts for the single weight vector $(1,1,1,2)$, which the layers determine.

Write a Keller potential of degree $\le d$ in four variables, after removing its affine part, as $\phi=\frac12x^{\mathsf T}Lx+\sum_{j=3}^dh_j$ with $L=\Hess\phi(0)$ invertible and $h_j$ forms of degree $j$. Put $H_j=\Hess h_j$.

\begin{lemma}[weights]\label{lem:weights}%
$\det\bigl(L+\sum_{j\ge3}s^{j-2}H_j(x)\bigr)=\det L$ in $\kk[s,x]$.
\end{lemma}
\begin{proof}
$\Hess\phi(sx)=L+\sum_js^{j-2}H_j(x)$, and $\det\Hess\phi$ is constant.
\end{proof}

\begin{lemma}[row weights]\label{lem:rowweights}%
Let $P(x)\in GL_4(\kk(x))$, $\tilde M(s)=P^{\mathsf T}\bigl(L+\sum_js^{j-2}H_j\bigr)P$, and let $w_1,\dots,w_4\ge0$ satisfy $\deg_s\tilde M_{ij}\le w_i+w_j$. Let $\tilde N_{ij}$ be the coefficient of $s^{w_i+w_j}$ in $\tilde M_{ij}$. If $\sum w_i>0$, then $\det\tilde N=0$.
\end{lemma}
\begin{proof}
Every term of the Leibniz expansion of $\det\tilde M$ has degree $\le2\sum w_i$ in $s$, and the coefficient of $s^{2\sum w_i}$ is $\det\tilde N$. By Lemma~\ref{lem:weights}, $\det\tilde M=\det(P)^2\det L$ does not depend on $s$.
\end{proof}

\begin{lemma}[leading coefficient]\label{lem:lead}%
Let $K_4\subset\kk^4$, $\dim K_4=4-e$, with $\partial_\kappa h_4=0$ for $\kappa\in K_4$. In a basis whose last $4-e$ vectors span $K_4$, let $D=(\partial_i\partial_jh_4)_{i,j\le e}$ and $\gamma=(\partial_i\partial_jh_3)_{i,j>e}$. Then $\det(L+sH_3+s^2H_4)$ has degree $\le e+4$ in $s$, with coefficient $\det D\cdot\det\gamma$ at $s^{e+4}$.
\end{lemma}
\begin{proof}
Lemma~\ref{lem:rowweights} with weights $w_i=1$ for $i\le e$ and $w_i=\frac12$ for $i>e$, after doubling the weights to make them integral and substituting $s^2$ for $s$. The rows and columns of $K_4$ meet $H_4$ only in zeros, so the leading block matrix is block diagonal with blocks $D$ and $\gamma$.
\end{proof}

\begin{theorem}[degree four]\label{thm:deg4}%
Let $\phi$ be a Keller potential of degree $\le4$ in four variables. There is $\xi\neq0$ with $D_\xi^2h_3=0$ and $\partial_\xi h_4=0$. In particular $D^2_\xi\phi=\xi^{\mathsf T}L\xi$ is constant, and $\xi$ is an affine direction or a constant pivot.
\end{theorem}
\begin{proof}
If $\deg\phi\le2$ every $\xi$ works, so let $h_3+h_4\neq0$. Let $K_4=\{\kappa:\partial_\kappa h_4=0\}$ and $e=4-\dim K_4$.

If $h_4=0$, the coefficient of $s^4$ in Lemma~\ref{lem:weights} is $\det H_3=0$, and \ref{E:GN} gives $\xi\neq0$ with $\partial_\xi h_3=0$.

Let $h_4\neq0$, so that $e\ge1$, and choose coordinates as in Lemma~\ref{lem:lead}. Then $h_4$ is a form in the first $e$ coordinates without a vertex direction among them, so $\det D\neq0$ by \ref{E:GN}. Lemmas~\ref{lem:weights} and~\ref{lem:lead} give $\det D\cdot\det\gamma=0$, since $e+4>0$, hence $\det\gamma=0$. In particular $e\neq4$, since for $e=4$ the matrix $\gamma$ is empty and has determinant $1$. For $\xi\in K_4$ one has $\partial_\xi h_4=0$ and $D_\xi^2h_3=\xi^{\mathsf T}\gamma\,\xi$, so it suffices to find $\xi\in K_4\setminus\{0\}$ with $\xi^{\mathsf T}\gamma\,\xi=0$.

For $e=3$ the matrix $\gamma$ is $1\times1$, so $\gamma=0$ and $\xi$ spanning $K_4$ works.

For $e=2$ write $\gamma=\begin{psmallmatrix}\alpha&\beta\\ \beta&\delta\end{psmallmatrix}$ with linear forms $\alpha,\beta,\delta$ and $\alpha\delta=\beta^2$. If $\alpha=0$, then $\beta=0$ and $\xi=(0,0,1,0)$ works. Otherwise unique factorization in $\kk[x]$ gives a linear form $\ell$ and constants $a,b,d$ with $\alpha=a\ell$, $\beta=b\ell$, $\delta=d\ell$ and $b^2=ad$. Then $\gamma=\ell\,\gamma_0$ with $\gamma_0$ a constant symmetric matrix of rank at most one, and a nonzero $\xi\in K_4$ in the kernel of $\gamma_0$ works.

For $e=1$ write $t$ for the first coordinate and $y\in\kk^3$ for the coordinates on $K_4$, so that $h_4=at^4$ with $a\neq0$ and $h_3=c_0(y)+tc_1(y)+t^2c_2(y)+c_3t^3$ with forms $c_j$ of degree $3-j$. Then $\gamma=\Hess_yc_0+tC_1$ with $C_1=\Hess_yc_1$ constant, and $\det(\Hess_yc_0+tC_1)=0$ identically in $t$. The term of degree $0$ in $t$ gives $\det\Hess_yc_0=0$, and by \ref{E:GN} the space $V_0=\{\xi\in K_4:\partial_\xi c_0=0\}$ is nonzero. A vector $\xi\in V_0$ satisfies $\xi^{\mathsf T}\gamma\,\xi=\xi^{\mathsf T}C_1\xi$. If $\dim V_0\ge2$, the quadratic form $C_1$ has a nonzero isotropic vector in $V_0$, since $\kk$ is algebraically closed. If $\dim V_0=1$, choose coordinates on $K_4$ with $V_0=\kk e_3$. Then $c_0\in\kk[y_1,y_2]$, and $\det\Hess_{(y_1,y_2)}c_0\neq0$, since otherwise the binary cubic $c_0$ would be the cube of a linear form and $\dim V_0\ge2$. The term of degree one in $t$ of $\det(\Hess_yc_0+tC_1)$ is $\operatorname{tr}\bigl(\adj(\Hess_yc_0)\,C_1\bigr)=\det\Hess_{(y_1,y_2)}c_0\cdot(C_1)_{33}$, so $(C_1)_{33}=0$ and $\xi=e_3$ works.
\end{proof}

\begin{corollary}\label{cor:deg4}%
$\Rc_4$ contains no potential of degree $\le4$, and $\grad\phi$ is an automorphism for every Keller potential $\phi$ of degree $\le4$ in four variables, with $\deg(\grad\phi)^{-1}\le27$.
\end{corollary}
\begin{proof}
By Theorem~\ref{thm:deg4}, $\phi$ has a constant pivot, and Theorem~\ref{thm:pivot} applies, or an affine direction. In the latter case Theorem~\ref{thm:affine} reduces injectivity either to \ref{E:dB} or to a planar Keller map $(C,A)$ of degree $\le3$, which is invertible by \cite{Moh83}.
The degree bound is \cite{BCW82}.
\end{proof}

\begin{remark}\label{rem:GNdim}
Result \ref{E:GN} is the only place in the proof of Theorem~\ref{thm:deg4} where $n\le4$ is used. In five variables the cubic $z_1z_4^2+z_2z_4z_5+z_3z_5^2$ of Perazzo \cite{Per1900} has vanishing Hessian determinant and no vertex direction.%
\end{remark}

\begin{remark}\label{rem:ni}
Corollary~\ref{cor:deg4} is \cite[Theorem~1.2, Corollary~7.1]{Ni26}. Ni shows that every quartic Keller potential in four variables becomes, after a linear change, $P(x_1,x_2,x_3)+x_4Q(x_1,x_2,x_3)+ax_4^2$ with $\deg Q\le2$ \cite[Propositions~4.1 and~5.3]{Ni26}. For such a potential $D^2_{e_4}\phi=2a$, so $e_4$ is a constant pivot or an affine direction. Theorem~\ref{thm:deg4} therefore has the same content as the structural part of \cite{Ni26}. The proof given here differs. It uses the weight identity of Lemma~\ref{lem:weights} and the leading coefficient of Lemma~\ref{lem:lead} instead of the individual homogeneous layers of $\det\Hess\phi$, and the same tools apply in degree five. The inversion step \cite[Proposition~6.1]{Ni26} is the case $\deg Q\le2$ of Theorems~\ref{thm:affine} and~\ref{thm:pivot}. In that case the planar Keller map of Theorem~\ref{thm:affine}(i) has a component of degree at most two, and \cite{Ni26} inverts it directly.
\end{remark}

\subsection{Degree five}

\begin{lemma}\label{lem:homlambda}%
Let $\psi\in\kk[y_1,\dots,y_r,u_1,\dots,u_n]$ be a form of degree $m$, and let $\lambda\in\kk(y)^n\setminus\{0\}$, $\gamma\in\kk(y)$ satisfy $\sum_i\lambda_i\partial_{u_i}\psi=\gamma$. Then there are $s\ge0$, forms $\lambda_i'\in\kk[y]$ of degree $s$ without a common factor, not all zero, and a form $\gamma'\in\kk[y]$ of degree $s+m-1$ with $\sum_i\lambda_i'\partial_{u_i}\psi=\gamma'$. If $r=1$, then $s=0$, that is $\lambda'\in\kk^n$.
\end{lemma}
\begin{proof}
Clear denominators so that $\lambda_i,\gamma\in\kk[y]$. Write $\lambda_i=\sum_s\lambda_{i,s}$ with $\lambda_{i,s}$ homogeneous of degree $s$ in $y$, and choose $s$ with $\lambda_{\bullet,s}\neq0$. Then $\sum_i\lambda_{i,s}\partial_{u_i}\psi$ is homogeneous of total degree $s+m-1$, so taking the component of that degree gives $\sum_i\lambda_{i,s}\partial_{u_i}\psi=\gamma_{s+m-1}$. Divide by the greatest common divisor of the $\lambda_{i,s}$, which then divides $\gamma_{s+m-1}$. If $r=1$, forms in one variable without a common factor are constants.
\end{proof}

\begin{theorem}\label{thm:deg5}%
Let $\deg\phi=5$ and $\det\Hess\phi\in\kk^\times$. Then $D_v^3\phi=0$ for some $v\in\kk^4\setminus\{0\}$.
\end{theorem}
\begin{proof}
Write $\phi=\phi_5+\dots+\phi_0$ with $\phi_j$ homogeneous of degree $j$, $\Delta=\det\Hess\phi$, and let $[\Delta]_r$ be the homogeneous part of degree $r$ of $\Delta$. All $[\Delta]_r$ with $r\ge1$ vanish. We use the following observation. If $U\subseteq\kk^4$ is a subspace of dimension $\ge2$ with $D_v^3\phi_5=D_v^3\phi_4=0$ for $v\in U$, then $D_v^3\phi=6\phi_3(v)$ for $v\in U$, and the restriction of the cubic form $\phi_3$ to a plane in $U$ has a nontrivial zero. In that situation the theorem holds.

\emph{Step 1.} The top layer $[\Delta]_{12}=\det\Hess\phi_5$ vanishes, and \ref{E:GN} gives $\phi_5=a_5(x_1,x_2,x_3)$ after a linear change. With $t=x_4$ and $x=(x_1,x_2,x_3)$,
\[
\phi=a+bt+\tfrac12ct^2+\tfrac16dt^3+\tfrac1{24}et^4,\qquad\deg a\le5,\ \deg b\le3,\ \deg c\le2,\ \deg d\le1,\ e\in\kk .
\]
\emph{Step 2.} One has $[\Delta]_{11}=(c_2+d_1t+\frac12et^2)\det\Hess_xa_5$, where $c_2,d_1$ are the parts of degree $2$, $1$ of $c,d$. If $(c_2,d_1,e)\neq0$, then $\det\Hess_xa_5=0$. If $c_2=d_1=e=0$, then $d=d_0\in\kk$ and
\[
[\Delta]_{10}=\det\begin{pmatrix}\Hess_xa_5&\grad b_3\\ \grad b_3^{\mathsf T}&c_1+d_0t\end{pmatrix},
\]
whose coefficient of $t$ is $d_0\det\Hess_xa_5$. So either $d_0=0$, in which case $\deg_t\phi\le2$ and $v=e_4$ works, or $\det\Hess_xa_5=0$. In the remaining cases \ref{E:GN} in three variables gives $a_5=a_5(x_1,x_2)$ after a linear change of $x$. Put $y=(x_1,x_2)$, $z=(x_3,x_4)$ and $U=\operatorname{span}(e_3,e_4)$. Then $D_v^3\phi_5=0$ for $v\in U$.

\emph{Step 3.} Now $[\Delta]_{11}=0$ and $[\Delta]_{10}=\det\Hess_ya_5\cdot\det\Hess_z\phi_4$.

\emph{Branch A, $\det\Hess_ya_5=0$.} A binary form with vanishing Hessian determinant is a power of a linear form, so after a linear change of $y$ one has $a_5=x_1^5$. Then $[\Delta]_9=20x_1^3\det\Hess_{(x_2,x_3,x_4)}\phi_4$, so $\det\Hess_{(x_2,x_3,x_4)}\phi_4=0$ in $\kk(x_1)[x_2,x_3,x_4]$. By \ref{E:dB62} over $\kk(x_1)$ there are two cases. In the first, there is a two-dimensional subspace $W$ of $\kk(x_1)^3$ with $D_vD_w\phi_4=0$ for $v,w\in W$. For a general $c\in\kk$, $W$ specializes at $x_1=c$ to a two-dimensional subspace $U$ of $\operatorname{span}(e_2,e_3,e_4)$ with $D_vD_w\phi_4(c,\cdot)=0$ for $v,w\in U$. Since $\phi_4(x_1,y)=(x_1/c)^4\phi_4(c,cy/x_1)$ for $y=(x_2,x_3,x_4)$, one gets $D_vD_w\phi_4=0$ for $v,w\in U$. As also $D_v^3\phi_5=0$ for $v\in U$, the observation applies. In the second case $\sum_i\lambda_i\partial_{x_{i+1}}\phi_4\in\kk(x_1)$ for some $\lambda\in\kk(x_1)^3\setminus\{0\}$, and Lemma~\ref{lem:homlambda} with $r=1$ gives $u\in\operatorname{span}(e_2,e_3,e_4)\setminus\{0\}$ and $\gamma'\in\kk$ with $\partial_u\phi_4=\gamma'x_1^3$. After a linear change of $(x_2,x_3,x_4)$ we may take $u=e_4$, so $\phi_4=\tilde a_4(x_1,x_2,x_3)+\gamma'x_1^3x_4$. Thus $b_3=\gamma'x_1^3$, $c_2=0$, $d\in\kk$ and $e=0$. If $d=0$, then $\deg_{x_4}\phi\le2$. If $d\neq0$, then
\[
[\Delta]_8=\det\Hess_{(x_2,x_3)}\tilde a_4\cdot\bigl(20x_1^3(c_1+dt)-9\gamma'^2x_1^4\bigr).
\]
The second factor involves $t$, so $\det\Hess_{(x_2,x_3)}\tilde a_4=0$. By \ref{E:dBvdE} over $\kk(x_1)$ and Lemma~\ref{lem:homlambda} there is $w\in\operatorname{span}(e_2,e_3)\setminus\{0\}$ with $\partial_w\tilde a_4=\gamma''x_1^3$. After a linear change of $(x_2,x_3)$, $w=e_3$ and $\phi_4=\hat a_4(x_1,x_2)+x_1^3(\gamma''x_3+\gamma'x_4)$. Hence $D_v^3\phi_4=0$ for $v\in\operatorname{span}(e_3,e_4)$, and the observation applies.

\emph{Branch B, $\det\Hess_ya_5\neq0$ and $\det\Hess_z\phi_4=0$.} Apply \ref{E:dBvdE} over $\kk(y)$ to $\phi_4\in\kk(y)[z]$ and then Lemma~\ref{lem:homlambda}. There are coprime forms $\lambda_3,\lambda_4\in\kk[y]$ of degree $s\ge0$ and a form $\gamma'\in\kk[y]$ of degree $s+3$ with $\lambda_3\partial_3\phi_4+\lambda_4\partial_4\phi_4=\gamma'$. Put $\zeta=\lambda_4x_3-\lambda_3x_4$ and choose $\eta$ such that $(\zeta,\eta)$ are $\kk(y)$-linear coordinates in $z$ with $\partial_\eta=\lambda_3\partial_3+\lambda_4\partial_4$. Then $\partial_\eta(\phi_4-\gamma'\eta)=0$, so $\phi_4=\sum_{j=0}^4P_j(y)\zeta^j+\gamma'(y)\eta$ with $P_j\in\kk(y)$. For $j\ge2$ the part of $\phi_4$ of degree $j$ in $z$ is $P_j\zeta^j$. Its coefficients at $x_3^j$ and $x_4^j$ are $P_j\lambda_4^j$ and $P_j(-\lambda_3)^j$, both polynomials, and $\gcd(\lambda_3,\lambda_4)=1$. Hence $P_j\in\kk[y]$ is homogeneous of degree $4-j(s+1)$.

(B2) $s\ge1$. Then $4-j(s+1)<0$ for $j\ge3$, so $P_3=P_4=0$ and $\deg_z\phi_4\le2$. Thus $D_v^3\phi_4=0$ for $v\in U$, and the observation applies.

(B1) $s=0$. Then $\lambda$ is constant and, after a linear change of $z$, $\partial_4\phi_4=\gamma'(y)$, that is $\phi_4=\Phi(y,x_3)+\gamma'(y)x_4$ with $\gamma'$ a binary cubic. For $v\in U$ one has $D_v^3\phi_4=v_3^3\,\partial_3^3\Phi$. If $\deg_{x_3}\Phi\le2$, the observation applies. Suppose $\deg_{x_3}\Phi\ge3$, so $q=\partial_3^2\Phi\neq0$. The term of $\phi$ in $t^3$ is $\frac16dt^3$ with $d\in\kk$, because $\phi_5,\phi_4$ are at most linear in $t$. If $d=0$, then $\deg_t\phi\le2$. If $d\neq0$, then
\[
[\Delta]_9=q\cdot\det\begin{pmatrix}\Hess_ya_5&\grad_y\gamma'\\ \grad_y\gamma'^{\mathsf T}&c_1+dt\end{pmatrix}.
\]
The coefficient of $t$ in the second factor is $d\det\Hess_ya_5\neq0$, so $[\Delta]_9\neq0$, a contradiction.
\end{proof}

\begin{proof}[Proof of Theorem~\ref{thm:main-low}]
Part (i) is Theorem~\ref{thm:deg4}, and part (ii) is Theorem~\ref{thm:deg4} for $\deg\phi\le4$, since an affine direction or a constant pivot $v$ satisfies $D_v^3\phi=0$, and Theorem~\ref{thm:deg5} for $\deg\phi=5$.
\end{proof}

\begin{corollary}\label{cor:deg5HC}
$\HC_4$ holds for Keller potentials of degree at most five.
\end{corollary}
\begin{proof}
By Theorem~\ref{thm:main-low}(ii) and Theorem~\ref{thm:main-dir}, either $\grad\phi$ is an automorphism, or $\phi$ is a doubling of a non-injective Keller map of degree $\le4$, which contradicts \cite{Moh83}.
\end{proof}

\subsection{Degree six, top form of rank three}

Let $\deg\phi=6$ and $\det\Hess\phi\in\kk^\times$, and write $\phi=\phi_6+\dots+\phi_0$ with $\phi_j$ homogeneous of degree $j$. The top layer $\det\Hess\phi_6$ of $\det\Hess\phi$ vanishes, and \ref{E:GN} gives $\phi_6\in\kk[x']$, $x'=(x_1,x_2,x_3)$, after a linear change. Put $\delta=\det\Hess_{x'}\phi_6$. The case $\delta\neq0$ is the case in which $\Hess\phi_6$ has rank three.

\begin{proposition}\label{prop:deg6}%
Let $\deg\phi=6$, $\det\Hess\phi\in\kk^\times$, $\phi_6\in\kk[x']$ and $\delta\neq0$.
\begin{enumerate}[label=\textup{(\alph*)}]
\item $\phi=F+x_4G+x_4^2Q+\gamma x_4^3$ with $F,G,Q\in\kk[x']$, $\deg G\le4$, $\deg Q\le2$ and $\gamma\in\kk$. The parts $b_4$ and $c_2$ of degrees $4$ and $2$ of $G$ and $Q$ satisfy
\[
2\delta\,c_2=\nabla_{x'}b_4^{\mathsf T}\,\adj\bigl(\Hess_{x'}\phi_6\bigr)\,\nabla_{x'}b_4 .
\]
\item If $\gamma=0$, then $D_{e_4}^3\phi=0$.
\item For the weights $(1,1,1,2)$ the polynomial $\Psi=\phi_6+x_4b_4+x_4^2c_2+\gamma x_4^3$ is the part of weight six of $\phi$, all other terms of $\phi$ have weight less than six, $\det\Hess\Psi=0$ and $\rk\Hess\Psi=3$.
\item If $\gamma\neq0$, then $\Hess\Psi\,v\neq0$ for every $v\in\kk^4\setminus\{0\}$.
\end{enumerate}
\end{proposition}

\begin{proof}
Put $H_j=\Hess\phi_j$ and $A=\Hess_{x'}\phi_6$, and let $[\Delta]_r$ be the homogeneous part of degree $r$ of $\Delta=\det\Hess\phi$, so that $[\Delta]_r=0$ for $r\ge1$. The entries of $H_j$ are forms of degree $j-2$. Since the last row and column of $H_6$ vanish, $\adj H_6=\delta\,e_4e_4^{\mathsf T}$.

(a) A term of degree $15$ of the Leibniz expansion of $\Delta$ has three entries from $H_6$ and one from $H_5$, so $[\Delta]_{15}=\operatorname{tr}(\adj H_6\,H_5)=\delta\,\partial_4^2\phi_5$. Hence $\phi_5=a_5+x_4b_4$ with $a_5,b_4\in\kk[x']$. A term of degree $14$ has either three entries from $H_6$ and one from $H_4$, and these terms give $\delta\,\partial_4^2\phi_4$, or two entries from $H_6$ and two from $H_5$. In the second case the entries in the last row and in the last column come from $H_5$, and they are distinct, since the entry $(4,4)$ of $H_5$ is $\partial_4^2\phi_5=0$. These terms are the terms of the expansion of
\[
\det\begin{pmatrix}A&v\\ v^{\mathsf T}&0\end{pmatrix}=-v^{\mathsf T}\adj(A)\,v,\qquad v=\nabla_{x'}b_4 .
\]
Therefore $[\Delta]_{14}=\delta\,\partial_4^2\phi_4-v^{\mathsf T}\adj(A)\,v$. The second term lies in $\kk[x']$, so $\partial_4^2\phi_4=2c_2$ with $c_2\in\kk[x']$ as stated, and $\phi_4=a_4+x_4b_3+x_4^2c_2$ with $a_4,b_3\in\kk[x']$. Thus $\phi_6,\phi_5,\phi_4$ have degree at most $0,1,2$ in $x_4$. The coefficient of $x_4^3$ in $\phi_3$ is a constant $\gamma$, and $\phi_2,\phi_1,\phi_0$ have degree at most two in $x_4$. This gives the form of $\phi$, with $b_4$ and $c_2$ the parts of $G$ and $Q$ coming from $\phi_5$ and $\phi_4$.

(b) $D_{e_4}^3\phi=6\gamma$.

(c) A term of $F$, $x_4G$, $x_4^2Q$ of degree $d$ in $x'$ has weight $d$, $d+2$, $d+4$, which is at most six, with equality exactly for the terms of $\phi_6$, $x_4b_4$, $x_4^2c_2$. Every term of $\partial_i\partial_j\phi$ has weight at most $6-w_i-w_j$, and the part of weight $6-w_i-w_j$ is $\partial_i\partial_j\Psi$. So every product of the Leibniz expansion of $\Delta$ has weight at most $4\cdot6-2\cdot5=14$, and the part of weight $14$ of $\Delta$ is $\det\Hess\Psi$. Since $\Delta$ is constant, $\det\Hess\Psi=0$. The principal minor of $\Hess\Psi$ in $x'$ is $\delta\ne0$ at $x_4=0$, so $\rk\Hess\Psi=3$.

(d) Let $\Hess\Psi\,v=0$ with $v=(v',v_4)\in\kk^4$. Then $D_v\Psi$ is constant. Its parts of weights five and four are $D_{v'}\Psi$ and $v_4\partial_4\Psi$, so both vanish. The coefficient of $x_4^2$ in $\partial_4\Psi$ is $3\gamma\ne0$, so $v_4=0$. Then $D_{v'}\phi_6=(D_{v'}\Psi)|_{x_4=0}=0$, hence $Av'=0$ and $v'=0$.
\end{proof}

By Proposition~\ref{prop:deg6}(b), if $\gamma=0$, then $e_4$ is a direction of degree at most two. For $\gamma\neq0$ we use the weighted forms of Proposition~\ref{prop:deg6}(c), and we write $(x,y,z,t)$ for $(x_1,x_2,x_3,x_4)$. Linear changes of $(x,y,z)$ preserve the weights $(1,1,1,2)$ and the statements of Proposition~\ref{prop:deg6}.

\begin{lemma}\label{lem:deg6-dir}%
Let $h,b,c\in\kk[x,y,z]$ be forms of degrees $6$, $4$, $2$, let $\gamma\in\kk^\times$ and $\Psi=h+tb+t^2c+\gamma t^3$. If $\det\Hess_{(x,y,z)}h\neq0$ and $\det\Hess\Psi=0$, then $D_u^2\Psi=0$ for some $u\in\kk^3\setminus\{0\}$.
\end{lemma}
\begin{proof}
Put $H=\Hess h$, $B=\Hess b$ and $p=\grad b$, all in $(x,y,z)$, and let $[\,\cdot\,]_{t^j}$ denote the coefficient of $t^j$. The formulas for these coefficients below follow by expansion of $\det\Hess\Psi$. We distinguish the rank of the constant matrix $\Hess c$.

\emph{$c=0$.} Then $\det\Hess\Psi=6\gamma t\det(H+tB)-p\T\adj(H+tB)p$, and $[\det\Hess\Psi]_{t^4}=6\gamma\det B$, so $\det B=0$. Suppose $\rk B\le1$. By Lemma~\ref{lem:app-low}, $B$ vanishes on a two-dimensional subspace $U$. For $u\in U$ the cubic form $\partial_ub$ has zero gradient, hence vanishes, and $b=\alpha\ell^4$ with $\alpha\in\kk$ and a linear form $\ell$. With $v=\grad\ell$ one has $p=4\alpha\ell^3v$ and $B=12\alpha\ell^2vv\T$. Since $v\T\adj(M+\sigma vv\T)v=v\T\adj(M)v$ for every symmetric $M$, by the bordered determinant, $p\T\adj(H+tB)p$ does not depend on $t$, and $[\det\Hess\Psi]_t=6\gamma\det H\neq0$, a contradiction. Hence $\rk B=2$, and by \ref{E:GN}, after a linear change, $b\in\kk[x,y]$ and $\adj B=\delta_be_ze_z\T$ with $\delta_b=\det\Hess_{(x,y)}b\neq0$. As $\adj(H+tB)$ has degree at most two in $t$, $[\det\Hess\Psi]_{t^3}=6\gamma\operatorname{tr}(\adj(B)H)=6\gamma\delta_bh_{zz}$. So $h_{zz}=0$, and $u=e_z$.

\emph{$\rk\Hess c=2$.} After a linear change $c=xy$. Then $[\det\Hess\Psi]_{t^6}=-6\gamma b_{zz}$, and if $b_{zz}=0$, then $[\det\Hess\Psi]_{t^5}=6\gamma(2b_{xz}b_{yz}-h_{zz})$. Hence $b=b_0+z\beta$ and $h=z^2\beta_x\beta_y+z\alpha_5+\alpha_6$ with forms $b_0,\beta,\alpha_5,\alpha_6\in\kk[x,y]$. For such $b,h$ the coefficient of $z$ in $[\det\Hess\Psi]_{t^4}$ is $18\gamma(\beta_{xx}\beta_y^2+2\beta_{xy}\beta_x\beta_y+\beta_{yy}\beta_x^2)$. Write $\beta=\beta_0x^3+\beta_1x^2y+\beta_2xy^2+\beta_3y^3$. The coefficients of $x^5$, $x^4y$, $xy^4$, $y^5$ of the last factor are, up to nonzero constants,
\[
\beta_0(\beta_0\beta_2+\beta_1^2),\quad 27\beta_0^2\beta_3+42\beta_0\beta_1\beta_2+5\beta_1^3,\quad 27\beta_0\beta_3^2+42\beta_1\beta_2\beta_3+5\beta_2^3,\quad \beta_3(\beta_1\beta_3+\beta_2^2).
\]
If $\beta_0\beta_3\neq0$, then $\beta_0\beta_2=-\beta_1^2$ and $\beta_1\beta_3=-\beta_2^2$. For $\beta_1=0$ this gives $\beta_2=0$ and $27\beta_0^2\beta_3=0$, and for $\beta_1\neq0$ it gives $\beta_2=-\beta_1^2/\beta_0$, $\beta_3=-\beta_1^3/\beta_0^2$ and $-64\beta_1^3=0$. Therefore $\beta_0\beta_3=0$, and the second and third expressions give $\beta\in\kk x^3\cup\kk y^3$. So $\beta_x\beta_y=0$, $h_{zz}=0$, and $u=e_z$.

\emph{$\rk\Hess c=1$.} After a linear change $c=x^2$. Then $[\det\Hess\Psi]_{t^5}=12\gamma\det\Hess_{(y,z)}b$. By \ref{E:dBvdE} over $\kk(x)$ and Lemma~\ref{lem:homlambda} with $r=1$, after a linear change of $(y,z)$, $b=b_0(x,y)+\kappa x^3z$ with $\kappa\in\kk$.

Let $\partial_y^2b_0\neq0$. As $[\det\Hess\Psi]_{t^4}=12\gamma b_{yy}(h_{zz}-\frac12b_{xz}^2)$, one gets $h=\frac94\kappa^2x^4z^2+z\alpha_5+\alpha_6$ with $\alpha_5,\alpha_6\in\kk[x,y]$. For such $h$ the coefficient of $z$ in $[\det\Hess\Psi]_{t^3}$ is $-486\gamma\kappa^3x^5\partial_y^2b_0$. Hence $\kappa=0$, $h_{zz}=0$, and $u=e_z$.

Let $\partial_y^2b_0=0$. Then $b=ax^4+x^3\ell$ with $a\in\kk$ and $\ell=Ay+Bz$, $A,B\in\kk$. Put $\tilde h=h-\frac94x^4\ell^2$. Then $[\det\Hess\Psi]_{t^3}=12\gamma\det\Hess_{(y,z)}\tilde h$, and \ref{E:dBvdE} over $\kk(x)$ and Lemma~\ref{lem:homlambda} give, after a linear change of $(y,z)$ which preserves the form of $b$, $\tilde h=f(x,y)+\rho x^5z$ with $\rho\in\kk$. If $B=0$, then $u=e_z$. If $B\neq0$, the coefficient of $z$ in $[\det\Hess\Psi]_{t^2}$ is $-486\gamma B^3x^5f_{yy}$. Hence $f_{yy}=0$, $\tilde h$ is affine in $(y,z)$, and $u=(0,B,-A)$ satisfies $D_u\ell=0$, hence $D_u^2\tilde h=D_u^2(x^4\ell^2)=D_u^2b=D_u^2c=0$.
\end{proof}

\begin{lemma}\label{lem:deg6-nf}%
Let $\phi$ be as in Proposition~\ref{prop:deg6} with $\gamma\neq0$. After a linear change of $(x,y,z)$ there are $c,\kappa_0\in\kk^\times$ and $a_0,c_1,c_2,c_3\in\kk$ with
\[
\Psi=a_0xy^5+R(y,s)+czy^5,\qquad s=t+\kappa_0xy,\qquad R=\gamma s^3+c_1y^2s^2+c_2y^4s+c_3y^6 .
\]
\end{lemma}
\begin{proof}
By Proposition~\ref{prop:deg6}(c) and Lemma~\ref{lem:deg6-dir}, after a linear change, $\partial_z^2\Psi=0$, that is $\Psi=P+zQ$ with $P,Q\in\kk[X]$, $X=(x,y,t)$, weighted homogeneous of weights $6$ and $5$. By Proposition~\ref{prop:deg6}(d), $\Hess\Psi$ has no nonzero constant kernel vector. In particular $Q\neq0$. Put $q=\grad_XQ$ and $A=\Hess_XP+z\Hess_XQ$. Then $\det\Hess\Psi=-q\T\adj(A)q$, so the coefficients of $z^2$ and $z$ of $q\T\adj(A)q$ vanish. The first is $\BH(Q)$, and Theorem~\ref{thm:main-BH}(ii) shows that $N=\{\kappa\in\kk^3:\partial_\kappa Q=0\}$ is nonzero. As $Q$ is weighted homogeneous, $N$ is the sum of its intersections with $\kk e_x+\kk e_y$ and $\kk e_t$.

Suppose $\partial_tQ\neq0$. Then $N\cap(\kk e_x+\kk e_y)\neq0$, and after a linear change of $(x,y)$, $Q=Q(x,t)$. The coefficient of $z$ is then $P_{yy}\BH_{x,t}(Q)$, where $\BH_{x,t}$ is the bordered Hessian in $(x,t)$. If $\BH_{x,t}(Q)=0$, Proposition~\ref{prop:BHdet} and \ref{E:dBvdE} give $\lambda\in\kk^2\setminus\{0\}$ and $\gamma'\in\kk$ with $\lambda_1Q_x+\lambda_2Q_t=\gamma'$. The three terms have the weights $4$, $3$, $0$, so $\lambda_2=0$, $Q_x=0$, and $Q\in\kk[t]$ has weight five, which is impossible. So $P_{yy}=0$ and $P=P_0(x,t)+yP_1(x,t)$. In the order $(x,t,y,z)$ the Hessian of $\Psi$ has a zero lower right block of size two, so $\det\Hess\Psi=\det\bigl(\partial(P_1,Q)/\partial(x,t)\bigr)^2$ and $dP_1\wedge dQ=0$. Since $x\partial_x+2t\partial_t$ multiplies $P_1$ and $Q$ by $5$, $d(P_1/Q)=0$ and $P_1=\lambda Q$ with $\lambda\in\kk$. Then $e_y-\lambda e_z$ is a constant kernel vector, which is excluded.

Hence $Q=f(x,y)$ is a binary quintic. The coefficient of $z$ is $P_{tt}\BH(f)=\frac54P_{tt}f\det\Hess f$, and $P_{tt}\neq0$ contains $6\gamma t$. So $\det\Hess f=0$, and $f=cy^5$ after a linear change of $(x,y)$, with $c\neq0$. Now $q=5cy^4e_y$ and $\det\Hess\Psi=-25c^2y^8\det\Hess_{(x,t)}P$. By \ref{E:dBvdE} over $\kk(y)$ there are $\lambda\in\kk(y)^2\setminus\{0\}$ and $\eta\in\kk(y)$ with $\lambda_1P_x+\lambda_2P_t=\eta$. As $P_t$ contains $3\gamma t^2$, $\lambda_1\neq0$, and we may take $P_x-\kappa P_t=\eta$. Comparing the coefficients of $t^2$ shows that $\kappa=\kappa_0y$ with $\kappa_0\in\kk$, and then $\eta=a_0y^5$ by weights. In the polynomial coordinates $(x,y,s)$, $s=t+\kappa_0xy$, the derivation $\partial_x-\kappa_0y\partial_t$ is the derivative in $x$, so $P=a_0xy^5+R(y,s)$ with $R$ as stated. If $\kappa_0=0$, then $ce_x-a_0e_z$ is a constant kernel vector. Therefore $\kappa_0\neq0$.
\end{proof}

We keep the notation of Lemma~\ref{lem:deg6-nf}, and use $(x,y,s,z)$ as polynomial coordinates on $\kk^4$. Subscripts $x,y,s,z$ denote derivatives in these coordinates, while $\Hess_X$ and $\grad_X$ are taken in $X=(x,y,t)$. Put
\[
r=R_{ss}=6\gamma s+2c_1y^2,\qquad T=5a_0y^4+\kappa_0R_s=\tau_0y^4+\chi r^2,\qquad \chi=\frac{\kappa_0}{12\gamma},
\]
with $\tau_0\in\kk$, and $L=(1,0,-\kappa_0y)\T$, $n=\grad_Xs=(\kappa_0y,\kappa_0x,1)\T$. The polynomial $r$ is irreducible, and $y$ divides neither $r$ nor $T$. For $g\in\kk[x,y,s,z]$ one has
\begin{equation}\label{eq:deg6-L}
\grad_Xg=g_xe_x+g_ye_y+g_sn,\qquad L\T\grad_Xg=g_x,\qquad \Hess_Xg\cdot L=\grad_X(g_x)+\kappa_0g_se_y,
\end{equation}
and $L\T\Hess_Xg\,L=g_{xx}$. Let $A=\Hess_X\Psi$ and $q=\grad_X\Psi_z=5cy^4e_y$. Then $AL=Te_y$ and $Ae_t=rn+R_{sy}e_y$. In the basis $(L,e_y,e_t)$, of determinant one, the matrix of the form $A$ has the entries $0,T,0$ in its first row and $r$ in its last diagonal entry, so $\det A=-rT^2\neq0$. For $K=A^{-1}$ it follows that
\begin{equation}\label{eq:deg6-K}
Ke_y=\frac LT,\qquad Kn=\frac{e_t}r-\frac{R_{sy}}{rT}L,\qquad e_y\T Ke_y=n\T Ke_y=0,\qquad n\T Kn=\frac1r .
\end{equation}
Put $k=5cy^4/T$ and $\omega=e_z-kL$, so that $Kq=kL$ and $q\T Kq=0$. By \eqref{eq:deg6-L}, for every polynomial $g$,
\begin{equation}\label{eq:deg6-omega}
\omega\T\Hess g\,\omega=g_{zz}-2kg_{xz}+k^2g_{xx}.
\end{equation}

\begin{lemma}\label{lem:deg6-div}%
If $g_1,g_2,g_3\in\kk[x,y,s,z]$ are weighted homogeneous of weight less than four and $T^2g_3-10cy^4Tg_2+25c^2y^8g_1=0$, then $g_1=g_2=g_3=0$.
\end{lemma}
\begin{proof}
$T$ divides $y^8g_1$, hence $g_1$, which has smaller weight. Then $T$ divides $g_2$, and $g_3=0$.
\end{proof}

Let $\iota$ be the automorphism of $\kk(x,y,s,z)$ over $\kk(x,y,z)$ with $\iota(r)=-r$. It fixes $T$ and $k$. An element $g$ is \emph{even} if $\iota(g)=g$ and \emph{odd} if $\iota(g)=-g$. A nonzero element is not both.

Write $\phi=\Psi+\sum_{j\le5}\Phi_j$ with $\Phi_j$ of weight $j$, and put $H_j=\Hess_X\Phi_j$, $a_j=\grad_X\partial_z\Phi_j$ and $d_j=a_j-kH_jL$.

\begin{lemma}\label{lem:deg6-schur}%
With $K_1=-KH_5K$ and $K_2=KH_5KH_5K-KH_4K$,
\begin{align*}
\textup{(W13)}\quad&\omega\T\Hess\Phi_5\,\omega=0,\\
\textup{(W12)}\quad&\omega\T\Hess\Phi_4\,\omega=d_5\T Kd_5,\\
\textup{(W11)}\quad&\omega\T\Hess\Phi_3\,\omega=2d_5\T Kd_4+d_5\T K_1d_5,\\
\textup{(W10)}\quad&\omega\T\Hess\Phi_2\,\omega=d_4\T Kd_4+2d_5\T Kd_3+2d_5\T K_1d_4+d_5\T K_2d_5 .
\end{align*}
\end{lemma}
\begin{proof}
For an indeterminate $\varepsilon$ put $\phi^\varepsilon=\varepsilon^6\phi(\varepsilon^{-1}x,\varepsilon^{-1}y,\varepsilon^{-1}z,\varepsilon^{-2}t)=\Psi+\sum_{j\le5}\varepsilon^{6-j}\Phi_j$. Then $\det\Hess\phi^\varepsilon=\varepsilon^{14}\det\Hess\phi$, a nonzero constant times $\varepsilon^{14}$. In the order $(x,y,t,z)$ let $S=\begin{psmallmatrix}I&-kL\\0&1\end{psmallmatrix}$, of determinant one. Since $Kq=kL$ and $q\T Kq=0$,
\[
S\T\Hess\phi^\varepsilon S=\begin{pmatrix}A_\varepsilon&c_\varepsilon\\ c_\varepsilon\T&\sigma_\varepsilon\end{pmatrix},\qquad A_\varepsilon=A+\sum_j\varepsilon^{6-j}H_j,\]\[c_\varepsilon=\sum_j\varepsilon^{6-j}d_j,\quad \sigma_\varepsilon=\sum_j\varepsilon^{6-j}\omega\T\Hess\Phi_j\,\omega .
\]
Hence $(\sigma_\varepsilon-c_\varepsilon\T A_\varepsilon^{-1}c_\varepsilon)\det A_\varepsilon=\varepsilon^{14}\det\Hess\phi$ in $\kk(x,y,z,t)[[\varepsilon]]$, and $\det A_\varepsilon\equiv\det A\neq0$ modulo $\varepsilon$. So the coefficients of $\varepsilon,\dots,\varepsilon^4$ of $\sigma_\varepsilon-c_\varepsilon\T A_\varepsilon^{-1}c_\varepsilon$ vanish, and $A_\varepsilon^{-1}=K+\varepsilon K_1+\varepsilon^2K_2+O(\varepsilon^3)$.
\end{proof}

\begin{theorem}\label{thm:deg6}%
Let $\phi\in\kk[x_1,\dots,x_4]$ be a Keller potential of degree six whose top form has Hessian rank three. Then $D_v^3\phi=0$ for some $v\neq0$, and $\grad\phi$ is a polynomial automorphism.
\end{theorem}
\begin{proof}
By \ref{E:GN} we are in the situation of Proposition~\ref{prop:deg6}. If $\gamma=0$, Proposition~\ref{prop:deg6}(b) gives $v=e_4$. Let $\gamma\neq0$, and choose the coordinates of Lemma~\ref{lem:deg6-nf}. We show that, after a linear change of coordinates, $\partial_z^2\phi=0$.

\emph{Weight $13$.} By (W13), \eqref{eq:deg6-omega} and Lemma~\ref{lem:deg6-div}, the derivatives $\Phi_{5,xx}$, $\Phi_{5,xz}$, $\Phi_{5,zz}$ vanish, that is $\Phi_5=x\mathsf A+\mathsf B+z\mathsf C$ with $\mathsf A,\mathsf B,\mathsf C\in\kk[y,s]$ of weights $4$, $5$, $4$.

\emph{Weight $12$ and a shear.} By \eqref{eq:deg6-L}, $d_5=Un+Ve_y$ with $U=\mathsf C_s-k\mathsf A_s$ and $V\in\kk(x,y,s,z)$ affine in $x$ with $V_x=-k\kappa_0\mathsf A_s$. By \eqref{eq:deg6-K}, $d_5\T Kd_5=U^2/r$, and (W12) for $F=\Phi_4$ reads
\begin{equation}\label{eq:deg6-w12}
r\bigl(T^2F_{zz}-10cy^4TF_{xz}+25c^2y^8F_{xx}\bigr)=\bigl(T\mathsf C_s-5cy^4\mathsf A_s\bigr)^2 .
\end{equation}
So $r$ divides $T\mathsf C_s-5cy^4\mathsf A_s\equiv y^4(\tau_0\mathsf C_s-5c\mathsf A_s)$ modulo $r$, and $\tau_0\mathsf C_s-5c\mathsf A_s\in\kk r$ by weights. Write $\mathsf A_s=A_*y^2+\rho_Ar$ and $\mathsf C_s=C_*y^2+\rho_Cr$ with constants. Then $\tau_0C_*=5cA_*$. The substitution $t\mapsto t-\rho_Ax-\rho_Cz$ preserves $\Psi$ and replaces $\Phi_5$ by $\Phi_5-(\rho_Ax+\rho_Cz)R_s$, and a direction for the new potential gives one for $\phi$. So we may assume $\rho_A=\rho_C=0$, $\mathsf A_s=A_*y^2$ and $\mathsf C_s=C_*y^2$.

\emph{Weight $11$, $C_*=0$.} Suppose $C_*\neq0$. Then $U=\chi C_*y^2r^2/T$, the right side of \eqref{eq:deg6-w12} is $\chi^2C_*^2y^4r^4$, and by Lemma~\ref{lem:deg6-div} the unique solution of \eqref{eq:deg6-w12} is
\[
F_{xx}=\alpha r,\quad F_{xz}=\beta r,\quad F_{zz}=0,\qquad \alpha=-\frac{\tau_0\chi C_*^2}{25c^2},\quad \beta=-\frac{\chi C_*^2}{10c}.
\]
By \eqref{eq:deg6-K}, $Kd_5=pe_t+hL$ with $p=U/r$ and $h=V/T-UR_{sy}/(rT)$, so $h_x=-k\kappa_0\mathsf A_s/T$. By \eqref{eq:deg6-L}, $e_t\T d_4=F_{zs}-kF_{xs}$, $L\T d_4=F_{xz}-kF_{xx}=r(\beta-k\alpha)$, $e_t\T H_5e_t=\mathsf B_{ss}$, $e_t\T H_5L=\mathsf A_s$ and $L\T H_5L=0$, and the right side of (W11) is
\[
2p\,(F_{zs}-kF_{xs})+2h\,r(\beta-k\alpha)-p^2\mathsf B_{ss}-2ph\mathsf A_s .
\]
It is affine in $x$, since $F_{xs}-6\gamma\alpha x$ and $F_{zs}-6\gamma\beta x$ do not depend on $x$. Its coefficient of $x$ is
\[
12\gamma p(\beta-k\alpha)+2h_x\bigl(r(\beta-k\alpha)-p\mathsf A_s\bigr)=\frac{6\gamma\chi^2C_*^3\,y^2r\,(3\tau_0y^4-\chi r^2)}{5c\,T^2},
\]
which is odd and nonzero. On the left of (W11), $\Phi_{3,xx}$, $\Phi_{3,xz}$, $\Phi_{3,zz}$ have weight one and are linear forms in $(x,y,z)$, so by \eqref{eq:deg6-omega} the coefficient of $x$ lies in $\kk[k]$ and is even. This contradiction gives $C_*=0$, hence $A_*=0$, and $\Phi_5=\mathsf B(y,s)+A_0xy^4+C_0zy^4$ with $A_0,C_0\in\kk$.

\emph{Weights $12$, $11$, $10$.} Now $d_5=me_y$ with $m\in\kk(x,y,s,z)$, and $H_5L=(4A_0y^3+\kappa_0\mathsf B_s)e_y$ by \eqref{eq:deg6-L}. By \eqref{eq:deg6-K}, $d_5\T Kd_5=0$, and (W12) with Lemma~\ref{lem:deg6-div} gives $F_{xx}=F_{xz}=F_{zz}=0$, that is $F=xF_1+zF_2+F_0$ with $F_i\in\kk[y,s]$. Then $L\T d_4=0$ and $Kd_5=(m/T)L$, so $d_5\T Kd_4=0$ and $d_5\T K_1d_5=-(m/T)^2L\T H_5L=0$. By (W11) and Lemma~\ref{lem:deg6-div}, $\Phi_{3,xx}=\Phi_{3,xz}=\Phi_{3,zz}=0$, so $L\T d_3=0$. The same identities give $d_5\T Kd_3=d_5\T K_1d_4=0$ and $d_5\T K_2d_5=-(m/T)^2F_{xx}=0$, since $(H_5L)\T K(H_5L)$ is a multiple of $e_y\T Ke_y=0$. Finally $d_4=(F_{2,s}-kF_{1,s})n+V_4e_y$ with $V_4\in\kk(x,y,s,z)$, and $F_{1,s}=f_1y$, $F_{2,s}=u_1y$ with $f_1,u_1\in\kk$ by weights. The second derivatives $a,b,d$ of $\Phi_2$ in \eqref{eq:deg6-omega} are constants, and (W10) becomes
\[
r\bigl(dT^2-10cby^4T+25c^2ay^8\bigr)=y^2\bigl(u_1T-5cf_1y^4\bigr)^2 .
\]
The left side is odd and the right side even, so both vanish, and Lemma~\ref{lem:deg6-div} gives $a=b=d=0$. In particular $\Phi_{2,zz}=0$.

Since $\Psi$, $\Phi_5$, $\Phi_4$, $\Phi_3$, $\Phi_2$ are affine in $z$, and $\Phi_1,\Phi_0$ have weight at most one, $\partial_z^2\phi=0$, and $\phi$ has an affine direction. In both cases Theorem~\ref{thm:main-dir} applies. In its second alternative the planar Keller map has degree at most five and is not injective, which contradicts \cite{Moh83}.
\end{proof}

\begin{remark}\label{rem:deg6}
The forms of Lemma~\ref{lem:deg6-nf} occur. The polynomial $\Psi=\frac16(xy+t)^3+y^5z$ is of this form with $\kappa_0=1$, $\gamma=\frac16$, $c=1$ and $a_0=c_1=c_2=c_3=0$. It has $\det\Hess\Psi=0$ and $\delta=-25xy^{11}$ for $\phi_6=\Psi|_{t=0}$. Here $e_z$ is an affine direction of $\Psi$, and $e_t$ is not. The proof of Theorem~\ref{thm:deg6} uses the homogeneous parts of $\phi$ for the single weight vector $(1,1,1,2)$, and no paths of weights. For top forms of Hessian rank two and one see Theorems~\ref{thm:deg6r2} and~\ref{thm:deg6r1}, and for rank three in higher degree Proposition~\ref{prop:rk3-weights}.
\end{remark}

\subsection{Degree six, top form of rank two}

Let $\deg\phi=6$, $\det\Hess\phi\in\kk^\times$, and let $\Hess\phi_6$ have rank two. By \ref{E:GN}, applied twice, $\phi_6=h(y)$ after a linear change, with $y=(x_1,x_2)$ and $A=\Hess_yh$, $\det A\neq0$. Put $z=(x_3,x_4)$, $P_1=\Hess_y\phi_5$, $Q_k=(\partial_{y_i}\partial_{z_j}\phi_{6-k})_{i,j}$ and
\[
S_1=\Hess_z\phi_5,\qquad S_2=\Hess_z\phi_4-Q_1\T A^{-1}Q_1 .
\]
Let $\phi_\varepsilon=\sum_j\varepsilon^{6-j}\phi_j$. Then $\det\Hess\phi_\varepsilon=\sum_k\varepsilon^k[\Delta]_{16-k}$, and the Schur complement of the block $\Hess_y\phi_\varepsilon=A+\varepsilon P_1+O(\varepsilon^2)$ is $\varepsilon S_1+\varepsilon^2S_2+O(\varepsilon^3)$. Hence
\begin{equation}\label{eq:deg6r2-layers}
[\Delta]_{14}=\det A\det S_1,\qquad [\Delta]_{13}=\det A\operatorname{tr}(\adj(S_1)S_2)+\operatorname{tr}(\adj(A)P_1)\det S_1,
\end{equation}
and $[\Delta]_{12}=\det A\det S_2$ if $S_1=0$.

\begin{lemma}\label{lem:deg6r2-P}%
Let $P\in\kk[y_1,y_2,t]$ be homogeneous of weight $d\ge3$ for the weights $(1,1,2)$, with $\det\Hess_yP(y,0)\neq0$ and $\det\Hess P=0$. Then $\partial_tP=0$, or after a linear change of $y$, $P=\alpha_1y_1^d+\alpha_2y_1^{d-1}y_2+\alpha y_1^{d-2}t$ with $\alpha_1,\alpha_2,\alpha\in\kk$. In particular $\deg_tP\le1$.
\end{lemma}
\begin{proof}
Put $h=P(y,0)$ and $V_0=\{v\in\kk^3:D_v^2P=0\}$. By \ref{E:dB62}, $\partial_vP\in\kk$ for some $v=(v_y,v_t)\neq0$, or $P(Tx)=a_1(x_1)+a_2(x_1)x_2+a_3(x_1)x_3$ for some $T\in\mathrm{GL}_3(\kk)$. In the first case $v_y\T\grad_yP$ and $v_t\partial_tP$ are homogeneous of the positive weights $d-1$ and $d-2$, so both vanish, and $v_y\T\grad h=0$. If $v_y\neq0$, then $h$ is a polynomial in one linear form and $\det\Hess h=0$. Hence $\partial_tP=0$. In the second case we may assume that $a_2',a_3'$ are linearly independent over $\kk$, since otherwise $\partial_vP\in\kk$ for some $v\in T\operatorname{span}(e_2,e_3)\setminus\{0\}$. Then $a_2',a_3'$ are not both constant, and $D_v^2(P\circ T)=v_1^2(a_1''+a_2''x_2+a_3''x_3)+2v_1(v_2a_2'+v_3a_3')$ gives $V_0=T\operatorname{span}(e_2,e_3)$. Since $P\circ\sigma=2^dP$ for $\sigma=\operatorname{diag}(2,2,4)$, $\sigma V_0=V_0$, and $V_0$ is the sum of its intersections with the eigenspaces $\kk^2\times\{0\}$ and $\kk e_t$ of $\sigma$. If $V_0=\kk^2\times\{0\}$, then $\Hess h=0$. Therefore $V_0=\kk(u,0)\oplus\kk e_t$ with $u\in\kk^2\setminus\{0\}$. If $\ell_1,\ell_2$ are linear forms in $y$ with $\ell_1(u)=0$ and $\ell_2(u)=1$, then $P\in\kk[\ell_1]+\kk[\ell_1]\ell_2+\kk[\ell_1]t$, and homogeneity gives the claim.
\end{proof}

We give $y_1,y_2,x_3$ the weight $1$ and $t=x_4$ the weight $2$, and write $[\,\cdot\,]_w$ for the part of weight $w$. If $\deg_t\phi_{6-m}\le m$ for $m\le2$, then every term of $\phi$ has weight at most six, every term of $\partial_i\partial_j\phi$ has weight at most $6-w_i-w_j$, the parts of weight $6-w_i-w_j$ are the entries of $\Hess\Psi$, $\Psi=[\phi]_6$, and the parts of $\det\Hess\phi$ of positive weight vanish. If moreover $x_3$ does not occur in $\Psi$, the row and column of $x_3$ in $\Hess\Psi$ vanish, and with $\Phi_w=[\phi]_w$, $H'=\Hess_{(y,t)}\Psi$, $q=\grad_{(y,t)}\partial_3\Phi_5$,
\begin{equation}\label{eq:deg6r2-w}
[\det\Hess\phi]_{13}=\partial_3^2\Phi_5\det H',\qquad [\det\Hess\phi]_{12}=\partial_3^2\Phi_4\det H'-q\T\adj(H')q\ \text{ if }\ \partial_3^2\Phi_5=0 .
\end{equation}

\begin{theorem}\label{thm:deg6r2}%
Let $\phi\in\kk[x_1,\dots,x_4]$ be a Keller potential of degree six whose top form has Hessian rank two. Then $D_v^3\phi=0$ for some $v\neq0$, and $\grad\phi$ is a polynomial automorphism.
\end{theorem}
\begin{proof}
By \eqref{eq:deg6r2-layers}, $\det S_1=0$. As in Branch~B of the proof of Theorem~\ref{thm:deg5}, \ref{E:dBvdE} over $\kk(y)$ and Lemma~\ref{lem:homlambda} give coprime forms $\lambda_3,\lambda_4$ of degree $s$ and a form $\gamma'$ with $\lambda_3\partial_3\phi_5+\lambda_4\partial_4\phi_5=\gamma'$, and $\phi_5=\sum_jP_j(y)\zeta^j+r_3(y)x_3+r_4(y)x_4$ with $\zeta=\lambda_4x_3-\lambda_3x_4$ and forms $P_j$ of degree $5-j(s+1)$. Thus $S_1=G''\,ww\T$ with $w=(\lambda_4,-\lambda_3)$ and $G=\sum_jP_j\zeta^j$.

\emph{$S_1\neq0$, $s\ge1$.} Then $P_3=0$, $P_2\neq0$, $s=1$, and $P_2$ is linear. After a linear change of $y$, $\lambda=y$ and $w=(y_2,-y_1)$. Since $\adj(S_1)=2P_2\lambda\lambda\T$, \eqref{eq:deg6r2-layers} gives $\lambda\T\Hess_z\phi_4\lambda=v\T A^{-1}v$ with $v=Q_1\lambda=\grad\gamma'-\grad_z\phi_5$. The part of degree two in $z$ of the right side is $4P_2^2\,w\T A^{-1}w\,\zeta^2=\rho\zeta^2$ with $\rho=120P_2^2h/\det A\neq0$, since $w\T\adj(A)w=y\T Ay=30h$. On the left it is $y\T\Hess f(z)\,y$, where $f\in\kk[z]$ is the part of $\phi_4$ of degree four in $z$. Setting $z=e_3$ and $z=e_4$ shows $\rho\in\kk$. Hence $\Hess f=\rho\,uu\T$ with $u=(x_4,-x_3)$, and $12f=z\T\Hess f\,z=\rho(u\T z)^2=0$, so $\rho=0$, a contradiction.

\emph{$S_1\neq0$, $s=0$.} After a linear change of $z$, $\phi_5=\Phi(y,x_3)+\gamma'(y)t$ with $a=\partial_3^2\Phi\neq0$, and \eqref{eq:deg6r2-layers} gives $\partial_t^2\phi_4=(\grad\gamma')\T A^{-1}\grad\gamma'\in\kk[y]$. So the weights apply, and $\Psi=h+t\gamma'+\frac12t^2\partial_t^2\phi_4+\gamma t^3$, where $\gamma$ is the coefficient of $t^3$ in $\phi_3$, does not contain $x_3$. Here $\partial_3^2\Phi_5=a+t\,\partial_3^2\beta$, with $\beta$ the coefficient of $t$ in $\phi_4$, is nonzero, and \eqref{eq:deg6r2-w} gives $\det H'=0$.

\emph{$S_1=0$.} Then $\phi_5=P_0(y)+r_3x_3+r_4x_4$, $Q_1$ does not depend on $z$, and $f=\phi_4-\frac12z\T Nz$, $N=Q_1\T A^{-1}Q_1$, satisfies $\det\Hess_zf=\det S_2=0$. By \ref{E:dBvdE} over $\kk(y)$ there are $\kk(y)$-linear coordinates $(\ell_1,\ell_2)$ on $z$ with $f=G_1(\ell_1)+\kappa\ell_2$, $G_1\in\kk(y)[T]$, $\kappa\in\kk(y)$. If $\deg G_1\le2$, then $\deg_z\phi_4\le2$, so $D_v^3\phi=6\phi_3(v)$ for $v\in\operatorname{span}(e_3,e_4)$, and the binary cubic $\phi_3(0,0,x_3,x_4)$ has a zero $v\neq0$. Let $\deg G_1\ge3$. The parts of $\phi_4$ of degrees $4$ and $3$ in $z$ are $e_4\ell_1^4$ and $e_3\ell_1^3$, $e_i\in\kk(y)$, and their coefficients are forms in $y$ of degrees $0$ and $1$. If $e_4\neq0$, $\ell_1$ is proportional to a linear form over $\kk$. If $e_4=0$, write $e_3\ell_1^3=\sum_ja_jx_3^jx_4^{3-j}$. If $a_3\neq0$, then $3a_1a_3=a_2^2$, so $a_2\in\kk a_3$ and again $\ell_1$ is proportional to a linear form over $\kk$, and if $a_3=0$, then $\ell_1\in\kk(y)x_4$. After a linear change of $z$ we may take $\ell_1=x_3$. Then $\phi_4=\frac12z\T Nz+G_2(y,x_3)+\kappa'(y)t$ with $\deg_{x_3}G_2\ge3$, and $N_{34}=\partial_3\partial_t\phi_4$, $N_{44}=\partial_t^2\phi_4$ lie in $\kk[y]$. The weights apply, $\Psi=h+tr_4+\frac12t^2N_{44}+\gamma t^3$ does not contain $x_3$, $\Phi_5=P_0+r_3x_3+t(N_{34}x_3+\kappa')+t^2m$ with $m$ linear, $\partial_3^2\Phi_5=0$, and $q$ does not depend on $x_3$. Since $\partial_3^3\Phi_4=\partial_3^3G_2\neq0$, the derivative in $x_3$ of \eqref{eq:deg6r2-w} gives $\det H'=0$.

In the last two cases Lemma~\ref{lem:deg6r2-P} with $d=6$ gives $\deg_t\Psi\le1$, so $\gamma=0$ and $D_t^3\phi=6\gamma=0$. The second assertion follows from Theorem~\ref{thm:main-dir} and \cite{Moh83} as in the proof of Theorem~\ref{thm:deg6}.
\end{proof}

\begin{remark}\label{rem:deg6r2}
In the case $S_1\neq0$, $s\ge1$, the kernel of $S_1$ is spanned by $\lambda=y$, and $\lambda_3\partial_3+\lambda_4\partial_4$ plays the role of $\partial_t$. The argument shows that this case does not occur, already by $[\Delta]_{13}$. In the other cases the weights $(1,1,1,2)$ are used only after the layers $[\Delta]_{14}$, $[\Delta]_{13}$ or $[\Delta]_{14}$, $[\Delta]_{12}$ have shown $\deg_t\phi_{6-m}\le m$ for $m\le2$. The case of a top form of Hessian rank one is Theorem~\ref{thm:deg6r1}.
\end{remark}

\subsection{Degree six, top form of rank one}

Let $\deg\phi=6$, $\det\Hess\phi=c\in\kk^\times$, and let $\Hess\phi_6$ have rank one. By Lemma~\ref{lem:app-low}, $\phi_6=x^6$ after a linear change, with $x=x_1$. In this subsection $Y=(x_2,x_3,x_4)$, $z=(x_3,x_4)$, $H_j=\Hess_Y\phi_j$, $q_j=\partial_x\grad_Y\phi_j$, $p_5=\partial_x^2\phi_5$, and
\[
S_1=H_5,\qquad S_2=H_4-\frac{q_5q_5\T}{30x^4},\qquad S_3=H_3-\frac{q_5q_4\T+q_4q_5\T}{30x^4}+\frac{p_5\,q_5q_5\T}{900x^8}.
\]
With $\phi_\varepsilon$ as in the previous subsection, the Schur complement of the entry $\partial_x^2\phi_\varepsilon=30x^4+O(\varepsilon)$ of $\Hess\phi_\varepsilon$ is $\varepsilon S_1+\varepsilon^2S_2+\varepsilon^3S_3+O(\varepsilon^4)$, and its determinant is $O(\varepsilon^{16})$. In particular
\begin{equation}\label{eq:deg6r1-layers}
[\Delta]_{13}=30x^4\det S_1,\qquad [\Delta]_{12}=30x^4\operatorname{tr}(\adj(S_1)S_2)\ \text{ if }\ \det S_1=0 .
\end{equation}

\begin{lemma}\label{lem:deg6r1-plane}%
Let $U\subset\operatorname{span}(e_2,e_3,e_4)$ be a plane with $D_uD_{u'}\phi_5=0$ and $D_uD_{u'}D_{u''}\phi_4=0$ for $u,u',u''\in U$. Then $D_v^3\phi=0$ for some $v\in U\setminus\{0\}$.
\end{lemma}
\begin{proof}
For $v\in U$ one has $D_v^3\phi=D_v^3\phi_3=6\phi_3(v)$, and the binary cubic $\phi_3|_U$ has a zero $v\neq0$.
\end{proof}

\begin{lemma}\label{lem:deg6r1-P}%
Let $P=x^6+ax^4t+bx^2t^2+\gamma t^3$ with $a,b,\gamma\in\kk$, and $\delta_P=\det\Hess_{(x,t)}P$. If $\delta_P=0$, then $a=b=\gamma=0$. If $b=\frac4{15}a^2$, then
\[
\delta_P=\bigl(180\gamma-\tfrac{32}{15}a^3\bigr)x^4t+\bigl(72a\gamma-\tfrac{64}{75}a^4\bigr)x^2t^2+\tfrac{16}5a^2\gamma\,t^3 ,
\]
and the first two coefficients vanish for $\gamma=\frac8{675}a^3$.
\end{lemma}
\begin{proof}
The coefficients of $x^6$, $t^3$, $x^2t^2$, $x^4t$ in $\delta_P$ are $-4(4a^2-15b)$, $12b\gamma$, $12(6a\gamma-b^2)$, $4(45\gamma-2ab)$. The second assertion is a direct computation.
\end{proof}

We use the weights of the previous subsection, $1$ for $x,x_2,x_3$ and $2$ for a coordinate $t$, and the parts $\Phi_w$ of weight $w$.

\begin{lemma}\label{lem:deg6r1-W}%
Let $\deg_t\phi_{6-m}\le m$ for $m\le2$ and $\Psi=\Phi_6=P(x,t)$ as in Lemma~\ref{lem:deg6r1-P}. Put $u=(x_2,x_3)$, $H'=\Hess_{(x,t)}\Psi$, $U_w=\Hess_u\Phi_w$ and $m_i=\grad_{(x,t)}\partial_{x_i}\Phi_5$. Then $[\det\Hess\phi]_{14}=[\det\Hess\phi]_{13}=0$ and $[\det\Hess\phi]_{12}=\delta_P\det U_5$. If $\partial_3^2\Phi_5=\partial_2\partial_3\Phi_5=0$, then
\[
[\det\Hess\phi]_{11}=\partial_2^2\Phi_5\,\bigl(\delta_P\,\partial_3^2\Phi_4-m_3\T\adj(H')m_3\bigr).
\]
If $u$ does not occur in $\Phi_5$, then $[\det\Hess\phi]_{10}=\delta_P\det U_4$.
\end{lemma}
\begin{proof}
As before the parts of positive weight of $\det\Hess\phi$ vanish. Multiply the entries of $\Hess\Phi_w$ by $\varepsilon^{6-w}$. The blocks of $u$, of $(u,(x,t))$ and of $(x,t)$ are $\varepsilon U_5+\varepsilon^2U_4+O(\varepsilon^3)$, $\varepsilon M_5+O(\varepsilon^2)$ with rows $m_i\T$, and $H'+O(\varepsilon)$, since $u$ does not occur in $\Psi$. The Schur complement of the block of $(x,t)$ is $\varepsilon U_5+\varepsilon^2(U_4-M_5H'^{-1}M_5\T)+O(\varepsilon^3)$, which gives the formulas.
\end{proof}

\begin{theorem}\label{thm:deg6r1}%
Let $\phi\in\kk[x_1,\dots,x_4]$ be a Keller potential of degree six whose top form has Hessian rank one. Then $D_v^3\phi=0$ for some $v\neq0$, and $\grad\phi$ is a polynomial automorphism.
\end{theorem}
\begin{proof}
By \eqref{eq:deg6r1-layers}, $\det\Hess_Y\phi_5=0$ in $\kk(x)[Y]$. Apply \ref{E:dB62} over $\kk(x)$. In its second case, as in Branch~A of the proof of Theorem~\ref{thm:deg5}, the plane of \ref{E:dB62} specializes to a constant plane, and after a linear change of $Y$, $\phi_5=a_1+a_2x_3+a_3x_4$ with $a_i\in\kk[x,x_2]$. In the first case Lemma~\ref{lem:homlambda} with $r=1$ gives, after a linear change of $Y$, $\partial_4\phi_5\in\kk x^4$.

\emph{$\phi_5=a_1+a_2x_3+a_3x_4$, and $\partial_v\phi_5\notin\kk x^4$ for every constant $v\in\operatorname{span}(e_3,e_4)\setminus\{0\}$.} Put $\kappa=(0,\partial_2a_3,-\partial_2a_2)=\rho\,(0,k_3,k_4)$ with coprime forms $k_3,k_4\in\kk[x,x_2]$ of degree $d$. By hypothesis $d\ge1$. Since $S_1\kappa=0$ and $\adj(S_1)=-\kappa\kappa\T$, \eqref{eq:deg6r1-layers} gives
\[
D_\kappa^2\phi_4=\frac{(\kappa\T q_5)^2}{30x^4},\qquad \kappa\T q_5=\partial_2a_3\,\partial_xa_2-\partial_2a_2\,\partial_xa_3\in\kk[x,x_2].
\]
Let $F_k$ be the part of $\phi_4$ of degree $k$ in $z$. Then $D_\kappa^2F_k=0$ for $k=3,4$, so $F_k$ has degree at most one along $(k_3,k_4)$ over $\kk(x,x_2)$, and $\ell^{k-1}$ divides $F_k$ for $\ell=k_4x_3-k_3x_4$. As $\ell$ is primitive, Gauss's lemma gives $F_k=\ell^{k-1}L_k$ with $L_k\in\kk[x,x_2][z]$, and since the coefficients of $F_k$ are forms of degree $4-k<(k-1)d$, $F_k=0$. So Lemma~\ref{lem:deg6r1-plane} applies to $U=\operatorname{span}(e_3,e_4)$.

\emph{$\phi_5=\psi(x,x_2,x_3)+\gamma'x^4x_4$, $\det\Hess_{(x_2,x_3)}\psi\neq0$.} Then $\adj(S_1)=\det\Hess_{(x_2,x_3)}\psi\cdot e_4e_4\T$, and \eqref{eq:deg6r1-layers} gives $\partial_4^2\phi_4=\frac8{15}\gamma'^2x^2$. With $t=x_4$ the weights apply, $\Psi=P$ with $a=\gamma'$ and $b=\frac4{15}\gamma'^2$, and the rows of $x_2,x_3$ in $\Hess\Psi$ vanish. Hence $[\det\Hess\phi]_{12}=\det\Hess_{(x_2,x_3)}\Phi_5\cdot\delta_P$, and the first factor equals $\det\Hess_{(x_2,x_3)}\psi$ at $t=0$. So $\delta_P=0$, Lemma~\ref{lem:deg6r1-P} gives $\gamma=0$, and $D_t^3\phi=6\gamma=0$.

\emph{$\phi_5=\psi(x,x_2,x_3)+\gamma'x^4x_4$, $\det\Hess_{(x_2,x_3)}\psi=0$.} By \ref{E:dBvdE} over $\kk(x)$ and Lemma~\ref{lem:homlambda}, after a linear change of $Y$,
\[
\phi_5=\psi_1(x,x_2)+x^4\ell(z),\qquad f=\phi_4-\tfrac4{15}x^2\ell^2,
\]
with $\ell$ a linear form.

\emph{(i) $\eta=\partial_2^2\psi_1\neq0$.} Let $y=(x,x_2)$. Here $S_1=\eta\,e_2e_2\T$ and $(S_2)_{zz}=\Hess_zf$. With $T_2=(S_2)_{zz}$ and $T_3=(S_3)_{zz}-(S_2)_{z2}(S_2)_{2z}/\eta$, the Schur complement of the entry $\eta\varepsilon+O(\varepsilon^2)$ gives $[\Delta]_{12}=0$, $[\Delta]_{11}=30x^4\eta\det T_2$, and $[\Delta]_{10}=30x^4\eta\operatorname{tr}(\adj(T_2)T_3)$ if $\det T_2=0$. By \ref{E:dBvdE} over $\kk(y)$ and Lemma~\ref{lem:homlambda} there are coprime forms $\lambda_3,\lambda_4\in\kk[y]$ of degree $s$, forms $P_j\in\kk[y]$ of degree $4-j(s+1)$ and $r_3,r_4\in\kk[y]$ with $f=\sum_jP_j\zeta^j+r_3x_3+r_4x_4$, $\zeta=\lambda_4x_3-\lambda_3x_4$. If $\Hess_zf=0$, which holds for $s\ge2$, then $\deg_z\phi_4\le2$, and Lemma~\ref{lem:deg6r1-plane} applies to the plane of $z$. Let $s=1$ and $P_2\neq0$. Then $P_2\in\kk$, $\adj(T_2)=2P_2\lambda\lambda\T$, and $\lambda\T(S_3)_{zz}\lambda$ has degree at most one in $z$, since $\deg_z\phi_j\le6-j$ for $j=5,4,3$. Further $\lambda\T(S_2)_{z2}=J(2P_2\zeta+P_1)+\lambda\T\partial_2r-\frac2{15x}(\lambda\T\grad\ell)\,\partial_x\partial_2\psi_1$ with $J=\lambda_3\partial_2\lambda_4-\lambda_4\partial_2\lambda_3\in\kk^\times x$. So the part of $[\Delta]_{10}$ of degree two in $z$ is $-240x^4P_2^3J^2\zeta^2\neq0$, a contradiction. Let $s=0$ and $\Hess_zf\neq0$. After a linear change of $z$, $f=G(y,x_3)+\kappa(y)t$ with $t=x_4$ and $\partial_3^2G\neq0$. Put $\alpha=\ell(e_4)$. Then $\partial_t^2\phi_4=\frac8{15}\alpha^2x^2$, the weights apply, and $\Psi=P$ with $a=\alpha$, $b=\frac4{15}\alpha^2$. If $\alpha\neq0$, choose $x_3$ with $\ell(e_3)=0$. Then $\Phi_5=\psi_1+\kappa t+nt^2$ with $n$ linear, $m_3=(0,2n_3t)$ with $n_3=\partial_3n$, $\partial_2^2\Phi_5=\eta+t\,\partial_2^2\kappa\neq0$, $\Phi_4|_{t=0}=G$, and Lemma~\ref{lem:deg6r1-W} gives $\delta_P\,\partial_3^2\Phi_4=4n_3^2t^2\,\partial_x^2\Psi$. By Lemma~\ref{lem:deg6r1-P} either the coefficient of $t$ in $\delta_P$ is nonzero, and then the coefficient of $t$ gives $\partial_3^2G=0$, or $\delta_P=\frac{16}5\alpha^2\gamma t^3$ with $\gamma\neq0$, and then the coefficients of $t^2$ and $t^3$ give $n_3=0$ and $\partial_3^2G=0$. If $\alpha=0$, then $\ell=\gamma_3x_3$ and $\Psi=x^6+\gamma t^3$. For $\gamma\neq0$, $\delta_P=180\gamma x^4t$, $m_3=(4\gamma_3x^3,2n_3t)$, and Lemma~\ref{lem:deg6r1-W} gives $\partial_3^2\Phi_4=\frac8{15}\gamma_3^2x^2+\frac{2n_3^2}{3\gamma}t$. At $t=0$, $\Phi_4=G+\frac4{15}\gamma_3^2x^2x_3^2$, so $\partial_3^2G=0$. In all cases this contradicts $\partial_3^2G\neq0$, unless $\alpha=\gamma=0$, and then $\deg_t\phi\le2$.

\emph{(ii) $\partial_2^2\psi_1=0$.} Then $\phi_5=c_5x^5+x^4\ell(Y)$ with $\ell$ a linear form on $Y$, $S_1=0$, $S_2=\Hess_Yf$, $[\Delta]_{13}=[\Delta]_{12}=[\Delta]_{11}=0$, $[\Delta]_{10}=30x^4\det S_2$, and $[\Delta]_9=30x^4\operatorname{tr}(\adj(S_2)S_3)$ if $\det S_2=0$. Apply \ref{E:dB62} over $\kk(x)$ to $f$. In its second case $f$ is affine along a constant plane $U$, and Lemma~\ref{lem:deg6r1-plane} applies to $U$. In the first case Lemma~\ref{lem:homlambda} gives a constant $v_0$ with $\partial_{v_0}f\in\kk x^3$. In coordinates $(u,t)$ of $Y$ with $v_0=e_t$, $f=\psi'(x,u)+\beta x^3t$. If $\rk\Hess_u\psi'\le1$, \ref{E:dBvdE} over $\kk(x)$ and Lemma~\ref{lem:homlambda} give a constant $w$ with $\partial_w\psi'\in\kk x^3$, and Lemma~\ref{lem:deg6r1-plane} applies to $\operatorname{span}(w,v_0)$. Let $\delta'=\det\Hess_u\psi'\neq0$. Then $\adj(S_2)=\delta'e_te_t\T$ and $(S_3)_{tt}=0$. Put $\alpha=\ell(v_0)$. If $\alpha=0$, the $t$-component of $q_5$ vanishes, $(S_3)_{tt}=\partial_t^2\phi_3=0$, and $\deg_t\phi\le2$. If $\alpha\neq0$, the change $t\mapsto t+\alpha^{-1}\ell|_{t=0}$ keeps $v_0$ and the form of $f$, so that $\ell=\alpha t$. Then $(S_3)_{tt}=0$ gives
\[
\partial_t^2\phi_3=\tfrac{16}{225}\alpha^3t+\bigl(\tfrac45\alpha\beta-\tfrac{16}{45}\alpha^2c_5\bigr)x ,
\]
so $\gamma=\frac8{675}\alpha^3$. The weights apply, $\Psi=P$ with $a=\alpha$, $b=\frac4{15}\alpha^2$, $\delta_P=\frac{16}5\alpha^2\gamma t^3\neq0$, and $\Phi_5=c_5x^5+\beta x^3t+\frac12t^2\,\partial_t^2\phi_3|_{t=0}$ does not contain $u$. Lemma~\ref{lem:deg6r1-W} gives $\det\Hess_u\Phi_4=0$. At $t=0$, $\Phi_4=\psi'$, so $\delta'=0$, a contradiction.

The second assertion follows from Theorem~\ref{thm:main-dir} and \cite{Moh83} as in the proof of Theorem~\ref{thm:deg6}.
\end{proof}

\begin{corollary}\label{cor:deg6}%
Let $\kk$ be algebraically closed. Every Keller potential $\phi\in\kk[x_1,\dots,x_4]$ of degree at most six satisfies $\deg_v\phi\le2$ for some $v\neq0$, and $\grad\phi$ is a polynomial automorphism.
\end{corollary}

\begin{proof}
For $\deg\phi\le5$ this is Theorem~\ref{thm:main-low}(ii) with Theorem~\ref{thm:main-dir} and \cite{Moh83}. For $\deg\phi=6$, $[\Delta]_{16}=\det\Hess\phi_6=0$, and $\Hess\phi_6\neq0$ in characteristic zero, so Theorems~\ref{thm:deg6}, \ref{thm:deg6r2} and~\ref{thm:deg6r1} apply.
\end{proof}

By Lemma~\ref{lem:descent}, the gradient of a Keller potential of degree at most six in four variables is a polynomial automorphism over every field of characteristic zero. For degree at most four the bound $\deg(\grad\phi)^{-1}\le27$ of Corollary~\ref{cor:deg4} holds over every such field as well, since the inverse over the algebraic closure is defined over the field itself. More generally, if $\grad\phi$ is a polynomial automorphism of $\kk^4$, then $\deg(\grad\phi)^{-1}\le(\deg\phi-1)^3$ by \cite{BCW82}, which gives the bounds $64$ for $\deg\phi\le5$ and $125$ for $\deg\phi\le6$.

\begin{remark}\label{rem:deg6r1}
The proof of Theorem~\ref{thm:deg6r1} uses the weights $(1,1,1,2)$ only after the layers have shown $\deg_t\phi_{6-m}\le m$ for $m\le2$. In the last case the layer $[\Delta]_9$ forces $\gamma=\frac8{675}\alpha^3$, the value for which $t^3$ divides $\delta_P$, and the contradiction comes from the part of weight ten. In the case of a twisted kernel, as in Theorem~\ref{thm:deg6r2}, the contradiction comes from the part of highest degree in $z$ of a single layer.
\end{remark}

\subsection{Top form of rank three in higher degree}

The first step of the proof of Theorem~\ref{thm:deg6}, the passage to the weights $(1,1,1,2)$, holds in every degree. Let $\deg\phi=d\ge3$ and $\det\Hess\phi=c\in\kk^\times$, and as above $\phi_d\in\kk[x']$ and $\delta=\det\Hess_{x'}\phi_d\neq0$. Put $t=x_4$, $A=\Hess_{x'}\phi_d$, and let $b\in\kk[x']$ be the coefficient of $t$ in $\phi_{d-1}$.

\begin{proposition}\label{prop:rk3-weights}%
\textup{(a)} $\deg_t\phi_{d-m}\le m$ for $1\le m\le d-2$. Hence, for the weights $(1,1,1,2)$, every term of $\phi$ has weight at most $d$. The part of weight $d$ is $\Psi=\sum_{m\ge0}t^mc_m$, where $c_m\in\kk[x']$ is the coefficient of $t^m$ in $\phi_{d-m}$, and $\det\Hess\Psi=0$, $\rk\Hess\Psi=3$.

\textup{(b)} With $w=A^{-1}\grad_{x'}b$,
\[
2c_2=\grad_{x'}b\T A^{-1}\grad_{x'}b,\qquad 6c_3=3\,w\T\Hess_{x'}b\;w-2\,D^3\phi_d[w,w,w].
\]
\end{proposition}
\begin{proof}
Write $z=x'$, $h=\phi_d$, $P_m=\phi_{d-m}$, and $\Phi_s=s^d\phi(x/s)=\sum_{m\ge0}s^mP_m\in\kk[s,z,t]$. Then $\Hess\Phi_s(x)=s^{d-2}\Hess\phi(x/s)$ and $\det\Hess\Phi_s=cs^\nu$, $\nu=4(d-2)$. Let $L=\kk(z)$ and work in $L[t][[s]]$. There is a unique $\xi=z+\sum_{k\ge1}s^k\xi_k$ with $\xi_k\in L[t]^3$ and $\grad_z\Phi_s(\xi,t)=\grad h(z)$. Indeed, the coefficient of $s^k$ of this equation reads $A\xi_k=-\grad_zP_k(z,t)+G'_k$, where $G'_k$ is a sum of products of derivatives in $z$ of $h$ and of $P_j$, $j<k$, at $(z,t)$ with components of $\xi_i$, $i<k$, and $A\in\mathrm{GL}_3(L)$. Put
\[
u=\vv{\grad h(z),\xi}-\Phi_s(\xi,t)=\sum_{m\ge0}s^mu_m .
\]
Since $\grad_z\Phi_s(\xi,t)=\grad h(z)$, one has $\partial_tu=-(\partial_t\Phi_s)(\xi,t)$, and differentiation of $\grad_z\Phi_s(\xi,t)=\grad h(z)$ in $t$ gives $\partial_t\xi=-A_s^{-1}\beta_s$, where $A_s=\Hess_z\Phi_s$, $\beta_s=\grad_z\partial_t\Phi_s$ and $\eta_s=\partial_t^2\Phi_s$ are evaluated at $(\xi,t)$. By the Schur complement, as in the proof of Lemma~\ref{lem:legendre}(b),
\[
\partial_t^2u=-\bigl(\eta_s-\beta_s\T A_s^{-1}\beta_s\bigr)=-\frac{cs^\nu}{\det A_s(\xi,t)} .
\]
Here $\det A_s(\xi,t)\equiv\det A$ modulo $s$ is a unit of $L[t][[s]]$. Hence $\partial_t^2u_m=0$ for $m<\nu$.

(a) The Taylor expansion of $h$ and of the $P_j$ at $z$ gives $u_m=-P_m(z,t)+G_m$ for $m\ge1$, where every term of $G_m$ is a product of a derivative of $h$ with $\xi_{k_1},\dots,\xi_{k_r}$, $r\ge2$, $k_1+\dots+k_r=m$, or of a derivative of some $P_j$, $j<m$, with $\xi_{k_1},\dots,\xi_{k_r}$, $j+k_1+\dots+k_r=m$. Suppose $\deg_tP_j\le j$ and $\deg_t\xi_j\le j$ for $j<m\le d-2$. Then $\deg_tG_m\le m$, and since $m<\nu$, $\deg_tu_m\le1$. So $\deg_tP_m\le m$, and the equation for $\xi_m$ gives $\deg_t\xi_m\le m$. This proves the first claim by induction. A term of $\phi_{d-m}$ of degree $k$ in $t$ has weight $d-m+k$, and $\phi_1,\phi_0$ have weight at most two. As in the proof of Proposition~\ref{prop:deg6}(c), the part of weight $4d-10>0$ of $\det\Hess\phi$ is $\det\Hess\Psi$, and the principal minor of $\Hess\Psi$ in $x'$ is $\delta\neq0$ at $t=0$.

(b) One has $\xi_1=-A^{-1}\grad_zP_1$, $u_2=-P_2+\frac12\grad_zP_1\T A^{-1}\grad_zP_1$ and $u_3=-P_3+\grad_zP_2\T w_1-\frac12w_1\T(\Hess_zP_1)w_1+\frac16D^3h[w_1,w_1,w_1]$ with $w_1=A^{-1}\grad_zP_1$. By (a), $P_1=a+tb$ and $w_1=A^{-1}\grad_za+tw$. The coefficient of $t^2$ in $u_2$ gives the first formula. The coefficient of $t^3$ in $u_3$ is $-c_3+\grad_zc_2\T w-\frac12w\T(\Hess_zb)w+\frac16D^3h[w,w,w]$, and the first formula gives $\grad_zc_2\T w=w\T(\Hess_zb)w-\frac12D^3h[w,w,w]$.
\end{proof}

\begin{remark}\label{rem:rk3-weights}
For $d=6$ part (a) is Proposition~\ref{prop:deg6}(a) and (c), and (b) gives in addition $\gamma=c_3=\frac12w\T(\Hess_{x'}b_4)w-\frac13D^3\phi_6[w,w,w]$. For the polynomial of Remark~\ref{rem:deg6} this gives $\gamma=\frac16$. The same argument applied to $\Psi$, for which $\det\Hess\Psi=0$, gives $\partial_t^2u_m=0$ for every $m$, so all $c_m$, $m\ge2$, are determined by $\phi_d$ and $b$. In degree $d\ge7$ a continuation of this method would need a classification of these weighted forms analogous to Lemmas~\ref{lem:deg6-dir} and~\ref{lem:deg6-nf}. For the dichotomy it is replaced by Theorem~\ref{thm:main-aff}.

Over $\R$, at a point where $\Hess_{x'}\phi_d$ is invertible, $u$ is the partial Legendre transform of $\Phi_s$ in the variables $x'$, written in the coordinates $p=\grad_{x'}\phi_d$, and $\Hess_pu=A_s(\xi,t)^{-1}$. The identity of the proof then reads $\partial_t^2u+cs^\nu\det\Hess_pu=0$, a Monge--Amp\`ere equation in $p$ with $t$ as evolution variable, degenerate at $s=0$.
\end{remark}

In Sections~\ref{sec:translation} to~\ref{sec:tree} the field $\kk$ is algebraically closed of characteristic zero.

\section{Translation directions of gradient images}
\label{sec:translation}

Three notions of direction occur in the proof. A kernel direction of $\Hess P$ lies in the source space and may vary with the point. A translation direction of the closure of the gradient image of $P$ lies in the target space and is constant. An affine direction of $\phi$ is a constant source vector along which $\phi$ is affine. The quasi-translations of Section~\ref{ssec:qt} connect the first two notions (Proposition~\ref{prop:qt-hess} and Example~\ref{ex:apex}), and the weighted analysis of Sections~\ref{sec:walls} and~\ref{sec:tree} produces the third.

The proof of Theorem~\ref{thm:main-aff} for $\deg\phi\ge3$ in Sections~\ref{sec:translation} to~\ref{sec:tree} has the following steps. (1) The highest homogeneous part $\phi_D$ of $\phi$ has vanishing Hessian determinant, and by \ref{E:GN} it depends, in suitable linear coordinates, on $r\le3$ variables, with nonsingular Hessian in these variables (proof of Proposition~\ref{prop:tree}). (2) The closure of the gradient image of a polynomial in four variables with Hessian of rank three has a nonzero translation direction (Theorem~\ref{thm:apex}). This is the content of the present section. (3) Starting from $\phi_D$, the weights are changed along paths. At a wall the minors and cofactors of the previous leading form (Lemma~\ref{lem:survive}) and the apex directions (Lemma~\ref{lem:apexbasic}) restrict the wall form. (4) Every wall form which occurs leads to one of the exit lemmas and tools of Section~\ref{sec:walls}, listed in Table~\ref{tab:exit}, or to a weight vector with positive entries and $W=0$, where Corollary~\ref{cor:endpoint} applies. The cases are listed at the beginning of Section~\ref{sec:tree}. (5) Every branch of this procedure is finite (Lemma~\ref{lem:termination}), and it ends in an affine direction of $\phi$ or in a contradiction.

\begin{remark}[attribution]\label{rem:attribution}
The dichotomy \textup{(D)} for all Keller potentials in four variables is the statement of \cite[Theorem~1.1]{Tan26}, and the strategy of the proof of Theorem~\ref{thm:main-aff} in Sections~\ref{sec:translation} to~\ref{sec:tree} is that of \cite{Tan26}. That preprint refers to \cite[Theorem~4.17]{AHL}, a criterion of invertibility by the Legendre potential, and to \cite[Remark~4.5]{AHL}, where it is observed that the argument of Meng does not give $\JC_2\Rightarrow\HC_4$. Complete proofs are given for the statements taken from that preprint. The preprint is distributed under the Creative Commons Attribution~4.0 International license. Where our text follows its arguments, the table below gives the source and the changes made.
\begin{center}\small
\begin{tabular}{>{\raggedright}p{4.3cm}>{\raggedright}p{3.7cm}>{\raggedright\arraybackslash}p{5.6cm}}
\hline
This paper & \cite{Tan26} & Relation\\
\hline
Corollary~\ref{cor:equiv}(a) & Theorem~1.1 & the dichotomy of Theorem~1.1, here over every field of characteristic zero, with exactly one alternative and a degree bound\\
Theorem~\ref{thm:main-aff} & none & not contained in \textup{(D)}, whose first alternative carries no information on directions, obtained from the strategy of \cite{Tan26} by removing the inversion arguments\\
Theorem~\ref{thm:cong}, Lemma~\ref{lem:slice}, Theorem~\ref{thm:apex} & Theorem~2.1, Lemma~2.2, Theorem~2.3 & same statements, rank-two case of Theorem~\ref{thm:cong} proved without \cite[Theorem~0.1]{DP04} (Lemma~\ref{lem:cong-dir})\\
Proposition~\ref{prop:qt-hess} & beginning of the proof of Lemma~2.2, proofs of Theorem~2.3 and Corollary~2.4 & same argument, formulated through quasi-translations\\
Lemmas~\ref{lem:cong-rk1}, \ref{lem:cong-lines} & proof of Theorem~2.1 & same argument\\
Lemmas~\ref{lem:cong-ker}, \ref{lem:cong-par} & proof of Theorem~2.1 (orthogonality of $\ker N$ and $\operatorname{im}N$, lines through a point at infinity) & same argument\\
Corollary~\ref{cor:dBproblem} & Corollary~2.4 & same statement and argument\\
Lemmas~\ref{lem:E6}, \ref{lem:E7} & uses \cite[Theorem~7.2, Corollary~7.3]{dB16} & proved here\\
Lemmas~\ref{lem:apexbasic}(a)--(c), \ref{lem:survive}, \ref{lem:shear}, \ref{lem:euler} & preliminaries of Section~3, beginning of the proof of Proposition~4.1 & same arguments, stated as lemmas\\
Lemmas~\ref{lem:D1}, \ref{lem:D2}, \ref{lem:D3}, \ref{lem:D4}--\ref{lem:D7} & Lemmas~3.2 to~3.8 & in this order\\
Lemma~\ref{lem:D1-walls} & proof of Lemma~3.2 & stated as a lemma\\
Theorem~\ref{thm:touch} & Lemma~3.1, end of Lemma~3.4 & replaces them, no inversion of the gradient\\
Lemmas~\ref{lem:P1}, \ref{lem:tree-rk2}, \ref{lem:tree-C1}, \ref{lem:tree-binary}, Proposition~\ref{prop:tree} & parts of Proposition~4.1 and of its proof & case analysis organized by Lemmas~\ref{lem:restrict} to~\ref{lem:K}, termination by Lemma~\ref{lem:termination}\\
\hline
\end{tabular}
\end{center}
Every branch of the weighted analysis ends in an affine direction of $\phi$ or in a contradiction, and invertibility enters only through Theorem~\ref{thm:affine}. Lemma~\ref{lem:2d} is used in Lemma~\ref{lem:euler}, and Lemma~\ref{lem:det2} in the determinant identities. The external inputs are the published results \ref{E:dB} to \ref{E:dBvdE}, and \ref{E:dB62} is proved in Appendix~\ref{app:dB62}.
\end{remark}

The statements of Theorem~\ref{thm:cong}, Lemma~\ref{lem:slice} and
Theorem~\ref{thm:apex} are due to Tang \cite[Theorem~2.1, Lemma~2.2, Theorem~2.3]{Tan26}, and so is Corollary~\ref{cor:dBproblem} \cite[Corollary~2.4]{Tan26}. We
formulate them through the quasi-translations of
Section~\ref{ssec:qt}. A polynomial in $n$ variables with Hessian of
rank $n-1$ defines a quasi-translation $Z+Q$ whose lines are contained
in the fibers of the gradient (Proposition~\ref{prop:qt-hess}), and a
translation direction of the gradient image is the same as a linear
relation among the components of $Q$. Theorem~\ref{thm:apex} says that
such a relation exists in four variables. In the rank-two case of
Theorem~\ref{thm:cong} the proof given here does not use the
classification of congruences of lines \cite[Theorem~0.1]{DP04}, on
which Tang's proof rests.
Ranks of polynomial and rational matrices are taken over the rational
function field in the source coordinates. We write
$X=(X_1,X_2,X_3)$ or $(X_1,\dots,X_4)$, $e_j$ for the standard basis,
and $a\cdot b=a\T b$ for the standard bilinear pairing. A
\emph{translation direction} of a Zariski closed set $Y\subseteq K^m$
is a constant vector $d$ with $Y+\tau d=Y$ for all $\tau\in K$.

\begin{proposition}\label{prop:qt-hess}
Let $h\in\kk[Z_1,\dots,Z_n]$ with $\rk\Hess h=n-1$, let
$R\in\kk[Y_1,\dots,Y_n]$ be an irreducible polynomial vanishing on the
image of $\nabla h$, and put $Q=(\nabla R)(\nabla h)$. Then $Q\ne0$,
$\Hess h\cdot Q=0$, and for every $\tau\in\kk$ the map $Z+\tau Q$ is
a quasi-translation, Pascal finite of Pascal depth at most two, with
\[
\nabla h\big(Z+\tau Q(Z)\big)=\nabla h(Z),\qquad
Q\big(Z+\tau Q(Z)\big)=Q(Z).
\]
A vector $d\in\kk^n$ is a translation direction of the closure of the
image of $\nabla h$ if and only if $d\T Q=0$.
\end{proposition}

\begin{proof}
The closure of the image of $\nabla h$ is an irreducible hypersurface,
since in characteristic zero its dimension is the generic rank of
$\Hess h$. Its ideal is generated by $R$, so $R$ has least degree
among the nonzero relations of $\nabla h$. Lemma~\ref{lem:graph}(a)
gives $Q\ne0$, $\Hess h\cdot Q=0$ and $D_Q(Q_i)=0$, and the remark
after it gives the two identities. In particular $JQ\cdot Q=0$, and
Proposition~\ref{prop:qt-pf} applies to $\tau Q$. Finally
$d\T Q=(D_dR)(\nabla h)$. If $d$ is a translation direction, then
$R(Y+ad)$ generates the ideal of the hypersurface for every $a\in\kk$,
so $R(Y+ad)=R(Y)$, as the homogeneous parts of highest degree agree,
and $D_dR=0$. Conversely, if $d\T Q=0$, then $D_dR$ vanishes on the
hypersurface, is divisible by $R$, and is zero by degree, so
$R(Y+ad)=R(Y)$.
\end{proof}

\begin{example}\label{ex:apex}
Let $s=Z_1+Z_4Z_2+Z_4^2Z_3$ and $h=\frac12s^2\in\kk[Z_1,\dots,Z_4]$, the polynomial $g(Z_4,s)$ of Remark~\ref{rem:cyl} with $g=\frac12s^2$. With $t=Z_4$ and $w=Z_2+2tZ_3$ one has
\[
\grad h=(s,\ ts,\ t^2s,\ ws).
\]
The map $Z\mapsto(s,t,w)$ is dominant, so the closure of the image of $\grad h$ is the hypersurface $\{R=0\}$ with $R=Y_2^2-Y_1Y_3$, and $\rk\Hess h=3$. Proposition~\ref{prop:qt-hess} gives
\[
Q=(\grad R)(\grad h)=-s\,(t^2,\ -2t,\ 1,\ 0).
\]
The functions $s$, $t$ and $w$ are annihilated by $D_Q$, hence constant along the lines $Z+\tau Q(Z)$. So $\grad h(Z+\tau Q)=\grad h(Z)$, $Q(Z+\tau Q)=Q(Z)$, and $Z+Q$ is a quasi-translation of Pascal depth at most two. The translation directions of $\{R=0\}$ are the vectors $d$ with $D_dR=2d_2Y_2-d_1Y_3-d_3Y_1=0$, that is $d\in\kk e_4$, in accordance with $e_4\T Q=0$. The kernel of $\Hess h$ is spanned by $(t^2,-2t,1,0)$, so $\Hess h$ has no nonzero constant kernel vector. The kernel direction in the source varies with $t$, while the closure of the gradient image has the constant translation direction $e_4$ in the target, as Theorem~\ref{thm:apex} asserts for every polynomial in four variables with Hessian of rank three.
\end{example}

\begin{theorem}[congruences of lines]
\label{thm:cong}
Let $V\colon\kk^3\dashrightarrow\kk^3$ be a rational map with
nilpotent Jacobian matrix $JV$. Suppose $H=\Hess f$ for a polynomial
$f\in\kk[X_1,X_2,X_3]$ with $\det H\ne0$ as a polynomial, and that
$H\,JV$ is symmetric. Then the image of $V$ is contained in an affine
plane.
\end{theorem}

The condition that $H\,JV$ is symmetric says that $N(X)=JV(X)$ is self-adjoint for the nondegenerate symmetric bilinear form $(u,v)\mapsto u\T H(X)v$. Pointwise linear algebra then gives isotropic subspaces, $\operatorname{im}N$ if $\rk N=1$ and $\ker N=\operatorname{im}N^2$ if $\rk N=2$ (Lemma~\ref{lem:cong-ker}). The content of the theorem is that the integrability $N=JV$, together with the fact that $H$ is a Hessian matrix, turns this moving structure into a constant affine relation among the components of $V$. The proof of Theorem~\ref{thm:cong} is divided into five lemmas. Lemmas~\ref{lem:cong-rk1} and~\ref{lem:cong-lines} follow the proof of \cite[Theorem~2.1]{Tan26}. Lemma~\ref{lem:cong-dir} replaces the use of \cite[Theorem~0.1]{DP04} in that proof. In these lemmas $V\colon\kk^3\dashrightarrow\kk^3$ is a rational map with nilpotent Jacobian matrix $N=JV$, $H=\Hess f$ for a polynomial $f\in\kk[X_1,X_2,X_3]$ with $\det H\ne0$, and $HN$ is symmetric.

\begin{lemma}[rank one]\label{lem:cong-rk1}
If $\rk N=1$, then the image of $V$ is contained in an affine plane.
\end{lemma}

\begin{proof}
The image closure is a point or an irreducible curve.
A point is contained in a plane, so let it be a curve. Its function field, regarded as the subfield $\kk(V)$ of $\kk(X)$, has transcendence degree one, and by Lemma~\ref{lem:luroth} it is generated by one rational function $t=p/q$ of the source
coordinates, with $p,q$ coprime polynomials, $q\ne0$. Write $V=A(t)$
with $A$ a vector of rational functions of one indeterminate $T$.
Assume, for a contradiction, that the image lies in no affine plane.
Equivalently, the components of $A'(T)$ are linearly independent over
$\kk$ (a constant relation $\lambda\cdot A'=0$ integrates to
$\lambda\cdot A=\mathrm{const}$).

Nilpotence gives $\tr JV=A'(t)\cdot\nabla t=0$. Since
$q^2\nabla t=q\nabla p-p\nabla q$, we get
$A'(t)\cdot\nabla(p-Tq)\big|_{T=t}=0$, where the gradient is taken
with $T$ fixed. Clearing the denominators of $A'(T)$, one obtains a
polynomial vector $P(T)$ with $P(T)\cdot\nabla(p-Tq)=0$ at
$T=p/q$, so $p-Tq$, which is primitive as a polynomial in $T$ over
$\kk[X]$ because $p,q$ are coprime, divides
$P(T)\cdot\nabla(p-Tq)$ in $\kk[X][T]$ by Gauss's lemma. Every entry
of $\nabla(p-Tq)$ has source degree smaller than the source degree of
$p-Tq$, and multiplication by polynomials in $T$ alone does not raise
source degrees. A nonzero multiple of $p-Tq$ would have source degree
at least that of $p-Tq$. Hence
\begin{equation}\label{eq:traceT}
A'(T)\cdot\nabla(p-Tq)=0\qquad\text{in }\kk(T)[X_1,X_2,X_3].
\end{equation}

Write $\nabla(p-Tq)=\sum_m m(X)\,c_m(T)$, the sum over the monomials
$m$ of $X$, with coefficient vectors $c_m(T)\in\kk^3$ affine in $T$,
and let $M(T)$ be the matrix with rows $c_m(T)\T$. By
\eqref{eq:traceT} the row space of $M$ is orthogonal to $A'(T)\ne0$,
so $\rk M\le2$. Rank zero would make $p,q$, hence $t$, constant.

Suppose $\rk M=2$. The cross product $n(T)$ of two independent rows is
a nonzero polynomial vector of degree at most two in $T$ with
$A'(T)\parallel n(T)$. Writing $n(T)=n_0+n_1T+n_2T^2$, the vectors
$n_0,n_1,n_2$ span $\kk^3$, since a common annihilator $\lambda$ would give
$\lambda\cdot A'(T)=0$, a constant relation. Conjugating $V$ by the
constant linear map taking $n_0,n_1,n_2$ to the standard basis (which
preserves the hypotheses, since $JV$ transforms by conjugation and $H$ by
congruence) makes $n(T)=(1,T,T^2)\T$. The orthogonality
$c_m(T)\cdot n(T)=0$, for every monomial $m$, spells out to the four
identities
\[
p_x=0,\qquad p_y=q_x,\qquad p_z=q_y,\qquad q_z=0
\]
in coordinates $x,y,z$. The first and last give $p=p(y,z)$,
$q=q(x,y)$. The second gives $\partial_yp=\partial_xq=h(y)$ for some
one-variable polynomial $h$, so $p=P(y)+B(z)$ and $q=xP'(y)+C(y)$ with
$P'=h$. The third gives $B'(z)=xP''(y)+C'(y)$, so $P''=0$ and $B',C'$
are constants. Then $p,q$ are affine, $\nabla(p-Tq)$ is constant in
$X$, and $M$ has a single row, contradicting rank two.

Hence $\rk M=1$, and all rows of $M$ are proportional over $\kk(T)$ to a
single covector, defining a nonzero linear form $L(T)$ in the source
variables over $\kk(T)$ whose kernel annihilates all derivatives of
$p-Tq$. Thus $p-Tq\in\kk(T)[L(T)]$. If the source degree $D$ of
$p-Tq$ is at least two, its top homogeneous part is $c(T)L(T)^D$ with
$c\ne0$. This part is affine in $T$, so its second $T$-derivative
vanishes. Dividing that derivative by $L^{D-2}$ and reducing modulo
the prime $L$ gives $D(D-1)c(T)(L')^2\equiv0\pmod L$, so $L'$ is
proportional to $L$ and the ratios of the coefficients of $L$ are
constant. Therefore $L$ is, up to a factor in $\kk(T)$, a constant linear form
$\ell$. Evaluating $p-Tq$ at two constants $T$ shows
$p,q\in\kk[\ell]$. The same holds for $D=1$ when the linear parts of
$p,q$ are proportional. In either case \eqref{eq:traceT} becomes a
fixed nonzero linear relation on $A'(T)$, a contradiction. So $p,q$
are affine with independent linear parts. After the affine source
change $x=p$, $y=q$ (again preserving the hypotheses) we have
$t=x/y$ and $\nabla t=(1,-t,0)\T/y$.

The identity $\tr JV=0$ now reads $A_1'=tA_2'$. Independence gives
$A_2'\ne0$, and we put $k=A_3'/A_2'\in\kk(t)$, so that
$A'=A_2'\,(t,1,k)\T$. Symmetry of
$H\,JV=H\,A'(t)\,\nabla t\T$ forces $HA'(t)\parallel\nabla t$, that
is, $H\,(t,1,k)\T\parallel(1,-t,0)\T$. Contracting with the two
vectors $(t,1,0)$ and $(0,0,1)$, orthogonal to $(1,-t,0)$, and using
the chain rule for $F(t,y,z):=f(ty,y,z)$ at fixed $t$ (so that
$F_y+kF_z=(t,1,k)\cdot\nabla f(ty,y,z)$) gives
\[
\partial_y\big(F_y+kF_z\big)=\partial_z\big(F_y+kF_z\big)=0 ,
\]
so $G:=F_y+kF_z$ depends only on $t$. The coefficient of $y^jz^\ell$
in $F$ has $t$-degree at most $j$, because $x^ay^bz^\ell$ becomes
$t^ay^{a+b}z^\ell$. Write $F=\sum_{\ell=0}^NF_\ell(t,y)z^\ell$ with
$N=\deg_zF$. If $N\ge1$, the coefficient of $z^N$ in $G$ is
$\partial_yF_N=0$, and the degree restriction at $j=0$ makes
$F_N=c\in\kk^\times$. If moreover $N\ge2$, the coefficient of
$z^{N-1}$ is $\partial_yF_{N-1}+Nck=0$. Thus
$\partial_yF_{N-1}=-Nck(t)$ is independent of $y$, it is the
coefficient of $y$ in $F_{N-1}$ up to higher $y$-terms, and the degree
restriction at $j=1$ makes $k=at+b$ affine. Then
$A_3'=kA_2'=aA_1'+bA_2'$, a constant relation, contradiction. If
$N\le1$, then $F=F_0(t,y)+cz$ and $\partial_yF_0=G-kc$ depends only on
$t$. The degree restrictions make $F=(at+b)y+cz+d$, so $f$ is affine
and $\det H=0$, a contradiction.
\end{proof}

\begin{lemma}[involutivity and lines]\label{lem:cong-lines}
Let $\rk N=2$ and suppose that the image of $V$ is contained in no affine plane. Then the distribution $\operatorname{im}N$ is closed under the Lie bracket, and the closure of every component of a general fiber of $V$ is an affine line with direction $\ker N$.
\end{lemma}

\begin{proof}
Choose a point where $N$ has rank two
and $H$ is invertible, and normalize at that point. A constant
conjugation makes $N$ equal to the regular nilpotent Jordan block $J$
($J_{12}=J_{23}=1$). The symmetry $HJ=J\T H$ forces, at the point,
\[
H=\begin{pmatrix}0&0&a\\0&a&b\\a&b&c\end{pmatrix}=Ra(I+U),\qquad
a\ne0,\quad U=\tfrac1a\big(bJ+cJ^2\big),\quad U^3=0,
\]
where $R$ is the antidiagonal identity and $RJ=J\T R$, hence
$Rg(U)=g(U)\T R$ for polynomials $g$. The congruence by
$a^{-1/2}\big(I-\tfrac12U+\tfrac38U^2\big)$, a polynomial in $J$ and
hence commuting with $J$, replaces $H$ by $R$ at the point and leaves
$N=J$ there. Indeed, for this congruence $M$ one has $M\T R=RM$, and
$M(I+U)M=a^{-1}I$ since $U^3=0$, so $M\T HM=aRM(I+U)M=R$.

Let $K_{i,j\ell}=\partial_j\partial_\ell V_i$ and
$T_{ij\ell}=\partial_i\partial_j\partial_\ell f$ at the point.
Differentiating $\tr N^3=\tr N^2=\tr N=0$ and dividing by $3,2,1$, one obtains
\[
K_{3,1\ell}=0,\qquad K_{2,1\ell}+K_{3,2\ell}=0,\qquad
K_{1,1\ell}+K_{2,2\ell}+K_{3,3\ell}=0 .
\]
The $(1,2)$ and $(1,3)$ entries of the symmetry of
$\partial_\ell(HN)=(\partial_\ell H)J+R\,\partial_\ell N$ give
\[
T_{11\ell}+K_{3,2\ell}=K_{2,1\ell},\qquad
T_{12\ell}+K_{3,3\ell}=K_{1,1\ell}.
\]
Combining these identities gives $T_{11\ell}=-2K_{3,2\ell}$ and $K_{3,31}=K_{3,13}=0$. From
$T_{112}=T_{121}$, $K_{1,11}=-2K_{3,22}$. From the second trace
identity, $K_{2,21}=-K_{3,22}$. Substituting into the third at
$\ell=1$ gives $K_{3,22}=0$. Let $\omega$ be a nonzero rational one-form annihilating $\operatorname{im}N$, equal to $dX_3$ at the point.
Differentiating $\omega N=0$ gives
$\partial_\ell\omega_1=-K_{3,2\ell}$ and
$\partial_\ell\omega_2=-K_{3,3\ell}$, and the coefficient of
$dX_1\wedge dX_2\wedge dX_3$ in $\omega\wedge d\omega$ at the point is
$-K_{3,31}+K_{3,22}=0$. The point is general, so $\omega\wedge d\omega=0$. For rational
vector fields $L,L'$ with values in $\operatorname{im}N$ and a rational
vector field $T$ with $\omega(T)=1$ this gives
$\omega([L,L'])=-d\omega(L,L')=-(\omega\wedge d\omega)(T,L,L')=0$, so
the rank-two distribution $\operatorname{im}N$ is closed under the Lie
bracket.

Locally $V=(u,v,w(u,v))$ with
$u,v$ algebraically independent components and $w$ algebraic over
$\kk(u,v)$. Derivations in $u,v$ extend uniquely to the image function
field. Let $p=w_u$, $q=w_v$, let $M$ be the $3\times2$ matrix with
columns $(1,0,p)\T,(0,1,q)\T$ and $J_0=D(u,v)$, so $N=MJ_0$. From
$N^2=M(J_0M)J_0$ of rank one and $N^3=0$ we get
$\rk(J_0M)=1$ and $(J_0M)^2=0$. A $2\times2$ matrix of rank one with square zero has the form $b\,\mu\nu\T$ with $\nu\T\mu=0$, so $\nu$ is proportional to $(-\mu_2,\mu_1)\T$. After exchanging $u$ and $v$ if necessary, $\mu=(1,\beta)\T$, and

\[
J_0M=b\begin{pmatrix}-\beta&1\\-\beta^2&\beta\end{pmatrix},
\qquad b\ne0 .
\]
The columns of $M$ define the vector fields $L_1=\partial_1+p\partial_3$,
$L_2=\partial_2+q\partial_3$ with
$[L_1,L_2]=(L_1(q)-L_2(p))\partial_3$, and the chain rule with the
entries of $J_0M$ gives
$[L_1,L_2]=-b\,(w_{uu}+2\beta w_{uv}+\beta^2w_{vv})\partial_3$. Since
$\partial_3\notin\operatorname{im}N$, involutivity makes the
coefficient vanish. If all second derivatives of $w$ vanished, then
$p,q$, and after subtracting an affine function also $w$, would be
constant (an element of the image field with vanishing derivatives is
constant, by differentiating its minimal polynomial), contradicting
nonplanarity. Hence $\beta$ is algebraic over the image function
field, and on each component of the geometric generic fiber of $V$ the
functions $p,q,\beta$ are constant. The tangent line of such a
component is $\ker N=\operatorname{im}N^2$, generated by
$(1,\beta,p+q\beta)\T$, a constant vector on the component. Two
independent linear forms vanishing on it are constant on the
component, whose closure is therefore an affine line.
\end{proof}

\begin{lemma}\label{lem:cong-ker}
Let $\rk N=2$ and let $G\in\kk(X)^3\setminus\{0\}$ span $\ker N$. Then $G\T HG=0$ and $N\T HG=0$.
\end{lemma}

\begin{proof}
Since $\ker N=\operatorname{im}N^2\subseteq\operatorname{im}N$ and
$N$ is self-adjoint for the pairing $a\T Hb$, one has
$G\T HNb=(NG)\T Hb=0$ for every $b$, hence
\[
G\T HG=0,\qquad N\T HG=HNG=0 .
\]
\end{proof}

\begin{lemma}[pencils of planes]\label{lem:pencil}
Let $p,q\in\kk[X_1,X_2,X_3]$ be coprime with $p/q\notin\kk$, and let $n\in\kk(T)^3\setminus\{0\}$. Suppose that for general $c\in\kk$ every irreducible component of $\{p=cq\}$ is an affine plane with normal vector $n(c)$. Then $n\in\kk(T)\,n_1$ for a constant vector $n_1\in\kk^3$, or there is $e\in\kk^3\setminus\{0\}$ with $n\T e=0$ in $\kk(T)$.
\end{lemma}

\begin{proof}
Suppose that $n$ is not a multiple of a constant vector. Fix a general $c_0$ and a component $Z_0$ of $\{p=c_0q\}$. The $2\times2$ minors of $(n(c_0),n(T))$ are rational functions of $T$, not all zero, so for all but finitely many $c$ the vectors $n(c)$ and $n(c_0)$ are not proportional. For such a general $c$ a component $Z$ of $\{p=cq\}$ is an affine plane not parallel to $Z_0$, and $Z\cap Z_0$ is a line. On it $p=cq$ and $p=c_0q$, so $p=q=0$. Since $p$ and $q$ are coprime, $\{p=q=0\}$ has dimension at most one and contains finitely many lines. So one of them, of direction $e\neq0$, lies in $Z$ for infinitely many $c$. Then $n(c)\T e=0$ for infinitely many $c$, and $n\T e=0$ in $\kk(T)$.
\end{proof}

\begin{lemma}[constant direction]\label{lem:cong-dir}
Let $\rk N=2$, suppose that the image of $V$ is contained in no affine plane, let $G\in\kk(X)^3\setminus\{0\}$ span $\ker N$, and put $E=\sum_iG_i\partial_i$. Then $G\in\kk(X)\,g$ for a constant vector $g\in\kk^3$.
\end{lemma}

\begin{proof}
By Lemma~\ref{lem:cong-lines} the direction of $G$ is constant along the components of the general fibers of $V$, and by Lemma~\ref{lem:cong-ker} one has $G\T HG=0$ and $N\T HG=0$. Let $\mathcal R\in\kk[Y_1,Y_2,Y_3]$ be an irreducible polynomial
vanishing on the image of $V$, and put $n=(\nabla\mathcal R)(V)$. Then
$n\ne0$ and $n\T N=0$. The left kernel of $N$ is a line and $HG\ne0$,
so $HG=\mu n$ with $\mu\in\kk(X)^\times$, and $n\T G=\mu^{-1}G\T HG=0$.

Choose a linear form $\pi$ with $\pi(G)\ne0$, and put $g=G/\pi(G)$ and $X_0=X-\pi(X)g$. Hence $E(g)=0$, and then $E(X_0)=G-\pi(G)g=0$. Moreover
$E(V)=NG=0$, so $E(n)=0$. Let $F_0=\{r\in\kk(X):E(r)=0\}$. Then
$X_0,g,n\in F_0^3$ and $X=X_0+\pi g$, so $\kk(X)=F_0(\pi)$. Since $E$ acts
on $F_0[\pi]$ as $\pi(G)\,d/d\pi$, the element $\pi$ is transcendental
over $F_0$, and every derivation $D$ of $F_0$ extends to $\kk(X)$ with
$D(\pi)=0$.

\emph{Step 1, $n\T D(g)=0$ for every derivation $D$ of $F_0$.} The
polynomial $f(X_0+\pi g)$ lies in $F_0[\pi]$, and its second derivative
in $\pi$ is $g\T Hg=0$. Hence $f=\alpha_0+\pi\alpha_1$ with $\alpha_0,\alpha_1\in F_0$. The derivative of $\nabla f(X_0+\pi g)\in
F_0[\pi]^3$ in $\pi$ is $Hg$, a nonzero multiple of $n\in F_0^3$, so
$\nabla f=\gamma+\theta n$ with $\gamma\in F_0^3$ and $\theta\in F_0[\pi]$ of positive
degree. Applying $D$ to $f=\alpha_0+\pi\alpha_1$ and using
$D(X)=D(X_0)+\pi D(g)$ gives
\[
\gamma\T D(X_0)+\pi\,\gamma\T D(g)+\theta\,\big(n\T D(X_0)+\pi\,n\T D(g)\big)
=D(\alpha_0)+\pi D(\alpha_1).
\]
The right side has degree at most one in $\pi$ and $\theta$ has positive
degree, so $n\T D(g)=0$.

\emph{Step 2, the field $L=\kk(g_1,g_2,g_3)$ has transcendence degree at
most one.} After a linear change of coordinates $\pi=X_3$, so
$g=(g_1,g_2,1)$ and $L=\kk(g_1,g_2)\subseteq F_0$. If $g_1,g_2$ were
algebraically independent, then, as $\operatorname{trdeg}F_0=2$, there would be
derivations $D_1,D_2$ of $F_0$ with $D_i(g)=e_i$, and Step 1 together
with $n\T g=0$ would give $n=0$.

\emph{Step 3, $\operatorname{trdeg}L=1$ is impossible.} By Lemma~\ref{lem:luroth}, $L=\kk(s)$ with $s\in L\subseteq F_0$. Write $g=g(s)$, let a
prime denote $d/ds$, and put $n_0=g\times g'$. Since $\pi(g)=1$, $\pi(g')=0$
and $g'\ne0$, the vectors $g,g'$ are linearly independent. Taking
in Step 1 a derivation with $D(s)=1$ gives $n\T g'=0$, so $n=\rho n_0$
with $\rho\in\kk(X)^\times$. If the direction of $n_0$ were constant,
the linear form $n_0\T Y$ would have zero differential on the image
of $V$, since $n_0\T N=0$, and the image would lie in a plane. Hence
$n_0$ and $n_0'$ are linearly independent over $\kk(s)$.

Put $h=n_0(s)\T V$ and $k=n_0'(s)\T V$. Since $n_0\T N=0$, every
derivation $D$ of $\kk(X)$ satisfies $D(h)=k\,D(s)$. Thus $h$ is
annihilated by the derivations vanishing on $s$, and is algebraic over
$\kk(s)$ because the characteristic is zero. Applying a derivation with
$D(s)=1$ to the minimal polynomial of $h$ over $\kk(s)$ shows that
$k=D(h)\in\kk(s,h)$, so $k$ is algebraic over $\kk(s)$ as well.

Write $s=p/q$ with coprime $p,q\in\kk[X]$. For a general $c\in\kk$ let
$Z$ be an irreducible component of $\{p=cq\}$. On $Z$ the functions
$h$ and $k$ are constant, so $V(Z)$ lies in the line
$\{n_0(c)\T Y=h,\ n_0'(c)\T Y=k\}$. Its direction is $g(c)$, since
$n_0\T g=n_0\T g'=0$ and hence $(n_0')\T g=0$. Let $T$ be
the tangent plane of $Z$ at a general point. It contains
$\ker N=\kk g(c)$, because $E(s)=0$, and $N(T)\subseteq\kk g(c)$. As
$\dim T=2>\dim\ker N$, one gets $N(T)=\kk g(c)$, so $N^2(T)=0$ and
$T\subseteq\ker N^2=\operatorname{im}N=n_0(c)^\perp$. Hence
$n_0(c)\T X$ is constant on $Z$, and $Z$ is an affine plane with
normal $n_0(c)$.

Since the direction of $n_0$ is not constant, Lemma~\ref{lem:pencil} gives $e\in\kk^3\setminus\{0\}$ with $n_0\T e=0$ in $\kk(s)$. Hence $(n_0')\T e=0$ and $e\in\kk(s)\,g$. Since $\pi(g)=1$, the vector $g$ is constant, a
contradiction. Hence $L=\kk$, and $g$ is constant.
\end{proof}

\begin{lemma}[parallel lines]\label{lem:cong-par}
Let $\rk N=2$ and suppose that $\ker N=\kk(X)\,g$ for a constant vector $g\in\kk^3\setminus\{0\}$. Then the image of $V$ is contained in an affine plane.
\end{lemma}

\begin{proof}
After a linear change of coordinates $g=e_3$. Then
$V=(A(x,y),B(x,y),C(x,y))$, and $g\T Hg=0$ by Lemma~\ref{lem:cong-ker}, that is $f_{zz}=0$,
gives $f=zP(x,y)+Q(x,y)$ with $P$ nonconstant. The planar Jacobian of
$(A,B)$ is nilpotent. If $(A,B)$ is constant, the image lies in a
plane. Otherwise the nilpotent Jacobian of $(A,B)$ has rank one, $\kk(A,B)$ has transcendence degree one, and Lemma~\ref{lem:luroth} gives $(A,B)=(a(t),b(t))$ for a nonconstant rational $t$. With $r=(a'(t),b'(t))\T$ the trace gives
$\nabla t\cdot r=0$, and $He_3=(P_x,P_y,0)\T$ annihilates
$\operatorname{im}N$ by Lemma~\ref{lem:cong-ker}, so also $\nabla P\cdot r=0$. Hence
$dt\wedge dP=0$ and $t$ is algebraic over $\kk(P)$. On each component
of a general level of $P$, $t$ is constant and the tangent direction
is the constant vector $r$, so the component is dense in a line lying
in the level ($P$ is a polynomial). By Lemma~\ref{lem:linfib} with $n=2$ and $k=1$ these lines are parallel, and the common direction is independent of
the component. A constant linear form annihilating it annihilates
$(a',b')$, giving a constant affine relation between $A$ and $B$, and the
image lies in a plane.
\end{proof}

\begin{proof}[Proof of Theorem~\ref{thm:cong}]
If $\rk JV=0$, then $V$ is constant. If $\rk JV=1$, Lemma~\ref{lem:cong-rk1} applies. Let $\rk JV=2$ and suppose that the image of $V$ is contained in no affine plane. By Lemma~\ref{lem:cong-dir} the kernel of $JV$ is spanned by a constant vector, and Lemma~\ref{lem:cong-par} gives a contradiction.
\end{proof}

\begin{lemma}[slices]\label{lem:slice}
Let $h\in\kk[X_1,X_2,X_3,w]$ have Hessian rank three, and suppose
$A(X)=\Hess_Xh(X,0)$ is nonsingular over $\kk(X_1,X_2,X_3)$, where
$\Hess_X$ differentiates in $X_1,X_2,X_3$ alone. Set
$v(X)=A(X)^{-1}\nabla_Xh_w(X,0)$. Then
\[
\nabla_Xh\big(X-wv(X),\,w\big)=\nabla_Xh(X,0),
\]
$Dv$ is nilpotent, and $A\,Dv$ is symmetric.
\end{lemma}

\begin{proof}
Let $R$ and $Q$ be as in Proposition~\ref{prop:qt-hess} for $h$, with
$Z=(X,w)$. Then $\nabla h(Z+\tau Q(Z))=\nabla h(Z)$ for all $\tau$.

The slice $w=0$ maps dominantly to the hypersurface, because the first
three gradient coordinates have Jacobian $A$ there. Hence
$Q(X,0)\ne0$ generically. On the slice the first three Hessian rows
are $(A\mid\nabla_Xh_w(X,0))$, whose kernel is spanned by $(-v,1)$.
The fourth row leaves this kernel unchanged because the full Hessian
has rank three. Therefore $Q(X,0)=Q_4(X,0)\,(-v(X),1)$, and substituting
$\tau=w/Q_4(X,0)$ in the quasi-translation identity gives the claimed
formula. Differentiation in $X$ gives
\[
\Hess_Xh\big(X-wv(X),w\big)\,(I-w\,Dv)=A(X).
\]
Both factors on the left have entries in $\kk(X)[w]$ and the product
of their determinants is the $w$-independent $\det A$. Each
determinant is therefore constant in $w$, and at $w=0$ the second one
is $1$. So $\det(I-wDv)=1$, the characteristic polynomial of $Dv$ is
$T^3$, and $(Dv)^3=0$ by Cayley--Hamilton. Finally, differentiating
$Av=\nabla_Xh_w(X,0)$ shows that $A\,Dv$ is the difference of
$\Hess_Xh_w(X,0)$ and the matrix with entries
$\sum_kv_k\partial_jA_{ik}$, both symmetric since
$\partial_jA_{ik}=\partial_iA_{jk}$.
\end{proof}

\begin{theorem}[translation directions]
\label{thm:apex}
If $h\in\kk[X_1,\dots,X_4]$ has Hessian rank three, then the closure
of its gradient image has a nonzero translation direction.
Equivalently, the components of the quasi-translation $Q$ of
Proposition~\ref{prop:qt-hess} are linearly dependent over $\kk$.
\end{theorem}

\begin{proof}
At a point where the (symmetric) Hessian has rank three, it restricts
nondegenerately to some coordinate three-plane. Indeed, its adjugate there is $c\,\eta\eta\T$ with $c\ne0$ and $\eta$ spanning the kernel, and for $\eta_i\ne0$ the diagonal entry $c\eta_i^2$ is the principal minor complementary to $e_i$. A constant linear
change and a translation of the remaining coordinate put us in the
situation of Lemma \ref{lem:slice} with slice $w=0$. The map $v$ of
that lemma satisfies the hypotheses of Theorem \ref{thm:cong} with
$H=\Hess(h(\,\cdot\,,0))=A$, since $Dv$ is nilpotent, $A\,Dv$ is symmetric,
and $\det A\ne0$. Hence there are $b\in\kk^3\setminus\{0\}$ and
$d\in\kk$ with $b\cdot v\equiv d$. With $R,Q$ as in Lemma \ref{lem:slice}, on the slice
$Q(X,0)=Q_4(X,0)(-v,1)$ and $(b,d)\cdot(-v,1)=d-b\cdot v=0$. Since
$(b,d)\T Q=(D_{(b,d)}R)(\nabla h)$ and the slice maps dominantly to
the gradient image, $(b,d)\T Q=0$, and $(b,d)\ne0$ is a translation
direction by Proposition~\ref{prop:qt-hess}. An affine source change moves the gradient image by an
invertible affine target map, so the conclusion also holds in the
original coordinates.
\end{proof}

\begin{corollary}[{\cite[Corollary~2.4]{Tan26}}]\label{cor:dBproblem}
Let $h\in\kk[X_1,\dots,X_4]$ have identically vanishing Hessian
determinant. Then there are a nonzero $R\in\kk[Y_1,\dots,Y_4]$ with
$R(\nabla h)=0$ and a nonzero $d\in\kk^4$ with $D_dR=0$. Thus
Problem~1 of de Bondt \cite[Section 6]{dB18} has an affirmative
answer.
\end{corollary}

\begin{proof}
For rank three take an irreducible equation $R$ of the gradient-image
closure and the direction $d$ of Theorem \ref{thm:apex}, and apply
Proposition~\ref{prop:qt-hess}. For rank at most two the first three gradient
coordinates are algebraically dependent, and any nonzero relation $R$
in $Y_1,Y_2,Y_3$ works with $d=e_4$.
\end{proof}

\begin{remark}\label{rem:qt-lindep}
Theorem~\ref{thm:apex} is a statement on linear dependence of
quasi-translations. In dimension at most three every quasi-translation
$Z+Q$ has linearly dependent components \cite[Corollary~5.6(i)]{dB18}.
Theorem~\ref{thm:apex} gives the same for the quasi-translations in
dimension four which arise from Hessians of rank three by
Proposition~\ref{prop:qt-hess}. By Proposition~\ref{prop:qt-pf} these quasi-translations are Pascal finite of Pascal depth at most two.
\end{remark}

\section{Weighted forms, walls, and exit lemmas}\label{sec:walls}

\subsection*{Standing hypotheses}
Throughout Sections \ref{sec:walls} and \ref{sec:tree},
$\phi\in K[x_1,\dots,x_4]$ is a Keller potential,
$\det\Hess\phi=\kappa\in K^\times$, the affine part of $\phi$ has been removed
($\phi(0)=0$, $\nabla\phi(0)=0$, which changes neither
$\det\Hess\phi$ nor the existence of affine directions), and
$D:=\deg\phi\ge3$. ``A direction'' always means an affine direction,
$v\ne0$ with $D_v^2\phi=0$. Hessian ranks are taken over
$K(x_1,\dots,x_4)$.

Table~\ref{tab:exit} lists the exit lemmas of this section and the tools
which call them. For each lemma it gives the leading form at the start
of the path, the apex configurations which are possible at a wall, the
minors or cofactors which survive to the wall form, the quantity which
bounds the number of walls, and the way in which the lemma ends. The
notions used in the table are defined below. In the column on minors,
$H$ denotes the leading form of the lemma, $A'$ and $B'$ are the derivatives
of $aT^r$ and $bT^s$ at $T=S$, and $F$ is a fiber variable.

\begin{table}[tp]
\footnotesize
\setlength{\tabcolsep}{3pt}
\begin{tabular}{@{}>{\raggedright\arraybackslash}p{0.075\textwidth}*{5}{>{\raggedright\arraybackslash}p{0.161\textwidth}}@{}}
Lemma & Leading form at the start & Apex at a wall & Surviving minor or cofactor & Bound on the walls & End\\
\hline
\ref{lem:D1} & $R+zA+wB$, $A=aS^r$, $B=bS^s$, $S$ genuinely binary, $w_z\ne w_w$ & a plane, graph over the base block & cofactors $-(B')^2(zA'+wB')\BH_2(S)$ and $-(A')^2(zA'+wB')\BH_2(S)$, minors $-(D_vA)^2$, $-(D_vB)^2$ & at most four base and fiber weight equalities, Lemma~\ref{lem:D1-walls} & affine direction $e_z$\\
\ref{lem:D2} & $a(x)(Y_1+bx^rY_2+cx^sY_3)^m$ & the active axis & passive cofactors \eqref{eq:phase-cof} & phase degree $M$, strictly increasing, $M\le D_Y$ & contradiction\\
\ref{lem:D3} & $ax^\ell(w+x^Ay)^m+bx^pz$ & the active axis & diagonal entries of $-F_{LL}\,nn\T$ & phase degree, at most two shears & affine direction $e_z$, or Lemma~\ref{lem:D2}\\
\ref{lem:P1} & $G(x,y)+c_1x^pz+c_2x^qw$ & no wall, one small step & none & none & Lemma~\ref{lem:D4}\\
\ref{lem:D4} & $ax^\ell y^m+bx^pz+cx^qw$ & the active axis, or a plane & cofactors $-H_{yy}(cqx^{q-1})^2$, $-H_{yy}(bpx^{p-1})^2$, minors on $(x,F)$ & phase degree, at most two shears, $W$ decreasing & Lemmas \ref{lem:D1}, \ref{lem:D2}, \ref{lem:D3}, or $W=0$\\
\ref{lem:D5} & rank three, one apex in a block of smallest weight & given & as in the lemma called & one small step & Lemmas \ref{lem:D2}, \ref{lem:D3}, \ref{lem:D4}\\
\ref{lem:D6} & $\sum_ib_ix^{p_i}Y_i$ & the active axis, or a plane containing it & minors $-(b_ip_ix^{p_i-1})^2$ on $(x,Y_i)$ & $W$ decreasing & Lemmas \ref{lem:D1}, \ref{lem:D4}, \ref{lem:D5}, or $W=0$\\
\ref{lem:D7} & $ay^M+bx^pz+cx^qw$ & a plane, a line is excluded & fiber cofactors as in Lemma~\ref{lem:D4}, minor $-(qcx^{q-1})^2$ & three weight equalities, $W$ decreasing & Lemma~\ref{lem:D1}, or $W=0$\\
\ref{lem:T} & $H(u_1,u_2,u_3)$, $v$ unused & one or two apices in the ternary block & cofactor deleting $v$, restricting to $\det\Hess_3H$ & $W$ decreasing & the caller\\
\ref{lem:B} & $G(x,y)+zA+wB$, rank three & a plane & constant kernel test & none & Lemma~\ref{lem:D1}, or Lemmas \ref{lem:P1} and \ref{lem:D4}\\
\ref{lem:K} & $as^D+bts^p+czs^q$, $u$ unused & $e_s$, or a plane containing $e_s$ & $M_{st}$, $M_{sz}$ & $W$ decreasing & Lemmas \ref{lem:D5}, \ref{lem:D6}, \ref{lem:B}\\
\hline
\end{tabular}
\caption{Exit lemmas and the tools which call them. ``$W=0$'' refers to
Corollary~\ref{cor:endpoint}.}\label{tab:exit}
\end{table}

\subsection{Weights, walls, and the touching criterion}

For a weight vector $w\in\Q^4$ let
$E(w)=\max\{\vv wm:m\in\supp\phi\}$ be the support function of the
Newton polytope of $\phi$, let $\phi_w$ be the \emph{leading form}
(the sum of the terms of weight $E(w)$), and put
\[
W(w)\;=\;4E(w)-2\sum_{i=1}^4w_i .
\]
Algebraically, $w$ defines a filtration of $\kk[x_1,\dots,x_4]$ by the weights of monomials, and $\phi_w$ is the initial form of $\phi$. Lemma~\ref{lem:W} below compares the initial form of $\det\Hess\phi$ with $\det\Hess\phi_w$, so the weighted argument is a sequence of initial degenerations of the equation $\det\Hess\phi=\kappa$. A wall is a weight vector at which the initial form changes, and the exit lemmas determine which initial forms can occur.

\begin{lemma}\label{lem:W}
$W(w)\ge0$ for every $w\in\R^4$. If $W(w)>0$, then $\det\Hess\phi_w=0$. If $W(w)=0$, then
$\det\Hess\phi_w=\kappa$.
\end{lemma}

\begin{proof}
Grade monomials by $m\mapsto\vv wm$. Every monomial of
$\partial_i\partial_j\phi$ has weight at most $E(w)-w_i-w_j$, so every
product in the Leibniz expansion of $\det\Hess\phi$ has weight at most
$W(w)$, with the component of weight exactly $W(w)$ coming from
$\phi_w$ alone and equal to the corresponding component of
$\det\Hess\phi_w$. Since $\kappa\ne0$ is a monomial of weight $0$, we
get $W(w)\ge0$. Every product in the Leibniz expansion of $\det\Hess\phi_w$ has weight exactly $W(w)$, so $\det\Hess\phi_w$ is the component of weight $W(w)$ of $\det\Hess\phi=\kappa$. It is therefore zero when $W(w)>0$ and equal to
$\kappa$ when $W(w)=0$.
\end{proof}

\begin{remark}\label{rem:centre}
Let $P$ be the Newton polytope of $\phi$ and $b=\frac12(1,1,1,1)$. The function $E$ is the support function of $P$, and $W(w)=4\bigl(E(w)-\vv{w,b}\bigr)$. Hence the inequality $W\ge0$ of Lemma~\ref{lem:W} says that $b\in P$, and $W(w)=0$ holds if and only if $b$ lies on the face of $P$ on which $w$ attains its maximum. By Lemma~\ref{lem:W}, the faces which do not contain $b$ have leading forms with singular Hessian, and the faces which contain $b$ have leading forms with Hessian determinant $\kappa$. By Theorem~\ref{thm:touch} below, the weighted analysis ends at once if, in some linear coordinates, $b$ lies on a face of $P$ with a nonnegative normal vector. This hypothesis is not invariant under linear changes of coordinates.
\end{remark}

\begin{theorem}[touching criterion]\label{thm:touch}
Suppose $E(w)=\tfrac12(w_1+\dots+w_4)$ for some
$w\in\R_{\ge0}^4\setminus\{0\}$. Then, after a permutation of the
variables, one of the following holds. \textup{(i)} $\deg\phi\le2$.
\textup{(ii)} $\partial_1^2\phi=0$. \textup{(iii)}
$\phi=R(x_3,x_4)+x_1A(x_3,x_4)+x_2B(x_3,x_4)$. In each case $\phi$ has
an affine direction.
\end{theorem}

\begin{proof}
Put $B:=\tfrac12\sum_iw_i$ and permute so that
$w_1\ge w_2\ge w_3\ge w_4\ge0$. Every $m\in\supp\phi$ satisfies
$\vv wm\le B$, and since $w,m\ge0$ also $m_iw_i\le B$. If $w_1>B$
then $\phi$ does not involve $x_1$ and $\det\Hess\phi=0$. Hence
$w_1\le B$. Let $S=\{i:w_i>0\}$. $|S|=1$ would give $B=w_1/2<w_1$, so
$|S|\ge2$. If $|S|=2$, then $w_1\le B=\tfrac12(w_1+w_2)$ forces
$w_1=w_2=B$ and every monomial has $m_1+m_2\le1$, which is case (iii). If
$|S|=3$, then $w_1\ge\tfrac23B>\tfrac12B$, so $m_1\le B/w_1<2$, which is case (ii).
If $|S|=4$, then $w_1\ge\tfrac12B$, with equality forcing all weights equal
and $\deg\phi\le2$ (case (i)). Otherwise $m_1<2$ (case (ii)). In case
(ii), $e_1$ is an affine direction. In case (iii), $e_1$ (and $e_2$)
is. In case (i), $D_v^2\phi=v\T\Hess\phi\,v$ is a nondegenerate
quadratic form on $K^4$ and has a nonzero isotropic vector over the
algebraically closed $K$.
\end{proof}

\begin{corollary}[endpoints]\label{cor:endpoint}
If $W(w)=0$ for some $w$ with positive entries, then $\phi$ has an
affine direction. \qed
\end{corollary}

Along an affine path of weights each monomial weight is an affine
function of the path parameter. We choose paths on which the monomials
of the current leading form stay tied, and call the first parameter at
which another monomial joins them a \emph{wall}. The finite support of
$\phi$ gives a strict gap until that parameter. At a wall, a further
small step in a specified direction selects from the tied monomials
those of greatest slope. Every path in Sections \ref{sec:walls} and
\ref{sec:tree} is stopped at its first wall or at $W=0$, whichever
comes first. In the latter case Corollary \ref{cor:endpoint} produces
a direction and the branch terminates. For this reason walls and
endpoints never need to be located explicitly.

\begin{example}\label{ex:wall}
For illustration let $p=x^4+x^2y+y^3\in K[x,y]$ and $w(\tau)=(1,\tau)$, $\tau\ge1$. The three monomials have the weights $4$, $2+\tau$ and $3\tau$. For $1\le\tau<\frac43$ the leading form is $x^4$. At $\tau=\frac43$ the monomial $y^3$ joins it, the leading form is $x^4+y^3$, and $w(\frac43)$ is orthogonal to the edge of the Newton polygon joining $(4,0)$ and $(0,3)$. This parameter is the wall. For $\tau$ slightly larger than $\frac43$ the leading form is $y^3$, the monomial of greatest slope, and the monomial $x^2y$ is never leading along this path. In general, at a wall the weight vector lies in the normal cone of a face of the Newton polytope which contains the previous leading face, and the leading form acquires the monomials of that face.
\end{example}

\subsection{Apices of leading forms}

An \emph{apex direction} of a polynomial $P$ with $\rk\Hess P\ge2$ is
a vector $d\in K^4$ such that the Zariski closure $Z_P$ of the image
of $\nabla P$ is invariant under translation by $\lambda d$ for all
$\lambda\in K$, that is, a translation direction of $Z_P$ in the sense of Section~\ref{sec:translation}. Thus $Z_P$ is a cylinder over a variety of smaller dimension, with rulings of direction $d$. Apex directions form a vector space, the \emph{apex
space}. An apex direction is a vector in the target space of $\nabla P$, whereas an affine direction is a vector in the source space. By Theorem \ref{thm:apex} the apex space of a rank-three
polynomial is nonzero.

\begin{lemma}\label{lem:apexbasic}
Let $P$ be weighted homogeneous for weights $w\in\Q^4$, of weight
$E$, with no affine part. The weights need not be positive.
\begin{enumerate}[label=\textup{(\alph*)}]
\item \textup{(Equivariance and blocks.)} $\nabla P$ is equivariant
for the source scaling with weights $w_i$ and the target scaling with
weights $E-w_i$. The apex space and the constant Hessian kernel split
into their intersections with the equal-weight blocks.
\item \textup{(Contraction.)} If $\rk\Hess P=3$, then
$\adj(\Hess P)=\rho\,nn\T$ for a nonzero $\rho$ and
$n=\nabla R(\nabla P)$, where $R$ is an irreducible equation of the
hypersurface $Z_P$. Every apex direction $d$ satisfies
$d\T\adj(\Hess P)\,d=0$.
\item \textup{(Three apices.)} If $\rk\Hess P=3$, the apex space has
dimension three if and only if $\Hess P$ has a nonzero constant
kernel vector. In that case $Z_P$ is a hyperplane.
\item \textup{(Constant kernel and variables.)} A constant kernel
vector $k$ of $\Hess P$ lying in one weight block satisfies
$\partial_kP=0$, and after a linear change within the block, $P$ does not
involve one variable.
\item \textup{(Rank one.)} If $\rk\Hess P=1$, then $P=f(\ell)$ for a
constant linear form $\ell$ and $f\in K[T]$. In particular $\Hess P$
has a three-dimensional constant kernel.
\item \textup{(Two variables.)} If $P$ involves at most
two variables and its Hessian in them is singular, then $\Hess P$
has a nonzero constant kernel vector in their span.
\item \textup{(Three variables.)} If $P$ involves exactly three variables and its Hessian in them is singular, then $\rk\Hess P\le2$.
\end{enumerate}
\end{lemma}

\begin{proof}
(a) Clearing denominators of $w$, the action
$x\mapsto(t^{w_i}x_i)$, $t\in K^\times$, multiplies $P$ by $t^E$ and $\partial_iP$ by
$t^{E-w_i}$. Hence $Z_P$ is invariant under the diagonal target action
$\sigma_t$, and with $d$ also $\sigma_t(d)$ is an apex direction.
Expanding into the eigenspaces of the action, which are the
equal-weight blocks, gives the splitting. The same argument applies to
the constant kernel, using the source action.

(b) The image of $\nabla P$ is dense in the irreducible hypersurface
$Z_P=\{R=0\}$, and differentiating $R(\nabla P)=0$ gives
$\Hess P\cdot n=0$ with $n=\nabla R(\nabla P)\ne0$ over the rational
function field ($R$ is reduced, and its gradient is nonzero on the
smooth locus which meets the dense image). As $\rk\Hess P=3$, the
kernel is spanned by $n$ and $\adj(\Hess P)$, whose columns lie in the
kernel and which is symmetric of rank one, equals $\rho\,nn\T$,
$\rho\ne0$. If $d$ is an apex, $D_dR$ vanishes on $Z_P$, is divisible
by the irreducible $R$, and has smaller degree, hence $D_dR=0$. Then
$n\cdot d=(D_dR)(\nabla P)=0$ and
$d\T\adj(\Hess P)d=\rho\,(n\cdot d)^2=0$.

(c) If $k$ is a constant kernel vector, $\nabla(\partial_kP)=
\Hess P\cdot k=0$, so $\partial_kP$ is a constant, which vanishes
because $P$ has no affine part. Thus $\vv k{\nabla P}=0$, so
$Z_P$ lies in the hyperplane $\vv kY=0$. Having dimension three it
equals that hyperplane, whose translation directions form the
three-dimensional space $k^\perp$. Conversely, if the apex space has
dimension three, $R$ is annihilated by three independent constant
derivations, hence is a polynomial in one linear form $\ell$, and
irreducibility gives $R=\ell-c$. Then $\vv{\ell}{\nabla P}\equiv c$
and $\Hess P\cdot\ell^\vee=0$ for the constant vector $\ell^\vee$ of
coefficients of $\ell$.

(d) As in (c), $\partial_kP=0$. A linear change within the block (which
preserves weighted homogeneity) moves $k$ to a coordinate vector.

(e) Since $\rk\Hess P=1$, the partial derivatives $\partial_1P,\dots,\partial_4P$
generate a $K$-algebra of transcendence degree one (in characteristic
zero the transcendence degree equals the generic rank of the Jacobian
matrix). By \ref{E:Zaks} there are a nonconstant $g\in K[x]$ and
$h_i\in K[T]$ with $\partial_iP=h_i(g)$. Then
$\partial_j\partial_iP=h_i'(g)\,g_j$, and symmetry of $\Hess P$ gives
$h_i'(g)g_j=h_j'(g)g_i$. Fix $i$ with $h_i'\ne0$ and put
$c_j=h_j'/h_i'\in K(T)$. Then $g_j=c_j(g)\,g_i$ for all $j$, that is,
$\nabla g(x)$ is a multiple of the vector $c(g(x))$ wherever $g(x)$ is
not a pole of the $c_j$. For all but finitely many $s\in K$ the level
set $\{g=s\}$ is smooth with $\nabla g\ne0$ on it, and $s$ is not a
pole of the $c_j$. Let $Z$ be an irreducible component of such a level
set. Its normal direction is the constant covector $N=c(s)$, so the
linear form $\vv Nx$ has zero differential on $Z$ and is constant on
$Z$ in characteristic zero. As $\dim Z=3$, $Z$ is a hyperplane with
normal $c(s)$. By Lemma~\ref{lem:linfib} with $k=n-1$ these hyperplanes are
parallel. Hence the direction of $c(s)$ is a fixed covector $v$ for
generic $s$, and $\nabla g\parallel v$ on a dense set. Thus
$g\in K[\ell]$ for $\ell=\vv vx$, and $\partial_iP=H_i(\ell)$ for
polynomials $H_i$. Now $\Hess P=H'(\ell)\,v\T$ is symmetric, so
$H'(\ell)\parallel v$, that is, $H=vF'+a$ for a polynomial $F$ and a
constant vector $a$. Therefore $P=F(\ell)+\vv ax+P(0)$, and the absence of
terms of degree at most one gives $\vv ax=-F'(0)\ell$, so
$P\in K[\ell]$. Then $\Hess P=P''(\ell)vv\T$ has the three-dimensional
constant kernel $v^\perp$.

(f) By \ref{E:dBvdE} with $F=K$ there are $\lambda\ne0$ and
$\gamma$ with $\sum_i\lambda_i\partial_iP=\gamma$. The left side has
no constant term because $P$ has no monomial of degree at most one, so
$\gamma=0$ and $\Hess P\cdot\lambda=\nabla(\sum_i\lambda_i\partial_iP)=0$.
The theorem \ref{E:GN} does not apply here, because wall forms are
weighted homogeneous but in general not homogeneous.

(g) The Hessian of $P$ in its three variables has rank at most two, and the fourth row and column of $\Hess P$ vanish. Unlike (f), a singular Hessian in three variables need not have a nonzero constant kernel vector. By the example after \ref{E:dBvdE}, $x_1^2x_2+x_1^3x_3$ is weighted homogeneous for the weights $(1,2,1)$, with singular Hessian and no nonzero constant kernel vector.
\end{proof}

\begin{lemma}[leading parts of minors]\label{lem:survive}
Let $P$ be a polynomial, $w\in\Q^4$ a weight vector, called the old
weights, and $H$ the leading form of $P$ for $w$, of old weight $E$, so
that $P-H$ consists of terms of old weight less than $E$. In the main
application $P$ is a wall form and $H$ is the leading form just before
the wall.
\begin{enumerate}[label=\textup{(\alph*)}]
\item Each product in entry $(i,j)$ of $\adj(\Hess P)$ has old weight
at most $3E-2\sum_kw_k+w_i+w_j$, and its component at that bound is
the corresponding entry for $H$. For constant vectors $u,v$ that are
homogeneous for the old weights, the old Hessian minor of $H$ on
$(u,v)$ survives as the component of maximal old weight of the minor
of $P$. In particular, if the old minor is nonzero, so is the new one.
\item \textup{(Contraction criterion.)} If a constant vector $d$ has a
unique nonzero coordinate of greatest old weight and the corresponding
principal cofactor of $\Hess H$ is nonzero, then
$d\T\adj(\Hess P)\,d\ne0$. By Lemma \textup{\ref{lem:apexbasic}(b)},
$d$ is not an apex direction of a rank-three $P$.
\end{enumerate}
\end{lemma}

\begin{proof}
(a) In the deleted-row-and-column expansion of an entry of the
adjugate, the omitted row and column carry the weights $w_j,w_i$, so
each product of three Hessian entries has old weight at most
$3E-2\sum_kw_k+w_i+w_j$. Terms of $P-H$ only contribute below that
bound, and the component at the bound is computed from $H$. The
$2\times2$ minor of the Hessian bilinear form on homogeneous constant
vectors $u,v$ is treated by the same comparison.
(b) In $d\T\adj(\Hess P)d=\sum_{i,j}d_id_j\adj(\Hess P)_{ij}$ the term
of greatest old weight is contributed by the diagonal entry of the
unique heaviest coordinate alone, and equals (the leading part of) the
corresponding principal cofactor of $\Hess H$, which is nonzero.
\end{proof}

\begin{lemma}[graph shears]\label{lem:shear}
Under a constant invertible source change $X'=SX$ the gradient becomes
$S^{-\mathsf T}(\nabla\phi)(S^{-1}X')$, so an apex $d$ becomes
$S^{-\mathsf T}d$. Suppose source and target are split into
corresponding base and fiber blocks and the apex plane of a wall form
is the graph of a matrix $U$ from the target base block to the target
fiber block. Then the source shear
$\mathrm{base}'=\mathrm{base}+U\T\,\mathrm{fiber}$ carries the apex
plane to the target base block. At a wall it preserves weights, and if
every fiber variable mixed into a base variable has strictly smaller
old weight, the substitution of the inverse shear preserves the old
leading form.
\end{lemma}

\begin{proof}
The transformation rule is immediate from the chain rule. The graph of
$U$ is invariant under the diagonal target scaling at the wall, so $U$
pairs only equal new weights and the shear preserves them. In the
substituted expansion of any monomial, each occurrence of a
replacement variable lowers the old weight strictly, so the old
leading form is unchanged.
\end{proof}

\begin{lemma}[common polynomial]\label{lem:euler}
Let $A,B$ be nonzero weighted homogeneous polynomials in two base
variables of positive weights, of weighted degrees $p,q>0$, with
$dA\wedge dB=0$. Then $A=aS^r$, $B=bS^s$ for nonzero constants $a,b$,
positive integers $r,s$, and a common weighted homogeneous polynomial
$S$. Call $S$ \emph{genuinely binary} if it is not a polynomial in a
constant linear form. Equivalently \textup{(Lemma~\ref{lem:2d})}
$\BH_2(S)\neq0$, equivalently $D_vS\ne0$ for every nonzero constant
base vector $v$. If $S$ is not genuinely binary and the base weights
differ, then $S$ is, up to a constant, a power of a coordinate.
\end{lemma}

\begin{proof}
Contracting $dA\wedge dB=0$ with the weighted Euler field gives
$pA\,dB=qB\,dA$. Clearing weight denominators to make $p,q$ integral,
$d(A^q/B^p)=0$, so $A^q/B^p\in K^\times$ and for each irreducible
factor the multiplicities in $A$ and $B$ are in ratio $p:q$. Dividing
by $\gcd(p,q)$ yields $A=aS^r$, $B=bS^s$. Comparison of the extreme
weights in $A=aS^r$ makes $S$ weighted homogeneous. If $S\in K[\ell]$
for a linear form $\ell$ and the base weights differ, the highest
ordinary power of $\ell$ in $S$ would otherwise contain monomials of
different weights and the same degree, so $\ell$ is proportional to a
coordinate. The equivalence with $\BH_2(S)\ne0$ is Lemma~\ref{lem:2d} (if
$S=f(\ell)$, then $\BH_2(S)=0$ by direct computation), and with the
nonvanishing of all $D_vS$ it is immediate.
\end{proof}

\subsection*{Phases}
Let $P$ be a wall form of rank three with one-dimensional apex space
aligned, by Lemma \ref{lem:apexbasic}(a) and a constant change within
its weight block, with a target coordinate axis. Call the
corresponding source variable $x$ \emph{active}, of weight $w_x>0$,
and the remaining variables $Y=(Y_1,Y_2,Y_3)$ \emph{passive}. By
\ref{E6} below, the part of $P$ of passive degree at least two is a
$K(x)$-linear combination of powers $L^k$, $k\ge2$, of a single
\emph{phase} $L=\ell(x)\T Y$, and the passive Hessian is
$\big(\partial_L^2\big)\,\ell\ell\T$ of rank one. A \emph{phase
variable} is a passive variable with nonzero coefficient in $L$. When
one nonzero entry $\ell_i$ is normalized to $1$, the ratio
$\ell_j=\big(\Hess_YP\big)_{ij}/\big(\Hess_YP\big)_{ii}$ is a rational
function of $x$, homogeneous of weight $w_{Y_i}-w_{Y_j}$ for the
scaling of Lemma \ref{lem:apexbasic}(a). A nonzero homogeneous
rational function of one variable has no zero or pole in
$K^\times$, so $\ell_j$ is zero or a constant times
$x^{(w_{Y_i}-w_{Y_j})/w_x}$ with integral exponent. For $k\ge2$ the
coefficient of $L^k$ equals the coefficient of $Y_i^k$ in the
restriction to the other passive variables set to zero, hence is a
polynomial in $x$.

\subsection*{The rank classification of de Bondt}
The following two statements are due to de Bondt
\cite[Theorem 7.2, Corollary 7.3]{dB16}. They are also used in the
quartic classification of Ni \cite{Ni26}. We prove them in Lemmas \ref{lem:E6} and \ref{lem:E7}
by the methods of this paper. We state them for a polynomial $P$ in four variables over
$K$ with the apex space aligned with target coordinate axes.

\begin{enumerate}[label=\textup{(B\arabic*)}]
\item\label{E6} If $\rk\Hess P=3$ and the apex space is
$\spn(e_1)$, then the Hessian of $P$ in the complementary variables
$x_2,x_3,x_4$ has rank one, and
$P=F(x_1,L)+M$, where $L=\ell_2(x_1)x_2+\ell_3(x_1)x_3+\ell_4(x_1)x_4$ with $\ell_2,\ell_3,\ell_4\in K[x_1]$ without common factor, $F\in K[x_1][T]$ contains only powers $T^j$ with $j\ge2$ and $\partial_T^2F\ne0$, and
$M\in K[x_1][x_2,x_3,x_4]$ has degree at most one in $x_2,x_3,x_4$. If the apex space is
$\spn(e_1,e_2)$, then $P=G(x_1,x_2)+x_3A(x_1,x_2)+x_4B(x_1,x_2)$, and the Jacobian determinant of $(A,B)$ vanishes.
\item\label{E7} If $\rk\Hess P=2$, then either $\Hess P$ has a
two-dimensional constant kernel, or the apex space is one-dimensional
and, aligned with $e_1$, the polynomial $P$ is jointly affine in
$x_2,x_3,x_4$, that is $P=G(x_1)+x_2A(x_1)+x_3B(x_1)+x_4C(x_1)$.
\end{enumerate}

\noindent
The three-apex case of rank three is Lemma \ref{lem:apexbasic}(c).

Let $Z=\{R=0\}$ be the closure of the image of $\nabla P$, with $R$
irreducible, and let $d$ be an apex direction. For $a\in K$ the
polynomial $R(Y+ad)$ generates the ideal of $Z$, so
$R(Y+ad)=c_aR(Y)$ with $c_a\in K^\times$. The homogeneous parts of
highest degree agree, so $c_a=1$, and differentiation at $a=0$ gives
$D_dR=0$. Hence, if $e_1,\dots,e_r$ are apex directions, then
$R\in K[Y_{r+1},\dots,Y_m]$, where $m$ is the number of variables.

\begin{lemma}\label{lem:E6}
Statement \ref{E6} holds.
\end{lemma}

\begin{proof}
Since $\rk\Hess P=3$, the closure $Z$ of the image of $\nabla P$ is a
hypersurface $\{R=0\}$.

Let the apex space be $Ke_1$. Then $R\in K[Y_2,Y_3,Y_4]$, so
$Z=K\times S$ with $S=\{R=0\}\subseteq K^3$ an irreducible surface.
Every translation direction $w$ of $S$ gives the apex direction
$(0,w)$ of $Z$, so $S$ has none. Put $z=(x_2,x_3,x_4)$ and
$C=\Hess_zP$. For every $t\in K$ the values of $\nabla_zP(t,z)$ lie in
$S$, so $\rk C\le2$. Suppose $\rk C=2$. For general $t_0\in K$ the
polynomial $P_0=P(t_0,z)$ has Hessian of rank two, and the closure of
the image of $\nabla P_0$ is $S$. By \ref{E:dB62} (Appendix
\ref{app:dB62}) either $c\T\nabla P_0$ is constant for some
$c\in K^3\setminus\{0\}$, and $S$ is an affine plane, or after a
linear change of coordinates $P_0=a_1(u_1)+a_2(u_1)u_2+a_3(u_1)u_3$
and, in the dual coordinates of the target, $S=K\times\Gamma$ with
$\Gamma$ the closure of the curve $(a_2,a_3)(K)$. In both cases $S$
has a translation direction, a contradiction. If $C=0$, then
$\rk\Hess P\le2$. Hence $\rk C=1$. By Lemma \ref{lem:app-low} over the
algebraic closure of $K(x_1)$, the matrix $C$ vanishes on a
two-dimensional subspace, which is spanned by vectors over $K(x_1)$ by Lemma~\ref{lem:lindescent}. Let $L=\ell_2x_2+\ell_3x_3+\ell_4x_4$ be a linear form
vanishing on it, with $\ell_2,\ell_3,\ell_4\in K[x_1]$ without common factor. Then $\partial_wP\in K(x_1)$ for every $w$ in this
subspace, and $P=\sum_jf_jL^j+M_0$ with $f_j\in K(x_1)$ and $M_0$ affine in $z$ over $K(x_1)$. For $j\ge2$ the part of $P$ of degree $j$ in $z$ is $f_jL^j$. It lies in $K[x_1][z]$, and $L^j$ is primitive over $K[x_1]$ by Gauss's lemma, so $f_j\in K[x_1]$. Put $F=\sum_{j\ge2}f_jT^j$. Then $M=P-F(x_1,L)$ lies in $K[x_1][z]$ and has degree at most one in $z$, and $\partial_T^2F\ne0$ since $C\ne0$.

Let the apex space be $\spn(e_1,e_2)$. Then $R\in K[Y_3,Y_4]$, so
$Z=K^2\times\Gamma$ with $\Gamma=\{R=0\}\subseteq K^2$ an irreducible
curve. A translation along $\Gamma$ would give a third apex
direction, so $\Gamma$ is not a line. The values of
$(P_{x_3},P_{x_4})$ lie in $\Gamma$, so $B=\Hess_{(x_3,x_4)}P$ has
rank at most one. Suppose $\rk B=1$. By Lemma \ref{lem:app-binary}
over $K(x_1,x_2)$ there are $c\in K(x_1,x_2)^2\setminus\{0\}$ and
$\gamma\in K(x_1,x_2)$ with $c_1P_{x_3}+c_2P_{x_4}=\gamma$. For general
$(x_1,x_2)\in K^2$ the image of $(x_3,x_4)\mapsto(P_{x_3},P_{x_4})$
is then a curve contained in $\Gamma$ and in a line, so $\Gamma$ is a
line, a contradiction. Hence $B=0$, and $P$ is affine in
$(x_3,x_4)$ with coefficients in $K[x_1,x_2]$. For this form $\det\Hess P$ is the square of the Jacobian determinant of $(A,B)$, which therefore vanishes.
\end{proof}

\begin{lemma}\label{lem:E7}
Statement \ref{E7} holds.
\end{lemma}

\begin{proof}
Put $H=\Hess P$ and $F=\nabla P$. The closure $Z$ of the image of $F$
is an irreducible surface. Choose $R_1,R_2$ in the ideal of $Z$ with
$\nabla R_1,\nabla R_2$ linearly independent at a general point of
$Z$, and put $G_i=(\nabla R_i)(F)$. Differentiating $R_i(F)=0$ gives
$HG_i=0$, and then $D_{G_j}G_i=\Hess R_i(F)\,HG_j=0$ and $D_{G_j}F=0$.
As in the remark after Lemma \ref{lem:graph}, applied first to $G_1$
and then to $G_2$, this gives
$F(x+sG_1(x)+tG_2(x))=F(x)$ for all $s,t\in K$. For general $x$ the
vectors $G_1(x),G_2(x)$ span $\ker H(x)$, and the affine plane
$x+\ker H(x)$ lies in the fiber of $F$ through $x$. The general fibers
have dimension two. By Lemma~\ref{lem:linfib} the irreducible components of the general fibers are these affine
planes, two of them are equal or disjoint, and the direction of the component through a
general point $x$ is $\ker H(x)$. Hence any two directions $U,U'$ of such components satisfy $U+U'\neq K^4$, that is, they meet
nontrivially. Then either all directions contain a fixed line, or all
lie in a fixed hyperplane. Indeed, let $U\ne U'$, $L=U\cap U'$ and
$W=U+U'$. A direction $U_0$ not containing $L$ meets $U$ and $U'$ in
distinct lines, so $U_0\subseteq W$. If such a $U_0$ exists, a
direction $U_1\supseteq L$ with $U_1\not\subseteq W$ satisfies
$U_1\cap W=L$, hence $U_1\cap U_0=0$, which is excluded. So all
directions lie in $W$.

\emph{All directions in a hyperplane $W$.} Let $d$ span $W^\perp$.
Since $H$ is symmetric, $\operatorname{im}H(x)=\ker H(x)^\perp$
contains $d$. After a linear change of coordinates $d=e_1$. Write
$z=(x_2,x_3,x_4)$ and $C=\Hess_zP$. The projection of
$\operatorname{im}H$ to the $z$-coordinates has dimension at most one,
so $\rk C\le1$, and the column $\partial_{x_1}\nabla_zP$ lies in the
column space of $C$. If $C=0$, then $P$ is jointly affine in $z$, the
second alternative. Let $\rk C=1$ and $F_1=K(x_1)$. Apply
\ref{E:dB62} to $P\in F_1[z]$. In its second case
$P=a_1(z_1)+a_2(z_1)z_2+a_3(z_1)z_3$ in suitable $F_1$-linear
coordinates $z_1,z_2,z_3$, and $\rk C=1$ forces $a_2'=a_3'=0$. In its
first case $P=g(z_1,z_2)+\gamma z_3$ in suitable coordinates, with
$\gamma\in F_1$ and $\det\Hess g=0$, and Lemma \ref{lem:hess2F}
applies to $g$. In both cases
\[
P=Q(x_1,\ell)+c\T z+c_0,\qquad \ell=a\T z,
\]
with $a\in F_1^3\setminus\{0\}$, $c\in F_1^3$, $c_0\in F_1$,
$Q\in F_1[T]$ and $Q_{TT}\ne0$, and $C=Q_{TT}(x_1,\ell)\,aa\T$. Let a
prime denote $\partial_{x_1}$. The vector
\[
\partial_{x_1}\nabla_zP=\big(Q_{Tx_1}(x_1,\ell)+Q_{TT}(x_1,\ell)\,a'^{\mathsf T}z\big)a+Q_T(x_1,\ell)\,a'+c'
\]
lies in $F_1(z)\,a$, so $Q_T(x_1,\ell)\,a'+c'\in F_1[z]\,a$. Since
$Q_{TT}\ne0$, a positive power of $\ell$ occurs in $Q_T(x_1,\ell)$,
and comparison of degrees in $z$ gives $a'\in F_1a$ and $c'\in F_1a$.
Hence $a=\mu a_0$ with $a_0\in K^3$ and $\mu\in F_1^\times$, and
$c\T v\in K$ for every $v\in K^3$ with $a_0\T v=0$. For such $v$ the
derivative $\partial_{(0,v)}P=c\T v$ is constant, and $\Hess P$ has
a two-dimensional constant kernel.

\emph{All directions contain a fixed line $Kk$.} Then $Hk=0$, so
$\partial_kP$ is constant. After a linear change of coordinates $k=e_4$
and $P=P_0(x_1,x_2,x_3)+\gamma x_4$ with $\gamma\in K$ and
$\rk\Hess P_0=2$. By \ref{E:dB62} either $\partial_cP_0\in K$ for some
$c\in K^3\setminus\{0\}$, and $(c,0)$ and $e_4$ span a constant
kernel, or after a linear change
$P_0=a_1(x_1)+a_2(x_1)x_2+a_3(x_1)x_3$, and $P$ is jointly affine in
$x_2,x_3,x_4$.

\emph{The apex space.} Let $P=G(x_1)+x_2A(x_1)+x_3B(x_1)+x_4C(x_1)$ as
in \ref{E7} and $\alpha=(A,B,C)\T$. Then
$\nabla P=(G'+z\T\alpha',\,\alpha)$, and the closure of its image is
$K\times\Gamma$, where $\Gamma$ is the closure of $\alpha(K)$, a curve
since $\rk H=2$. A vector $(\lambda,w)$ is an apex
direction if and only if $\Gamma$ is invariant under translation by
$w$. For $w\ne0$ this makes $\Gamma$ a line, then $\alpha'$ is proportional
to a fixed vector and $\Hess P$ has a two-dimensional constant
kernel. Hence the apex space is $Ke_1$ if there is no such kernel.
\end{proof}

\subsection*{Vanishing bordered Hessian in four variables}

\begin{proposition}\label{prop:BH4}
Let $A\in K[x_1,\dots,x_4]$. Then $\BH(A)=0$ if and only if, after a linear change of coordinates, one of the following holds.
\begin{enumerate}[label=\textup{(\roman*)}]
\item $A\in K[\lambda_1,\lambda_2,\lambda_3]$ for linear forms $\lambda_1,\lambda_2,\lambda_3$.
\item $A=g(x_1,L)$ with $g\in K[x_1,s]$ and $L=\ell_2(x_1)x_2+\ell_3(x_1)x_3+\ell_4(x_1)x_4$, $\ell_2,\ell_3,\ell_4\in K[x_1]$.
\end{enumerate}
\end{proposition}

\begin{proof}
If $\rk\Hess A\le2$, then $\adj\Hess A=0$. If $\rk\Hess A=3$ and $u$ spans the kernel of $\Hess A$, then $\adj\Hess A=\rho\,uu\T$ with $\rho\ne0$ and $\BH(A)=\rho\vv{\grad A,u}^2$.

Let (i) hold. There is $\kappa\in K^4\setminus\{0\}$ with $\partial_\kappa A=0$, so $\Hess A\,\kappa=0$ and $\vv{\grad A,\kappa}=0$, and $\BH(A)=0$. Let (ii) hold, write $\ell=(\ell_2,\ell_3,\ell_4)$ and put $u=(0,\ell\times\ell')$. If $u=0$, then $\ell\in K[x_1]\,c$ with $c\in K^3$, and (i) holds. Otherwise $x_1$, $L$ and $\ell_2'x_2+\ell_3'x_3+\ell_4'x_4$ are annihilated by $D_u$, the entries of $\grad A$ are polynomials in them, so $\Hess A\,u=0$ and $\vv{\grad A,u}=0$, and $\BH(A)=0$.

Conversely let $\BH(A)=0$. By Theorem~\ref{thm:main-BH}(i), $\det\Hess A=0$.

If $\rk\Hess A\le1$, Lemma~\ref{lem:app-low} gives a three-dimensional subspace $U$ with $\Hess A\,U=0$. Then $\partial_uA\in K$ for $u\in U$, and $A=f(\lambda_1)+\lambda_2$ with linear forms $\lambda_1,\lambda_2$, which is (i).

If $\rk\Hess A=2$, then by \ref{E7} either $\Hess A$ has a two-dimensional constant kernel $N$, and then $\partial_kA\in K$ for $k\in N$ and $A\in K[\lambda_1,\lambda_2,\lambda_3]$, or $A=G(x_1)+x_2A_2(x_1)+x_3A_3(x_1)+x_4A_4(x_1)$, which is (ii) with $g=G+s$ and $\ell=(A_2,A_3,A_4)$.

Let $\rk\Hess A=3$. Let $Z=\{R=0\}$ be the closure of the image of $\grad A$, with $R$ irreducible, and $Q=(\grad R)(\grad A)$. By Proposition~\ref{prop:qt-hess}, $Q$ spans the kernel of $\Hess A$. Since $\vv{\grad A,Q}=(\mathcal ER)(\grad A)$ with $\mathcal E=\sum_iY_i\partial_{Y_i}$, the condition $\BH(A)=0$ says that $R$ divides $\mathcal ER$, that is, that $R$ is homogeneous. By Theorem~\ref{thm:apex} the apex space $V$ of $Z$ is nonzero, and $\dim V\le3$, so that $Z$ is a cone with vertex $0$ which is invariant under the translations by $V$.

$\dim V=3$. Then $R$ is linear, and after a linear change $A=G(x_1,x_2,x_3)+cx_4$ with $c\in K$. Then $\BH(A)=c^2\det\Hess G$, and $\det\Hess G\ne0$ since the rank is three. Hence $c=0$, which is (i).

$\dim V=2$, aligned with $e_1,e_2$. By \ref{E6}, $A=G(y)+x_3B_1(y)+x_4B_2(y)$ with $y=(x_1,x_2)$, and the Jacobian determinant of $(B_1,B_2)$ vanishes. Choose $c\in K(y)^2\setminus\{0\}$ with $c_1\grad B_1+c_2\grad B_2=0$. Then $u=(0,0,c_1,c_2)$ spans the kernel of $\Hess A$, and $\BH(A)=0$ gives $c_1B_1+c_2B_2=0$. If $c_2=0$, then $B_1=0$. Otherwise $B_2=rB_1$ with $r=-c_1/c_2$, and the two relations give $c_2B_1\grad r=0$, so $B_1=0$ or $r\in K$. In both cases $B_1,B_2$ are linearly dependent over $K$, which is (i).

$\dim V=1$, aligned with $e_1$. By \ref{E6}, $A=F(x_1,L)+M_0(x_1)+\vv{b(x_1),z}$ with $z=(x_2,x_3,x_4)$, $L=\vv{\ell(x_1),z}$, $\ell\in K[x_1]^3$ with coprime entries, $F\in K[x_1][T]$ with $\partial_T^2F\ne0$, and $b\in K[x_1]^3$. Then $Z=K\times S$, where $S\subseteq K^3$ is the closure of the image of $(\partial_TF)(x_1,L)\,\ell+b$, and $S$ has no translation direction (proof of Lemma~\ref{lem:E6}). For general $t\in K$ the polynomial $(\partial_TF)(t,T)$ is nonconstant and $\ell(t)\ne0$, so $S$ contains the line $b(t)+K\ell(t)$. Since $R$ is homogeneous, $S$ contains $\spn(\ell(t),b(t))$. If this span were a plane, it would equal $S$, which has no translation direction. Hence $b\in K(x_1)\,\ell$, and $b=\beta\ell$ with $\beta\in K[x_1]$, since the entries of $\ell$ are coprime. So $A=g(x_1,L)$ with $g=F+M_0+\beta s$, which is (ii).
\end{proof}
\subsection{Exit lemmas}

Seven of the lemmas of this subsection, Lemmas~\ref{lem:D1}, \ref{lem:D2},
\ref{lem:D3}, \ref{lem:D4}, \ref{lem:D5}, \ref{lem:D6} and~\ref{lem:D7},
follow \cite[Lemmas 3.2--3.8]{Tan26}, in this order, with complete proofs
and with all inversion arguments removed. Lemma~\ref{lem:D1-walls} isolates
the count of walls from the proof of \cite[Lemma~3.2]{Tan26}, and
Lemma~\ref{lem:P1} is the weight change of the proof of
\cite[Proposition~4.1]{Tan26}. Each of the seven lemmas either derives a
contradiction or exhibits an affine direction of $\phi$. The standing hypotheses are in force.

\begin{lemma}[two fibers with a common binary polynomial]\label{lem:D1}
Suppose the leading form of $\phi$, for positive weights with the
fiber weights $w_z\ne w_w$ and $W>0$, is
$R(x,y)+zA(x,y)+wB(x,y)$ with $A=aS^r$, $B=bS^s$ as in Lemma
\textup{\ref{lem:euler}} and $S$ genuinely binary. Then $\phi$ has an
affine direction.
\end{lemma}

\begin{proof}
Raise both fiber weights at unit speed, fixing the base weights. After
a small step the leading form is $H=zA+wB$, and $W$ is constant (the
leading weight has slope $1$ and the weight sum slope $2$, so $W$ has
slope $4\cdot1-2\cdot2=0$). In particular $W>0$ along the path.

Put $U:=\BH_2(S)\ne0$ (Lemma \ref{lem:euler}). Writing $A',B'$ for
the derivatives of $aT^r,bT^s$ at $T=S$, the base Hessian of $H$ is
$(zA'+wB')\Hess S+(zA''+wB'')(\nabla S)(\nabla S)\T$. By the bordered
expansion (Lemma~\ref{lem:det2}) the cofactor of $\Hess H$ complementary
to $z$ equals $-(B')^2(zA'+wB')\,U$, and the one complementary to $w$
equals $-(A')^2(zA'+wB')\,U$. Both are nonzero. For every nonzero
constant base vector $v$, the minors on $(v,e_z)$ and $(v,e_w)$ are
$-(D_vA)^2$ and $-(D_vB)^2$, nonzero because $D_vS\ne0$.

At a first wall these cofactors and minors survive
(Lemma \ref{lem:survive}(a)), so the wall form has rank three and, by
Theorem \ref{thm:apex}, a nonzero apex space, which splits into blocks
(Lemma \ref{lem:apexbasic}(a)). The two fiber axes are excluded by the
contraction criterion with the two principal cofactors, and the fiber
weights are distinct, so the apex space meets the fiber plane
trivially. If it were a line, its weight block would contain at most
one fiber axis. Choosing the other fiber variable $F$ and a constant
base direction $v$, homogeneous for the old weights and annihilated by
the apex (possible since the base weights are equal, or the apex line
has at most one base component), the surviving minor on $(v,e_F)$ is
nonzero, contradicting \ref{E6}, since on the complement of a line apex
the Hessian has rank one and all its $2\times2$ minors vanish. In
detail, the vector $v$ is chosen as follows. If the two base weights
are equal, $v$ is a nonzero base vector annihilated by the base
component of the apex. If they differ, the apex line lies in one weight
block, its base component is a multiple of one base axis, and $v$ is
the other base axis. In both cases $v$ is homogeneous for the old
weights. The minor of $H$ on $(v,e_F)$ is $-(D_vA)^2$ or $-(D_vB)^2$,
which is nonzero since $D_vS\ne0$. By Lemma \ref{lem:survive}(a) it is
the component of maximal old weight of the minor of the wall form on
$(v,e_F)$, because every other product in that minor contains a factor
coming from a term of old weight less than $E$. Hence
the apex space is a plane, a graph over the target base block, and the
graph shear of Lemma \ref{lem:shear} puts the wall form in jointly
affine form in the two fibers.

At a wall a distance $\tau>0$ along the segment, a tied monomial of
total fiber degree $j$ has old weight $E+(1-j)\tau$. Hence the fiber-free
part of the wall form vanishes, and its fiber-linear part is exactly
$zA+wB$. The inverse shear fixes the fibers and is the identity on the
base at $z=w=0$, so it leaves the fiber-linear part of a form with
zero fiber-free part unchanged. The shear restores $H$ exactly, and
the path continues with the same speeds. By Lemma \ref{lem:D1-walls}
there are at most four genuine walls, and the cofactors and minors
above are the same at each of them. Beyond the last wall the path continues
indefinitely, and a monomial of $\phi$ of total fiber degree at least
two would eventually overtake $H$. Hence there is none, and $\phi$ is,
in the constant coordinates accumulated from the shears, jointly
affine in the two fiber variables. In particular
$\partial_z^2\phi=0$, so $e_z$ is an affine direction.
\end{proof}

The following lemma makes the bound on the walls in the proof of
Lemma \ref{lem:D1} explicit. It separates the count of walls, which is
a statement on the path of weights, from the analysis of a single wall.

\begin{lemma}\label{lem:D1-walls}
In the proof of Lemma \textup{\ref{lem:D1}} the path has at most four
genuine walls. At each of them the leading form just before the wall
is $H=zA+wB$ with the same $A$, $B$ and $S$, in the coordinates reached
at that wall, and the base weights keep their values and their order.
\end{lemma}

\begin{proof}
Let $\tau\ge0$ be the path parameter. The weights are attached to the
four coordinate positions, and a shear at a wall preserves them by
Lemma \ref{lem:shear}. The speeds are not changed after a shear.
Hence, along the whole path and across all shears, the base weights
$w_x$, $w_y$ are constant and the fiber weights are
\[
w_z(\tau)=w_z(0)+\tau,\qquad w_w(\tau)=w_w(0)+\tau .
\]
In particular the base weights keep their values and their order.

At a wall the graph shear pairs a target fiber coordinate with a target
base coordinate only if their new weights are equal (Lemma
\ref{lem:shear}). If no fiber weight equals a base weight at the wall,
the shear is the identity. Then the wall form is jointly affine in the
fibers in the original coordinates, and by the proof of Lemma
\ref{lem:D1} it equals $H$. No monomial joins $H$ at such a parameter,
so it is not a genuine wall. Therefore every genuine wall occurs at a
parameter $\tau$ which satisfies one of the four equations
\[
w_z(\tau)=w_x,\qquad w_z(\tau)=w_y,\qquad w_w(\tau)=w_x,\qquad
w_w(\tau)=w_y .
\]
Each of them is affine in $\tau$ with slope one, so it has at most one
solution. The parameter $\tau$ increases along the path, so each
equation holds at one point of the path at most. A shear does not
change the functions $w_z(\tau)$, $w_w(\tau)$, $w_x$, $w_y$, so it
cannot make an equation hold again at a later parameter. Hence there
are at most four genuine walls.

By the proof of Lemma \ref{lem:D1}, after the shear at a wall the wall
form is $H$ in the new coordinates, with the same $A$, $B$ and $S$. A
small step beyond the wall selects $H$ again, since all its monomials
have slope one. So the leading form just before the next wall is $H$,
and the cofactors and minors used in the proof of Lemma \ref{lem:D1}
are computed from the same $H$ at every wall.
\end{proof}

Lemma~\ref{lem:D2} excludes a leading form whose part of passive degree at least two is a power of a phase in all three passive variables. It is called in the proofs of Lemmas~\ref{lem:D3}, \ref{lem:D4} and~\ref{lem:D5}.

\begin{lemma}[no ternary phase]
\label{lem:D2}
The leading form of $\phi$, for weights with $w_x=1$, three distinct
passive weights greater than one, and $W>0$, cannot be
\[
a(x)\,\big(Y_1+bx^rY_2+cx^sY_3\big)^m,\qquad
a\ne0,\ m\ge2,\ bc\ne0,
\]
with $0,r,s\in\mathbb Z$ distinct \textup{(}the expression being
polynomial even if some phase coefficients are Laurent
monomials\textup{)}.
\end{lemma}

\begin{proof}
For a polynomial $F(x,L)$ over $K(x)$ in a passive linear form
$L=\ell(x)\T Y$, the passive Hessian is $F_{LL}\,\ell\ell\T$ and the
mixed Hessian row is $(F_{xL}+F_{LL}\partial_xL)\ell\T+F_L(\ell')\T$.
Expanding the principal cofactor that deletes a passive variable along
the active row and column, with the other two passive indices $i,j$,
gives
\begin{equation}\label{eq:phase-cof}
\text{cofactor}
=-F_{LL}\,F_L^2\,\big(\ell_i'\ell_j-\ell_j'\ell_i\big)^2 .
\end{equation}
For coefficients proportional to $1,x^r,x^s$ with distinct exponents,
all three passive cofactors are nonzero.

Raise all passive weights at unit speed. The leading weight has slope
$m$ and $W$ slope $4m-6\ge2$, so $W$ stays positive and a first wall
occurs, since a monomial of $\phi$ outside the phase powers gains on the
leading form, which is bounded by the passive degree of $\phi$. If no
wall occurs, the final-segment analysis below applies directly.
At the wall the three cofactors survive, the wall form has rank three,
and since all four weights are distinct, the apex space is spanned by
coordinate axes. The three passive axes are excluded by the
contraction criterion, so the apex space is the active axis. We write out the comparison at the wall. Let $\tau>0$ be the
distance of the wall from the start of the segment, and $E$ the old
leading weight. The weight of a monomial grows with slope equal to its
passive degree $j$, and the leading weight with slope $m$. So a
monomial which is tied at the wall has old weight $E+(m-j)\tau$. Old
weights are at most $E$, hence $j\ge m$, and $j=m$ holds exactly for
the monomials of the old form. By \ref{E6} the wall form is
$P=\sum_{k\ge2}f_kL^k+M_1$ with
$L=Y_1+\ell_2Y_2+\ell_3Y_3$, $f_k\in K[x]$ as in the paragraph on
phases, and $M_1$ of passive degree at most one. Here the $Y_1$-coefficient of $L$ is normalized to one,
which is possible since the term $a(x)Y_1^m$ of the old form occurs in
$P$, as shown next. The passive degree $k$ part of $P$ is $f_kL^k$ for $k\ge2$,
since $L^k$ is homogeneous of passive degree $k$ and $M_1$ has passive
degree at most one. Hence $f_k=0$ for $k<m$, $M_1=0$, and
\[
f_mL^m=a(x)\big(Y_1+bx^rY_2+cx^sY_3\big)^m .
\]
The coefficient of $Y_1^m$ on the left is $f_m$, and only the power
$L^m$ contributes to it, since the other powers have other passive
degrees. So $f_m=a(x)$. The coefficient of $Y_1^{m-1}Y_j$ is
$mf_m\ell_j$ on the left, again from $L^m$ alone, and $ma(x)bx^r$ or
$ma(x)cx^s$ on the right. Hence $\ell_2=bx^r$ and $\ell_3=cx^s$, and the
phase of the wall form is the old phase. Therefore
$P=\sum_{k=m}^Mf_kL^k$ with $f_M\ne0$. If $M=m$, then $P$ is the old
form, no monomial joins it, and the parameter is not a wall. So
$M>m$. A small common raise of the passive weights beyond the wall
gives every monomial the slope of its passive degree and selects
$f_ML^M$, which has the shape of the lemma with $m$ replaced by $M$ and
the same $b,r,c,s$. The weights $w_{Y_1}-w_{Y_j}$ are constant along
the path, so $r,s$ do not change. Each genuine wall raises the phase
degree strictly, and the phase degree is at most the maximal passive
degree $D_Y$ of $\phi$. So no earlier value of the phase degree
recurs, there are at most $D_Y-2$ genuine walls, and after the last of
them the current power equals $D_Y$.

Let $N$ be the sum of the terms of $\phi$ of passive degree $D_Y$.
For the weights $(0,1,1,1)$ the leading form of $\phi$ is $N$ and
$W=4D_Y-6>0$, so $\det\Hess N=0$ by Lemma \ref{lem:W}. At a point of the final
segment, the leading form of $\phi$ consists of passive-degree-$D_Y$
terms only, so it is also the leading form of $N$. Its three passive
cofactors survive in $N$ by Lemma \ref{lem:survive}(a), applied to $N$
and its leading form for a weight vector of the final segment.

We also record a count of passive degrees, which is used below. Let
$Q=\sum_dQ_d$ be a polynomial, with $Q_d$ of passive degree $d$, and let
$d_0$ be the largest $d$ with $Q_d\ne0$. An entry $\partial_p\partial_qQ_d$
of the Hessian has passive degree $d-\pi(p,q)$, where $\pi(p,q)$ is the
number of passive variables among $p,q$. The principal cofactor which
deletes one passive variable is a sum of products of three entries
indexed by $(p,\sigma(p))$, where $p$ runs over the active variable and
the other two passive variables and $\sigma$ is a permutation of these
three. Each of the two passive variables occurs once as a row index and
once as a column index, so in every product the numbers
$\pi(p,\sigma(p))$ add up to four. A product
of entries taken from $Q_{d_1},Q_{d_2},Q_{d_3}$ therefore has passive
degree $d_1+d_2+d_3-4\le3d_0-4$, with equality only for
$d_1=d_2=d_3=d_0$. Hence the component of passive degree $3d_0-4$ of
this cofactor of $Q$ is the same cofactor of $Q_{d_0}$. If the latter is
nonzero, so is the former. The
passive weights on the segment are distinct, so by the contraction
criterion (Lemma \ref{lem:survive}(b) for $N$ and its leading form)
and the splitting under $(0,1,1,1)$ (Lemma \ref{lem:apexbasic}(a)),
the apex space of the rank-three $N$ is the active axis. By \ref{E6} and passive
homogeneity, $N=a(x)\big(l(x)\T Y\big)^{D_Y}$ with $l$ rational. The
components of $l$ are linearly independent over $K$, since a constant
relation would give a constant passive kernel vector and, by
Lemma \ref{lem:apexbasic}(c), a three-dimensional apex space.

Since $l_1,l_2,l_3$ are linearly independent over $K$, their Wronskian is nonzero. Indeed, at a regular point,
constant row operations produce distinct leading orders
$n_1,n_2,n_3$ of the Taylor expansions, and the lowest-order
coefficient of the Wronskian is a nonzero multiple of a Vandermonde
determinant in the $n_i$. Choose $x_0$ regular for $a$ and $l$, with
$a(x_0)\ne0$ and $l(x_0),l'(x_0),l''(x_0)$ independent. Translate
$x=x_0+\xi$ and take passive coordinates
\[
u=l(x_0)\T Y,\qquad v=l'(x_0)\T Y,\qquad t=\tfrac12\,l''(x_0)\T Y .
\]
Then $l(x_0+\xi)\T Y=u+\xi v+\xi^2t+(\text{terms of order}\ge3$ in
$\xi)$, so every monomial $\xi^ku^iv^jt^h$ of $N$ satisfies
$k\ge j+2h$, with the terms of equality summing to
$P_0=a(x_0)(u+\xi v+\xi^2t)^{D_Y}$.

Subtract the affine part of $\phi$ in the new coordinates and assign
the weights $(-\eta,1,1+\eta,1+2\eta)$ to $(\xi,u,v,t)$, $\eta>0$. The
weight of $\xi^ku^iv^jt^h$ is $d+\eta(j+2h-k)$ with $d=i+j+h$ the
passive degree. The terms of $P_0$ have weight $D_Y$ for every
$\eta$. A monomial of $\phi$ of passive degree $d$ with $2\le d<D_Y$
and positive slope reaches weight $D_Y$ at some $\eta\le D_Y-2$. For
passive degree one, the only slope-two term is the removed multiple of
$t$, and slope at most one reaches weight $D_Y$ no sooner than
$\eta=D_Y-1$. Passive degree zero has nonpositive slope. We give the details. For small $\eta>0$ the leading form is $P_0$,
of weight $D_Y$, since the other monomials of $N$ have negative slope
and the monomials of passive degree $d<D_Y$ have weight close to $d$.
The leading weight stays $D_Y$ up to the first wall. A wall parameter
solves an equation $d+\eta\sigma=D_Y$ with integers $d,\sigma$, so it is
rational. All weights are rational numbers attached to monomials of the
polynomial ring $K[\xi,u,v,t]$, and every comparison below is a
comparison of such numbers. No fractional powers of $\xi$ occur.

A monomial of passive degree $d$ with $2\le d<D_Y$ and slope
$\sigma=j+2h-k\ge1$ reaches the weight $D_Y$ at
\[
\eta_\sigma=\frac{D_Y-d}{\sigma}\le D_Y-2 ,
\]
with equality only for $d=2$, $\sigma=1$. A monomial of passive degree
one has $\sigma\le2$, and $\sigma=2$ holds only for $t$ itself, which
has been removed with the affine part. For $\sigma\le1$ its weight is at
most $1+\eta$, which is less than $D_Y$ for $\eta<D_Y-1$. A monomial of
passive degree zero has weight $-\eta k\le0<D_Y$.

Suppose a wall occurs at some $\eta$ with $0<\eta\le D_Y-2$. The
endpoint $\eta=D_Y-2$ is included. There $W=4D_Y-6-4\eta\ge2>0$. Let
$P_w$ be the wall form. Its monomials have weight $D_Y$. Those of
passive degree $D_Y$ are the monomials of $N$ with $j+2h-k=0$, that is,
the monomials of $P_0$. Those of passive degree one or zero do not
occur, by the bounds above, since $\eta\le D_Y-2<D_Y-1$. The three
passive cofactors of $P_0$ are \eqref{eq:phase-cof} with coefficients
$1,\xi,\xi^2$, and they are nonzero. By the count of passive degrees
above, with $d_0=D_Y$, the same cofactors of $P_w$ are nonzero. All four
weights are distinct, so again the unique apex is the $\xi$-axis by
Lemmas \ref{lem:apexbasic}(a) and \ref{lem:survive}(b), which do not
require positive weights, and \ref{E6} applies. It gives
$P_w=\sum_{j\ge2}f_jL^j+M_1$ with $f_j\in K[\xi]$, $M_1$ of passive degree
at most one, and a phase $L$ whose $u$-coefficient is normalized to one.
As in the first part of the proof, the passive degree $D_Y$ part
$f_{D_Y}L^{D_Y}=P_0$ fixes $L=u+\xi v+\xi^2t$ and $f_{D_Y}=a(x_0)$, and
$M_1=0$ since $P_w$ has no monomials of passive degree at most one. For
$2\le j<D_Y$ the coefficient $f_j$ is the coefficient of $u^j$ in
$P_w|_{v=t=0}$. A monomial $\xi^ku^j$ of $P_w$ has weight $j-\eta k=D_Y$,
that is,
\[
k=\frac{j-D_Y}{\eta}.
\]
For $j<D_Y$ this exponent is negative, which is impossible for a
monomial, so $f_j=0$. For $j=D_Y$ it is zero, and the monomial is
$a(x_0)u^{D_Y}$ of $P_0$. A positive exponent would require $j>D_Y$,
which is impossible since $D_Y$ is the largest passive degree of $\phi$.
Hence $P_w=a(x_0)L^{D_Y}=P_0$. No monomial joins $P_0$ at this
parameter, so it is not a genuine wall. Hence no wall occurs for
$0<\eta\le D_Y-2$. If a monomial of $\phi$ with $k=0$ and passive degree
$d$, $2\le d<D_Y$, had $j+2h>0$, it would reach the weight $D_Y$ at
$\eta_\sigma\le D_Y-2$ with $\sigma=j+2h$, and the first wall would occur
at a parameter at most $\eta_\sigma$, which is excluded. A monomial of
passive degree $D_Y$ with $k=0$ belongs to $N$ and satisfies
$0=k\ge j+2h$. So every monomial of $\phi$ with $k=0$ and passive
degree at least two has $j+2h=0$. Therefore the
restriction $\phi|_{\xi=0}$ has all passive-degree-$\ge2$ monomials
pure powers of $u$. Its passive Hessian has rank at most one, so
$\Hess\phi$ restricted to $\xi=0$ has rank at most three and
$\det\Hess\phi$ vanishes on the hyperplane $\xi=0$, contradicting
$\kappa\ne0$.
\end{proof}

\begin{lemma}[two-variable phase]
\label{lem:D3}
Suppose the leading form of $\phi$ is
\[
H=ax^{\ell}\,(w+x^{A}y)^m+bx^pz,\qquad ab\ne0,\ \ell\ge0,\
A,p\ge1,\ m\ge2,
\]
with integer exponents, $w_x=1$, the other weights exceeding one, and
$W>0$. Then $\phi$ has an affine direction.
\end{lemma}

\begin{proof}
Raise the weights of $y,w$ at unit speed and that of $z$ at speed $m$.
The leading weight has slope $m$ and $W$ slope $2m-4\ge0$, so $W$
stays positive. For a form $F(x,L)+B(x)z$ with $L=w+r(x)y$, the vector
\[
n=\big(0,\,-B',\,F_Lr',\,rB'\big)\T
\]
(in the coordinates $x,y,z,w$, primes denoting $\partial_x$) lies in the
Hessian kernel. The passive rows annihilate it because their passive
contributions involve only the phase direction $(r,0,1)$ and $n$ has
zero $x$-component, while in the active row the combination
$(\text{entry}_{xy})-r\,(\text{entry}_{xw})=F_Lr'$ and
$(\text{entry}_{xz})=B'$ cancel on $n$. The cofactor deleting $y$ is
$-F_{LL}(B')^2\ne0$, so the Hessian has rank three and
$\adj=-F_{LL}\,nn\T$ by comparison of the $yy$ entries. With $r=x^A$,
$B=bx^p$, $F=ax^\ell L^m$, all three passive cofactors are nonzero.
Both differences of the $z$-weight from a phase weight have slope
$m-1>0$, so they vanish at finitely many parameters only. We take the
old weights at a point of the segment just before the wall at which
the three passive weights are distinct.
The contraction criterion with old weights excludes every passive apex
at the next wall, even when new weights tie, and the active weight is
strictly smallest, so the apex space of the wall form is the active
axis.

Let $\gamma,\delta$ be the weights of $z,w$ at the wall and $E$ the
leading weight. The pure $w^m$ term normalizes the $w$-coefficient of
the phase. For a monomial in $x,y,w$ of passive degree $j$, the old
weight of a tied monomial is $E_{\mathrm{old}}+(m-j)\tau$, so the
restriction of the wall form to $z=0$ has no passive degree below $m$
and its degree-$m$ part is $ax^\ell(w+x^Ay)^m$. Comparing the pure
$w^m$ and $yw^{m-1}$ coefficients fixes the $y$-entry of the phase,
and the homogeneous ratio formula gives the phase
$w+x^Ay+kx^{\delta-\gamma}z$. Phase powers of degree $j\ge2$ have
coefficients proportional to $x^{E-j\delta}$ and only degrees
$j\ge m$ occur. The monomials linear in $z$ alone have old weight
$E_{\mathrm{old}}$, so the affine part in $z$ is exactly $bx^pz$.

If $k\ne0$ and $\delta-\gamma\notin\{0,A\}$, raise the three phase
weights at unit speed for a small step. Then the affine $z$-term drops, the
leading form becomes a single phase power with three distinct
coefficient exponents, and Lemma \ref{lem:D2} gives a contradiction.
If $\delta-\gamma\in\{0,A\}$, remove the $z$-entry by the constant
shear $w'=w+kz$ or $y'=y+kz$. Each preserves the old leading form,
because the weight of $z$ crosses the weight of the affected phase
variable from below (Lemma \ref{lem:shear}). Restart the path with the
highest phase power $M$ in place of $m$ and $z$-speed $M$. At every
wall either $M$ increases or one of the two shears occurs. The passive
degree of $\phi$ bounds $M$, and each difference of the $z$-weight
from a phase weight strictly increases on all segments, so there are
at most two shears and finitely many restarts.

After the last wall the path continues indefinitely, and every
monomial of $\phi$ then satisfies
$\deg_y+\deg_w+m\deg_z\le m$ (a violating monomial would overtake the
leading form). Hence $\deg_z\phi\le1$, and a monomial with
$\deg_z=1$ has $\deg_y=\deg_w=0$. In the accumulated constant
coordinates, $\phi=F(x,y,w)+zB(x)$ and $\partial_z^2\phi=0$, an
affine direction.
\end{proof}

\begin{lemma}[isolating the highest power of $y$]\label{lem:P1}
Suppose the leading form of $\phi$, for positive weights with $x$ of
minimum weight and $W>0$, is
\[
G(x,y)+c_1x^pz+c_2x^qw,\qquad p>q\ge1,\ c_1c_2\ne0,
\]
and the largest exponent $m$ of $y$ in $G$ is at least two. Then for a
modified path the leading form becomes
$c_0x^{\ell}y^m+c_1x^pz+c_2x^qw$ with $c_0\ne0$, $\ell\ge0$, $W>0$,
and all weights of $y,z,w$ strictly greater than that of $x$. After
normalizing $w_x=1$, Lemma \textup{\ref{lem:D4}} applies.
\end{lemma}

\begin{proof}
Increase the weight of $y$ at speed one and those of $z,w$ at speed
$m$, fixing the weight of $x$. The maximal slope among the monomials
of the form is $m$, attained exactly at the displayed monomials
($G$'s terms with $y$-exponent $m$ share one $x$-exponent $\ell$ by
homogeneity, and the two fiber terms have slope $m$ by construction).
All other monomials of $\phi$ have strict weight gaps, so a
sufficiently small positive step selects the displayed form, keeps
$W>0$ by continuity, and makes the other weights strictly greater than
$w_x$. Divide the weights by $w_x$.
\end{proof}

\begin{lemma}[a monomial phase and two linear terms]\label{lem:D4}
Suppose the leading form of $\phi$ is
$ax^{\ell}y^m+bx^pz+cx^qw$ with nonzero coefficients and integer
exponents $\ell\ge0$, $m\ge2$, $p>q\ge1$, $w_x=1$, the other weights
exceeding one, and $W>0$. Then $\phi$ has an affine direction.
\end{lemma}

\begin{proof}
Let $\beta$ be the weight of $y$ and $E=\ell+m\beta$. The fiber
weights are $E-p,E-q$. Raise $\beta$ at unit speed and both fiber
weights at speed $m$. Then $W=2(p+q-1-\beta)$, which vanishes at
$\beta=p+q-1$. If the path reaches that value with no wall, all
weights are positive there and Corollary \ref{cor:endpoint} gives a
direction. For a wall at distance $\tau$, a tied monomial
$x^iy^jz^kw^h$ has old weight $E_{\mathrm{old}}+[m-j-m(k+h)]\tau$.

At a wall before the endpoint, the principal cofactors deleting $z$
and $w$, equal for the old form to $-H_{yy}(cqx^{q-1})^2$ and
$-H_{yy}(bpx^{p-1})^2$ by expansion along the $w$ row (respectively the $z$ row),
survive. The wall form has rank three and a nonzero apex space with
zero intersection with the fiber plane, whose weights are distinct. A
line apex other than the active axis would have $y$ in its weight
block together with at most one fiber axis. The surviving minor on $x$
and the other fiber contradicts \ref{E6}.

If the apex space is a plane, it is a graph over $(x,y)$, and only $y$
can mix with a fiber at a weight equality. After the graph shear the
wall form is $R(x,y)+z\tilde G(x,y)+w\tilde K(x,y)$ with
$\det D(\tilde G,\tilde K)=0$ (Lemmas \ref{lem:doubling-det} and
\ref{lem:W}) and coefficient degrees $p\ne q$. The restrictions to
$y=0$ are $bx^p$ and $cx^q$, because the corresponding monomials have
old weight $E_{\mathrm{old}}$ and the shear preserves the old leading
form. If the common polynomial of Lemma \ref{lem:euler} is genuinely
binary, Lemma \ref{lem:D1} gives a direction. Otherwise it is a power
of $x$ (distinct base weights), both coefficients are exactly
$bx^p,cx^q$, the old-weight formula with $k=h=0$ excludes $y$-degrees
below $m$ and fixes the degree-$m$ term, and we restart the path
retaining the largest $y$-degree $M$ and both fiber terms, with slopes
$1,M,M$ for $y,z,w$. If $M=m$, a genuine wall requires a nontrivial
graph shear.

If the apex space is the active axis, \ref{E6} gives the phase
$y+dx^{\beta-(E-p)}z+ex^{\beta-(E-q)}w$ at the wall, with the
$y$-coefficient normalized by the pure $y^m$ term. Setting $z=w=0$
and comparing old weights excludes all phase powers below $m$,
affine terms in $y$, and active-only terms. The monomials linear in a
single fiber have old weight $E_{\mathrm{old}}$, retaining exactly the
partners $bx^pz,cx^qw$. If a fiber occurs in the phase with exponent
zero, absorb it into $y$ by a constant shear preserving the old form
(the fiber weight crosses $\beta$ from below). If three phase entries
remain, their exponents are distinct, and a small common raise of the
phase weights selects the highest power, and Lemma \ref{lem:D2} gives a
contradiction. With two entries, raise the two phase weights at unit
speed and the absent variable at the speed of the highest power. This
selects the highest power and the absent variable's affine partner,
and Lemma \ref{lem:D3} applies after naming the phase variables by
increasing weight and absorbing a constant. With one entry, restart
with slopes $1,M,M$.

All restarts use the finite support and a small positive step. Both
fiber weights minus $\beta$ strictly increase on every segment, so at
most two shears occur, and every other continuation increases the
phase degree, bounded by $\deg\phi$. The finite sequence of segments
therefore exits through Lemma \ref{lem:D1}, \ref{lem:D2} or
\ref{lem:D3}, or reaches $W=0$ with positive weights, where
Corollary \ref{cor:endpoint} finishes.
\end{proof}

Lemma~\ref{lem:D5} treats a leading form of rank three whose apex space is a line in a block of smallest weight. According to the number of phase variables it passes to Lemma~\ref{lem:D2}, \ref{lem:D3} or~\ref{lem:D4}.

\begin{lemma}[one apex of minimal weight]
\label{lem:D5}
Suppose a leading form of $\phi$ has Hessian rank three, $W>0$,
positive weights, and a one-dimensional apex space lying in a
smallest-weight block. Then $\phi$ has an affine direction.
\end{lemma}

\begin{proof}
A constant change within the block makes the apex the active
coordinate $x$. Scale $w_x=1$. By \ref{E6} the complementary
rank-one form has a phase linear in the passive variables over
$K(x)$. Normalize one nonzero coefficient to one. The other
coefficients are constants times Laurent monomials in $x$ (phase
formalism), and constant changes within the passive weight blocks
leave at most one phase entry per block. A passive variable absent
from the phase has a nonconstant affine coefficient, since a constant
coefficient would give a zero Hessian column up to the active one,
that is, a constant kernel vector, hence three apices by
Lemma \ref{lem:apexbasic}(c), contradicting uniqueness. Two absent
variables cannot share a weight, since their affine coefficients would
be proportional monomials of the same weight and a constant
combination would again give a kernel vector.

Let $m\ge2$ be the phase degree ($F_{LL}\ne0$ by rank three). Raise
the phase-variable weights at speed one and each absent variable at
the speed equal to the highest phase power. After a small step the
leading form is that highest power together with the affine terms of
the absent variables, and all passive weights exceed one. With three
phase variables, Lemma \ref{lem:D2} gives a contradiction. With two,
normalize the coefficient of the larger-weight one. The other
coefficient is a nonzero constant times a positive integral power of
$x$, absorb the constant and apply Lemma \ref{lem:D3} (the
coefficient of the pure power of the greatest-weight phase variable is
a polynomial in $x$, providing the exponent $\ell\ge0$ there). With
one phase variable, the other two passive variables have nonconstant
affine coefficients with distinct positive integral exponents, and
Lemma \ref{lem:D4} applies.
\end{proof}

\begin{lemma}[three linear terms]
\label{lem:D6}
Suppose the leading form of $\phi$ is
$\sum_{i=1}^3b_ix^{p_i}Y_i$ with $b_i\ne0$, integers
$p_1>p_2>p_3\ge1$, $w_x=1$, passive weights greater than one, and
$W>0$. Then $\phi$ has an affine direction.
\end{lemma}

\begin{proof}
Raise the three passive weights at unit speed. The leading weight has
slope one and $W$ slope $-2$. Before $W=0$ the old minors on
$(x,Y_i)$, equal to $-(b_ip_ix^{p_i-1})^2$, survive, and the passive
weights $E-p_i$ stay distinct. A constant kernel vector would contain
a weight axis (Lemma \ref{lem:apexbasic}(a)), and no axis $e$ is a
kernel vector, since for each axis some surviving mixed entry
$e_x\T\Hess\,e_{Y_i}$, equal at leading order to
$b_ip_ix^{p_i-1}\ne0$, pairs $e$ nontrivially. Hence there is no
constant kernel.
At a rank-two wall, \ref{E7} gives a single apex. It must be $x$,
since a passive apex would make the form jointly affine in a
complement containing $x$ and another passive variable, killing a
surviving mixed minor. The affine normal form then fixes all three
partners by weight, and an active-only term would exceed the old
leading weight. Therefore the wall form equals the old one, not a genuine wall.
At a rank-three wall the apex space has dimension one or two. A single
apex is again $x$, and Lemma \ref{lem:D5} applies. A two-apex plane
contains $x$, since a plane spanned by two passive axes would make the
form jointly affine in $x$ and the third passive variable,
contradicting the surviving minor. Aligned with $(x,Y_i)$, the form is
$G(x,Y_i)+Y_jA+Y_kB$ with coefficient degrees two of the $p$'s. If the
common polynomial is genuinely binary, Lemma \ref{lem:D1} finishes.
Otherwise it is a polynomial in one linear form $\ell$ in $x,Y_i$. Since the weights of $x$ and $Y_i$ are distinct, weighted homogeneity makes $\ell$ proportional to $x$ or to $Y_i$, and the old nonzero $x$-powers in the restrictions of the two coefficients to $Y_i=0$ exclude $\ell=Y_i$. So it is a power of $x$, both coefficients are unchanged, and a
genuinely new term must be nonlinear in $Y_i$ (active-only terms would
exceed the old degree, passive-linear terms are fixed by the old
form). Retaining its highest power and the two affine partners by the
paired raise gives Lemma \ref{lem:D4}. If the first wall occurs at
$W=0$, or no wall occurs by then, Corollary \ref{cor:endpoint}
finishes.
\end{proof}

\begin{lemma}[a pure power of a variable of smaller weight]\label{lem:D7}
Suppose the leading form of $\phi$ is $ay^M+bx^pz+cx^qw$ with
$abc\ne0$, integers $0<q<p<M$, and weights
$1,\beta,M\beta-p,M\beta-q$ on $x,y,z,w$ with $p/M\le\beta<1$. Then
$\phi$ has an affine direction.
\end{lemma}

\begin{proof}
Increase $\beta$ at unit speed, fixing $w_x=1$. The three displayed
terms stay tied and $W=2(p+q-1-\beta)$. At a wall before $W=0$ the
two fiber cofactors survive, by the computation of
Lemma \ref{lem:D4}. The wall form has rank three and its apex space
avoids the fiber plane and has dimension one or two. The weight of $w$
exceeds the other three at the wall, since $M\beta-q>M\beta-p$,
$(M-1)\beta-q\ge p-q-p/M>0$, and $M\beta-q\ge p-q\ge1$ grows past
one after any positive step. The fiber weights differ by $p-q$, so the
only possible weight equalities along the path are
$M\beta-p=\beta$, $M\beta-p=1$, and $\beta=1$, at
$\beta=p/(M-1),(p+1)/M,1$.

Suppose the apex space is a line. Its generator has zero
$w$-component ($w$ is alone in the heaviest block and its axis is
excluded by a cofactor) and nonzero $x$-component (otherwise the
annihilator contains the $x,w$ axes, whose old restricted minor
$-(qcx^{q-1})^2$ survives, contradicting \ref{E6}). Scaled to
$(1,\lambda,\mu,0)\T$, the source shear $x'=x+\lambda y+\mu z$ aligns
the apex with $x$. The coefficients $\lambda,\mu$ can be nonzero only
for variables of wall weight one, which had strictly smaller old
weight, so the alignment preserves the old leading form
(Lemma \ref{lem:shear}). Complementary rank one gives a phase in
$y,z,w$ with nonzero $y$-coefficient, of highest power exactly $M$ by
the retained $ay^M$ and the weight bound. Normalize the
$y$-coefficient. The coefficients of $y^{M-1}z$ and $y^{M-1}w$ force
the remaining phase entries to be monomials of exponents
$p-(M-1)\beta\le p/M<1$ and $q-(M-1)\beta<0$, or zero. Only a
constant $z$-entry can occur, at the first listed equality, and is
removed by a constant shear of $y$ by $z$ preserving the old form.
Lower powers $y^j$ are excluded by old weights, and the affine terms
in $z,w$ are fixed by the old comparison, so the wall form is the
original leading form. But its gradient image closure is the product
of $K^2$ with the closure of the curve $(bx^p,cx^q)$, since for fixed
$x\ne0$ the first output varies freely with the fibers and the second
freely with $y$. This product has at least two independent translation
directions, contradicting the one-dimensional apex space.

If the apex space is a plane, it is transverse to the fiber plane,
and the graph shear puts the form into
$R(x,y)+zG(x,y)+wK(x,y)$, preserving the old leading form since it
mixes only a fiber of smaller old weight into a base variable. The
coefficients have degrees $p,q$ and vanishing Jacobian
(Lemmas \ref{lem:doubling-det} and \ref{lem:W}), hence a common
weighted homogeneous polynomial, whose restrictions to $y=0$ are
$bx^p,cx^q$ by the old-weight comparison of fiber-linear,
$y$-independent monomials. If it is genuinely binary,
Lemma \ref{lem:D1} finishes. Otherwise it is a power of $x$ for
$\beta\ne1$, and at $\beta=1$ possibly of $x+\lambda y$, rectified by
one more old-weight-preserving shear. Both affine partners are
recovered, and the comparison of $y$-powers fixes $R=ay^M$. Hence
every continuing wall is a constant shear at one of the three listed
equalities, with the same leading form afterwards. There are finitely
many. If no earlier exit occurs, the path reaches $\beta=p+q-1$ with
positive weights ($M(p+q-1)-p\ge Mp-p>0$) and $W=0$, and
Corollary \ref{cor:endpoint} finishes.
\end{proof}

\subsection{Five tools}

The exit lemmas end every path whose leading form at the start has one of the shapes listed in Table~\ref{tab:exit}. The tools of this subsection and the lemmas of Section~\ref{sec:tree} show that the paths which start from the highest homogeneous part of $\phi$ reach one of these shapes.

The five lemmas of this subsection extract patterns which recur in
the case analysis of \cite[Proposition 4.1]{Tan26}, and they shorten
the case analysis of Section~\ref{sec:tree}.

\begin{lemma}[restriction principle]\label{lem:restrict}
Let $P$ be a wall form of rank three whose apex space is the single
axis of the variable $x_k$, and let $Y$ denote the other three
variables. Then the following hold.
\begin{enumerate}[label=\textup{(\alph*)}]
\item all $2\times2$ minors of $\Hess_YP$ vanish identically in
$K[x_1,\dots,x_4]$.
\item for every set $Z$ of variables, the matrix
$\Hess_{Y\setminus Z}\big(P|_{Z=0}\big)$ has rank at most one.
\item in particular, if $h=P|_{Z=0}$ and
$\Hess_{Y\setminus Z}h$ has rank at least two, the apex $e_k$ is
impossible.
\end{enumerate}
\end{lemma}

\begin{proof}
(a) is the rank-one statement of \ref{E6}, since the minors are
polynomials vanishing identically. They remain zero after the
substitution $Z=0$, and differentiation in the variables of
$Y\setminus Z$ commutes with that substitution, so
$\Hess_{Y\setminus Z}(P|_{Z=0})=\big(\Hess_{Y\setminus Z}P\big)|_{Z=0}$,
giving (b) and (c).
\end{proof}

\begin{lemma}[kernel step]\label{lem:kernelstep}
Suppose the leading form of $\phi$ is a form $h$ in $r\le3$ of the
variables with nonsingular Hessian in those variables, and the
remaining \emph{unused} variables share a common weight at least equal
to the maximal weight of the used ones. Raise the common unused weight.
Then one of the following occurs.
\begin{enumerate}[label=\textup{(\alph*)}]
\item $W=0$ is reached with positive weights, and
Corollary \textup{\ref{cor:endpoint}} gives a direction.
\item the first wall form has a nonzero constant Hessian kernel. Then,
after a constant linear change among the unused variables, it is a
form in $r'>r$ variables, the remaining unused variables again sharing the maximal weight, and either its Hessian in those $r'$ variables is nonsingular, or $r'=3$ and $\rk\Hess P=2$.
\item the first wall form $P$ has no constant Hessian kernel,
satisfies $P|_{\mathrm{unused}=0}=h$, and $\rk\Hess P\in\{2,3\}$.
\end{enumerate}
\end{lemma}

\begin{proof}
The unused variables form one weight block, strictly heavier than all
used ones after any positive step. The terms of $h$ keep their
weights, so $P|_{\mathrm{unused}=0}=h$ at the first wall, and $P$
involves the unused block. A constant kernel vector splits into
blocks (Lemma \ref{lem:apexbasic}(a)). A component $k$ in the used
block satisfies, after restriction, $\Hess h\cdot k=0$, hence
vanishes. So the kernel lies in the unused block, and it is proper
there because $P$ involves the block. In case (b), a linear change
within the unused block makes $P$ independent of
$\dim$-of-kernel many variables (Lemma \ref{lem:apexbasic}(d)). The
resulting form $P'$ in $r'=4-\dim\ker\le3$ variables, $r'>r$, has no further constant kernel vector. For $r'=2$ its Hessian in the two variables is nonsingular by Lemma~\ref{lem:apexbasic}(f). For $r'=3$ it is nonsingular, or $\rk\Hess P\le2$ by Lemma~\ref{lem:apexbasic}(g), and then $\rk\Hess P=2$, since the restriction $h$ has rank at least two when $r\ge2$, and since a form of rank one would have a three-dimensional constant kernel by Lemma~\ref{lem:apexbasic}(e) when $r=1$. In case (c), $\rk\Hess P\le3$ by Lemma \ref{lem:W}
and $W>0$. For $r\ge2$ the restriction $h$ gives $\rk\Hess P\ge2$. For
$r=1$ a wall form of rank one would have a three-dimensional constant
kernel by Lemma~\ref{lem:apexbasic}(e), which case (c) excludes, so
again $\rk\Hess P\ge2$.
\end{proof}

\begin{lemma}[unused variable over a ternary form]\label{lem:T}
Suppose the leading form of $\phi$ is $H(u_1,u_2,u_3)$, a form with
$\det\Hess_3H\ne0$, with positive weights, $W>0$, and the fourth
variable $v$ unused, of weight at least $\max_iw_{u_i}$. Raise only
$w_v$. Then either $W=0$ is reached with positive weights
\textup{(}Corollary \textup{\ref{cor:endpoint}}\textup{)}, or the
first wall form $Q$ satisfies the following.
\begin{enumerate}[label=\textup{(\alph*)}]
\item $Q|_{v=0}=H$ and $Q_v\ne0$.
\item $\rk\Hess Q=3$.
\item $\Hess Q$ has no constant kernel.
\item the apex space of $Q$ lies in $\spn(e_{u_1},e_{u_2},e_{u_3})$
and has dimension one or two.
\end{enumerate}
\end{lemma}

\begin{proof}
(a) Monomials without $v$ keep their weights, so those of $H$ remain
in the wall form and no other $v$-free monomial joins. The wall is
genuine, so $Q$ contains a monomial with $v$. (b) The principal
cofactor of $\Hess Q$ deleting the $v$ row and column restricts at
$v=0$ to $\det\Hess_3H\ne0$, so the rank is at least three, and $W>0$
makes $\det\Hess Q=0$. (c) After a positive step $v$ is its own
weight block. A kernel component in the ternary block vanishes after
restriction ($\det\Hess_3H\ne0$), and $e_v$ in the kernel would make
$Q_v$ constant, hence zero (no affine part), contradicting (a). (d)
The nonzero restriction of the $v$-cofactor shows
$e_v\T\adj(\Hess Q)e_v\ne0$, so $e_v$ is not an apex
(Lemma \ref{lem:apexbasic}(b)). The apex space splits and lies in the
ternary block, is nonzero by Theorem \ref{thm:apex}, and has
dimension at most two by (c) and Lemma \ref{lem:apexbasic}(c).
\end{proof}

\begin{lemma}[two apices]\label{lem:B}
Let a wall form of $\phi$, for positive weights with $W>0$, have
Hessian rank three and a two-dimensional apex space, and, after the
alignment of Lemma \textup{\ref{lem:shear}} and \ref{E6}, the shape
\[
P=G(x,y)+zA(x,y)+wB(x,y),
\]
with $x$ of minimal weight among the four variables. Assume
\begin{enumerate}[label=\textup{($*$)}]
\item if $w_x<w_y$, then $A|_{y=0}\ne0$ or $B|_{y=0}\ne0$.
\end{enumerate}
Then $\phi$ has an affine direction.
\end{lemma}

\begin{proof}
By Lemmas \ref{lem:doubling-det} and \ref{lem:W},
$\det D(A,B)=0$. The coefficients $A,B$ are nonconstant and not
proportional, since a constant one, or a constant combination, would give a
constant kernel vector and a three-dimensional apex space
(Lemma \ref{lem:apexbasic}(c)). By Lemma \ref{lem:euler}, $A=aS^r$,
$B=bS^s$. Their weighted degrees differ. Indeed, $S$ has positive
weighted degree, so equal weighted degrees of $A$ and $B$ would give
$r=s$, and $A,B$ would be proportional, which was just excluded. If $S$
is genuinely binary, Lemma \ref{lem:D1} gives a direction. Otherwise
$S$ is a polynomial in a constant linear form $\ell$. If $w_x=w_y$, a
constant change in the base block makes $\ell=x$, and if
$w_x<w_y$, $\ell$ is a coordinate and ($*$) excludes $y$. Hence
$A=ax^p$, $B=bx^q$ with $p>q\ge1$ after possibly exchanging $z,w$. If
$G_{yy}=0$, then $P$ is jointly affine in $y,z,w$, its Hessian
vanishes on the corresponding principal $3\times3$ block, and the rank
is at most two, contradicting rank three. So the largest $y$-exponent
of $G$ is at least two, and Lemma \ref{lem:P1} followed by
Lemma \ref{lem:D4} gives a direction.
\end{proof}

Lemma~\ref{lem:K} treats the trinomial leading forms which arise from wall forms of rank two in Lemmas~\ref{lem:tree-rk2} and~\ref{lem:tree-binary}. Its proof lists the possible apex lines and apex planes at the next wall, excludes some of them, and passes the others to Lemma~\ref{lem:D5}, \ref{lem:D6} or~\ref{lem:B}.

\begin{lemma}[trinomial faces]\label{lem:K}
Suppose the leading form of $\phi$ is
\[
F=as^D+bts^p+czs^q,\qquad abc\ne0,\quad D>p>q\ge1,
\]
with weights $w_s=1\le w_t=D-p<w_z=D-q$, the fourth variable $u$
unused of weight $w_u=w_z$, and $W>0$. Raise only $w_u$. Then $\phi$
has an affine direction.
\end{lemma}

\begin{proof}
If $W=0$ is reached, Corollary \ref{cor:endpoint} finishes. Let $Q$
be the first wall form, so $Q|_{u=0}=F$ and $Q_u\ne0$.

\emph{Identities.} The principal $2\times2$ minors of $\Hess F$ on
$(s,t)$ and $(s,z)$ are
\[
M_{st}=-\big(bp\,s^{p-1}\big)^2,\qquad
M_{sz}=-\big(cq\,s^{q-1}\big)^2 ,
\]
and they survive to $Q$ by restriction to $u=0$
(Lemma \ref{lem:survive}).

\emph{No kernel.} The only constant kernel direction of $F$, viewed in
four variables, is $e_u$, since the $s$-derivative of $F$ has nonzero $t$-
and $z$-terms, and the $t$- and $z$-derivatives are independent powers
of $s$. On the wall $Q_u\ne0$, so $Q$ has no constant kernel. Its
rank is two or three, and by Lemma \ref{lem:apexbasic}(c) and
\ref{E7} it has one or two apices.

\emph{One apex.} The axes $e_z,e_u$ are excluded, since their annihilators
contain $s,t$, and the surviving minor $M_{st}$ contradicts
\ref{E6} (rank three) or joint affineness (\ref{E7}, rank two).
If $w_t>1$, the weight blocks are axes and the minor $M_{sz}$
likewise excludes $e_t$, so the apex is $e_s$. If $w_t=1$, an apex
$d_se_s+d_te_t$ has $d_te_s-d_se_t$ and $e_z$ in its annihilator, and
the surviving minor of $Q$ on these two directions restricts to
$-\big(d_t\,cq\,s^{q-1}\big)^2$, so $d_t=0$ and again the apex is
$e_s$. At rank three, Lemma \ref{lem:D5} applies ($s$ has minimal
weight). At rank two, \ref{E7} and weighted homogeneity give the
normal form
\[
Q=as^D+bts^p+czs^q+du\,s^{\rho},\qquad d\ne0,\ 0<\rho<q ,
\]
and a small common increase of the weights of $t,z,u$ selects the
three linear partners, of slope one, while $as^D$ has slope zero.
Hence Lemma \ref{lem:D6} applies.

\emph{Two apices.} The apex plane contains $e_s$. Indeed, a plane
whose annihilator contains $e_s$ is excluded, since by \ref{E6} the
form $Q$ would then be affine in $s$, while $Q_{ss}$ contains
$aD(D-1)s^{D-2}\ne0$. For $w_t>1$ the weight blocks are axes, and this
excludes every plane spanned by two of $e_t,e_z,e_u$. The plane
spanned by $e_z,e_u$ is also excluded by $M_{st}$. A plane spanned by
$d_se_s+d_te_t$ and $e_u$ reduces to the restricted-minor computation
above, forcing $d_t=0$. And for a plane spanned by $d_se_s+d_te_t$
and $e_z$, joint affineness in the complementary plane makes the
second derivative of $Q$ along $d_te_s-d_se_t$ vanish, while its
$t$-coefficient restricts to $d_t^2\,bp(p-1)s^{p-2}$ with
$p\ge2$ when $w_t=1$ (and the case $w_t>1$ has axes for blocks), so
$d_t=0$. The remaining planes are $(s,t)$, $(s,z)$, $(s,u)$. For
$(s,t)$ one has $Q=G(s,t)+czs^q+uB(s,t)$ with $G=as^D+bts^p$ fixed by the
old form. The vanishing Jacobian of the coefficient pair gives
$cqs^{q-1}B_t=0$, so $B=B(s)$ and $Q$ is jointly affine in $t,z,u$, of
rank at most two, so this plane is impossible at rank three. For
$(s,z)$ one gets symmetrically $bps^{p-1}B_z=0$ with the same conclusion. For
$(s,u)$ the $t,z$ coefficients restrict at $u=0$ to $bs^p,cs^q$. A
common polynomial in a constant linear form is a power of $s$
(condition ($*$) of Lemma \ref{lem:B} holds by these restrictions),
and Lemma \ref{lem:B} gives a direction, its genuinely binary branch
included.
\end{proof}

\section{The highest homogeneous part}\label{sec:tree}

As in step~(1) of the outline in Section~\ref{sec:translation}, the rank $r$ of $\Hess\phi_D$ is $4-\dim T_D$, where $T_D=\{v:D_v\phi_D=0\}$ is the space of translation symmetries of the top form, and in each case below the weights of the variables spanning $T_D$ are raised. In this section all paths start from the weights $(1,1,1,1)$ and stop at the
first wall or at $W=0$, whichever comes first. In the latter case
Corollary \ref{cor:endpoint} finishes the branch, so only walls are
discussed. At every wall form, a nonzero constant vector of the
Hessian kernel gives, because affine parts have been removed, a zero
directional derivative, and rank-three forms have nonzero apex spaces
(Theorem \ref{thm:apex}) splitting into weight blocks
(Lemma \ref{lem:apexbasic}).

The cases of the proof of Proposition~\ref{prop:tree} are the following. Let $r=\rk\Hess\phi_D$ be the Hessian rank of the top homogeneous part $\phi_D$. In the items below, the rank of a wall form $P$ always means $\rk\Hess P$, which in general differs from $r$. Each item names the steps which end the branch, by an affine direction or by a contradiction.
\begin{enumerate}[label=\textup{(\arabic*)}]
\item $r=3$. Raise the weight of the fourth variable (Lemma~\ref{lem:T}). A single apex is treated by Lemma~\ref{lem:D5}, two apices by Lemma~\ref{lem:B}.
\item $r=2$, Lemma~\ref{lem:tree-rk2}. Raise the two unused weights (Lemma~\ref{lem:kernelstep}). At the wall the following occur.
\begin{enumerate}[label=\textup{(\alph*)}]
\item A ternary form with nonsingular Hessian. Lemma~\ref{lem:T}, then Lemma~\ref{lem:D5} for one apex, and Lemmas~\ref{lem:D1}, \ref{lem:P1} with~\ref{lem:D4}, or~\ref{lem:B} for two apices.
\item A wall form $P$ with $\rk\Hess P=2$. The normal form \eqref{eq:form1} and Lemma~\ref{lem:K}.
\item A wall form $P$ with $\rk\Hess P=3$. Lemma~\ref{lem:D5} for one apex and Lemma~\ref{lem:B} for two apices. Three apices do not occur, since the wall form has no constant kernel (Lemma~\ref{lem:kernelstep}(c)).
\end{enumerate}
\item $r=1$, $\phi_D=a_0s^D$. Raise the three other weights (Lemma~\ref{lem:kernelstep}). A wall form $P$ with $\rk\Hess P=1$ does not occur. A wall form $P$ with $\rk\Hess P=2$ is a binary form with two unused variables, treated by Lemma~\ref{lem:tree-binary}. A wall form $P$ with $\rk\Hess P=3$ is treated by Lemma~\ref{lem:tree-C1}.
\item Lemma~\ref{lem:tree-binary}. Raise the two unused weights (Lemma~\ref{lem:kernelstep}). A wall form $P$ with $\rk\Hess P=2$ leads to a trinomial face and Lemma~\ref{lem:K}. A ternary form is treated by Lemma~\ref{lem:T} and Lemma~\ref{lem:tree-C1}, a wall form $P$ with $\rk\Hess P=3$ by Lemma~\ref{lem:tree-C1}.
\item Lemma~\ref{lem:tree-C1}. A single apex of minimal weight is treated by Lemma~\ref{lem:D5}, and a single heavier apex does not occur. For two apices the common polynomial of Lemma~\ref{lem:euler} gives Lemma~\ref{lem:D1}, Lemmas~\ref{lem:P1} and~\ref{lem:D4}, or Lemma~\ref{lem:D7}. Three apices occur only in item (3), where the three heavy variables have a common weight, and they lead by Lemma~\ref{lem:T} to a wall form with one or two apices.
\end{enumerate}
A path which reaches $W=0$ before a wall ends by Corollary~\ref{cor:endpoint}. Lemmas~\ref{lem:D2}, \ref{lem:D3} and~\ref{lem:D6} are used inside the proofs of Lemmas~\ref{lem:D3}, \ref{lem:D4}, \ref{lem:D5} and~\ref{lem:K}.

\emph{Local and global parts of the argument.} The local part is the
analysis of a single wall. It determines the wall form, its Hessian
rank, its apex space and the surviving minors, and it is carried out
in the proofs of the lemmas of Section~\ref{sec:walls} and of this
section. The global part is the termination of the procedure, which is
the content of the following lemma.

\begin{lemma}[termination]\label{lem:termination}
Every branch of the procedure of Proposition~\textup{\ref{prop:tree}} calls finitely many of the lemmas of Sections~\textup{\ref{sec:walls}} and~\textup{\ref{sec:tree}}, passes through finitely many paths of weights, and meets finitely many walls. It ends at a weight vector with positive entries and $W=0$, where Corollary~\textup{\ref{cor:endpoint}} gives an affine direction, or in an affine direction, or in a contradiction.
\end{lemma}

\begin{proof}
By the local analysis each lemma of the case list ends in one of the three outcomes of the statement or passes to another lemma of the list. It remains to bound the number of steps. We use the following facts, which concern only the bookkeeping of the procedure and can be checked independently of the local analysis.
\begin{enumerate}[label=\textup{(F\arabic*)}]
\item The lemmas call each other in a fixed order without cycles.
Lemma~\ref{lem:D3} calls Lemma~\ref{lem:D2}. Lemma~\ref{lem:D4} calls
Lemmas \ref{lem:D1}, \ref{lem:D2}, \ref{lem:D3}. Lemma~\ref{lem:D5}
calls Lemmas \ref{lem:D2}, \ref{lem:D3}, \ref{lem:D4}.
Lemma~\ref{lem:D6} calls Lemmas \ref{lem:D1}, \ref{lem:D4},
\ref{lem:D5}. Lemma~\ref{lem:D7} calls Lemma~\ref{lem:D1}.
Lemma~\ref{lem:B} calls Lemmas \ref{lem:D1}, \ref{lem:P1},
\ref{lem:D4}. Lemma~\ref{lem:K} calls Lemmas \ref{lem:D5},
\ref{lem:D6}, \ref{lem:B}. Lemmas \ref{lem:D1} and \ref{lem:D2} call
no other lemma. Proposition~\ref{prop:tree} calls Lemmas \ref{lem:T},
\ref{lem:D5}, \ref{lem:B}, \ref{lem:tree-rk2}, \ref{lem:kernelstep},
\ref{lem:tree-binary}, \ref{lem:tree-C1}. Lemma~\ref{lem:tree-rk2} calls
Lemmas \ref{lem:kernelstep}, \ref{lem:T}, \ref{lem:D1}, \ref{lem:D5},
\ref{lem:P1}, \ref{lem:D4}, \ref{lem:B}, \ref{lem:K}.
Lemma~\ref{lem:tree-binary} calls Lemmas \ref{lem:kernelstep},
\ref{lem:T}, \ref{lem:tree-C1}, \ref{lem:K}. Lemma~\ref{lem:tree-C1}
calls Lemmas \ref{lem:D1}, \ref{lem:D5}, \ref{lem:P1}, \ref{lem:D4},
\ref{lem:D7}, \ref{lem:T}. The only repetition is the call of
Lemma~\ref{lem:tree-C1} through Lemma~\ref{lem:T} in item (5). It
occurs at most once, since the second wall form has no constant kernel
by Lemma~\ref{lem:T}(c).
\item Inside a lemma every restart either raises the phase degree $M$
or comes with one of the shears counted in (F3). At a wall the
monomials of the old leading form remain in the wall form, so the
phase degree after the restart is at least $M$. Every phase degree is
the passive degree of a monomial of $\phi$ in the current coordinates,
hence at most the total degree $D=\deg\phi$, which invertible linear
changes of coordinates do not alter. So there are at most $D-M_0$
restarts which raise the phase degree, where $M_0$ is the phase degree
at the start of the lemma. In Lemma~\ref{lem:D2} every genuine wall
raises the phase degree. In Lemmas~\ref{lem:D3} and~\ref{lem:D4} a
restart which does not raise it requires one of the at most two shears
of (F3), and in Lemma~\ref{lem:D7} every continuing wall is a shear at
one of the at most three weight equalities of (F3).
\item Shears occur only at equalities of the weights of two
variables whose difference is a nonconstant affine function of the
path parameter. Each such equality holds at one parameter at most, and
shears do not change the weights. So there are at most four shears in
Lemma~\ref{lem:D1} (Lemma~\ref{lem:D1-walls}), at most two in Lemmas
\ref{lem:D3} and \ref{lem:D4}, and at most three parameters with shears
in Lemma~\ref{lem:D7}.
\item Along every segment $W$ is an affine function of the parameter.
Where it decreases, the path stops at $W=0$ at the latest, and
Corollary~\ref{cor:endpoint} applies, since all weights are positive
there. The path with the negative weight $-\eta$ in the proof of
Lemma~\ref{lem:D2} has $W\ge2$ throughout and never reaches $W=0$. Where $W$ is constant or increasing, $W>0$ is preserved, and the
number of walls is bounded by (F2) and (F3).
\item Each kernel step of Lemma~\ref{lem:kernelstep} in case (b)
increases the number of variables of the leading form. In items (3)
and (4) this number passes from one to two and from two to three, so
there are at most two kernel steps on a branch.
\end{enumerate}
By (F1) the calls between the lemmas form a finite tree, and by (F5) a branch contains at most two kernel steps. Inside a single lemma the number of restarts is bounded by (F2) and the number of shears by (F3), and by (F4) every path of weights stops at $W=0$ or meets finitely many walls. Hence every branch is finite. The quantities which control the termination are thus the position of a lemma in the order of calls (F1), the phase degree (F2), the equalities of weights available for shears (F3), the slope of $W$ along a segment (F4), and the number of variables of the leading form (F5).
\end{proof}

\begin{lemma}[top form of rank two]\label{lem:tree-rk2}
Suppose that in linear coordinates $(s,t,z,u)$ the top homogeneous part of $\phi$ is $h(s,t)$ with $\det\Hess_2h\ne0$. Then $\phi$ has an affine direction.
\end{lemma}

\begin{proof}
Raise the weights of $z,u$ together. This is the situation of
Lemma \ref{lem:kernelstep} with $r=2$. In its case (b) we reach,
after a constant change among $z,u$, a ternary form $H(s,t,z)$ with
$H(s,t,0)=h$, the fourth variable
unused at the maximal weight. If its ternary Hessian is nonsingular, it is treated under ``ternary forms''
below. Otherwise the wall form has rank two, with the one-dimensional constant kernel spanned by the unused axis, and the paragraph ``rank two at the wall'' below applies to it verbatim. Indeed, that paragraph uses \ref{E7}, whose first alternative, a two-dimensional constant kernel, is excluded here, and the restriction to $z=u=0$, which is $h$. In its case (c) the
wall form $P$ satisfies $P(s,t,0,0)=h$ and has rank two or three and
no constant kernel.

\emph{Rank two at the wall.} By \ref{E7} there is a single apex. It
lies in the $(s,t)$-block, since an apex on a $z$- or $u$-axis would
make $P$ jointly affine in a complement containing $s,t$, whence
$\Hess_2h=0$. Align the apex with $s$ by a constant change in the
block. By \ref{E7}, $P$ is jointly affine in $t,z,u$. Weighted
homogeneity and $P|_{z=u=0}=h$ give $h=as^D+bts^{D-1}$, $b\ne0$, and
the $z,u$ coefficients, of equal weighted degree, become after a
constant change in the $(z,u)$-block
\begin{equation}\label{eq:form1}
P=as^D+bts^{D-1}+czs^{D-\beta},\qquad bc\ne0,
\end{equation}
with $D-\beta$ a positive integer and $u$ unused. If $a=0$, substitute $t+\mu s$ for $t$ with $\mu\in\kk^\times$. This constant change within the $(s,t)$-block, whose two weights are equal, keeps $P$ jointly affine in $t,z,u$ and of the form \eqref{eq:form1}, with $a$ replaced by $b\mu\neq0$. So we may assume $a\neq0$. This is the situation of
Lemma \ref{lem:K} with $p=D-1$, $q=D-\beta$, up to raising $w_u$ to
$w_z$ first (their common weight at the wall is already equal). That
lemma finishes.

\emph{Rank three at the wall.} A single apex lies in the
$(s,t)$-block, since a $z$- or $u$-axis apex is excluded by the
restriction principle \ref{lem:restrict} with $Z=\{z,u\}$ and
$P|_{z=u=0}=h$, and it has minimal weight, so Lemma \ref{lem:D5}
finishes. Two apices cannot span the $(z,u)$-block (then $P$ would be
jointly affine in $s,t$, making $h$ affine) nor the $(s,t)$-block.
In the latter case $P=G(s,t)+zA(s,t)+uB(s,t)$, and $\det D(A,B)=0$ by
Lemmas \ref{lem:doubling-det} and \ref{lem:W}. The coefficients $A,B$
are nonzero, since otherwise $e_z$ or $e_u$ is a constant kernel
vector of $\Hess P$, and nonconstant, since $P$ has no affine part.
As $w_z=w_u$, they have the same weighted degree, which is positive.
By Lemma \ref{lem:euler}, applied with the base variables $s,t$ of
weight one, $A=aS_0^{r_1}$ and $B=bS_0^{r_2}$ with a common weighted
homogeneous polynomial $S_0$ of positive weighted degree, and equal
weighted degrees force $r_1=r_2$. So $A,B$ are proportional, which gives a
constant kernel vector and, by Lemma \ref{lem:apexbasic}(c), three
apices. So the apex plane splits as one direction in each block.
Constant changes align it with $(s,z)$, and
$P=G(s,z)+tA(s,z)+uB(s,z)$ with $A(s,0)$ a nonzero multiple of
$s^{D-1}$ (restrict to $z=u=0$) and coefficient degrees
$D-1>D-\beta$. Condition ($*$) of Lemma \ref{lem:B} holds by the
restrictions $A(s,0)\ne0$, and Lemma \ref{lem:B} finishes. Three
apices do not occur. By Lemma \ref{lem:apexbasic}(c) they would give a
nonzero constant kernel vector, which case (c) of
Lemma \ref{lem:kernelstep} excludes.

\emph{Ternary forms.} Whenever the leading form is a ternary
$H(u_1,u_2,u_3)$ with nonsingular ternary Hessian, positive weights,
$W>0$, and the unused variable of maximal weight, raise only the
unused weight. Then Lemma \ref{lem:T} applies, and the wall form $Q$ has
rank three, no kernel, $Q|_{v=0}=H$, and one or two apices in the
ternary block. Here
$H(s,t,z)$ has $w_s=w_t=1$ minimal and $H|_{z=u=0}=h$
nondegenerate binary. A single $z$-apex is excluded by the
restriction principle with $Z=\{z,u\}$ (a $\mathrm{GL}_2$ change in
the $(z,u)$-block does not affect the restriction $z=u=0$), so a
single apex lies in $(s,t)$ and Lemma \ref{lem:D5} finishes. If there are two
apices, a plane in the $(s,t)$-block gives, after the alignment and
Lemma \ref{lem:euler}, either Lemma \ref{lem:D1} or coefficients that
are powers of $s$ after a constant change of the equal-weight base.
In the latter case the pure base term $G$ satisfies
$G_{tt}\ne0$, because $\det\Hess_3H=-\big(\partial_s(\text{$z$-%
coefficient})\big)^2G_{tt}\ne0$, and Lemma \ref{lem:P1} with
Lemma \ref{lem:D4} finishes. This is Lemma \ref{lem:B} with
($*$) vacuous. A mixed plane aligns with $(s,z)$ and
$Q=G(s,z)+tA(s,z)+uB(s,z)$. The restriction $Q|_{z=u=0}=h$ forces
$h=G(s,0)+tA(s,0)$ with $A(s,0)$ a nonzero multiple of $s^{D-1}$, the
coefficient degrees $D-1>D-\gamma$ are distinct, ($*$) holds, and
Lemma \ref{lem:B} finishes (within its proof, if the common
polynomial is a power of $s$, the pure base term has
$G_{zz}\ne0$ by $\det\Hess_3H\ne0$, expanded along the $t$-row).

\emph{The trinomial continuation of \eqref{eq:form1}} is
Lemma \ref{lem:K}.
\end{proof}

\begin{lemma}[rank-three wall forms with a pure power]\label{lem:tree-C1}
Let $P$ be a wall form of $\phi$ of rank three and weight $D$ which contains $a_0s^D$, $a_0\ne0$, for weights with $w_s=1$, the other weights greater than one, and $W>0$. If $\Hess P$ has a nonzero constant kernel vector, assume moreover that the three variables other than $s$ have a common weight. Then $\phi$ has an affine direction.
\end{lemma}

\begin{proof}
A single apex of minimal weight is Lemma \ref{lem:D5}. Suppose the
unique apex is a heavier coordinate $t$ (after a constant change in
its block). The complementary form of \ref{E6} over $K(t)$ has a
phase with nonzero $s$-entry, since $P$ contains $a_0s^D$. Normalize
it to $s+a(t)u+b(t)v$. Setting $u=v=0$ shows the highest phase power
is $D$ with coefficient $a_0$. The coefficients of $s^{D-1}u$ and
$s^{D-1}v$ are $Da_0a(t)$ and $Da_0b(t)$, polynomials of weights
$1-w_u<0$ and $1-w_v<0$, hence zero. Then
$P=G(s,t)+uA(t)+vB(t)$, and the gradient-image closure lies in the
product of the $(s,t)$-target plane with the closure of the curve
$(A,B)$. Its dimension is $\rk\Hess P=3$, so it equals that product, which has at
least two independent translation directions, contradicting the
one-dimensional apex space. So heavy single apices do not occur.

A two-apex plane contains the $s$-axis, since by weight invariance it would
otherwise lie in the heavy block, and its complementary space would
contain $s$, forcing $P_{ss}=0$ against $a_0s^D$. Align the plane
with $(s,t)$ and write $P=G(s,t)+uA(s,t)+vB(s,t)$. By
Lemmas \ref{lem:doubling-det} and \ref{lem:W}, $\det D(A,B)=0$. As in
the proof of Lemma \ref{lem:B}, the coefficients $A,B$ are nonconstant
and not proportional, since otherwise $\Hess P$ has a constant kernel
vector and the apex space is three-dimensional
(Lemma \ref{lem:apexbasic}(c)). They have different weighted degrees,
since by Lemma \ref{lem:euler} they are $aS_0^{r_1}$ and $bS_0^{r_2}$
for a common polynomial $S_0$ of positive weighted degree, and equal
weighted degrees would force $r_1=r_2$ and proportionality. In particular
$w_u\ne w_v$, as Lemma \ref{lem:D1} requires. If the common
polynomial of Lemma \ref{lem:euler} is genuinely binary,
Lemma \ref{lem:D1} finishes. Otherwise the unequal weights of $s,t$
make it a power of $s$ or of $t$. For a power of $s$ one has
$G_{tt}\ne0$ (else $P$ is jointly affine in $t,u,v$ of rank at most
two), and Lemma \ref{lem:P1} with Lemma \ref{lem:D4} finishes. For a
power of $t$, order the fibers so that the coefficients are
$c_1t^p,c_2t^q$, $p>q\ge1$. Increase $w_s$ at speed one and
$w_u,w_v$ at speed $D$, fixing $w_t$. In $G$ only $a_0s^D$ has
$s$-exponent $D$ and none has a higher one, so the slope-$D$ terms
are exactly $a_0s^D+c_1t^pu+c_2t^qv$. A small step selects them, and
dividing all weights by $w_t$ gives a new $s$-weight
$\beta<1$ with fiber weights $D\beta-p,D\beta-q>0$, whence
$0<q<p<D$ and $p/D<\beta<1$, and Lemma \ref{lem:D7} finishes, with $t$ as
its weight-one variable. Three apices give a nonzero constant kernel vector by
Lemma~\ref{lem:apexbasic}(c). By hypothesis the three variables other
than $s$ then have a common weight, so they form one weight block,
which is the heaviest one. The kernel splits into blocks
(Lemma~\ref{lem:apexbasic}(a)), and $a_0s^D$ excludes an
$s$-component, so the kernel lies in that block. It is a line, since
$\rk\Hess P=3$. A change in that block yields a ternary form
$H(s,t,u)$ with nonsingular ternary Hessian, and the unused $v$ has the
common weight of the block, which is maximal. Lemma \ref{lem:T}
applies. Its wall form $Q$ satisfies $Q|_{v=0}=H$, so it contains
$a_0s^D$, it has rank three and no constant kernel, and it has one or
two apices. So $Q$ satisfies the hypotheses of this lemma, the
additional one being void, and the cases treated above apply to it.
\end{proof}

\begin{lemma}[binary forms with two unused variables]\label{lem:tree-binary}
Suppose that the leading form of $\phi$ is a form $H(s,t)$ with $\det\Hess_2H\ne0$ which contains $a_0s^D$, $a_0\ne0$, for weights with $w_s=1<w_t$ and $W>0$, and that the two other variables $u,v$ are unused, of equal weight at least $w_t$. Then $\phi$ has an affine direction.
\end{lemma}

\begin{proof}
Raise $w_u=w_v$ together. This is the kernel step
(Lemma \ref{lem:kernelstep}, $r=2$). Case (b) yields a ternary form
with unused fourth variable of maximal weight. If its ternary Hessian is nonsingular, Lemma \ref{lem:T} and
Lemma~\ref{lem:tree-C1} apply. The wall form of Lemma~\ref{lem:T} has no constant kernel by Lemma~\ref{lem:T}(c), so the additional hypothesis of Lemma~\ref{lem:tree-C1} is void. Otherwise the wall form has rank two and a one-dimensional constant kernel, and it is treated as the rank-two wall forms of case (c) below. In case (c) the wall form has no kernel. Rank three
is Lemma~\ref{lem:tree-C1}, whose additional hypothesis is again void. Rank two has a single apex, equal to
$s$ by $a_0s^D$ (the binary forms here contain $a_0s^D$) and
\ref{E7}, the $u,v$-coefficients are proportional (equal weight),
and a constant change gives the trinomial face
\[
a_0s^D+bts^p+cus^q,\qquad bc\ne0,\ p>q\ge1,
\]
with $p=D-w_t$, $q=D-w_u$, $q>0$ (no affine part) and the unused $v$
of weight $w_v=w_u$, and Lemma \ref{lem:K} finishes, after renaming
$(z,u)\to(u,v)$.
\end{proof}

\begin{proposition}\label{prop:tree}
Every Keller potential $\phi\in\kk[x_1,\dots,x_4]$ of degree $D\ge3$ has an affine direction.
\end{proposition}

\begin{proof}
Remove the affine part. Write $\phi_D$ for the top homogeneous part.
The degree-$4(D-2)$ component of $\det\Hess\phi$ is
$\det\Hess\phi_D$, which therefore vanishes, while some second
derivative of $\phi_D$ is nonzero. Iterating the vertex theorem \ref{E:GN}, by which
a form in at most four variables with vanishing Hessian
determinant loses a variable after a linear change, we reach
constant linear coordinates in which $\phi_D$ depends on $r\in\{1,2,3\}$
variables with nonsingular Hessian in those variables. $r$ is the
Hessian rank of $\phi_D$.

In the three cases below, $r=\rk\Hess\phi_D$ refers to the top
homogeneous part, and the rank of a wall form $P$ or $Q$ is the rank of
its own Hessian.

\emph{The case $\rk\Hess\phi_D=3$.} Here $\phi_D=h(x,y,z)$ with
$\det\Hess_3h\ne0$. Keep the weights of $x,y,z$ at one and raise that
of $w$. This is the situation of Lemma \ref{lem:T}. Its wall form $Q$
has $\rk\Hess Q=3$ and no constant kernel, and its apex space lies in
$\spn(e_x,e_y,e_z)$ and has dimension one or two (Lemma \ref{lem:T}).
At the wall the weights of $x,y,z$ are equal to one and the weight of
$w$ is larger, so $\spn(e_x,e_y,e_z)$ is a single weight block, the
block of smallest weight. The apex space splits into its intersections
with the weight blocks (Lemma \ref{lem:apexbasic}(a)), and its
intersection with the block $Ke_w$ is zero. Hence the apex space is
either a line or a plane in this block, and no other configuration
occurs. A three-dimensional apex space would give a constant kernel by
Lemma \ref{lem:apexbasic}(c). At the wall $Q$,
a single apex lies in the smallest-weight ternary block and
Lemma \ref{lem:D5} finishes. A two-dimensional apex space lies in the
ternary block of equal weights, a constant change aligns it with
$x,y$, and Lemma \ref{lem:B} finishes, condition ($*$) being vacuous
for $w_x=w_y$.

\emph{The case $\rk\Hess\phi_D=2$.} Lemma~\ref{lem:tree-rk2}. There
the two unused weights are raised together. At the first wall the wall
form $P$ either has a constant kernel, and a change among the unused
variables gives a ternary form treated by Lemma \ref{lem:T}, or it has
none, and $\rk\Hess P\in\{2,3\}$ (Lemma \ref{lem:kernelstep}). For
$\rk\Hess P=2$ the normal form \eqref{eq:form1} leads to
Lemma \ref{lem:K}. For $\rk\Hess P=3$ the apices are located by the
restriction principle, Lemma \ref{lem:restrict}, and Lemmas
\ref{lem:D5} and \ref{lem:B} finish.

\emph{The case $\rk\Hess\phi_D=1$.} Here $\phi_D=a_0s^D$, $a_0\ne0$. Keep $w_s=1$ and
raise the other three weights together (the situation of
Lemma \ref{lem:kernelstep} with $r=1$). At a wall $P$ before $W=0$, $P_{ss}\ne0$ by
$a_0s^D$. So $\rk\Hess P\ge1$, and $\rk\Hess P\le3$ by
Lemma \ref{lem:W}, since $W>0$. If $\rk\Hess P=2$, then by \ref{E7}
either $\Hess P$ has a two-dimensional constant kernel or $P$ has a
single apex, and these two subcases are treated below. If $\rk\Hess P=1$, then by Lemma
\ref{lem:apexbasic}(e) $P$ depends on one linear form only, so its
constant kernel is three-dimensional, splits, cannot contain the
$s$-axis, and is the whole heavy block, contradicting the new
terms at the wall. If $\rk\Hess P=2$ and $\Hess P$ has a two-dimensional constant
kernel, the kernel lies in the heavy block, since $P_{ss}\ne0$ and the kernel splits into blocks (Lemma \ref{lem:apexbasic}(a)). By Lemma \ref{lem:apexbasic}(d), applied to two independent kernel vectors, a linear change within the heavy block makes
$P$ independent of two of the three heavy variables. So $P$ is a binary form $H(s,t')$,
necessarily with nonsingular binary Hessian (a singular one would
have a further kernel vector by Lemma \ref{lem:apexbasic}(f)), and the two unused
variables have the common weight of the heavy block, which is maximal, so Lemma~\ref{lem:tree-binary} applies. If $\rk\Hess P=2$ and $P$ has a single apex, the apex is
$s$ ($P_{ss}\ne0$ excludes heavy apices by \ref{E7}), $P$ is
jointly affine in the heavy variables. Its coefficient at a heavy variable has weight $D-w$, where $w$ is the common heavy weight, and is a polynomial in $s$ alone, hence a multiple of $s^{D-w}$, and a change in the heavy block
leaves $P=a_0s^D+bts^p$, again a binary form with two unused
variables, and Lemma~\ref{lem:tree-binary} applies. If $\rk\Hess P=3$, Lemma~\ref{lem:tree-C1} applies. Its additional hypothesis holds, since the three variables other than $s$ were raised together and have a common weight at the wall.

By the case list at the beginning of this section and Lemma~\ref{lem:termination}, every branch terminates in an affine direction or a contradiction, and the proposition is proved.
\end{proof}

\begin{example}\label{ex:branch}
We follow one branch of the procedure for the potential $\phi=x_1x_4+(x_2+x_1x_3)\bigl(x_1x_2+(1+x_1^2)x_3\bigr)$ of Example~\ref{ex:case-ii}, which has the affine direction $e_4$. Its terms are
\[
x_1^3x_3^2,\quad 2x_1^2x_2x_3,\quad x_1x_2^2,\quad x_1x_3^2,\quad x_1x_4,\quad x_2x_3,
\]
so $\phi_5=x_1^3x_3^2$ is a binary form in $(x_1,x_3)$ with $\det\Hess_{(x_1,x_3)}\phi_5=-24x_1^4x_3^2\neq0$. This is the case $r=2$ of Lemma~\ref{lem:tree-rk2}.

(1) Raise the weights of $x_2$ and $x_4$ together, $w=(1,\omega,1,\omega)$ with $\omega\ge1$. The terms above have the weights $5$, $3+\omega$, $1+2\omega$, $3$, $1+\omega$, $1+\omega$, and $W=16-4\omega$ for $\omega\le2$. At $\omega=2$ the terms $2x_1^2x_2x_3$ and $x_1x_2^2$ join, with $W=8$, and the wall form is $P=x_1(x_2+x_1x_3)^2$. Its Hessian has rank three and the constant kernel vector $e_4$. This is case (b) of Lemma~\ref{lem:kernelstep}, and $P$ is a ternary form with the unused variable $x_4$ of maximal weight.

(2) Raise only the weight of $x_4$, $w=(1,2,1,\omega)$ with $\omega\ge2$, as in Lemma~\ref{lem:T}. The term $x_1x_4$ joins at $\omega=4$, where $W=12-2\omega=4$, and the wall form is $Q=x_1L^2+x_1x_4$ with $L=x_2+x_1x_3$, of Hessian rank three. Its gradient is $(\partial_1Q,\ 2x_1L,\ 2x_1^2L,\ x_1)$, so the closure of the gradient image is the hypersurface $Y_3=Y_2Y_4$, whose translation directions are the multiples of $e_1$. This is a single apex in the block $\{x_1,x_3\}$ of smallest weight, and Lemma~\ref{lem:D5} applies.

(3) In Lemma~\ref{lem:D5} the active variable is $x_1$, the phase is $L$ with the two phase variables $x_2$ and $x_3$, and $x_4$ is absent from the phase, with the affine coefficient $x_1$. After the small step of that lemma the weights of $x_2,x_3,x_4$ exceed one, and $Q$ is the leading form of Lemma~\ref{lem:D3} with $(w,y,z)=(x_2,x_3,x_4)$ and $(a,\ell,A,m,b,p)=(1,1,1,2,1,1)$. That lemma ends the branch in the affine direction $e_z=e_4$ of $\phi$, or, through Lemma~\ref{lem:D2}, in a contradiction. Here $e_4$ is indeed an affine direction of $\phi$.
\end{example}

\begin{example}[two first walls in the case $r=1$]\label{ex:branch-r1}
The potentials
\[
\phi_1=x_1^3+x_1x_3+x_2x_4,\qquad \phi_2=x_1^5+x_1^3x_4+x_1x_3+x_2x_4
\]
are doublings of the plane automorphisms $(x_1,x_2)$ and $(x_1,x_2+x_1^3)$, with $E=x_1^3$ and $E=x_1^5$, so $\det\Hess\phi_1=\det\Hess\phi_2=1$. Their top forms $x_1^3$ and $x_1^5$ have Hessian rank one, which is the case $r=1$ of Proposition~\ref{prop:tree}. Raise the weights of $x_2,x_3,x_4$ together, $w=(1,\omega,\omega,\omega)$ with $\omega\ge1$.

For $\phi_1$ the terms have the weights $3$, $1+\omega$ and $2\omega$. The first wall is at $\omega=\frac32$, where $W=1$, and the wall form is $P=x_1^3+x_2x_4$. Its Hessian has rank three and the constant kernel vector $e_3$, so the apex space is three-dimensional, and Lemma~\ref{lem:tree-C1} passes by Lemma~\ref{lem:T} to the ternary form $x_1^3+x_2x_4$, with nonsingular Hessian in $(x_1,x_2,x_4)$ and the unused variable $x_3$.

For $\phi_2$ the terms have the weights $5$, $3+\omega$, $1+\omega$ and $2\omega$. The first wall is at $\omega=2$, where $W=6$, and the wall form is $P=x_1^5+x_1^3x_4$. Its Hessian has rank two and the two-dimensional constant kernel spanned by $e_2,e_3$. Thus $P$ is a binary form in $(x_1,x_4)$ with nonsingular Hessian and two unused variables of equal weight, and Lemma~\ref{lem:tree-binary} applies.

So the same case $r=1$ leads at the first wall to a ternary and to a binary form.
\end{example}

\begin{proof}[Proof of Theorem~\textup{\ref{thm:main-aff}}]
Remove the affine part of $\phi$. If $\deg\phi\le2$, then $D_v^2\phi=v\T\Hess\phi\,v$ is a nondegenerate quadratic form on $\kk^4$, which has a nonzero isotropic vector since $\kk$ is algebraically closed. If $\deg\phi\ge3$, Proposition~\ref{prop:tree} applies.
\end{proof}

\section{Affine directions and the dichotomy}\label{sec:tang}

In this section we draw the consequences of Theorem~\ref{thm:main-aff}, which is proved in Sections~\ref{sec:translation} to~\ref{sec:tree}. Apart from this theorem they use only Sections~\ref{sec:prelim} to~\ref{sec:affine}.

\subsection{The dichotomy}\label{ssec:dichotomy}

\begin{corollary}\label{cor:equiv}
Let $\kk$ be a field of characteristic zero and let $\phi\in\kk[x_1,\dots,x_4]$ be a Keller potential.
\begin{enumerate}[label=\textup{(\alph*)}]
\item Exactly one of the following holds. Either $\grad\phi$ is a polynomial automorphism of $\kk^4$, or, after a linear change of coordinates over $\overline\kk$, $\phi=x_3C(x_1,x_2)+x_4A(x_1,x_2)+E(x_1,x_2)$ with $(C,A)$ a non-injective Keller map of the plane and $\deg C,\deg A\le\deg\phi-1$.
\item If $\kk$ is algebraically closed, then $\HC_4$ holds over $\kk$ if and only if $\JC_2$ holds over $\kk$.
\item If $\grad\phi$ is not a polynomial automorphism, then $\deg\phi\ge102$. Moreover $\deg\phi\ge109$, and if $\deg\phi\le125$, then the planar Keller map $(C,A)$ of \textup{(a)} has component degrees $72$ and $108$.
\end{enumerate}
\end{corollary}

\begin{proof}
(a) By Theorem~\ref{thm:main-aff} applied over $\overline\kk$, $\phi$ has an affine direction $v\in\overline\kk^{\,4}$. A linear change over $\overline\kk$ moves $v$ to $e_4$, and Theorem~\ref{thm:affine} applies. In its case (ii), and in its case (i) with $(C,A)$ injective, $\grad\phi$ is injective over $\overline\kk$, hence a polynomial automorphism of $\kk^4$ by Lemma~\ref{lem:descent}. In case (i) with $(C,A)$ not injective, $\grad\phi$ is not injective over $\overline\kk$, and the degree bound is read off the form of $\phi$ as in the proof of Theorem~\ref{thm:main-dir}. The two alternatives exclude each other.

(b) The implication $\HC_4\Rightarrow\JC_2$ is given by the doubling \eqref{eq:doubling} \cite{Men06}, \cite{dBvdE05}, see \cite[Remark~4.5]{AHL}. Conversely, under $\JC_2$ over $\kk=\overline\kk$ the second alternative of (a) cannot occur, since its planar Keller map would be an automorphism, hence injective.

(c) In the second alternative of (a) the map $(C,A)$ is a non-injective Keller map of the plane with $\deg C,\deg A\le\deg\phi-1$. By \cite{Moh83}, $\max(\deg C,\deg A)\ge101$, and \cite[Theorem~2.1]{GGHV22} gives $\max(\deg C,\deg A)\ge125$ or $\{\deg C,\deg A\}=\{72,108\}$.
\end{proof}

\begin{remark}\label{rem:commonplane}
In the second alternative of Corollary~\ref{cor:equiv}(a) the potential is affine jointly in the directions of the plane $U\subset\overline\kk^{\,4}$ spanned by $e_3,e_4$, that is, $u\T\Hess\phi(x)\,v=0$ for all $x$ and all $u,v\in U$. So $U$ is a totally isotropic plane, common to all the forms $\Hess\phi(x)$, as for every doubling. Theorem~\ref{thm:main-aff} gives a common isotropic line for every Keller potential in four variables. Geometrically, the quadrics $\{[v]\in\PP^3:v\T\Hess\phi(x)\,v=0\}$ all contain the line $\PP(U)$. Dually, $(\Hess\phi)^{-1}$ has a zero block in the variables $x_1,x_2$, so the covectors $(\xi_1,\xi_2,0,0)$ are characteristic at every point for the linearized operator $L_\phi u=c\,\tr\bigl((\Hess\phi)^{-1}\Hess u\bigr)$ of the equation $\det\Hess\phi=c$. Differentiating $\det\Hess\phi=c$ along a constant vector $v$ gives $L_\phi(D_v\phi)=0$, and for a doubling the solutions $D_{e_3}\phi=C$ and $D_{e_4}\phi=A$ of the linearized equation do not depend on $x_3,x_4$.
\end{remark}

\begin{remark}[Lagrangian interpretation]\label{rem:lagrangian}
The graph $\Lambda_\phi=\{(x,\grad\phi(x))\}$ of a Keller potential $\phi$ in $n$ variables is an exact Lagrangian subvariety of $T^*\A^n$. Its projection to the coordinates $x$ is an isomorphism, and $(\grad\phi)^*(dp_1\wedge\dots\wedge dp_n)=c\,dx_1\wedge\dots\wedge dx_n$, so its projection to the momentum coordinates $p$ is \'etale. By Proposition~\ref{prop:leg-global} the Legendre transform of $\phi$ is a polynomial if and only if $\Lambda_\phi$ is also a graph over the momentum space. For a map $F$ of the plane the doubling $\phi_F=\vv{y,F(x)}$ is a generating function of the cotangent lift $(x,p)\mapsto\bigl(F(x),JF(x)^{-\mathsf T}p\bigr)$ of $F$ \cite{Arn89}, with $p=JF(x)\T y$ and new momentum $y$, and adding $E(x)$ composes it with the translation $p\mapsto p+\grad E(x)$. By Corollary~\ref{cor:equiv}, for $n=4$ the graph $\Lambda_\phi$ fails to be a graph over the momentum space only if, in suitable linear coordinates, it arises in this way from a non-injective Keller map of the plane. The momentum projection has no critical points, so $\Lambda_\phi$ has no caustics, and the Hessian conjecture asks whether this everywhere regular projection is injective. By Remark~\ref{rem:lines}, $\Lambda_\phi$ contains the lines $\{(x+tv,\grad\phi(x)+t\Hess\phi(x)v):t\in\kk\}$ for an affine direction $v$, and for a doubling $E+x_3C+x_4A$ in four variables it is a union of affine planes over which $p_3=C$ and $p_4=A$ are constant.
\end{remark}

\begin{corollary}\label{cor:R4}
$\Rc_4=\varnothing$, while $\Rc_5\neq\varnothing$.
\end{corollary}

\begin{proof}
The first statement is Theorem~\ref{thm:main-aff}, the second is Theorem~\ref{thm:main-MY}.
\end{proof}

\subsection{Real solutions}\label{ssec:real}

Over $\R$ the results above are controlled by the signature of $\Hess\phi$, which is constant on $\R^4$. A negative Hessian determinant forces the gradient to be a polynomial automorphism, and for a positive one of degree at least three the signature is $(2,2)$, the signature of the doubling, which is then the only possible source of a non-injective gradient (Corollary~\ref{cor:sign}). For the linearized operator $L_\phi$ of Remark~\ref{rem:commonplane}, whose principal symbol $c\,\xi\T(\Hess\phi)^{-1}\xi$ has the signature of $\Hess\phi$ up to sign, these are the elliptic, the hyperbolic and the ultrahyperbolic case, and a non-injective gradient can occur only in the ultrahyperbolic case.

\begin{remark}\label{rem:realdir}
The assumption that $\kk$ is algebraically closed in Theorem~\ref{thm:main-aff} cannot be removed. For $\phi=\frac12(x_1^2+x_2^2+x_3^2+x_4^2)\in\R[x]$ one has $\det\Hess\phi=1$ and $D_v^2\phi=\vv{v,v}>0$ for every $v\in\R^4\setminus\{0\}$, so $\phi$ has no real affine direction. More generally, a real affine direction $v$ satisfies $v^{\mathsf T}\Hess\phi(x)\,v=0$ for all $x$, so it exists only if $\Hess\phi$ is indefinite. The conclusions on injectivity do not depend on the field, since an automorphism over $\overline\kk$ which is defined over $\kk$ is an automorphism over $\kk$ (Lemma~\ref{lem:descent}). This is why the conclusions on injectivity hold over $\kk$ and over $\R$, although the affine direction is found only over $\overline\kk$. By Corollary~\ref{cor:realaff} below, over $\R$ the quadratic example is the only exception. The passage from a complex to a real affine direction in its proof is not a formal descent. A non-real affine direction $a+ib$ gives real linear coordinates in which $\phi=\frac12f\,(s^2+t^2)+sA+tB+C$ with $f,A,B,C$ independent of $s,t$, and the constancy of the Hessian determinant then forces a real affine direction.
\end{remark}

\begin{corollary}\label{cor:realaff}
Let $\phi\in\R[x_1,\dots,x_4]$ with $\det\Hess\phi\in\R^\times$. Then $\phi$ has an affine direction $w\in\R^4\setminus\{0\}$, unless $\deg\phi=2$ and $\Hess\phi$ is definite.
\end{corollary}

\begin{proof}
If $\deg\phi\le2$, then $\deg\phi=2$, $\Hess\phi$ is a constant nondegenerate matrix, and a real affine direction is a nonzero real isotropic vector of $\Hess\phi$. It exists if and only if $\Hess\phi$ is indefinite.

Let $\deg\phi\ge3$. By Theorem~\ref{thm:main-aff} there is $v\in\C^4\setminus\{0\}$ with $D_v^2\phi=0$. Write $v=a+ib$ with $a,b\in\R^4$. If $a$ and $b$ are linearly dependent over $\R$, then $v$ is a complex multiple of a vector $w\in\R^4\setminus\{0\}$, and $D_w^2\phi=0$. Otherwise the real and imaginary parts of $D_v^2\phi=0$ give $D_a^2\phi=D_b^2\phi$ and $D_aD_b\phi=0$. In real linear coordinates $(s,t,x,y)$ with $a=e_s$ and $b=e_t$ one has $\phi_{st}=0$, so $\phi=P(s,x,y)+R(t,x,y)$, and $\phi_{ss}=\phi_{tt}$ depends only on $x,y$. Hence
\[
\phi=\tfrac12f\,(s^2+t^2)+sA+tB+C,\qquad f,A,B,C\in\R[x,y].
\]
If $f=0$, then $a$ is a real affine direction. Let $f\neq0$. For $g\in\R[x,y]$ write $\grad g$ and $\Hess g$ for the gradient and the Hessian in $(x,y)$. The Hessian of $\phi$ has the block $f\,I_2$ in $(s,t)$, and the Schur complement gives $\det\Hess\phi=f^2\det S$ with
\[
S=S_0+sS_1+tS_2+(s^2+t^2)S_3,
\]
where $S_0=\Hess C-f^{-1}(\grad A\grad A\T+\grad B\grad B\T)$, $S_1=\Hess A-f^{-1}(\grad f\grad A\T+\grad A\grad f\T)$, $S_2$ is obtained from $S_1$ by replacing $A$ with $B$, and $S_3=\frac12\Hess f-f^{-1}\grad f\grad f\T=-\frac{f^2}2\Hess(1/f)$. For symmetric $2\times2$ matrices put $\vv{M,N}=\det(M+N)-\det M-\det N$. Since $\det S=c/f^2$ does not depend on $s,t$, the coefficients of $s^4$, $s^3$, $t^3$, $s^2$, $t^2$, $st$, $s$ and $t$ in $\det S$ vanish, that is
\begin{gather*}
\det S_3=\vv{S_1,S_3}=\vv{S_2,S_3}=\vv{S_1,S_2}=\vv{S_0,S_1}=\vv{S_0,S_2}=0,\\
\det S_1=\det S_2=-\vv{S_0,S_3}.
\end{gather*}
Note that $\vv{e_1e_1\T,M}=M_{22}$ for every symmetric $M$.

Suppose $S_3\neq0$. Then $\det\Hess(1/f)=0$ and $f$ is not constant. By Theorem~\ref{thm:gamma}, $f$ is a polynomial in one linear form. The vectors $c\in\C^2$ with $D_cf=0$ form a line, which by Lemma~\ref{lem:lindescent} is spanned by a real vector. So in suitable real linear coordinates $(r,q)$ of the $(x,y)$-plane one has $f\in\R[r]$ and $S_3=\sigma e_1e_1\T$ with $\sigma\neq0$. From $\vv{S_i,S_3}=0$ the entries $(S_i)_{22}$ vanish for $i=1,2$. Write $m_i=(S_i)_{12}$. Then $\det S_i=-m_i^2$ and $\vv{S_1,S_2}=-2m_1m_2$, so $m_1m_2=0$ and $m_1^2=m_2^2$, hence $m_1=m_2=0$. Now $(S_1)_{22}=A_{qq}$ and $(S_1)_{12}=A_{rq}-f^{-1}f'A_q$, so $A_q=\alpha f$ with $\alpha\in\R$, and similarly $B_q=\beta f$ with $\beta\in\R$. Then $\det S_1=0$ gives $\vv{S_0,S_3}=\sigma(S_0)_{22}=0$, that is $C_{qq}=(\alpha^2+\beta^2)f$. For $w=-\alpha e_s-\beta e_t+e_q$ one obtains
\[
D_w^2\phi=(\alpha^2+\beta^2)f-2\alpha A_q-2\beta B_q+C_{qq}=0 .
\]

Suppose $S_3=0$. Then $1/f$ is affine, so $f\in\R^\times$, and $S_1=\Hess A$, $S_2=\Hess B$. Thus $\det\Hess A=\det\Hess B=0$, and by Lemma~\ref{lem:hess2F} over $\R$ each of $\Hess A$, $\Hess B$ is a multiple of $\grad\ell\grad\ell\T$ for a real linear form $\ell$. If both are nonzero, then $\vv{\Hess A,\Hess B}=0$ forces their linear forms to be proportional. If $\Hess A\neq0$ or $\Hess B\neq0$, choose real linear coordinates $(r,q)$ in which $\Hess A$ and $\Hess B$ are multiples of $e_1e_1\T$, with a nonzero multiple $\tau$, say for $A$. Then $A_q$ and $B_q$ are real constants, and $\vv{S_0,S_1}=\tau(S_0)_{22}=0$ gives $C_{qq}=f^{-1}(A_q^2+B_q^2)$. For $w=-f^{-1}A_qe_s-f^{-1}B_qe_t+e_q$ one obtains $D_w^2\phi=-f^{-1}(A_q^2+B_q^2)+C_{qq}=0$.

Finally let $\Hess A=\Hess B=0$. Then $A$ and $B$ are affine, and
\[
\phi=\tfrac f2\bigl((s+f^{-1}A)^2+(t+f^{-1}B)^2\bigr)+\widetilde C,\qquad\widetilde C=C-\tfrac1{2f}(A^2+B^2)\in\R[x,y].
\]
The coordinates $(s+f^{-1}A,t+f^{-1}B,x,y)$ are real and affine, so $\det\Hess\phi=f^2\det\Hess\widetilde C$ and $\deg\widetilde C=\deg\phi\ge3$. By \ref{E:Dil} there is $u\in\C^2\setminus\{0\}$ with $D_u^2\widetilde C=0$. Let $h$ be the leading form of $\widetilde C$, of degree $d\ge3$. The component of degree $2d-4>0$ of $\det\Hess\widetilde C$ is $\det\Hess h$, so $\det\Hess h=0$, and Lemma~\ref{lem:hess2F} over $\R$ gives $h=\lambda\ell^d$ with $\lambda\in\R^\times$ and a real linear form $\ell$. From $D_u^2h=\lambda d(d-1)\ell(u)^2\ell^{d-2}=0$ one gets $\ell(u)=0$, so $u$ may be taken in $\R^2$. Then $w=-f^{-1}D_uA\,e_s-f^{-1}D_uB\,e_t+u$ is real, the two affine coordinates $s+f^{-1}A$ and $t+f^{-1}B$ are constant along $w$, and $D_w^2\phi=D_u^2\widetilde C=0$.
\end{proof}

\begin{corollary}\label{cor:realdoubling}
Let $\phi\in\R[x_1,\dots,x_4]$ with $\det\Hess\phi\in\R^\times$, and suppose that $\grad\phi$ is not a polynomial automorphism of $\R^4$. Then there are linear coordinates over $\R$ in which $\phi=x_3C(x_1,x_2)+x_4A(x_1,x_2)+E(x_1,x_2)$ with $C,A,E\in\R[x_1,x_2]$, and $(C,A)$ is a Keller map of $\R^2$ which is not a polynomial automorphism.
\end{corollary}

\begin{proof}
By Corollary~\ref{cor:equiv}(a) there are linear coordinates over $\C$ in which $\phi$ is a doubling. Let $U\subset\C^4$ be the plane spanned by the directions of $x_3$ and $x_4$ in these coordinates. Then $D_uD_{u'}\phi=0$ for all $u,u'\in U$, and since $\phi$ has real coefficients, also for all $u,u'\in\bar U$.

If $U\cap\bar U=0$, choose linear coordinates $z=(z_1,z_2)$, $w=(w_1,w_2)$ over $\C$ whose directions span $U$ and $\bar U$. All second derivatives of $\phi$ in $z$ vanish, so $\phi=g(w)+\sum_iz_ih_i(w)$, and the vanishing of all second derivatives in $w$ makes $g$ and the $h_i$ affine. Hence $\phi$ has degree at most two, $\Hess\phi$ is a constant invertible matrix, and $\grad\phi$ is an affine automorphism, contrary to the assumption.

If $U=\bar U$, then $U\cap\R^4$ is a real plane, and a linear change over $\R$ takes it to the plane spanned by $e_3,e_4$. Then $\partial_3^2\phi=\partial_3\partial_4\phi=\partial_4^2\phi=0$, so $\phi=x_3C+x_4A+E$ with $C,A,E\in\R[x_1,x_2]$.

If $\dim(U\cap\bar U)=1$, then $U\cap\bar U$ is spanned by a vector $v\in\R^4$, and $v$ is an affine direction of $\phi$. After a linear change over $\R$ one has $v=e_4$ and $\phi=x_4A'+B'$ with $A',B'\in\R[x_1,x_2,x_3]$. The linear change of $(x_1,x_2,x_3)$ in Theorem~\ref{thm:affine} takes $x_3$ along a vector $c\neq0$ with $D_cA'=0$. By Lemma~\ref{lem:lindescent} such a vector $c$ can be taken in $\R^3$, so the change can be taken over $\R$. In case (ii) of Theorem~\ref{thm:affine} the gradient of $\phi$ is injective, hence a polynomial automorphism by Lemma~\ref{lem:descent}, contrary to the assumption. So case (i) holds, and $\phi=x_3C+x_4A+E$ with $C,A,E\in\R[x_1,x_2]$.

In both remaining cases $\det\Hess\phi=(\det J(C,A))^2$ by Lemma~\ref{lem:doubling-det}, so $(C,A)$ is a Keller map of $\R^2$. If it were a polynomial automorphism, it would be injective over $\C$, and $\grad\phi$ would be injective over $\C$ by Theorem~\ref{thm:affine}(i), hence a polynomial automorphism of $\R^4$ by Lemma~\ref{lem:descent}.
\end{proof}

\begin{corollary}\label{cor:lorentz}
Let $\phi\in\R[x_1,\dots,x_4]$ with $\det\Hess\phi\in\R^\times$, and let $\Hess\phi$ have signature $(3,1)$ or $(1,3)$. Then $\grad\phi$ is a polynomial automorphism of $\R^4$.
\end{corollary}

\begin{proof}
Suppose not. By Corollary~\ref{cor:realdoubling} there are linear coordinates over $\R$ in which $\phi=x_3C+x_4A+E$ with $C,A,E\in\R[x_1,x_2]$. Then $\Hess\phi(x)=\begin{psmallmatrix}*&N\T\\ N&0\end{psmallmatrix}$ with $N=J(C,A)(x)$ invertible, and the plane spanned by $e_3,e_4$ is totally isotropic for $\Hess\phi(x)$. A nondegenerate real quadratic form of signature $(3,1)$ or $(1,3)$ has no totally isotropic plane, a contradiction.
\end{proof}

\begin{corollary}\label{cor:sign}
Let $\phi\in\R[x_1,\dots,x_4]$ with $\det\Hess\phi=c\in\R^\times$.
\begin{enumerate}[label=\textup{(\alph*)}]
\item If $c<0$, then $\grad\phi$ is a polynomial automorphism of $\R^4$.
\item If $c>0$ and $\deg\phi\ge3$, then $\Hess\phi$ has signature $(2,2)$.
\item If $\grad\phi$ is not a polynomial automorphism, then $c>0$, the degree of $\phi$ satisfies the bounds of Corollary~\ref{cor:equiv}\textup{(c)}, and in suitable real linear coordinates $\phi=x_3C(x_1,x_2)+x_4A(x_1,x_2)+E(x_1,x_2)$, so that the plane spanned by $e_3,e_4$ is totally isotropic for every $\Hess\phi(x)$.
\end{enumerate}
\end{corollary}

\begin{proof}
Since $\det\Hess\phi$ has no zero on $\R^4$, the signature $(p,q)$ of $\Hess\phi(x)$ does not depend on $x$, and the sign of $c$ is $(-1)^q$. So $c<0$ means signature $(3,1)$ or $(1,3)$, and (a) is Corollary~\ref{cor:lorentz}. If $c>0$, the signature is $(4,0)$, $(2,2)$ or $(0,4)$. In the definite cases $\phi$ or $-\phi$ is convex on $\R^4$, and by the theorem of J\"orgens, Calabi and Pogorelov \cite{Jor54,Cal58,Pog72} $\phi$ is quadratic, which gives (b). In (c), $c>0$ by (a), the degree bounds are Corollary~\ref{cor:equiv}(c), and the normal form is Corollary~\ref{cor:realdoubling}. The isotropy of the plane follows from $\partial_3^2\phi=\partial_3\partial_4\phi=\partial_4^2\phi=0$.
\end{proof}

\begin{remark}
Corollary~\ref{cor:lorentz} also follows from Theorem~\ref{thm:main-aff} and the passage from complex to real constant pivots in Lorentzian signature of Liu \cite[Propositions~4.9 and~2.5]{Liu26}. By \cite[Theorem~A.1]{Liu26}, a real constant pivot gives $\phi=xt+F(x,y,z)$ in suitable affine coordinates $(t,x,y,z)$.
\end{remark}

\begin{remark}[signature $(2,2)$]\label{rem:22}
The signature $(2,2)$ does not decide between the two alternatives of the dichotomy. The potential of Example~\ref{ex:case-ii} has $\det\Hess\phi=1$ and degree five, so its Hessian has signature $(2,2)$ by Corollary~\ref{cor:sign}\textup{(b)}, and its gradient is a polynomial automorphism by Theorem~\ref{thm:affine}\textup{(ii)}. In the other direction, in signature $(2,2)$ neither positivity of the Hessian determinant nor analyticity of the solution suffices for injectivity of the gradient. For a real polynomial map $F$ of the plane with $\det JF>0$ the doubling $\phi_F$ of \eqref{eq:doubling} satisfies $\det\Hess\phi_F=(\det JF)^2>0$ by Lemma~\ref{lem:doubling-det}, its Hessian has signature $(2,2)$, and $\grad\phi_F$ is injective if and only if $F$ is. The non-injective polynomial map of Pinchuk \cite{Pin94} therefore gives a polynomial $\phi$ with $\det\Hess\phi>0$ on $\R^4$ and non-injective gradient. The real analytic map
\[
F(a,\theta)=\sqrt2\,e^{a/2}\bigl(\cos(\theta e^{-a}),\,\sin(\theta e^{-a})\bigr)
\]
of $\R^2$ satisfies $\det JF=1$ and $F(a,\theta+2\pi e^a)=F(a,\theta)$. Hence $x_3F_1(x_1,x_2)+x_4F_2(x_1,x_2)$ is a real analytic solution of $\det D^2\phi=1$ on $\R^4$ of signature $(2,2)$ with non-injective gradient. Thus there is no analog of the theorem of J\"orgens, Calabi and Pogorelov for entire solutions of signature $(2,2)$, and both the polynomiality and the constancy of $\det\Hess\phi$ enter the problem.
\end{remark}

\section{Five variables}\label{sec:MY}

We show that the five-variable potentials of \cite{MY26} have a three-dimensional space of directions along which they are quadratic, but neither an affine direction nor a constant pivot. Both facts are read off from a few coefficients of $D_\xi^2\psi$ and $D_\xi^3\psi$.

In the coordinates $(x,y,z)$ of $\kk^3$ let $u=1+xy$ and
\[
F_1=u^3z+y^2u(4+3xy),\qquad F_2=y+3xu^2z+3xy^2(4+3xy),\qquad F_3=2x-3x^2y-x^3z .
\]
The map $F=(F_1,F_2,F_3)$ is the non-injective Keller map of Alp\"oge \cite[Section~1]{MY26}, with $\det JF=-2$. Put $\phi=y_1F_1+y_2F_2+y_3F_3$, a polynomial in $(x,y,z,y_1,y_2,y_3)$ which is affine in $z$, write $\phi=B+zA$ with $A,B\in\kk[V]$, $V=(x,y,y_1,y_2,y_3)$, and put
\[
\psi_{\lambda,\mu}=B+\tfrac\lambda2A^2+\mu A,\qquad\lambda\in\kk^\times,\ \mu\in\kk .
\]
These are the Schur descents of $\phi$ of \cite[Lemma~A.1]{MY26}. The potential $2\psi_{1,13/2}=A^2+13A+2B$ is the counterexample to $\HC_5$ of \cite{MY26}.

\begin{lemma}\label{lem:MYkeller}%
$\det\Hess\psi_{\lambda,\mu}=4\lambda$.
\end{lemma}
\begin{proof}
Put $\Phi=\phi-(z-\mu)^2/(2\lambda)=\tilde B+z\Lambda+\frac12Sz^2$ with $\tilde B=B-\mu^2/(2\lambda)$, $\Lambda=A+\mu/\lambda$ and $S=-1/\lambda$. In the variables $(x,y,z;y_1,y_2,y_3)$ the Hessian of $\Phi$ has the block form $\begin{psmallmatrix}M&JF^{\mathsf T}\\ JF&0\end{psmallmatrix}$ with $3\times3$ blocks, so $\det\Hess\Phi=-(\det JF)^2=-4$ by Lemma~\ref{lem:doubling-det}. Lemma~\ref{lem:legendre}(c) with $n=5$, $m=1$ and $p=0$ gives $\det\Hess_V\psi_0=(-1)^5(-4)/S=-4\lambda$, where $\psi_0=-\tilde B+\frac12\Lambda^2W$ and $W=-\lambda$. A direct computation gives $\psi_0=-\psi_{\lambda,\mu}$, and the claim follows, since $V$ has five components.
\end{proof}

For a polynomial $\psi$ in $n$ variables $V$ write
\[
D_\xi^2\psi=\sum_\alpha q_\alpha(\xi)\,V^\alpha
\]
with quadratic forms $q_\alpha$ in $\xi\in\kk^n$. A vector $\xi\neq0$ is an affine direction or a constant pivot of $\psi$ if and only if $q_\alpha(\xi)=0$ for all $\alpha\neq0$.

\begin{theorem}\label{thm:MYcore}%
For all $\lambda\in\kk^\times$ and $\mu\in\kk$ the potential $\psi_{\lambda,\mu}$ has neither an affine direction nor a constant pivot.
\end{theorem}
\begin{proof}
Write $A=\sum_{j=1}^3y_ja_j$ and $B=\sum_{j=1}^3y_jb_j$ with
\begin{gather*}
a_1=u^3,\qquad a_2=3xu^2,\qquad a_3=-x^3,\\
b_1=y^2u(4+3xy),\qquad b_2=y+3xy^2(4+3xy),\qquad b_3=2x-3x^2y ,
\end{gather*}
and write $\xi=(\xi',\xi'')$, where $\xi'=(\xi_1,\xi_2)$ collects the components along $x,y$ and $\xi''=(\xi_3,\xi_4,\xi_5)$ those along $y_1,y_2,y_3$. Put
\[
S=D_{\xi''}A=\xi_3u^3+3\xi_4xu^2-\xi_5x^3 .
\]
The operator $D_{\xi''}$ lowers the degree in $(y_1,y_2,y_3)$ by one, $D_{\xi'}$ preserves it, and $D_\xi^2=D_{\xi'}^2+2D_{\xi'}D_{\xi''}+D_{\xi''}^2$. Since $\psi_{\lambda,\mu}=B+\mu A+\frac\lambda2A^2$, the part of $D_\xi^2\psi_{\lambda,\mu}$ of degree zero in $(y_1,y_2,y_3)$ is
\[
\lambda S^2+2\sum_{j=1}^3\xi''_jD_{\xi'}(b_j+\mu a_j),
\]
the coefficient of $y_1$ in the part of degree one is
\[
D_{\xi'}^2(b_1+\mu a_1)+2\lambda\bigl(D_{\xi'}(a_1)\,S+a_1D_{\xi'}S\bigr),
\]
and the coefficient of $y_1^2$ in the part of degree two is $\frac\lambda2D_{\xi'}^2(u^6)$.

The polynomials $a_j,b_j$ in $(x,y)$ have degree at most six, so the sum over $j$ in the part of degree zero has degree at most five. The terms of $u^3$, $xu^2$ and $x^3$ are multiples of $x^ky^k$, $x^{k+1}y^k$ and $x^3$ respectively. Hence in $S^2$ the monomial $x^6y^6$ arises only from $(\xi_3x^3y^3)^2$, the monomial $x^6y^4$ only from $(3\xi_4x\cdot x^2y^2)^2$, and $x^6$ only from $(\xi_5x^3)^2$. The coefficients of $D_\xi^2\psi_{\lambda,\mu}$ at $x^6y^6$, $x^6y^4$ and $x^6$ are therefore $\lambda\xi_3^2$, $9\lambda\xi_4^2$ and $\lambda\xi_5^2$.

The coefficient of $y_1y^4$ is the coefficient of $y^4$ in the restriction to $x=0$ of the coefficient of $y_1$ displayed above. At $x=0$ the polynomials $a_1$, $S$, $D_{\xi'}a_1$, $D_{\xi'}S$ and $D_{\xi'}^2a_1$ have degree at most two in $y$, while $b_1=4y^2+7xy^3+3x^2y^4$ gives $D_{\xi'}^2b_1|_{x=0}=6\xi_1^2y^4+42\xi_1\xi_2y^2+8\xi_2^2$. So the coefficient of $y_1y^4$ is $6\xi_1^2$. Finally $u^6=\sum_k\binom6kx^ky^k$, and the monomial $x^6y^4$ of $D_{\xi'}^2(u^6)$ arises only from $\xi_2^2\partial_y^2(x^6y^6)=30\xi_2^2x^6y^4$. So the coefficient of $y_1^2x^6y^4$ is $15\lambda\xi_2^2$.

If $\xi$ is an affine direction or a constant pivot of $\psi_{\lambda,\mu}$, the coefficients of $D_\xi^2\psi_{\lambda,\mu}$ at all nonconstant monomials vanish. Hence $\lambda\xi_3^2=9\lambda\xi_4^2=\lambda\xi_5^2=6\xi_1^2=15\lambda\xi_2^2=0$, and $\xi=0$ since $\lambda\neq0$.
\end{proof}

\begin{remark}
By Theorem~\ref{thm:MYcore} and Lemma~\ref{lem:MYkeller} the potentials $\psi_{\lambda,\mu}$ lie in $\Rc_5$, which is Theorem~\ref{thm:main-MY}. In particular Corollary~\ref{cor:hc5} does not produce a counterexample to $\HC_4$ from them.
\end{remark}

\begin{proposition}\label{prop:psi-L0}%
Let $\lambda\in\kk^\times$, $\mu\in\kk$ and $L_0=\{\xi\in\kk^5:\xi_1=\xi_2=0\}$. Then $D_\xi^3\psi_{\lambda,\mu}=0$ if and only if $\xi\in L_0$.
\end{proposition}
\begin{proof}
Since $A$ and $B$ are linear in $(y_1,y_2,y_3)$, the potential $\psi_{\lambda,\mu}$ has degree two in these variables, and $D_\xi^3\psi_{\lambda,\mu}=0$ for $\xi\in L_0$. Conversely, with the notation of the proof of Theorem~\ref{thm:MYcore}, the part of $D_\xi^3\psi_{\lambda,\mu}$ of degree two in $(y_1,y_2,y_3)$ is $\frac\lambda2D_{\xi'}^3(A^2)$, since $D_{\xi''}$ lowers this degree. Its coefficients at $y_1^2$ and $y_3^2$ are $\frac\lambda2D_{\xi'}^3(u^6)$ and $\frac\lambda2D_{\xi'}^3(x^6)=60\lambda\xi_1^3x^3$. In the first the monomial $x^6y^3$ arises only from $\xi_2^3\partial_y^3(x^6y^6)=120\,\xi_2^3x^6y^3$. Hence $D_\xi^3\psi_{\lambda,\mu}=0$ gives $60\lambda\xi_2^3=60\lambda\xi_1^3=0$, and $\xi_1=\xi_2=0$.
\end{proof}

\begin{remark}\label{rem:psi-L0}
For $\xi\in L_0$ one has $D_\xi^2\psi_{\lambda,\mu}=\lambda\,(\xi_3u^3+3\xi_4xu^2-\xi_5x^3)^2$, which is not constant for $\xi\neq0$. Hence $\deg_\xi\psi_{\lambda,\mu}=2$ for $\xi\in L_0\setminus\{0\}$, and $\deg_\xi\psi_{\lambda,\mu}\ge3$ for $\xi\notin L_0$ by Proposition~\ref{prop:psi-L0}. With a coordinate $t$ along $\xi\in L_0\setminus\{0\}$ the potential is $a+bt+\frac12ct^2$ with $c\notin\kk$, the case which in four variables is Theorem~\ref{thm:quad}. The potential $\psi_{1,13/2}$ has non-injective gradient \cite{MY26} and, by Theorem~\ref{thm:MYcore}, no affine direction. So in five variables a direction of degree at most two implies neither that the gradient is an automorphism nor that the potential has an affine direction, and Theorem~\ref{thm:main-aff} has no analog in five variables.
\end{remark}

The construction of \cite{MY26} applies the partial Legendre transform in one variable to a doubling, which is jointly affine in three variables (Lemma~\ref{lem:MYkeller}). In its free form this mechanism does not descend to four variables. Let $\Phi=tA(w)+B(w)$ with $w=(w_1,\dots,w_4)$, and suppose that after a linear change of $w$ the polynomial $\Phi$ is jointly affine in $w_2,w_3,w_4$. In the order $(w_2,w_3,w_4,t,w_1)$ the matrix $\Hess\Phi$ has a zero block of size three in the upper left corner, so its first three rows have nonzero entries only in the last two columns. Hence $\rk\Hess\Phi\le4$ and $\det\Hess\Phi=0$, and $\Phi$ is not a Keller potential.

\section{Concluding remarks}\label{sec:open}

By Theorem~\ref{thm:main-aff} every Keller potential in four variables has an affine direction, and Theorem~\ref{thm:affine} reduces the injectivity of its gradient to $\HC_2$ and $\JC_2$. By Theorem~\ref{thm:main-MY} the potentials of \cite{MY26} have neither an affine direction nor a constant pivot. Thus $\Rc_4=\varnothing$ and $\Rc_5\neq\varnothing$. The directions of degree at most two of the potentials of \cite{MY26} considered in Section~\ref{sec:MY} form exactly the three-dimensional space $L_0$ of Proposition~\ref{prop:psi-L0}. Thus the implication $\JC_2\Rightarrow\HC_4$, which the argument of Meng does not give \cite[Remark~4.5]{AHL}, holds, and the validity of $\HC_4$ is equivalent to $\JC_2$, which is not decided here.

The dichotomy can be read in several ways. Algebraically, the base locus of the linear system of quadrics $v\T\Hess\phi(x)\,v$ is nonempty, and it contains a line if and only if $\phi$ is a doubling. Geometrically, there is a vector which is isotropic for every $\Hess\phi(x)$, and a doubling has a common totally isotropic plane (Remark~\ref{rem:commonplane}). For a doubling the fibers of $\grad\phi$ are those of the planar map, and the extension $\kk(x_1,\dots,x_4)/\kk(\grad\phi)$ has the degree of the planar extension and a normal closure with the same Galois group (Remark~\ref{rem:fibers}). In the Lagrangian picture, $\Lambda_\phi$ is a graph over the momentum space exactly when the Legendre transform of $\phi$ is a polynomial, and otherwise it arises from the cotangent lift of a non-injective Keller map of the plane (Remark~\ref{rem:lagrangian}). Over $\R$ the second alternative requires the signature $(2,2)$ (Corollary~\ref{cor:sign}). In each of these settings the only possible obstruction to a polynomial Legendre duality is a non-injective Keller map of the plane.

The proof of Theorem~\ref{thm:main-aff} uses no classification quoted without proof. Its strategy, the analysis of the Newton polytope along paths of weights with exit lemmas at the walls, is that of \cite{Tan26}. For potentials with a direction of degree at most two, and in particular for potentials of degree at most five, Theorems~\ref{thm:main-dir} and~\ref{thm:main-low} give the conclusions of Corollary~\ref{cor:equiv} by homogeneous layers, without weighted leading forms. In degree six the homogeneous layers and the homogeneous parts for the weights $(1,1,1,2)$ give the same conclusions (Theorems~\ref{thm:deg6}, \ref{thm:deg6r2}, \ref{thm:deg6r1} and Corollary~\ref{cor:deg6}). In every degree a top form of Hessian rank three determines the weights $(1,1,1,2)$ and a weighted form $\Psi$ with $\det\Hess\Psi=0$ (Proposition~\ref{prop:rk3-weights}).

For a fixed bound $d$ on the degree, Theorem~\ref{thm:main-aff} is a first-order statement in the language of fields, since the Keller condition and the existence of $v\neq0$ with $D_v^2\phi=0$ are polynomial conditions on the coefficients of $\phi$ and on $v$. By the compactness theorem it therefore holds over every algebraically closed field of characteristic $p>p_0(d)$, for some $p_0(d)$ which the argument does not determine. In positive characteristic the dichotomy has a different meaning, since $(x_1-x_1^p,x_2)$ is a non-injective Keller map of the plane over $\overline{\F}_p$, and its doubling is a Keller potential with an affine direction and non-injective gradient.

\appendix
\section{\texorpdfstring{A proof of \ref{E:dB62}}{A proof of (E2)}}\label{app:dB62}
We give a proof of \ref{E:dB62} which uses Lemma~\ref{lem:graph}(a), Lemma~\ref{lem:hess2F}, generic smoothness in characteristic zero and the incidence of lines and planes in $\kk^3$. Theorem~\ref{thm:cylinder}, whose proof uses \ref{E:dB62}, is not used.

Let first $\kk$ be algebraically closed, $A\in\kk[x_1,x_2,x_3]$, $F=\grad A$ and $H=\Hess A$.

\begin{lemma}\label{lem:app-low}
Let $n\ge2$ and $A\in\kk[x_1,\dots,x_n]$ with $\rk\Hess A\le1$. Then $\Hess A$ vanishes on a subspace of $\kk^n$ of dimension $n-1$. In particular, for $n=3$ and $\rk H\le1$, $H$ vanishes on a two-dimensional subspace of $\kk^3$.
\end{lemma}

\begin{proof}
Put $F=\grad A$ and $H=\Hess A$, and let $\rk H=1$. The closure of the image of $F$ is a curve. Choose $R_1,\dots,R_{n-1}$ in its ideal with $\grad R_1,\dots,\grad R_{n-1}$ linearly independent at a general point of the curve, and put $G_i=(\grad R_i)(F)$. Then $HG_i=0$, $D_{G_j}G_i=\Hess R_i(F)\,HG_j=0$ and $D_{G_j}F=0$. By the remark after Lemma~\ref{lem:graph}, applied to $G_1,\dots,G_{n-1}$ in turn, $F(x+\sum_is_iG_i(x))=F(x)$ for all $s_i$. For general $x$ the affine hyperplane $x+\ker H(x)$ therefore lies in the fiber of $F$ through $x$, and by Lemma~\ref{lem:linfib} with $k=n-1$ all these hyperplanes have the same direction $U$, so $HU=0$.
\end{proof}

\begin{lemma}\label{lem:app-binary}
Let $L$ be a field of characteristic zero and $b\in L[u,v]$ with $\det\Hess b=0$. Then $(\Hess b)c=0$ for some $c\in L^2\setminus\{0\}$.
\end{lemma}

\begin{proof}
Let $\rk\Hess b=1$ and pass to $\overline L$. Lemma~\ref{lem:graph}(a) gives $U\ne0$ with $(\Hess b)U=0$ and $D(U_i)=0$ for $D=U_1\partial_u+U_2\partial_v$, so $\grad b$ is constant on the lines $z+sU(z)$. By Lemma~\ref{lem:linfib} over $\overline L$, with $n=2$ and $k=1$, these lines are parallel, of a fixed direction $c$, and $(\Hess b)c=0$. By Lemma~\ref{lem:lindescent} such a vector $c$ can be taken in $L^2$.
\end{proof}

Let now $\rk H=2$. The closure of the image of $F$ is a surface $\{R=0\}$, where $R$ has least degree among the nonzero relations of $F$. Put $G=(\grad R)(F)$. By Lemma~\ref{lem:graph}(a) and the remark after it, $HG=0$, $G\ne0$, and $G(x+tG(x))=G(x)$, $F(x+tG(x))=F(x)$ for all $t$.

\begin{lemma}\label{lem:app-twoplane}
If the values of $G$ lie in a plane of $\kk^3$, the conclusion of \ref{E:dB62} holds. The same holds over any field $\kk$ of characteristic zero if the components of $G$ satisfy a nontrivial linear relation over $\kk$.
\end{lemma}

\begin{proof}
If the values of $G$ lie on a line $\kk c$, then $Hc=0$. Otherwise, after a linear change of coordinates $G_1=0$. Then $(\partial_{Y_1}R)(F)=0$, and minimality of $\deg R$ gives $\partial_{Y_1}R=0$. Hence $A_{x_2}$ and $A_{x_3}$ are algebraically dependent, and $\det\Hess_{(x_2,x_3)}A=0$. Put $t=x_1$ and $L=\kk(t)$. Lemma~\ref{lem:app-binary} over $L$ and Lemma~\ref{lem:hess2F} give
\[
A=Q(t,\ell)+m_2x_2+m_3x_3+m_0,\qquad \ell=px_2+qx_3,
\]
with $p,q,m_i\in L$, $(p,q)\ne0$ and $Q\in L[T]$. With $\Delta=pq'-qp'$ and $\eta=pm_3'-qm_2'$, where the prime is $d/dt$, a direct computation gives
\[
\det\Hess A=-Q_{TT}(t,\ell)\big(\Delta\,Q_T(t,\ell)+\eta\big)^2 .
\]
If $Q_{TT}=0$, then $A$ is affine in $(x_2,x_3)$ with coefficients in $\kk[t]$, the second alternative. Otherwise $Q_T(t,\ell)$ is nonconstant in $\ell$, so $\Delta=\eta=0$. Then $(p,q)=\lambda(p_0,q_0)$ with $(p_0,q_0)\in\kk^2\setminus\{0\}$ and $\lambda\in L^\times$, and $p_0m_3-q_0m_2=\gamma\in\kk$. For $c=(0,q_0,-p_0)$ one gets $\partial_cA=-\gamma$. The argument uses only linear algebra over $\kk$ and Lemma~\ref{lem:app-binary} over $L$.
\end{proof}

Choose coordinates with $G_3\ne0$, let $\Pi=\{x_3=0\}$ and put $\Psi(y,t)=y+tG(y)$ for $y\in\Pi$. For general $x$ the line $x+\kk G(x)$ meets $\Pi$ in $y=x-(x_3/G_3(x))G(x)$, and $G(y)=G(x)$. Hence $\Psi$ is dominant, $F(\Psi(y,t))=F(y)$, and $F|_\Pi$ has rank two.

\begin{lemma}\label{lem:app-rank3}
Let $\nu\colon\Pi\dashrightarrow\kk^3$ be rational with $\nu=\rho G$ for a nonzero rational function $\rho$ on $\Pi$, and put $\Phi(y,t)=y+t\nu(y)$. For general $y\in\Pi$ and every $t\in\kk$ the differential of $\Phi$ at $(y,t)$ has rank three.
\end{lemma}

\begin{proof}
Since $G(\Phi(y,t))=G(y)$ and $HG=0$, the vector $\partial_t\Phi=\nu(y)$ lies in the kernel of $dF$ at $\Phi(y,t)$ for every $t$. As $F(\Phi(y,t))=F(y)$, one has $d(F|_\Pi)_y\circ\mathrm{pr}=dF_{\Phi(y,t)}\circ d\Phi_{(y,t)}$, where $\mathrm{pr}$ is the projection of the tangent space of $\Pi\times\kk$ onto that of $\Pi$. If the image of $d\Phi_{(y,t)}$ had dimension at most two, then, since it contains the nonzero vector $\nu(y)\in\ker dF_{\Phi(y,t)}$, its image under $dF_{\Phi(y,t)}$ would have dimension at most one, and $\rk d(F|_\Pi)_y\le1$. This contradicts the fact that $F|_\Pi$ has rank two.
\end{proof}

\begin{lemma}\label{lem:app-gauss}
The values of $G$ lie in a plane of $\kk^3$.
\end{lemma}

\begin{proof}
Let $\delta_1,\delta_2$ be the coordinate derivations of $\Pi$. By Lemma~\ref{lem:app-rank3} with $\nu=G$, the polynomial $\det(e_1+t\,\delta_1G,\,e_2+t\,\delta_2G,\,G)$ in $t$ has no zero for general $y$, so it is constant, and its coefficient of $t^2$ gives $\det(\delta_1G,\delta_2G,G)=0$. So $\delta_1G$, $\delta_2G$ and $G$ are linearly dependent, and the restriction of the map $x\mapsto[G(x)]\in\PP^2$ to $\Pi$ has rank at most one. Since $\Psi$ is dominant and $G(\Psi(y,t))=G(y)$, the map $[G]$ on $\kk^3$ is the composition of a projection onto $\Pi$ with its restriction to $\Pi$, so it has rank at most one. If it is constant, we are done. Otherwise its function field has transcendence degree one, so it is
$\kk(s)$ by Lemma~\ref{lem:luroth}, and $[G]=[n(s)]$ with $n\in\kk(s)^3$. Let a prime denote $d/ds$. The vectors $n,n'$ are linearly independent, and we put $E_s=\spn(n,n')$.

Let $\delta,\delta'$ be derivations of $\kk(\Pi)$ with $\delta(s)=1$, $\delta'(s)=0$ and $\delta'\ne0$. Lemma~\ref{lem:app-rank3} with $\nu=n(s)$ shows that $\det(\delta y+t\,n',\,\delta'y,\,n)$ is constant in $t$ and nonzero. Hence $\det(n',\delta'y,n)=0$ and $\delta'y\notin\kk n$, so the level curves of $s|_\Pi$ have their tangents in $E_s$ and not in $\kk n$. A component $Z$ of a general level surface $\{s=c\}$ is the union of the lines of direction $n(c)$ through a component of the corresponding level curve of $s|_\Pi$. Its tangent plane is $E_c$ at every general point, so $Z$ is an affine plane with direction $E_c$.

Write $s=p/q$ with coprime $p,q\in\kk[x]$. The plane $E_s$ has the normal vector $n\times n'$. If $E_c$ does not depend on $c$, the vectors $n(s)$ lie in a fixed plane. Otherwise Lemma~\ref{lem:pencil} gives $e\in\kk^3\setminus\{0\}$ with $(n\times n')\T e=0$, that is, $e\in E_s$ for all $s$. In coordinates with $e=e_3$ this reads $n_1n_2'-n_2n_1'=0$, so $(n_1,n_2)$ has constant direction and $n(s)$ lies in a fixed plane.
\end{proof}

\begin{proof}[Proof of \ref{E:dB62}]
Let $\kk$ be a field of characteristic zero. If $\rk H\le1$, Lemma~\ref{lem:app-low} over $\overline\kk$ applies, and the subspace on which $H$ vanishes contains some $c\in\kk^3\setminus\{0\}$ by Lemma~\ref{lem:lindescent}. Then $\partial_cA\in\kk$. Let $\rk H=2$. The relation ideal of $F$ over $\kk$ is generated by an irreducible $R$, and $R$ stays irreducible over $\overline\kk$, since $\kk[F]\otimes_\kk\overline\kk=\overline\kk[F]$ is a domain. Hence $G$ has coefficients in $\kk$. By Lemma~\ref{lem:app-gauss} over $\overline\kk$ its components satisfy a nontrivial linear relation over $\overline\kk$, hence over $\kk$ by Lemma~\ref{lem:lindescent}, and Lemma~\ref{lem:app-twoplane} concludes.
\end{proof}

\begin{remark}
A different proof of \ref{E:dB62} is given by de~Bondt and van~den~Essen \cite[Theorem~3.3]{dBvdE04}, see also \cite[Theorem~6.2]{dB18}. Lemma~\ref{lem:app-gauss} is a special case of \cite[Corollary~5.6(i)]{dB18}, by which every quasi-translation in dimension at most three has its image in a linear subspace of dimension at most two. 
\end{remark}

\section*{Acknowledgements}
The authors thank Qiyue Tang for his kind and helpful correspondence concerning the preprint \cite{Tan26}. This research was supported by the AGH University of Krakow within a subsidy of the Polish Ministry of Science and Higher Education (grant no. 16.16.420.054).

\end{document}